\documentclass[preprint, 11pt, authoryear]{elsarticle}
\usepackage{lineno}
\usepackage{mathtools}
\usepackage{amsthm}
\usepackage{empheq}
\usepackage{amsmath}
\usepackage{amssymb}
\usepackage{cases}
\usepackage{bm}
\usepackage{hyperref}
\usepackage{enumitem}
\usepackage{graphicx}
\usepackage[utf8]{inputenc}
\usepackage{listings}
\usepackage{booktabs}
\usepackage[margin=1.0in]{geometry}
\usepackage[title]{appendix}
\usepackage{subcaption}
\usepackage{booktabs} % tables
\usepackage{xcolor}
\usepackage{siunitx}

\usepackage{todonotes}

\usepackage{natbib}
\setcitestyle{square,numbers}
\newcommand{\dx}{\; dx}
\newcommand{\inner}[2]{\left(#1, \; #2\right)}
\newcommand{\duality}[2]{\left\langle#1, \; #2\right\rangle_{\mathring{H}^1}}
\newcommand{\bs}[1]{\boldsymbol{#1}}
\newcommand{\trial}{\mathcal{V}}
\newcommand{\trialh}{\mathcal{V}_h}
\newcommand{\hcirc}{\mathring{H}^1(\Omega)}
\newcommand{\hcircm}{\mathring{H}_m^1(\Omega)}
\newcommand{\hcircdual}{\mathring{H}^{-1}(\Omega)}
\newcommand{\erf}{\mathrm{erf}}
\newcommand{\beq}{\begin{linenomath} \begin{equation}}
\newcommand{\eeq}{\end{equation} \end{linenomath}}
\newcommand{\bcases}{\begin{linenomath} \begin{subnumcases}}
\newcommand{\ecases}{\end{subnumcases} \end{linenomath}}
\newcommand{\stepbound}{l_{max}}
\newcommand{\lx}{l_1}
\newcommand{\ly}{l_2}
\newcommand{\lz}{l_3}

\DeclareMathOperator*{\argmin}{arg\,min}

\begin{document}

\title{Optimal design of chemoepitaxial guideposts for the directed self-assembly of block copolymer systems using an inexact-Newton algorithm}

\author{Dingcheng Luo\corref{cor1}\fnref{fn1}}
\ead{dc.luo@utexas.edu}
\author{Lianghao Cao\fnref{fn1}}
\author{Peng Chen\fnref{fn1}}
% \ead{cvr@sayahna.org}
\author{Omar Ghattas\fnref{fn1,fn2}}
\author{J. Tinsley Oden \fnref{fn1}}
% \ead[url]{www.stmdocs.in}
\cortext[cor1]{Corresponding author}
\fntext[fn1]{Oden Institute for Computational Engineering and Sciences, The University of Texas at Austin}
\fntext[fn2]{Department of Mechanical Engineering, The University of Texas at Austin}

\begin{abstract}
Directed self-assembly (DSA) of block-copolymers (BCPs) is one of the most promising developments in the cost-effective production of nanoscale devices. The process makes use of the natural tendency for BCP mixtures to form nanoscale structures upon phase separation. The phase separation can be directed through the use of chemically patterned substrates to promote the formation of morphologies that are essential to the production of semiconductor devices. Moreover, the design of substrate pattern can formulated as an optimization problem for which we seek optimal substrate designs that effectively produce given target morphologies.

In this paper, we adopt a phase field model given by a nonlocal Cahn--Hilliard partial differential equation (PDE) based on the minimization of the Ohta--Kawasaki free energy, and present an efficient PDE-constrained optimization framework for the optimal design problem. The design variables are the locations of circular- or strip-shaped guiding posts that are used to model the substrate chemical pattern. To solve the ensuing optimization problem, we propose a variant of an inexact Newton conjugate gradient algorithm tailored to this problem. We demonstrate the effectiveness of our computational strategy on numerical examples that span a range of target morphologies. Owing to our second-order optimizer and fast state solver, the numerical results demonstrate five orders of magnitude reduction in computational cost over previous work. The efficiency of our framework and the fast convergence of our optimization algorithm enable us to rapidly solve the optimal design problem in not only two, but also three spatial dimensions. 

% where the computational cost associated with repeatedly solving the forward problem would otherwise render the problem computationally prohibitive. 
\end{abstract}

\maketitle

% \linenumbers
\section{Introduction}\label{sec:introduction}

Directed self-assembly (DSA) of block-copolymers (BCPs) is an effective method for the manufacturing of nanoscale structures. In the semiconductor industry, such self-assembly processes are investigated as promising means of extending the capabilities of conventional lithographic techniques to manufacture patterned structures for electronic circuitry at increasingly fine scales \cite{ChenXiong20}. DSA of BCPs has also seen applications in areas such as the manufacturing of photonic metamaterials \cite{StefikGuldinVignoliniEtAl15}, virus filtration membranes \cite{YangParkYoonEtAl08}, and functional materials \cite{KimYangLeeEtAl10}, to name a few. The process has demonstrated potential in being scalable and cost effective for large scale production of such structures. 

At the core of this process are BCPs, which are polymer chains consisting of blocks of distinct polymer species. DSA of BCPs takes advantage of the natural behavior of the BCPs. Typically, the process uses diblock copolymers, which are linear polymer chains composed of two bonded blocks of mutually repelling polymer species. Due to the repulsion, a mixture of such BCPs tends to separate into its distinct phases at lower temperatures. However, due to the bonding between the two polymer blocks, phase separation is able to occur only at scales corresponding to the length of the polymer. This leads to the formation of periodic nanoscale patterns in the form of lamellae, spheres, or cylinders, depending on the composition of the BCP \cite{BatesFredrickson99}. 

The phase separation can be further guided to form patterns, or morphologies, of nanomanufacturing interest. For example, guidance can be provided by chemically patterned substrates (chemoepitaxy) \cite{KimSolakStoykovichEtAl03}, topologically patterned substrates (graphoepitaxy) \cite{SegalmanYokoyamaKramer01}, external electric fields \cite{JeonJinKimEtAl17}, optical lasers \cite{YongJinKimEtAl18}, and thermal gradients \cite{SinghYagerSmilgiesEtAl12}.
In this study, we are concerned with the process of chemoepitaxy, in which a substrate patterned by a chemical attracting one of the two polymer species is used to guide the phase separation. This encourages the resulting BCP morphology to align with the prescribed chemical pattern, while the remaining regions are filled by the natural tendencies of the BCP mixture in a behavior known as pattern interpolation \cite{DetcheverryLiuNealeyEtAl10}. Thus, the resulting morphology can have finer detail than the prescribed chemical pattern. 

Typically, the substrate pattern is prescribed in terms of guideposts, which are units of the attractive substrate chemical with predefined shapes, such as circles or strips. As the guideposts can be manufactured by traditional lithographic techniques, DSA of BCPs can serve as a way of enhancing such lithographic techniques to further refine the feature sizes that they are capable of producing. This naturally leads to the problem of designing optimal substrate patterns to achieve given target morphologies, and the process by which this is performed is often termed inverse design. Moreover, one typically desires guidepost configurations that are well-spaced both to ensure manufacturability and to make efficient use of pattern interpolation to multiply the density of features.
% Moreover, one seeks optimal designs with sparse features to achieve manufacturing efficiency by pattern interpolation and density multiplication. 

Computational models have proven to be valuable in the study of the BCPs. Commonly used computational models of BCP phase separation include Monte Carlo methods \cite{WangYanNealeyEtAl99}, self-consistent field theory (SCFT) \cite{Fredrickson06}, and Cahn--Hilliard equations arising from the minimization of the Ohta--Kawasaki free energy \cite{OhtaKawasaki86}. This selection of models represent a range of accuracy and computational costs. The Monte Carlo methods are particle based methods, which are accurate but highly expensive. On the other hand, the phase field model with the Ohta--Kawasaki free energy represents a coarse graining of the BCP dynamics. Despite its more modest computational cost, the Ohta--Kawasaki phase field model captures the phase separation behavior of BCPs under the influence of the substrate, and has been successfully applied to the study of chemoepitaxial DSA \cite{LiQiuYangEtAl10, WuDzenis06}. 

In addition to simulation, these computational models have been employed to study the inverse design problem of chemoepitaxial DSA. In \cite{QinKhairaSuEtAl13}, the authors consider the inverse design problem using circular chemical guideposts on the substrate to produce desired morphologies. The work adopts the Cahn--Hilliard equations on a two-dimensional (2D) domain to model the phase separation of a thin BCP film. The effect of the substrate chemical is imposed through an additional linear polymer-substrate interaction energy, where the substrate chemical strength is represented by a spatially-varying function parameterized by the guidepost locations. Optimization for the guidepost locations is then carried out using the covariance-matrix-adaptation evolution strategy (CMA-ES) \cite{HansenOstermeier96}---a derivative-free, and thus slowly converging, optimization method---to achieve a range of target morphologies. The same algorithm is used in \cite{KhairaQinGarnerEtAl14} to design strip guideposts. There, the authors adopt a three-dimensional (3D) BCP model based on Monte-Carlo simulations using an analogous linear polymer-substrate interaction energy. Due to the 3D nature of the problem considered, the strength of the polymer-substrate interaction in the free energy is generally modeled to decay in the direction away from the substrate.

In \cite{HannonGotrikRossEtAl13, HannonDingBaiEtAl14, GadelrabHannonRossEtAl17}, the authors solve a guidepost design problem in 2D similar to \cite{QinKhairaSuEtAl13}. However, optimization for the guidepost locations is carried out by directly minimizing a free energy functional defined in terms of the target morphology and the guidepost locations using an evolutionary strategy, without directly computing the equilibrium states associated with the designs. The resulting guidepost configurations are then checked by simulation using an SCFT algorithm, and successful configurations are also tested experimentally. The study also considers a post-processing step in which closely spaced guideposts in the optimized configuration can be adjusted and removed to ensure manufacturability. 
A similar approach overall is taken in \cite{ZhangZhangLinEtAl19} with a modified evolutionary optimization algorithm.

We remark that in most of the existing literature, the optimization of the substrate design is conducted using evolutionary algorithms, which do not make use of derivative information. Such algorithms typically exhibit slow convergence, if they can be guaranteed to converge at all.
Furthermore, in the time-dependent Cahn--Hilliard model adopted in \cite{QinKhairaSuEtAl13}, a random initial state is evolved to the equilibrium state by gradient flow dynamics, while in the SCFT model adopted in \cite{HannonDingBaiEtAl14}, a random initial chemical potential is iteratively updated to solve for the equilibrium state. In both cases, different instances of the initial state may give rise to different equilibrium solutions, even for the same guidepost configuration. This presents a challenge for optimal guidepost configurations to be robust to different choices of initial conditions or initial guesses, respectively. However, this dependence on initialization of the model is not typically discussed in existing inverse design studies.

% Though in a different context, there are several studies on the optimization problem for topographically directed assembly of the BCPs such as \cite{OuakninLaachiDelaneyEtAl18}. In these studies, the authors formulate the problem as a shape optimization problem using SCFT as the underlying model. Due to the high dimensionality associated with the shape optimization, they make use of adjoint or linearization based shape derivatives and utilize derivative based algorithm. 

In this work, we study the problem of optimally designing the chemoepitaxial guideposts for the purpose of DSA of BCPs. We formulate this as a PDE-constrained optimization problem, adopting the phase field model for the BCPs based on minimization of the modified Ohta--Kawasaki free energy with linear polymer-substrate interactions, following \cite{QinKhairaSuEtAl13}. In this model, the guideposts define a spatially varying function over the substrate which characterizes the strength of the polymer-substrate interaction. For a specified number of guideposts, we optimize their locations to produce desired target morphologies. Additionally, we impose a penalty in the cost function on the closeness of guideposts to promote separation between them.
%For a specified number of guideposts with proper penalization on their closeness, the optimization problem then becomes one of selecting the guidepost locations to optimally produce given target morphologies while maintaining manufacturing efficiency. 

We propose an efficient PDE-constrained optimization computational framework for the optimal guidepost design problem.
In particular, we make use of an efficient strategy developed in \cite{CaoGhattasOden22} for finding the equilibrium states given by the nonlocal Cahn--Hilliard equations.
The formulation in \cite{CaoGhattasOden22} achieves significant computational cost reduction by minimizing a free energy functional via an energy-stable Newton algorithm to find the equilibrium state starting from an initial guess, instead of using the gradient flow that corresponds to time dependent nonlocal Cahn--Hilliard equations. Adopting this framework for the optimal design problem effectively leads to a bilevel optimization problem, in which there is an inner optimization problem for finding the minimal free-energy state, and an outer optimization problem for guidepost locations. We derive the corresponding gradient and Hessian action of the design objective with respect to the design variables by the Lagrangian formalism, based on the optimality conditions for the equilibrium states. We then tailor an inexact Newton conjugate gradient (CG) algorithm for the outer-level guidepost optimization problem to enhance the efficiency of the optimizer.  
% making use of derivatives which are available through the adjoint formulation. 
This allows us to improve computational efficiency by several orders of magnitude compared to existing procedures, such as those of \cite{QinKhairaSuEtAl13, KhairaQinGarnerEtAl14}, which use derivative-free optimization methods. Such methods require vast numbers of evaluations of the state problem due to their inability to efficiently exploit sensitivities of the design objective with respect to the design variables and curvatures of design space. Meanwhile, we are able to exploit this information scalably and efficiently using first order and second order adjoint PDEs. Importantly, we also discuss the influence of the initial guess for our state problem and investigate robustness of our solutions.

The rest of the paper is structured as follows. In Section \ref{sec:BCPs}, we describe the phase field model based on the Ohta--Kawasaki free energy and the fast solver developed in \cite{CaoGhattasOden22}. We then present the formulation of the optimal design problem in Section \ref{sec:design_problem}, followed by our proposed optimization procedure in Section \ref{sec:optimizer}. Finally, we present results from several 2D and 3D numerical experiments in Section \ref{sec:results} and provide our conclusions in Section \ref{sec:conclusion}.

% \pc{PC: complete this paragraph for the paper structure by briefly stating each section.}
% Directed self assembly of block-copolymers is a{}n effective method for manufacturing nanoscale structures of interest. In the semiconductor industry, such self-assembly processes are being investigated as promising alternatives for the manufacturing of patterned structures used electronic circuitry. This process makes use of the natural tendency of the phases of block copolymer mixtures to form periodic structures at equilibrium after phase separation. In addition, the process can be directed by various means to assist the formation of target morphologies, for example using a topologically (graphoepitaxy) or chemically (chemoepitaxy) pre-patterned substrate, thermal annealing, and external electrical fields. 

% Here, we study in particular the process of chemoepitaxy and its corresponding optimal design problem. In chemoepitaxy, a substrate chemical is placed at specified locations in an otherwise neutral substrate, which is to be then coated by the BCP. This substrate chemical attracts one of the BCP phases, and act to guide the phase separation of the BCP mixture. The design of the substrate pattern is therefore important for the success of the directed self assembly process. The pattern must be prescribed such that the phase separation leads to the formation of the desired morphology. It is also desirable for the substrate pattern to consist of features that are sparser spatially compared to the target morphology, as otherwise no manufacturing efficiency is gained. 

\section{Phase field model of BCPs}
\label{sec:BCPs}
We consider modeling of the equilibrium state of a BCP mixture under the influence of a substrate chemical field. In this study, we restrict ourselves to the case of diblock copolymers, which are BCPs consisting of two chemically distinct polymer species. We refer to one of the two species as phase A and the other as phase B. 

Given some domain $\Omega \subset \mathbb{R}^{d}$ in dimension $d=2$ or $3$, the state of the BCP mixture is represented by a function $u: \Omega \rightarrow [-1, 1]$, referred to as the order parameter field or phase field. The value of $u$ at $x \in \Omega$ gives an indication of the local composition of the mixture, taking a value of $u(x) = 1$ if locally the mixture is purely phase A, $u(x) = -1$ if purely phase B, and $-1 < u(x) < 1$ when interpolating between the pure phases. The equilibrium state of the BCP mixture is then given as the minimizer of a free energy functional defined in terms of the phase field $u$.

\subsection{Free-energy functional with substrate interaction}
We provide a definition of a free energy functional which accounts for the effects of the substrate chemical and derive the equations describing equilibrium state. The Ohta--Kawasaki energy functional $\mathcal{F}_{OK}$ is conventionally used for BCPs without substrate interactions, and is defined as 
\beq 
    \mathcal{F}_{OK}(u) = \int_{\Omega} W(u) \dx + \frac{\varepsilon^2}{2}\int_{\Omega} |\nabla u|^2 \dx + \frac{\sigma}{2} \int_{\Omega} (u-m)(-\Delta_N)^{-1}(u-m) \dx, 
\eeq
where $W(u) = \frac{1}{4}(u^2-1)^2$ is a double-well potential. Here, $\varepsilon$, $\sigma$ and $m$ are model parameters corresponding to the material properties. Specifically, $\varepsilon$ describes the strength of the repulsion between the two polymer phases, while $\sigma$ describes the strength of the nonlocal effect provided by the bonding of the polymer chain. The mass average $m$ represents the ratio between the masses of the two polymer phases in the overall mixture, and is defined such that $\int_{\Omega} (u - m)\dx$ = 0. Furthermore, $(-\Delta_N)^{-1}$ is an inverse Laplace operator, with $w = (-\Delta_N)^{-1}(u-m)$ satisfying
\bcases{}
    -\Delta w = u - m \quad
    &in $\Omega$,
    \\
    \nabla w \cdot n = 0 \quad 
    &on $\partial \Omega $,
\ecases
and $\int_{\Omega}w\dx = 0$. 
To account for the substrate chemical interactions, an additional energy term can be defined following \cite{QinKhairaSuEtAl13} as
\beq
    \mathcal{F}_{sub}(u) = \int_{\Omega} f u \dx.
\eeq
In this expression, $f(x)$ is a spatially varying function describing the effects of the substrate. We refer to $f(x)$ as the substrate interaction field. If $f(x) > 0$, the region locally attracts phase B ($u=-1$), while if $f(x) < 0$ the region locally attracts phase A ($u = 1$). When $f(x) = 0$, there is no preferential attraction for either of the two phases. Moreover, the magnitude of the function defines the strength of this attraction. The full free energy functional can then be written as 
\beq
    \mathcal{F}(u) = \mathcal{F}_{OK}(u) + \mathcal{F}_{sub}(u).
\eeq
The equilibrium state of the phase field is found by the minimization of the free energy functional while preserving a fixed mass-average of the phase field. To this end, we introduce the subspace of $H^1(\Omega)$ functions with zero average, i.e., 
\beq
    \hcirc := \{u \in H^1(\Omega) : \int_{\Omega}u \dx = 0\},
\eeq
with norm
\beq
\|u\|_{\hcirc} := (\nabla u, \nabla u)^{1/2},
\eeq
where $\inner{u}{v} = \int_{\Omega} uv\dx$ denotes the $L^2(\Omega)$ inner product. We also introduce the affine subspace of $H^1(\Omega)$ functions with average $m$,
\beq
\hcircm := \{u \in H^1(\Omega) : \int_{\Omega}(u - m) \dx = 0\}.
\eeq
We outline some basic properties of the function spaces in \ref{appendix:hcirc}. The equilibrium state is then given as a minimizer of the free energy functional, i.e.,
\beq
    u = \argmin_{u \in \hcircm} \mathcal{F}(u).
\eeq
The Gateaux derivative $D_u\mathcal{F}$ of this functional acting in any direction $\tilde{u}_0 \in \hcirc$ can be written as 
\beq
    \duality{D_u\mathcal{F}(u)}{\tilde{u}_0} = \inner{W'(u)}{\tilde{u}_0} + \inner{\varepsilon^2 \nabla u}{\nabla \tilde{u}_0} + \inner{\sigma(-\Delta_N)^{-1}(u-m)}{\tilde{u}_0} + \inner{f}{\tilde{u}_0}.
\eeq
Here we have used the $\duality{\cdot}{\cdot}$ to denote the duality pairing on $\hcirc$. Similarly, we have the second derivative for the functional, 
\beq
    \duality{D_u^2 \mathcal{F}(u)\hat{u}_0}{\tilde{u}_0} = \inner{W''(u)\hat{u}_0}{\tilde{u}_0}
    + \inner{\varepsilon^2 \nabla \hat{u}_0}{\nabla \tilde{u}_0}
    + \inner{\sigma(-\Delta_N)^{-1}\hat{u}_0}{\tilde{u}_0} \quad \forall \hat{u}_0, \tilde{u}_0 \in \hcirc.
\eeq
The first order optimality condition for the existence of a local minimizer is 
\beq
    \duality{D_u\mathcal{F}}{\tilde{u}_0} = 0 \quad \forall \tilde{u}_0 \in \hcirc,
\eeq
while the second order sufficient condition for the existence of a strict local minimizer is 
\beq
    \exists \; a > 0, \text{ such that } \duality{D_u^2\mathcal{F}\tilde{u}_0}{\tilde{u}_0} \geq a \|\tilde{u}_0\|^2_{\hcirc} \quad \forall \tilde{u}_0 \in \hcirc.
\eeq
The presence of the $\Delta_N^{-1}$ operator and the constraint of the mass average are inconvenient for computation. Instead, properties of the $\hcirc$ duality pairing can be used to derive an alternative representation of these conditions. We refer to \cite{CaoGhattasOden22} for a detailed discussion of this procedure. The transformed first order optimality condition for this system can be expressed in strong form as 
\bcases{}
    \Delta(W'(u) - \varepsilon^2 \Delta u + f) - \sigma(u - m) = 0 & $\text{in } \Omega$,
    \\
    \nabla u \cdot n = \nabla(W'(u) - \varepsilon^2 \Delta u + f) \cdot n = 0 & $\text{on } \partial \Omega$.
\ecases
We remark that when a gradient flow is used to reach the stationary point at equilibrium, one recovers the time dependent nonlocal Cahn--Hilliard equations
\bcases{}
    \partial_{t} u - \Delta(W'(u) - \varepsilon^2 \Delta u + f) + \sigma(u - m) = 0 & $\text{in } \Omega$,
    \\
    \nabla u \cdot n = \nabla(W'(u) - \varepsilon^2 \Delta u + f) \cdot n = 0 & $\text{on } \partial \Omega$.
\ecases
However, we instead directly seek a solution to the first order optimality condition, since minimizing the energy directly via a Newton method is several orders of magnitude faster than time-stepping to equilibrium \cite{CaoGhattasOden22}. Introducing a new variable $\mu = W'(u) - \varepsilon^2 \Delta u + f$, often referred to as the chemical potential, we can write the weak form of the first order optimality conditions as:
% Find $(u, \mu) \in H^1(\Omega) \times H^1(\Omega)$ such that 
\bcases{}
    (u, \mu) \in H^1(\Omega) \times H^1(\Omega) \nonumber \\
    (\nabla \mu, \nabla \tilde{u}) + (\sigma(u-m), \tilde{u}) = 0 & $\forall \tilde{u} \in H^1(\Omega) $, \label{eq:fwd_first}
    \\
    (\mu, \tilde{\mu}) - (W'(u), \tilde{\mu}) - (\varepsilon^2 \nabla u, \nabla \tilde{\mu}) - (f, \tilde{\mu}) = 0 &$\forall \tilde{\mu} \in H^1(\Omega)$.
\ecases
Notice that while both the trial and test spaces have been extended to include all of $H^1(\Omega)$, the mass constraint of $u \in \mathring{H}_m^{1}$ is still satisfied given these conditions. This can be seen by taking $\tilde{u} = 1$ as the test function, whereby \eqref{eq:fwd_first} reduces to $(u - m, 1) = 0$.

\subsection{Problem geometry and the thin film approximation}
The BCP phase separation process is physically a 3D problem. For the remainder of this paper, we restrict ourselves to rectangular domains $\Omega = (0, \lx) \times (0, \ly) \times (0, \lz)$. The substrate pattern is prescribed on the face at $x_3 = 0$ and is defined for the substrate domain $\Omega_s = (0, \lx) \times (0, \ly)$. In this set-up, the substrate interaction field $f$ is given by  
\beq
    f(x) = \tau(x_3) \eta(x_1, x_2).
\eeq
Here, the substrate pattern is described by $\eta : \Omega_s \rightarrow \mathbb{R}$, with $\eta < 0$ representing an A-attracting substrate, $\eta > 0$ representing a B-attracting substrate, and $\eta = 0$ representing a neutral substrate, as in the case of $f$. 
The function 
\beq \label{eq:tau-3}
\tau(x_3) = \exp(-x_3^2/2d_s^2)
\eeq 
extends the attraction of the substrate chemical normal to the plane $x_3 = 0$. This attractive force is assumed to decay as one moves away from the substrate at a rate characterized by a decay length scale $d_s$, which is defined in terms of the model parameters as
\beq \label{eq:ds}
d_s = \frac{12 \sqrt{3} c_s^2\varepsilon}{\sqrt{\sigma(1-m^2)}},
\eeq
where $c_s > 0$ is a material-dependent constant \cite{DetcheverryLiuNealeyEtAl10}.

For cases where the thickness of the mixture $\lz$ is on the order of the length scale of the polymer chain, we can approximate the problem as being purely two-dimensional and assume uniformity in the $x_3$ direction. This is often referred to as the thin film approximation. In this case, the PDE is solved on $\Omega = (0, \lx) \times (0, \ly)$. Furthermore, the substrate domain $\Omega_s$ coincides with the domain $\Omega$, corresponding to the fact that the polymer mixture is a thin film overlaying the substrate pattern. Thus, we may take $\tau = 1$, so that $f = \eta$ simply corresponds to the substrate pattern itself. We consider the cases of both $d=2$ and $d=3$ in our numerical experiments.

%%%%%% REWRITING THE NEWTON SOLVER PART TO BE BRIEF AND LEAVE DETAILS IN THE APPENDIX %%%%%% 

\subsection{Energy-stable Newton solver for the state equation}
% \dc{Condensed sections 2.3-2.6}
We employ the energy-stable Newton solver developed in \cite{CaoGhattasOden22} for the minimization of the free energy functional. Here, we briefly describe the elements of the algorithm and refer to \cite{CaoGhattasOden22} for a full presentation of the procedure. 

Given some initial guess for the state, $u_0$, we define an iterative procedure for $k \in \mathbb{N}$ by 
\beq\label{eq:StateUpdate}
    u^{(k+1)} = u^{(k)} + \beta_u^{(k)} \delta u^{(k)},
\eeq
where $u^{(k)}$ and $u^{(k+1)}$ refer to the states at the $k$th and $(k+1)$th iterations respectively, $\beta_u^{(k)}$ is the step size (for the state problem), and $\delta u^{(k)}$ defines the step direction.
% In the case of the standard Newton's method, the one solves
% \beq
%     \duality{D_u^2\mathcal{F}(u^{(k)}) \delta u^{(k)}}{\tilde{u}_0} = -\duality{D_u\mathcal{F}(u^{(k)})}{\tilde{u}_0} \quad \forall \tilde{u}_0 \in \hcirc
% \eeq 
% for the step direction $\delta u^{(k)}$.
% However, to ensure that $\delta u^{(k)}$ is a descent direction in the sense that,
% $$\duality{D_u\mathcal{F}(u^{(k)})}{\delta u^{(k)}} < 0,$$
% we instead solve a modified form of the Newton equations:
In the energy-stable Newton scheme, we find the step direction $\delta u^{(k)}$ by solving the modified Newton problem, 
% Find $(\delta u^{(k)}, \delta \mu^{(k)}) \in H^1(\Omega) \times H^1(\Omega)$ such that
\bcases{}
    (\delta u^{(k)}, \delta \mu^{(k)}) \in H^1(\Omega) \times H^1(\Omega) \nonumber \\
    \inner{\nabla \delta \mu^{(k)}}{\nabla \delta \tilde{u}} + \inner{\sigma \delta u^{(k)}}{\tilde{u}} = -\inner{\nabla \mu^{(k)}}{\nabla \tilde{u}} - \inner{\sigma(u^{(k)} - m)}{\tilde{u}} 
    & $\forall \tilde{u} \in H^1(\Omega)$,
    \\
    \inner{\delta \mu^{(k)}}{\tilde{\mu}} - \inner{W_\gamma''(u^{(k)})}{\tilde{\mu}} - \inner{\varepsilon^2\nabla \delta u^{(k)}}{\nabla \tilde{\mu}} = 0 
    & $\forall \tilde{\mu} \in H^1(\Omega)$.
\ecases
Here, we have introduced an auxiliary variable $\delta \mu^{(k)}$ in addition to the chemical potential $\mu^{(k)}$, which itself is given by
\beq
    \inner{\mu^{(k)}}{\lambda} - \inner{W'(u^{(k)})}{\lambda} - \inner{\varepsilon^2 \nabla u^{(k)}}{\nabla \lambda} - \inner{f}{\lambda} = 0 \quad \forall \lambda \in H^1(\Omega).
\eeq
Furthermore, the term
\beq
	W_{\gamma}''(u) := 2u^2 + \gamma(u^2-1)
\eeq
is an approximation for the second derivative of the double-well potential, and is parameterized by $\gamma \in [0,1]$ such that $W_1''(u) = W''(u)$ and $W_0''(u)$ is non-negative. In \cite{CaoGhattasOden22}, the existence of a $\gamma \in [0, 1]$ which ensures that $\delta u^{(k)}$ is a descent direction is established. Thus, at any step $k$, we perform a backtracking search for $\gamma
$ starting from $\gamma = 1$ down to $\gamma = 0$ until a descent direction is found, i.e. the step direction $\delta u^{(k)}$ satisfies
\beq
	\duality{D_u \mathcal{F}(u^{(k)})}{\delta u^{(k)}} < 0.
\eeq 
We then perform Armijo backtracking to find the step size $\beta_u^{(k)}$ that satisfies the descent condition
\beq
    \mathcal{F}(u^{(k)} + \beta_u^{(k)}\delta u^{(k)}) \leq \mathcal{F}(u^{(k)}) + c_{\text{Armijo}}\duality{D_u\mathcal{F}(u^{(k)})}{\beta_u^{(k)} \delta u^{(k)})},
\eeq
where $c_{\text{Armijo}}$ is the Armijo constant, typically chosen to be $c_{\text{Armijo}} = 10^{-4}$ \cite{NocedalWright06}.

\begin{sloppypar}
We assess the convergence of the Newton iterations using the residual function, $B: H^1(\Omega)^2 \rightarrow (H^1(\Omega)')^2$, which is given by
\begin{linenomath}
\begin{align}
\langle B(u, \mu), (\tilde{u}, \tilde{\mu})\rangle_{(H^1(\Omega))^2} := 
    & \left( (\nabla \mu, \nabla \tilde{u}) + (\sigma(u-m), \tilde{u}), \right. \nonumber \\
    &\left. (\mu, \tilde{\mu}) - (W'(u), \tilde{\mu}) - (\varepsilon^2 \nabla u, \nabla \tilde{\mu}) - (f, \tilde{\mu}) \right),
\end{align}
\end{linenomath}
and use the norm $\|B(u,\mu)\|_{(H^1(\Omega)')^2}$ as the convergence criterion, terminating when $\|B(u,\mu)\|_{(H^1(\Omega)')^2} \leq 10^{-8}$. Here $H^1(\Omega)'$ refers to the topological dual of $H^1(\Omega)$.
\end{sloppypar}
%%%%%%%% END OF NEW SECTION %%%%%%%%%
\subsection{Finite element discretization}
We employ a conventional finite element approach to discretize our function spaces. First, the domain $\Omega$ is discretized by either triangular elements in 2D or tetrahedral elements in 3D. The space $H^1(\Omega)$ is then approximated by finite dimensional subspaces $\trialh$ using continuous Lagrangian basis functions $\{N_i\}_{i=1}^{n}$, where $h$ corresponds to the mesh size and $n$ is the dimension of $\trialh$. 

For solving the Newton step problem, we use the same discretization $\trialh \subset \trial = H^1(\Omega)$ for the spaces of the state $u \in H^1(\Omega)$ and the chemical potential $\mu \in H^1(\Omega)$, as well as the step directions $\delta u^{(k)} \in H^1(\Omega)$ and $\delta \mu^{(k)} \in H^1(\Omega)$. Thus for $(u_h, \mu_h) \in \trialh \times \trialh$, we have
\beq
    u_h(x) = \sum_{i=1}^{n}(\bs{u})_i N_i(x) \quad \text{and} \quad \mu_h(x) = \sum_{i=1}^{n} (\bs{\mu})_i N_i(x),
\eeq
where $(\bs{u})_i$ is used to denote the $i$th component of a vector $\bs{u}$.  

In anticipation of the optimal design formulation, we also interpolate the field $f \in H^1(\Omega)$ into the subspace $\trialh$, which is discretized by the same basis functions. That is, we consider $f_h \in \mathcal{V}_h$ such that 
\beq
    f_h(x) = \sum_{i=1}^{n} (\bs{f})_i N_i(x).
\eeq

We implement our discretization and solvers using FEniCS, an open source finite element software package \cite{LoggMardalWells12}.

\subsection{Initial guesses and non-uniqueness of the equilibrium state}
\begin{figure}[htpb!]
\centering
\begin{subfigure}{0.32\textwidth}
    \centering
    \includegraphics[width=0.8\textwidth]{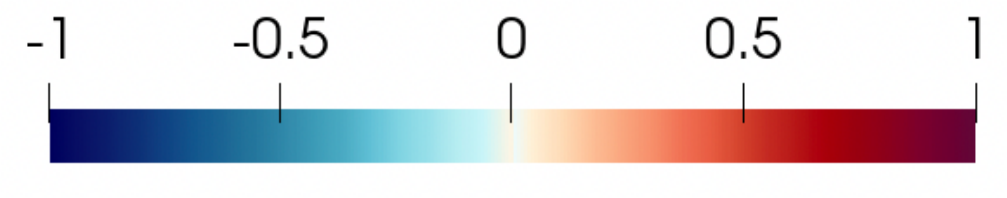}
\end{subfigure}

\begin{subfigure}{0.15\textwidth}
    \centering
    \includegraphics[width=\textwidth]{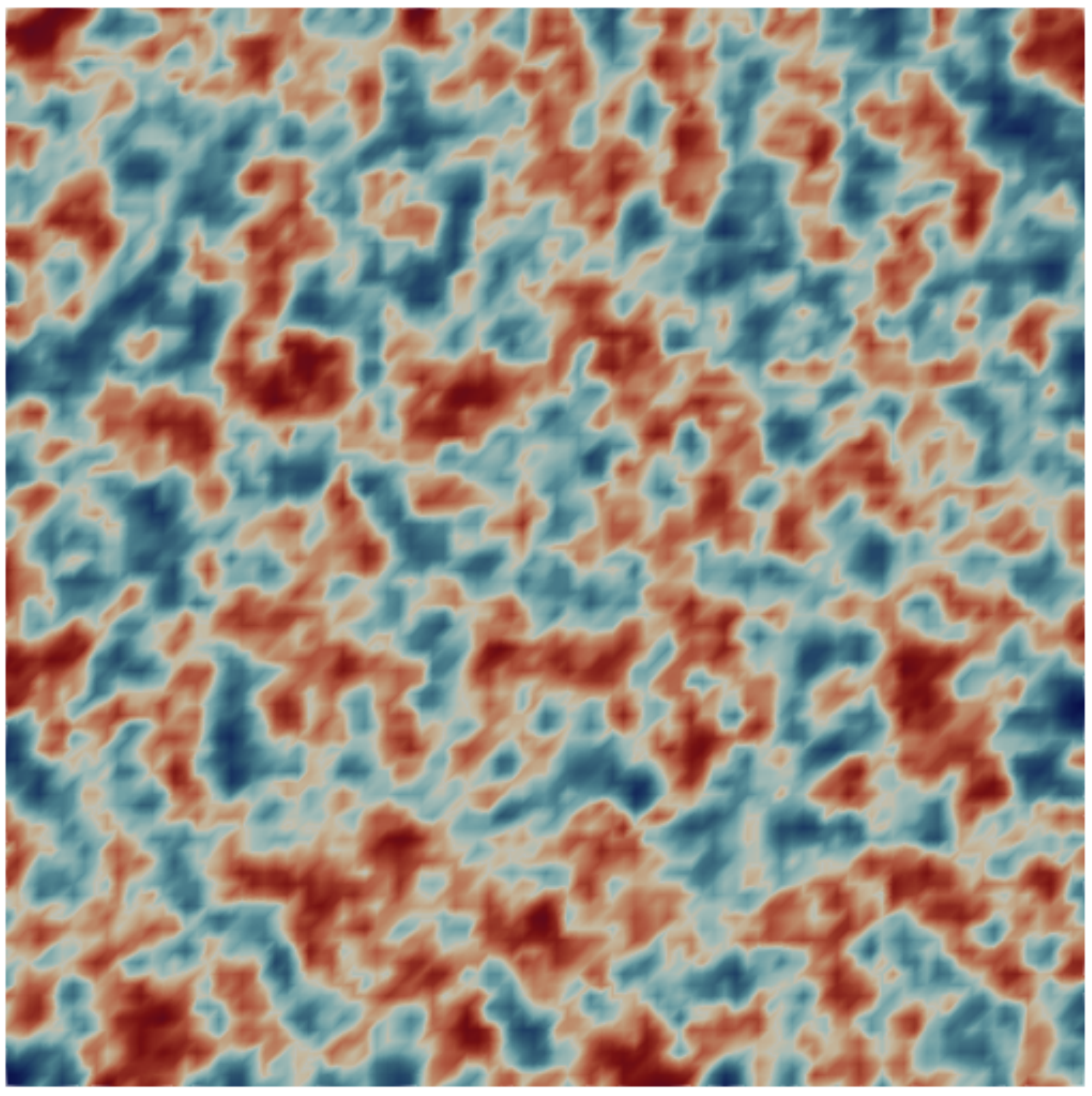}
\end{subfigure}
\begin{subfigure}{0.15\textwidth}
    \centering
    \includegraphics[width=\textwidth]{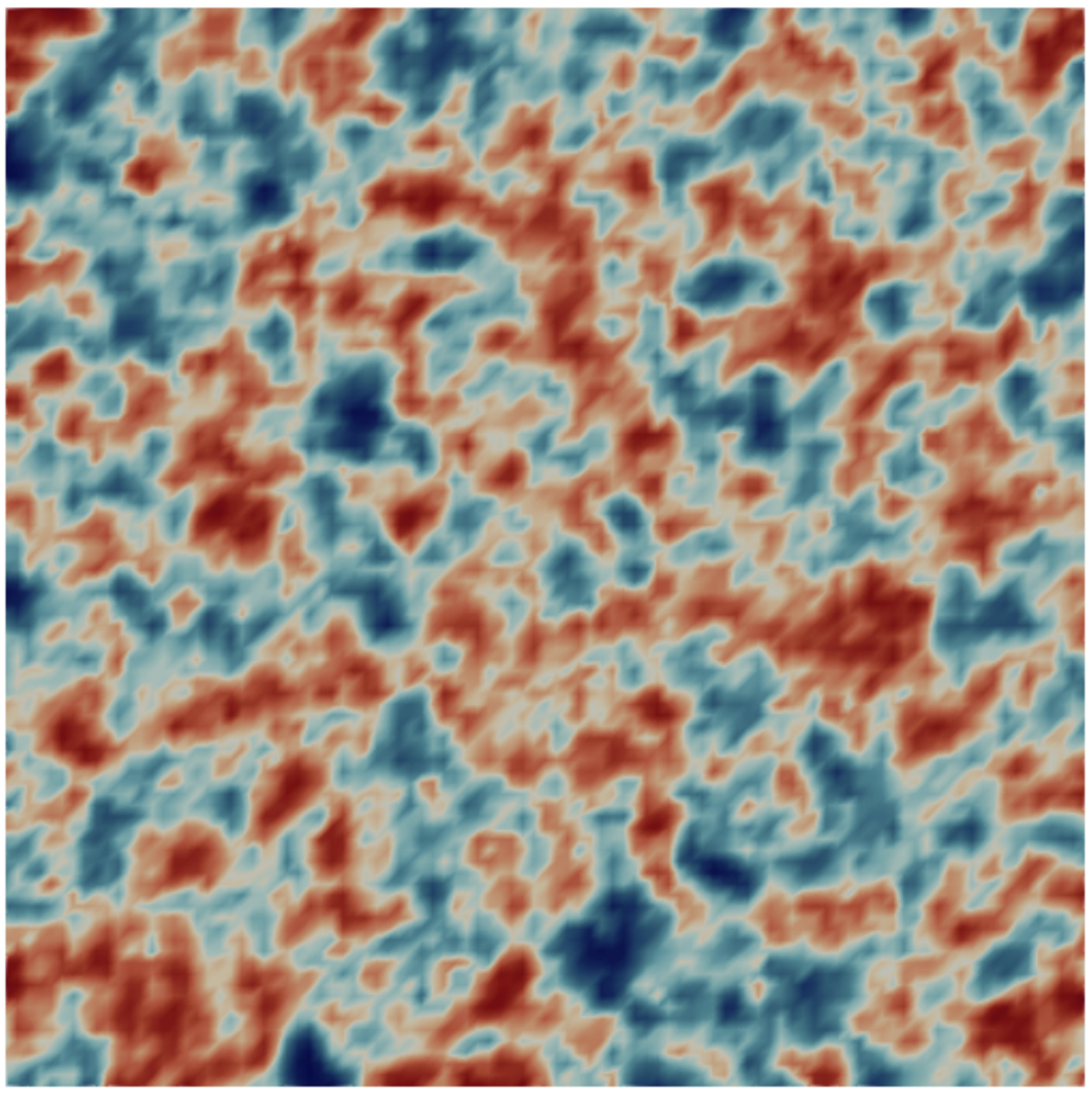}
\end{subfigure}
\begin{subfigure}{0.15\textwidth}
    \centering
    \includegraphics[width=\textwidth]{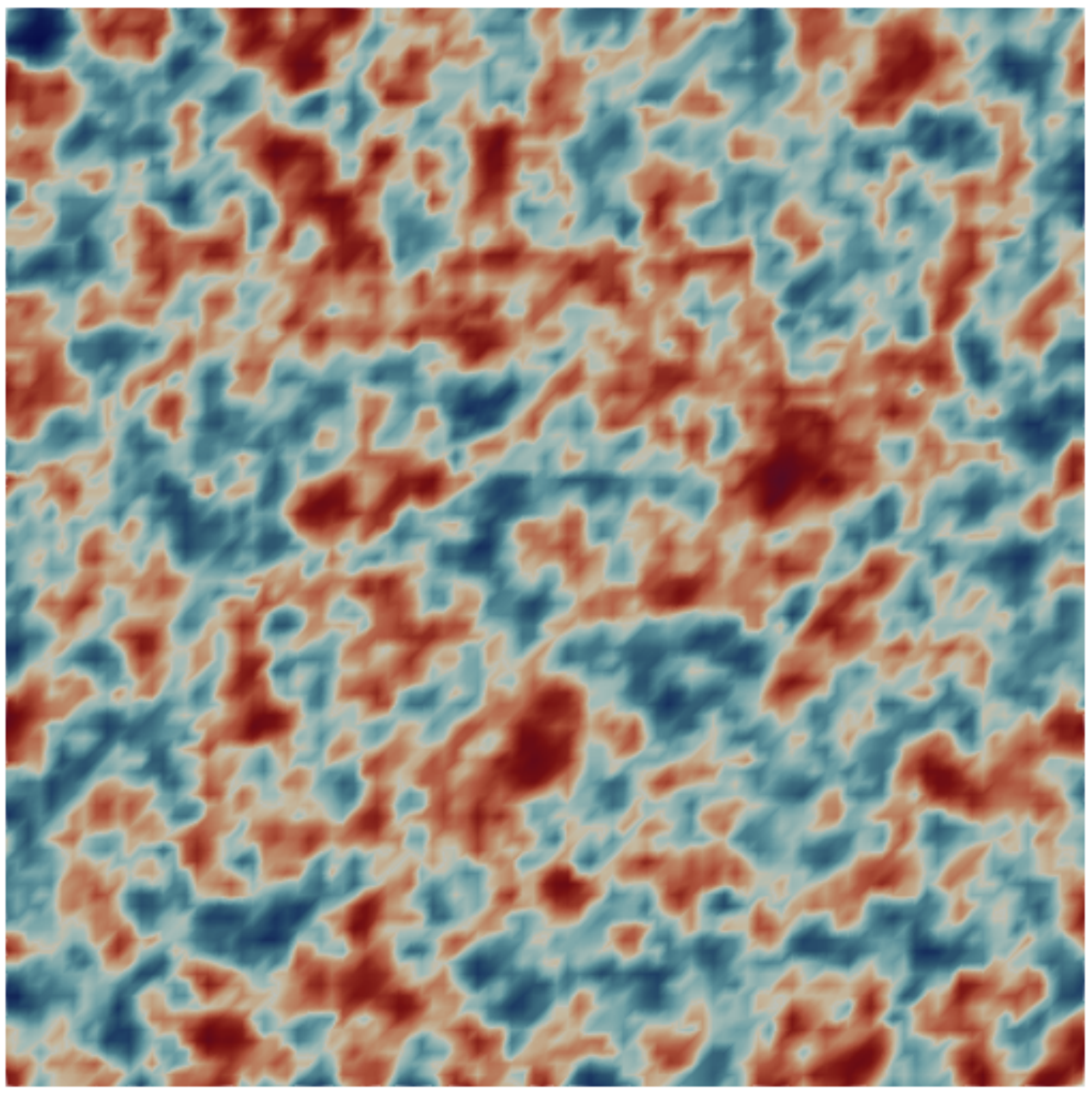}
\end{subfigure}

\begin{subfigure}{0.15\textwidth}
    \centering
    \includegraphics[width=\textwidth]{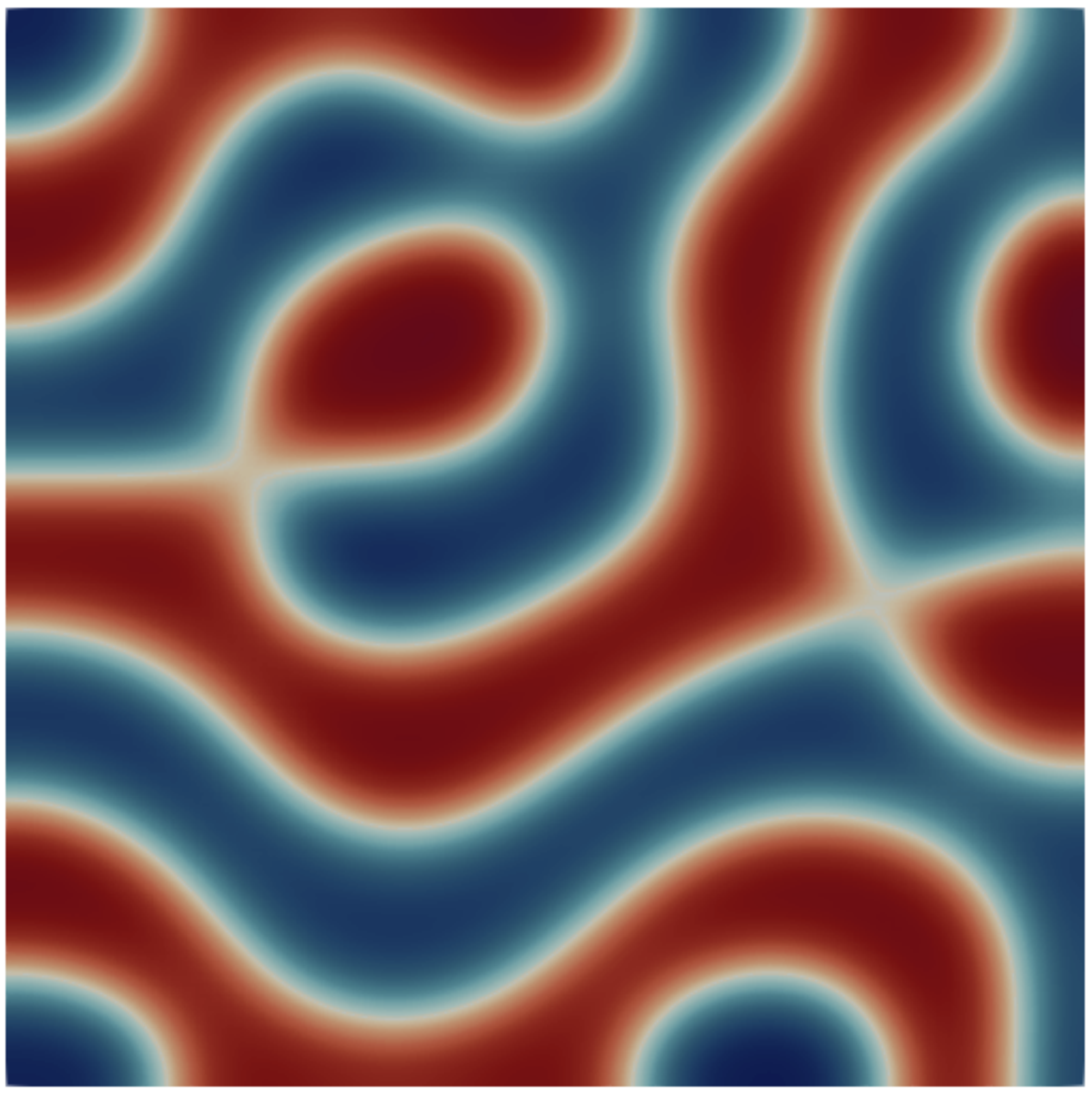}
\end{subfigure}
\begin{subfigure}{0.15\textwidth}
    \centering
    \includegraphics[width=\textwidth]{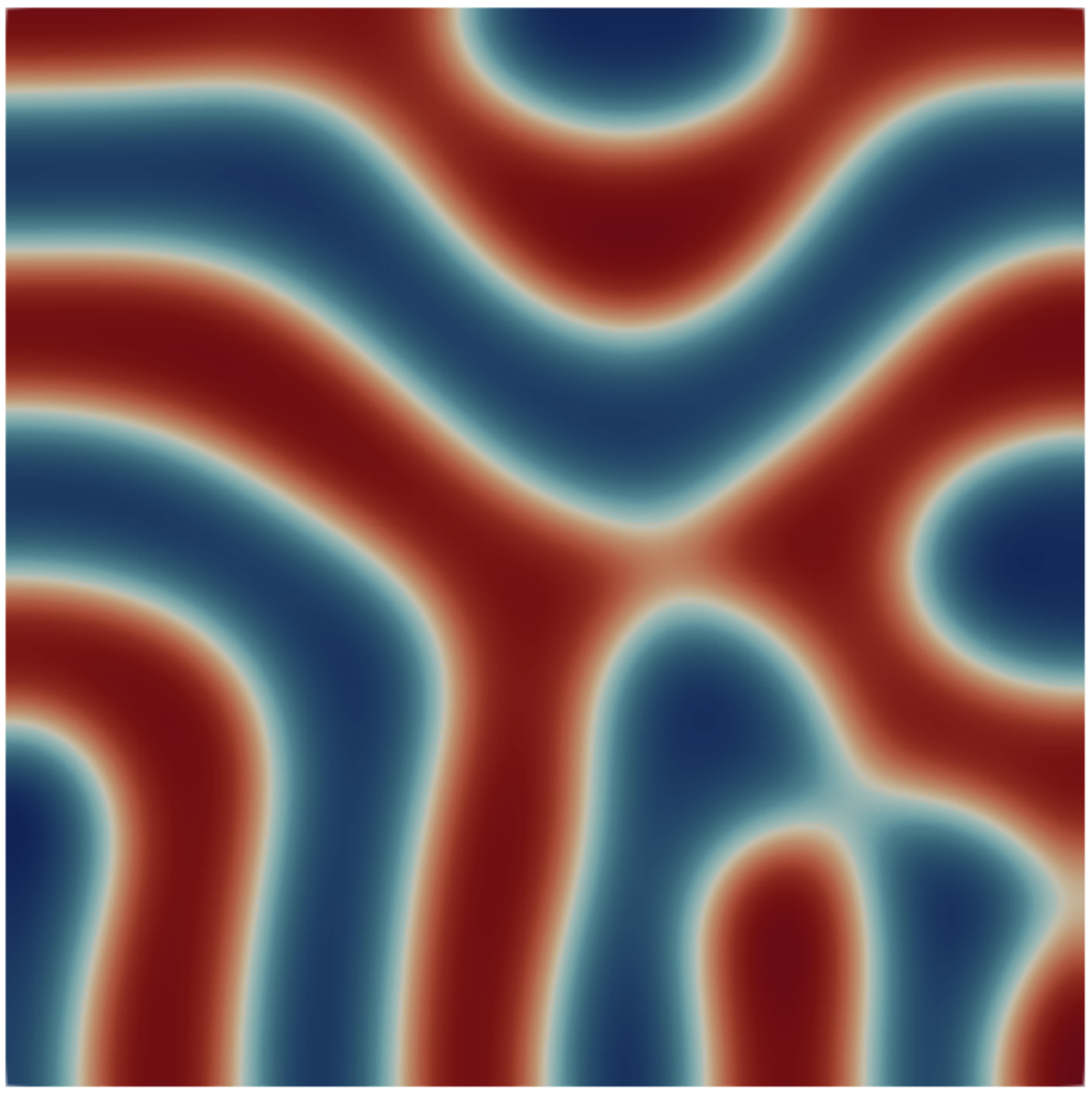}
\end{subfigure}
\begin{subfigure}{0.15\textwidth}
    \centering
    \includegraphics[width=\textwidth]{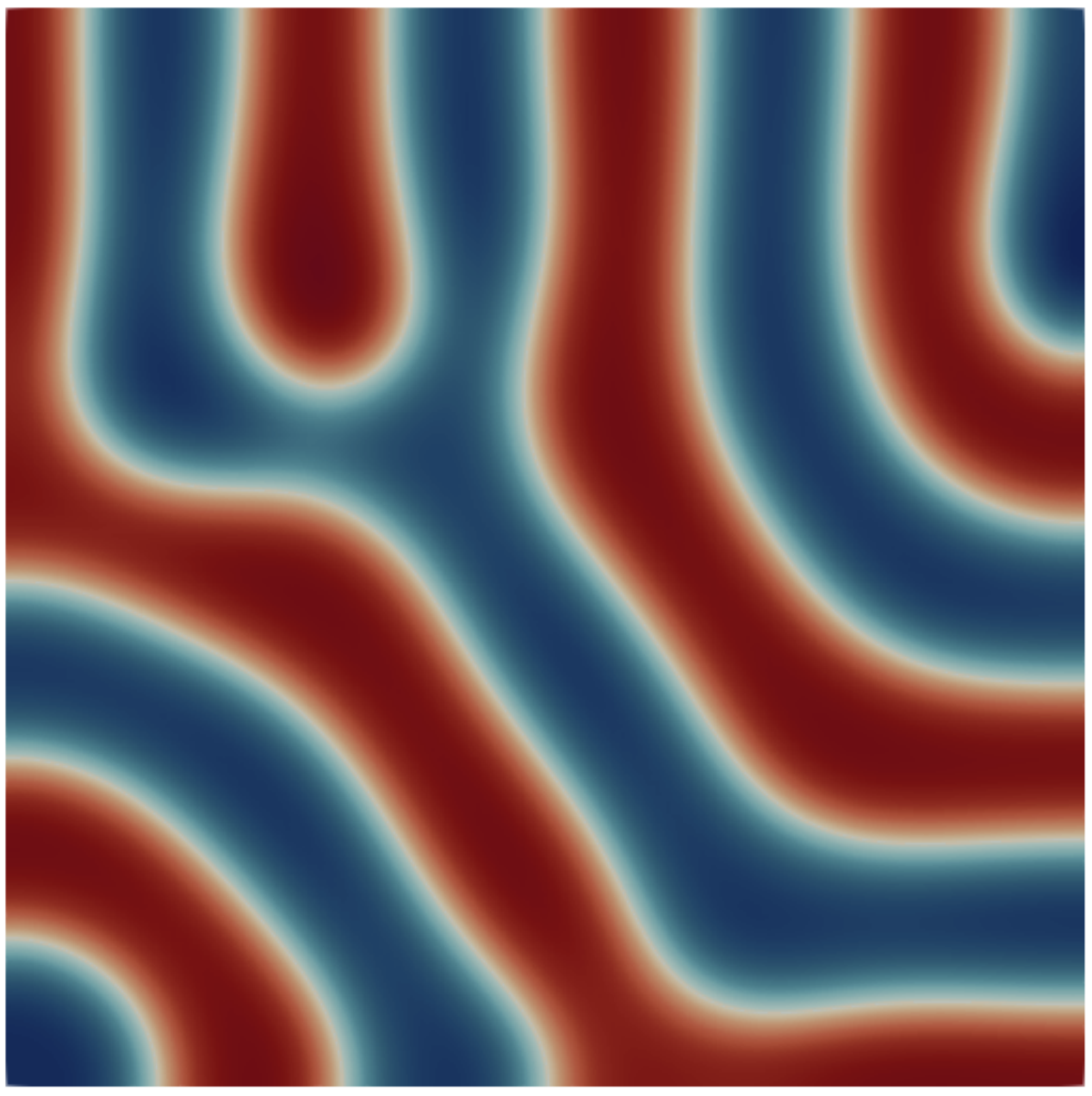}
\end{subfigure}
% \begin{subfigure}{0.15\textwidth}
%     \centering
%     \includegraphics[width=0.8\textwidth]{figures_pdf/equal_colorbar.pdf}
% \end{subfigure}
% \begin{subfigure}{0.32\textwidth}
%     \centering
%     \includegraphics[width=\textwidth]{figures_pdf/u1_equal.pdf}
% \end{subfigure}
% \begin{subfigure}{0.32\textwidth}
%     \centering
%     \includegraphics[width=\textwidth]{figures_pdf/u0_equal.pdf}
% \end{subfigure}
% % \begin{subfigure}{0.15\textwidth}
% %     \centering
% %     \includegraphics[width=0.8\textwidth]{figures_pdf/equal_colorbar.pdf}
% % \end{subfigure}
% \begin{subfigure}{0.32\textwidth}
%     \centering
%     \includegraphics[width=\textwidth]{figures_pdf/u1_equal.pdf}
% \end{subfigure}
\caption{Top: Initial guesses sampled from $\mathcal{D}$ with $(\delta_G, \gamma_G, s) = (0.8,0.02,1)$. Bottom: Solutions computed using the respective initial guesses by the energy-stable Newton algorithm for $\varepsilon = 0.08$, $\sigma = 12.8$, $m=0$, and $f\equiv0$ on a domain of $\lx = \ly = 3$.}

%\caption{Example equilibrium states from the energy stable Newton algorithm for $\varepsilon = 0.08$, $\sigma = 12.8$, $m=0$, $f\equiv0$ on a domain of $\lx = \ly = \pi$. The solutions are computed using two different initial guesses sampled from $\xi \sim \mathcal{U}(-0.5, 0.5)$ and interpolated into the finite element space.}
\label{fig:nonunique}
\end{figure}

% \dc{Updated this section.} 
One important property of the equilibrium states is that uniqueness of the solution is not guaranteed for any particular substrate field $f$. As with the time-dependent Cahn--Hilliard and SCFT models, our energy minimization algorithm may arrive at different equilibrium states corresponding to different basins of the free energy landscape when starting from different initial guesses $u^{(0)}$. 
For the time dependent Cahn--Hilliard equations, one typically uses a highly noisy random field with mean of $m$ as the initial condition, as is the case in \cite{QinKhairaSuEtAl13}. Thus, we introduce a distribution of initial guesses $\mathcal{D}$, defined, based on \cite{CaoGhattasOden22}, as a random perturbation from the homogeneous state $u(x) = m$ such that
\beq
    u^{(0)}(x) = m + s \cdot \erf(\xi(x)),
\eeq
where $s$ is a scalar governing the size of the random perturbation and $\xi \sim \mathcal{N}(0, \mathcal{C}_0)$ is a zero-mean Gaussian random field with covariance operator
$$ \mathcal{C}_0 = (\delta_G \mathcal{I} - \gamma_G \Delta)^{-2}.$$ 
Here, $\mathcal{I}$ is the identity operator, $\Delta$ is the Laplacian operator, and $\delta_G$ and $\gamma_G$ are scalars chosen to provide the desired variance and correlation length of the random field. 
As an example, we plot equilibrium state solutions arising from various initial guess samples in Figure \ref{fig:nonunique}.

% An example of this is seen for the case of a mixture comprising equal ratios of the two phases with a neutral substrate in Figure \ref{fig:nonunique}. We remark that this non-uniqueness is shared by the original time-dependent formulation of the nonlocal Cahn--Hilliard equations, where the role of the initial guess is played by the initial condition. The nonlocal Cahn--Hilliard equations produces different equilibrium solutions after time-evolution when starting from different initial conditions \dc{Is there a good citation?}.

% \todo{New sample and solutions picture}

% Though random initial conditions/guesses are adopted in \cite{QinKhairaSuEtAl13} and \cite{HannonDingBaiEtAl14}, sensitivity to the initial guesses are not explicitly considered or discussed.
However, as we have noted, the sensitivity to initial conditions or initial guesses is not typically discussed in existing studies on inverse design.
We postpone further discussion on the initial guesses to section \ref{sec:differentiability}, where we elaborate on its implications for our optimal design problem. We then present, in section \ref{sec:initial_guess_opt}, our approach to selecting the initial guess in the context of the optimal design problem.

\section{Optimal design of chemical guideposts}\label{sec:design_problem}
\subsection{Formulation of the PDE-constrained optimization problem}
We consider the following optimal design problem: Given predefined materials with fixed $(\varepsilon, \sigma, m)$ and some target morphology $u_d$, we seek an optimal substrate pattern that produces equilibrium states $u$ best matching the target $u_d$. In particular, we limit ourselves to the case where substrate pattern consists of either neutral or A-attracting chemicals, i.e. $f \leq 0$. Furthermore, we wish to obtain optimal designs which are manufacturable and have sparse features to attain manufacturing efficiency.

To formulate this as a PDE-constrained optimization problem, we first let $z \in \mathcal{Z}_{\text{ad}}$ be the design variable representing the guidepost configuration, where $\mathcal{Z}_{\text{ad}} \subset Z$ encapsulates the space of admissible designs in the design space $\mathcal{Z}$. The substrate interaction field $f$ is then a function of the design variable $z$. Thus, optimal design problem involves finding $z^* \in \mathcal{Z}_{\text{ad}}$ by solving
\beq
    z^* := \argmin_{z \in \mathcal{Z}_{\text{ad}}} \mathcal{J}(u,z; u_d) = Q(u; u_d) + \mathcal{P}(z)
    \label{eq:optimal_design}
\eeq
subject to the PDE constraints
\bcases{\label{eq:pde_constraint}}
    (\nabla \mu, \nabla \tilde{u}) + (\sigma(u-m), \tilde{u}) = 0 & $\forall \tilde{u} \in H^1(\Omega)$,
    \\
    (\mu, \tilde{\mu}) - (W'(u), \tilde{\mu}) - (\varepsilon^2 \nabla u, \nabla \tilde{\mu}) - (f(z), \tilde{\mu}) = 0 & $\forall \tilde{\mu} \in H^1(\Omega)$.
    
\ecases
Here $\mathcal{J}(u,z; u_d)$ is a cost functional consisting of two terms; the first is the design objective $Q(u; u_d) = \|u - u_d\|_{L^2(\Omega)}^2$ and the second is the penalization $\mathcal{P}(z)$ (to be specified later) that promotes certain features of the substrate, e.g., sparsity, boundedness away from the boundary, etc. 
% We can take simply $L^2(\Omega)$ misfit as the objective; that is, $Q(u; u_d) = \|u - u_d\|_{L^2(\Omega)}^2$. 

To compute the gradient of the cost functional with respect to the design variable, we use a Lagrangian approach with the multipliers/adjoint variables $(p, \lambda) \in H^1(\Omega) \times H^1(\Omega)$, i.e.,
\begin{linenomath}
\begin{align}
    \mathcal{L}(u, \mu, p, \lambda, z) :=  Q(u; u_d) + \mathcal{P}(z) + (\nabla \mu, \nabla p) + (\sigma(u-m), p) \nonumber \\
    + (\mu, \lambda) - (W'(u), \lambda) - (\varepsilon^2 \nabla u, \nabla \lambda) - (f(z), \lambda).
\end{align}
\end{linenomath}
The first-order optimality conditions can be written in terms of the Lagrangian through the state, adjoint, and gradient equations. That is, the state variables $(u, \mu)\in H^1(\Omega) \times H^1(\Omega)$, the adjoint variables $(p, \lambda) \in H^1(\Omega) \times H^1(\Omega)$, and the optimal design $z \in \mathcal{Z}_{\text{ad}}$ satisfy the state equation
\bcases{\label{eq:fwd}}
        \langle \partial_{p} \mathcal{L}, \tilde{p} \rangle
        = (\nabla \mu, \nabla \tilde{p}) + (\sigma(u-m), \tilde{p}) 
        = 0 & $\forall \tilde{p} \in H^1(\Omega)$, \\
        \langle \partial_{\lambda} \mathcal{L}, \tilde{\lambda}\rangle
        = (\mu, \tilde{\lambda}) - (W'(u), \tilde{\lambda}) - (\varepsilon^2 \nabla u, \nabla \tilde{\lambda}) - (f, \tilde{\lambda}) 
        = 0 & $\forall \tilde{\lambda} \in H^1(\Omega)$,
\ecases
the adjoint equation,
\bcases{\label{eq:adj}}
        % Adj1
        \langle \partial_{u} \mathcal{L}, \tilde{u}\rangle 
        = (\sigma p, \tilde{u}) - (W''(u)\lambda, \tilde{u}) - (\varepsilon^2 \nabla \lambda, \nabla \tilde{u})  + (\partial_{u}Q, \tilde{u}) 
        = 0 & $\forall \tilde{u} \in H^1(\Omega)$,
        \\
        % Adj2
        \langle \partial_{\mu} \mathcal{L}, \tilde{\mu}\rangle 
        = (\nabla p, \nabla \tilde{\mu}) + (\lambda, \tilde{\mu}) 
        = 0 & $\forall \tilde{\mu} \in H^1(\Omega)$,
\ecases
and the gradient equation
\beq\label{eq:grad-new}
        \langle \partial_{z} \mathcal{L}, \tilde{z} \rangle 
        = (\nabla_{z} \mathcal{P}, \tilde{z})  - (\lambda, \nabla_{z}f \tilde{z}) = 0 \quad \forall \tilde{z} \in \mathcal{Z}.
\eeq
% where $\mathcal{Z}$ is the test space for the design variable.

Newton-based optimization methods, which we pursue due to their superior convergence properties, require computation of the first and second derivatives of the cost with respect to\ $z$. It is possible to compute derivatives directly with respect to the design variable $z$ based on this Lagrangian formalism. However, we instead first compute the derivative with respect to the substrate interaction field $f$, and then proceed by chain rule to determine the $z$-derivatives. Though the two approaches are equivalent, the latter approach allows us to more conveniently write the derivatives for different parameterizations of $f(z)$, and also allows for more modular computer implementations.

In our implementation, we make use of the symbolic differentiation capabilities of FEniCS to obtain the adjoint and incremental equations. On the other hand, we cannot use FEniCS to directly compute variational forms of the derivatives with respect to $z$ when $z$ is itself not a function, but instead components of $z$ represent parameters defining the function $f$. In these cases, $f$ is defined through FEniCS's Expression objects, and currently FEniCS does not include the capability to differentiate through such objects to produce variational forms. This issue is circumvented by the approach we adopt here for computing derivatives.

\subsection{Derivatives with respect to the substrate interaction field}
Since the design variable is used to define the substrate interaction field $f$, we first aim to compute the derivative of the design objective $Q$ with respect to $f$. To do so, we define a Lagrangian using $\mathcal{P} = 0$. We then take the same definition of the adjoint variables $(p, \lambda) \in H^1(\Omega) \times H^1(\Omega)$ and write
\begin{linenomath}
\begin{align}
    \mathcal{L}(u, \mu, p, \lambda, f) :=  Q(u; u_d) + (\nabla \mu, \nabla p) + (\sigma(u-m), p) \nonumber \\
    + (\mu, \lambda) - (W'(u), \lambda) - (\varepsilon^2 \nabla u, \nabla \lambda) - (f, \lambda).
\end{align}
\end{linenomath}
The derivative then comes from solving the state equation (\ref{eq:fwd}), the adjoint equation (\ref{eq:adj}), and the new gradient equation
\beq\label{eq:grad}
    \langle D_{f} Q, \tilde{f} \rangle
    = \langle \partial_{f} \mathcal{L}, \tilde{f} \rangle 
    = - (\lambda, \tilde{f}) 
    \quad \forall \tilde{f} \in H^1(\Omega).
\eeq
% Given a solution $u$ satisfying the equation (\ref{eq:fwd}), we can also compute the corresponding gradient $\langle D_H\mathcal{J}, \tilde{z} \rangle$ by solving the adjoint equation for $(p, \lambda)$. We then have the gradient 
% \beq
%     \langle D_{f}\mathcal{J}, \tilde{f} \rangle = \langle \partial_{f} \mathcal{L}, \tilde{f} \rangle  = (\partial_{f} \mathcal{P}, \tilde{f})  - (\lambda, \tilde{f})
% \eeq

% As multiple forward solutions may exist for the same $f$, the gradient $D_f Q$ is local to the forward solution $(u, \mu)$ used to compute $D_f Q$, assuming that this gradient exists. That is, we assume that for a solution satisfying the forward problem $(u_0, \mu_0)$ corresponding to $f_0$, the solutions to the forward problem in some neighborhood of $f_0$ and $(u_0, \mu_0)$ can be represented by a mapping $u=u(f), \mu=\mu(f)$ that is differentiable in $f$, with $u_0 = u(f_0)$ and $\mu_0 = \mu(f)$. This leads to a functional $Q = Q(u(f))$ with the derivative $D_f Q$. 

Furthermore, we derive an expression for the Hessian acting on a direction $\hat{f}$ by defining a (Hessian) Lagrangian that captures the first order (state, adjoint, and gradient) equations 
\begin{linenomath}
\begin{align}
    \mathcal{L}^{H}(u, \mu, p, \lambda, f, \hat{u}, \hat{\mu}, \hat{p}, \hat{\lambda}, \hat{f})
    := \langle \partial_{p} \mathcal{L}, \hat{p} \rangle
    + \langle \partial_{\lambda} \mathcal{L}, \hat{\lambda} \rangle    
    + \langle \partial_{u} \mathcal{L}, \hat{u} \rangle
    + \langle \partial_{\mu} \mathcal{L}, \hat{\mu} \rangle 
    + \langle \partial_{f} \mathcal{L}, \hat{f} \rangle.
\end{align}
\end{linenomath}
By variation with respect to $(p, \lambda)$, we obtain $(\hat{u}, \hat{\mu}) \in H^1(\Omega) \times H^1(\Omega)$ from the incremental state equation
\bcases{\label{eq:incr_fwd}}
        \langle \partial_{p} \mathcal{L}^H, \tilde{p} \rangle 
        = (\nabla \hat{\mu}, \nabla \tilde{p}) + (\sigma \hat{u}, \tilde{p}) 
        = 0 & $\forall \tilde{p} \in H^1(\Omega)$, \\
        \langle \partial_{\lambda} \mathcal{L}^H, \tilde{\lambda} \rangle 
        = (\hat{\mu}, \tilde{\lambda}) - (W''(u)\hat{u}, \tilde{\lambda}) - (\varepsilon^2 \nabla \hat{u}, \nabla \tilde{\lambda}) - (\hat{f}, \tilde{\lambda}) 
        = 0 & $\forall \tilde{\lambda} \in H^1(\Omega)$.
\ecases
% \beq\label{eq:incr_fwd}
%     \begin{cases}
%         \langle \partial_{p} \mathcal{L}^H, \tilde{p} \rangle 
%         = (\nabla \hat{\mu}, \nabla \tilde{p}) + (\sigma \hat{u}, \tilde{p}) 
%         = 0 & \forall \tilde{p} \in H^1(\Omega) \\
%         \langle \partial_{\lambda} \mathcal{L}^H, \tilde{\lambda} \rangle 
%         = (\hat{\mu}, \tilde{\lambda}) - (W''(u)\hat{u}, \tilde{\lambda}) - (\varepsilon^2 \nabla \hat{u}, \nabla \tilde{\lambda}) - (\hat{f}, \tilde{\lambda}) 
%         = 0 & \forall \tilde{\lambda} \in H^1(\Omega)
%     \end{cases}
% \end

By variation with respect to $(u, \mu)$, we obtain $(\hat{p}, \hat{\lambda}) \in H^1(\Omega) \times H^1(\Omega)$ from the incremental adjoint equation
% Find $(\hat{p}, \hat{\lambda}) \in H^1(\Omega) \times H^1(\Omega)$ such that 
\bcases{\label{eq:incradj} \kern-2.5em}
        % Adj1
        \hspace{-0.25em} \langle \partial_{u} \mathcal{L^H}, \tilde{u} \rangle 
        = (\sigma \hat{p}, \tilde{u}) - (W''(u)\hat{\lambda}, \tilde{u}) - (\varepsilon^2 \nabla \hat{\lambda}, \nabla \tilde{u}) 
        + (\partial_{u}^2Q \hat{u} - W'''(u)\lambda \hat{u}, \tilde{u}) 
        = 0 %\nonumber 
        & $\hspace{-1.0em} \forall \tilde{u} \in H^1(\Omega)$,  \\
        % Adj2
        \hspace{-0.25em} \langle \partial_{\mu} \mathcal{L^H}, \tilde{\mu} \rangle 
        = (\nabla \hat{p}, \nabla \tilde{\mu}) + (\hat{\lambda}, \tilde{\mu}) 
        = 0  %\nonumber
        & $\hspace{-1.0em} \forall \tilde{\mu} \in H^1(\Omega)$.
\ecases
We can then compute the Hessian action as 
\beq
    \langle D_{f}^2Q \hat{f}, \tilde{f} \rangle 
    = \langle \partial_{f}\mathcal{L}^H, \tilde{f} \rangle 
    = - (\hat{\lambda}, \tilde{f}).
\eeq

The functional derivatives $D_f Q$ and $D_f^2Q$ can then be used as intermediate variables to compute the derivatives of $Q$ with respect to the design variable $z$ that parametrizes $f = f(z)$.

\subsection{Differentiability with respect to the substrate interaction field}\label{sec:differentiability}
Here we provide remarks on the differentiability of the state with respect to the substrate interaction field. Given $f^*$, we consider an equilibrium state $u^* \in \hcirc$ satisfying the first and second order optimality conditions of the energy minimization,
\beq
    \duality{D_u\mathcal{F}(u^*)}{\tilde{u}_0} = 0 \quad \forall \tilde{u}_0 \in \hcirc,
\eeq
and 
\beq
    \duality{D_u^2\mathcal{F}(u^*)\hat{u}_0}{\hat{u}_0} \geq a \|\hat{u}_0\|^2_{\hcirc} \quad \forall \hat{u}_0 \in \hcirc.
\eeq
The implicit function theorem (see, for example, Theorem 1.41, \cite{HinzePinnauUlbrichEtAl09}) guarantees the existence of a branch of equilibrium states $u(f)$ defined within a neighborhood of $(u^*,f^*)$, with $u^* = u(f^*)$, that is Fr\'{e}chet-differentiable with respect to $f$. This follows from the fact that the operator $D_u^2\mathcal{F}(u) : \hcirc \rightarrow \hcirc'$ is invertible, as it generates a continuous and coercive bilinear form when $u$ satisfies the second order optimality conditions. A proof of the continuity of bilinear form generated by $D_u^2 \mathcal{F}(u)$ is presented in \ref{appendix:continuity}.

However, we recall that the state problem is non-unique and multiple equilibrium states can exist for the same $f$. Therefore, the information given by the derivative $D_f u(f)$ is limited to the solution branch corresponding to the particular state $u$ at which it is computed. In particular, we have numerically observed that a small change in $f$ can result in large changes to the solution $u$, even if the same initial guess is used in both cases. For convenience, we introduce the notion of the solution operator $\mathcal{S}: L^2(\Omega) \times H^1(\Omega) \rightarrow H^1(\Omega)$, such that $u = \mathcal{S}(f, u^{(0)})$ is the solution for the state equation given the substrate interaction field variable $f \in L^2(\Omega)$ and an initial guess $u^{(0)} \in H^1(\Omega)$, using the energy stable Newton algorithm defined above.

\begin{figure}
\centering
    \begin{subfigure}{0.15\textwidth}
    \centering
    \includegraphics[width=0.95\textwidth]{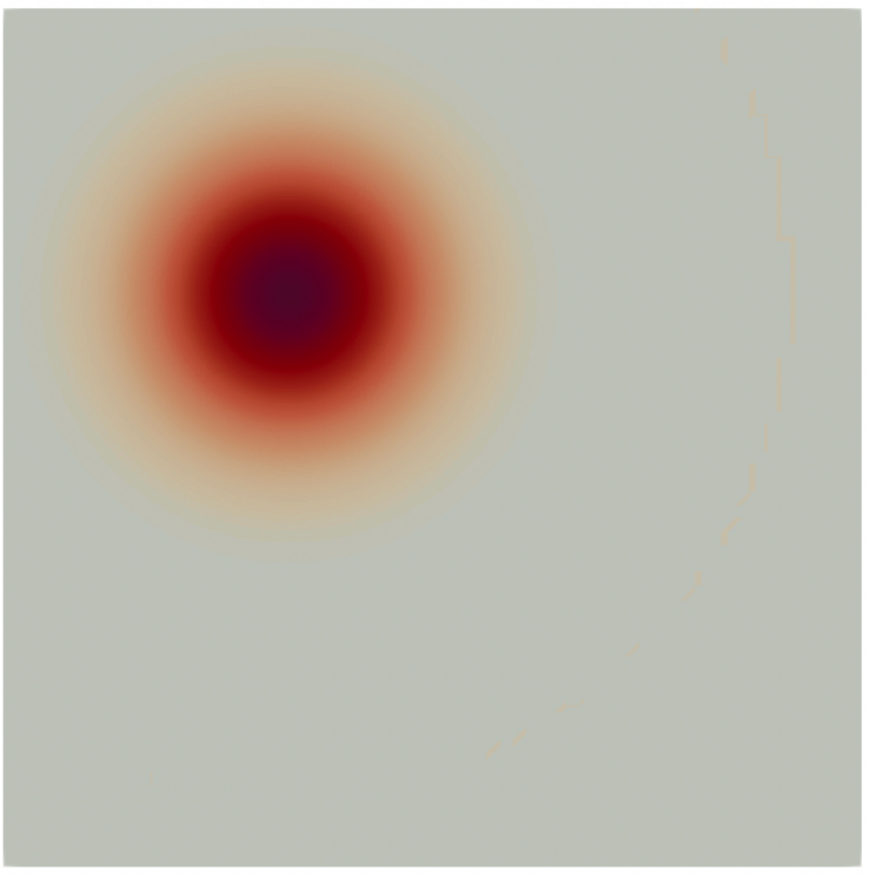}
    \caption{$f_1$}
    \end{subfigure}
    \begin{subfigure}{0.15\textwidth}
    \centering
    \includegraphics[width=0.95\textwidth]{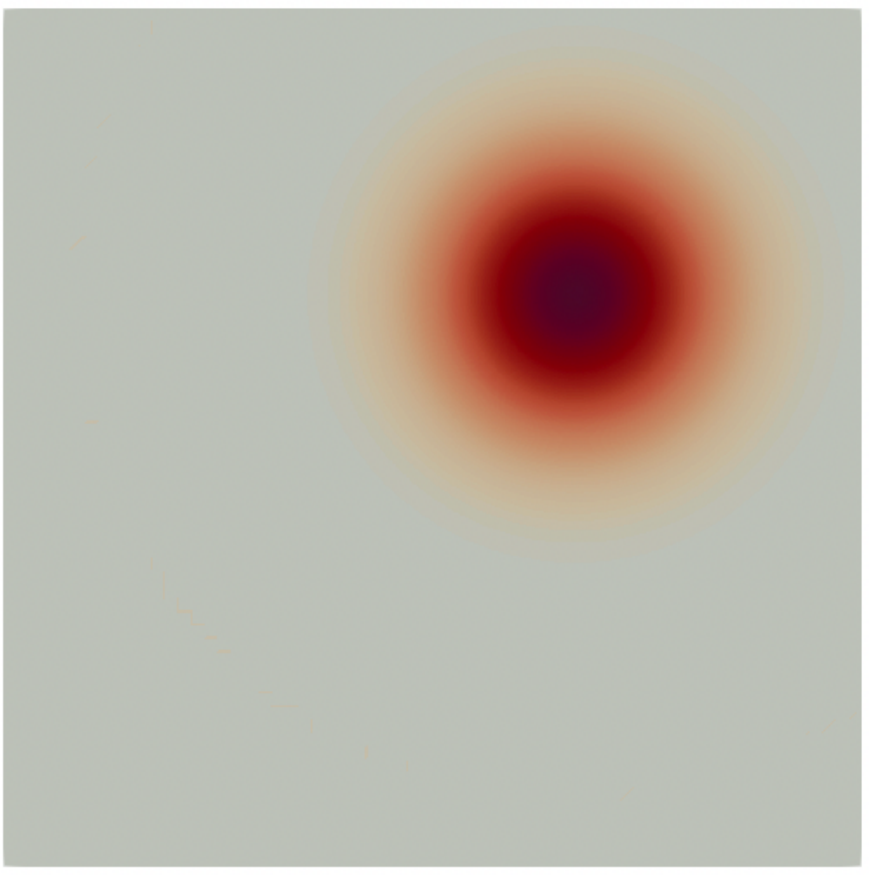}
    \caption{$\hat{f}$}
    \end{subfigure}
    \begin{subfigure}{0.15\textwidth}
    \centering
    \includegraphics[width=0.95\textwidth]{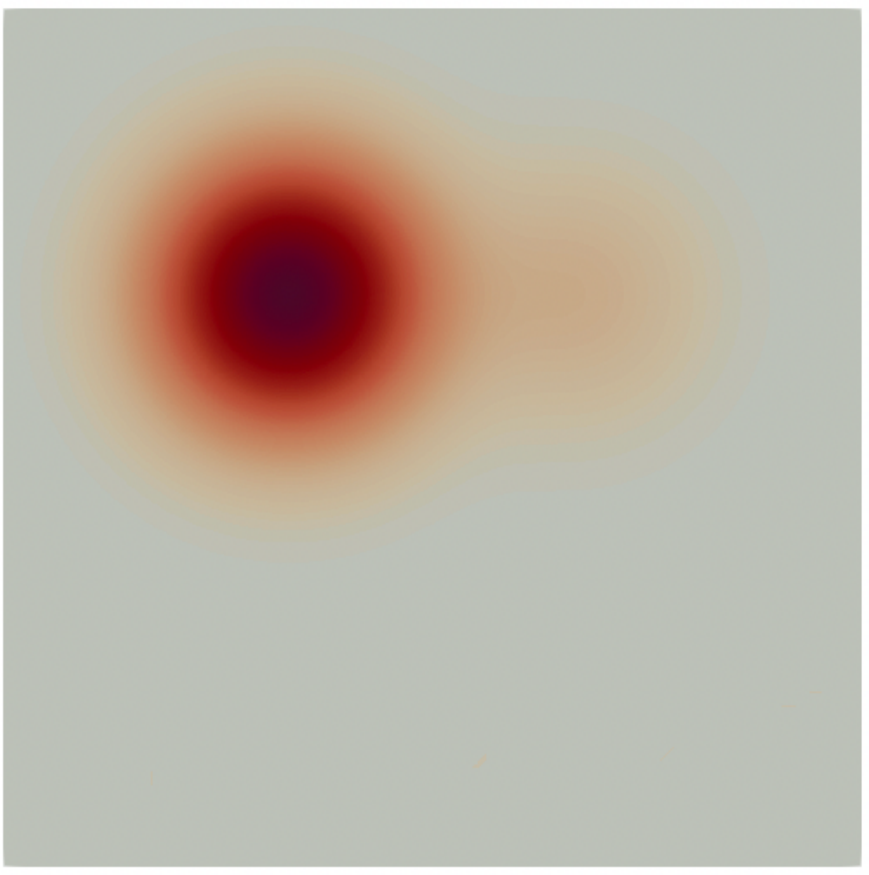}
    \caption{$f_2 = f_1 + 0.1 \hat{f}$}
    \end{subfigure}
    \begin{subfigure}{0.05\textwidth}
    \centering
    \includegraphics[width=0.95\textwidth]{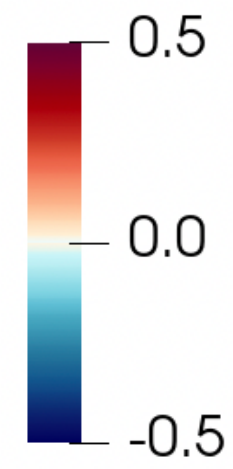}
    \end{subfigure}

    \begin{subfigure}{0.15\textwidth}
    \centering
    \includegraphics[width=0.95\textwidth]{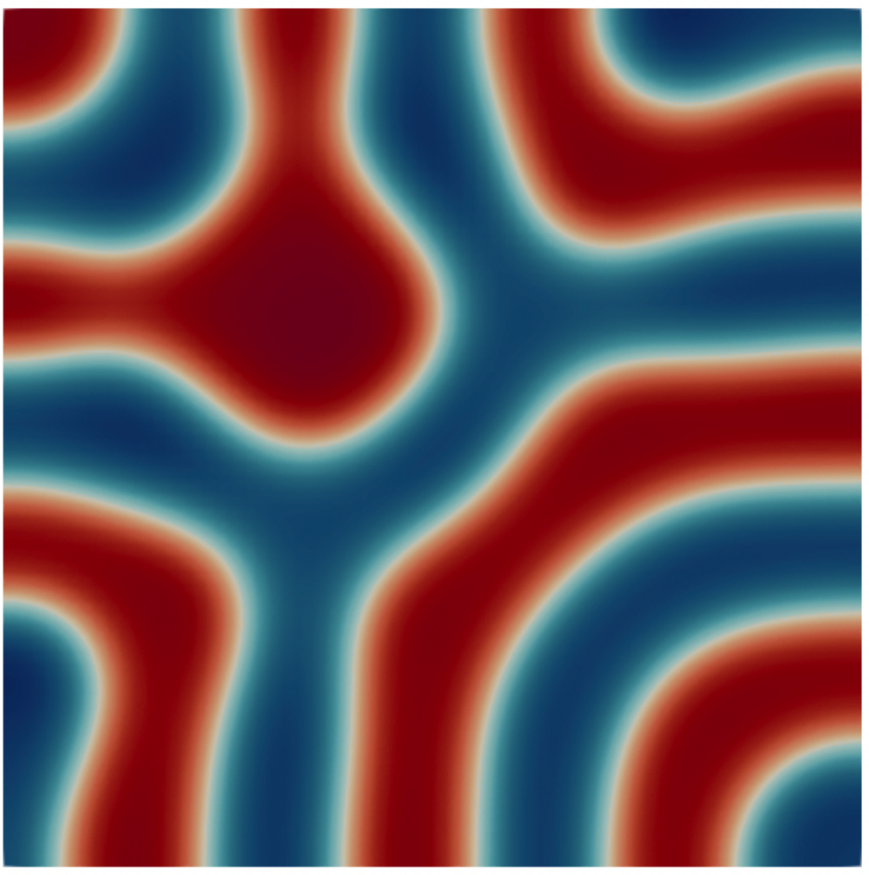}
    \caption{$u_1 = \mathcal{S}(f_1, 0)$}
    \end{subfigure}
    \begin{subfigure}{0.15\textwidth}
    \centering
    \includegraphics[width=0.95\textwidth]{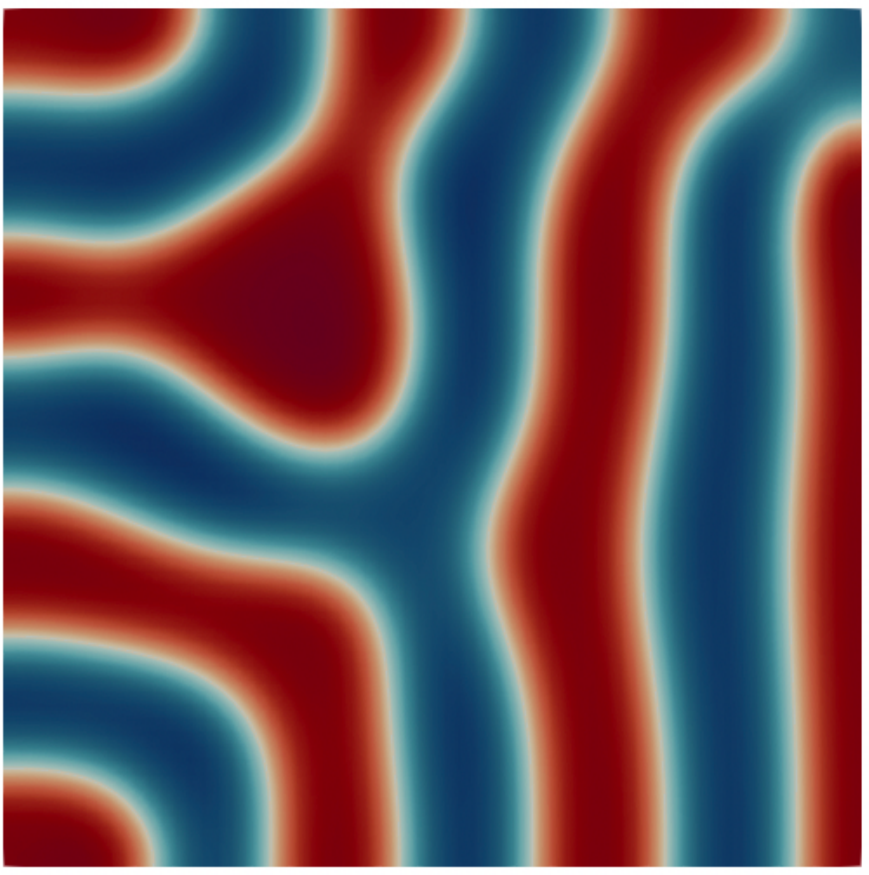}
    \caption{$u_2 = \mathcal{S}(f_2, 0)$}    
    \end{subfigure}
    \begin{subfigure}{0.15\textwidth}
    \centering
    \includegraphics[width=0.95\textwidth]{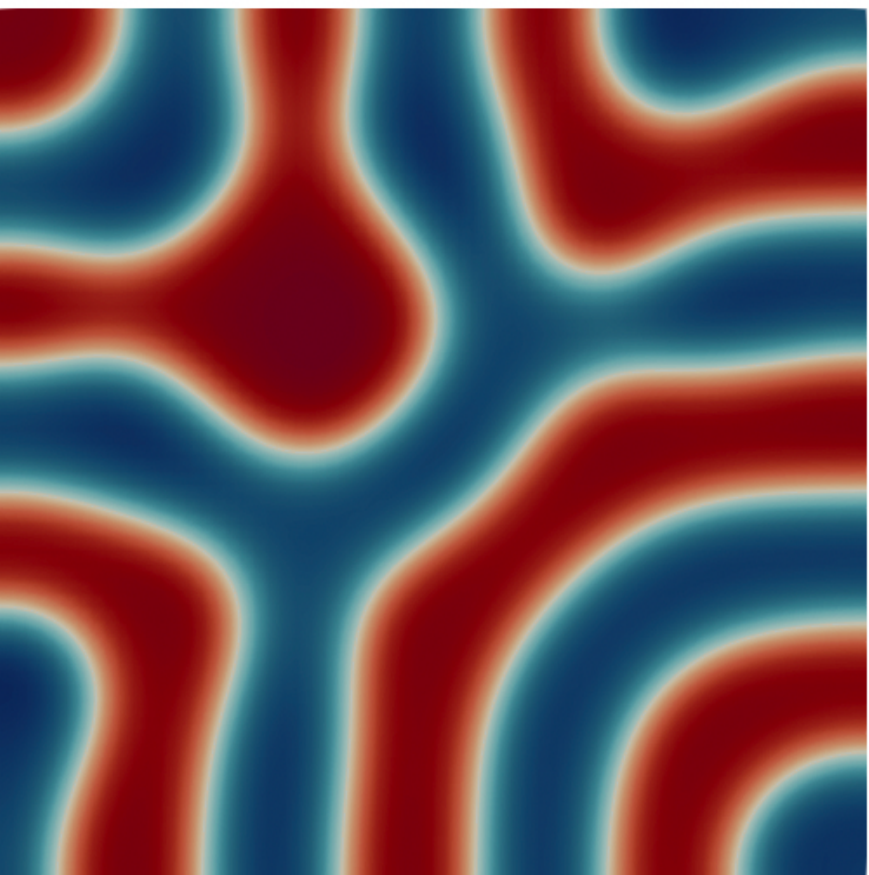}
    \caption{$\tilde{u}_2 = \mathcal{S}(f_2, u_1)$}
    \end{subfigure}
    \begin{subfigure}{0.05\textwidth}
    \centering
    \includegraphics[width=0.95\textwidth]{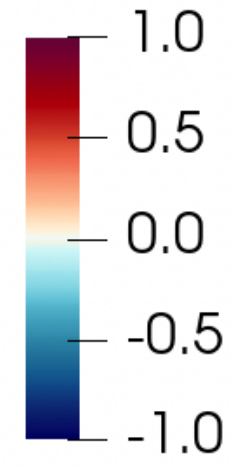}
    \end{subfigure}

    \caption{Numerical example illustrating the solutions for two slightly different substrate interaction fields, with the model parameters $(m, \varepsilon, \sigma) = (0, 0.08, 12.8)$. Top row: first substrate $f_1$ (left), perturbation direction $\hat{f}$ (middle), and the perturbed substrate $f_2 = f_1 + \delta \hat{f}$ with $\delta = 0.1$ (right). Bottom row: solution $u_1$ for first substrate using the initial guess $u^{(0)} = 0$ (left), solution $u_2$ for perturbed substrate using the same initial guess $u^{(0)} = 0$ (middle), and solution $\tilde{u}_2$ for perturbed substrate using $u_1$ as the initial guess (right).}
    \label{fig:differentiability}
\end{figure}

We demonstrate our observation by a numerical example on the 2D domain $\Omega = (0,3)^2$.  In our example, we consider a substrate interaction field $f_1$ and a perturbed substrate $f_2 = f_1 + \theta \hat{f}$, where $\hat{f}$ is the direction of the perturbation and $\theta \in (0,1)$ is small. We compute solutions $u_1 = \mathcal{S}(f_1, u^{(0)})$ and $u_2 = \mathcal{S}(f_2, u^{(0)})$, using the same initial guess of $u^{(0)} \equiv m = 0$ for both substrates. We compare this with an alternative procedure, where given the first solution $u_1 = \mathcal{S}(f_1, u^{(0)})$, we compute the solution $\tilde{u}_2 = \mathcal{S}(f_2, u_1)$ for the perturbed substrate using the first solution as the initial guess. The substrate interaction fields and resulting states are shown in Figure \ref{fig:differentiability}.

In the figure, we see that for a small perturbation of $\delta = 0.1$, the solutions $u_1$ and $u_2$ found using the same initial guess $u^{(0)} = 0$ differ significantly in the regions where the substrate is weak ($f$ is small). This suggests that when starting with an initial guess far from a particular solution branch, the energy-stable Newton algorithm can converge to other solution branches even for substrate fields $f_1$ and $f_2$ which are similar. When this happens, the derivative $D_f u(f_1)$ computed at $u_1$ does not reflect the behavior of $u_2$ corresponding to the perturbed substrate $f_2$, as it lies in a separate solution branch. On the other hand, when using $u_1$ as the initial guess for the perturbed substrate, the solution $\tilde{u}_2$ is only slightly perturbed and appears to remain on the same solution branch as $u_1$ where the derivative $D_f u(f_1)$ is meaningful.

\subsection{Guidepost based design}
\begin{figure}
\centering
    \begin{subfigure}{0.4\textwidth}
    \centering
    \includegraphics[trim={2cm 1cm 1cm 1cm}, width=\textwidth]{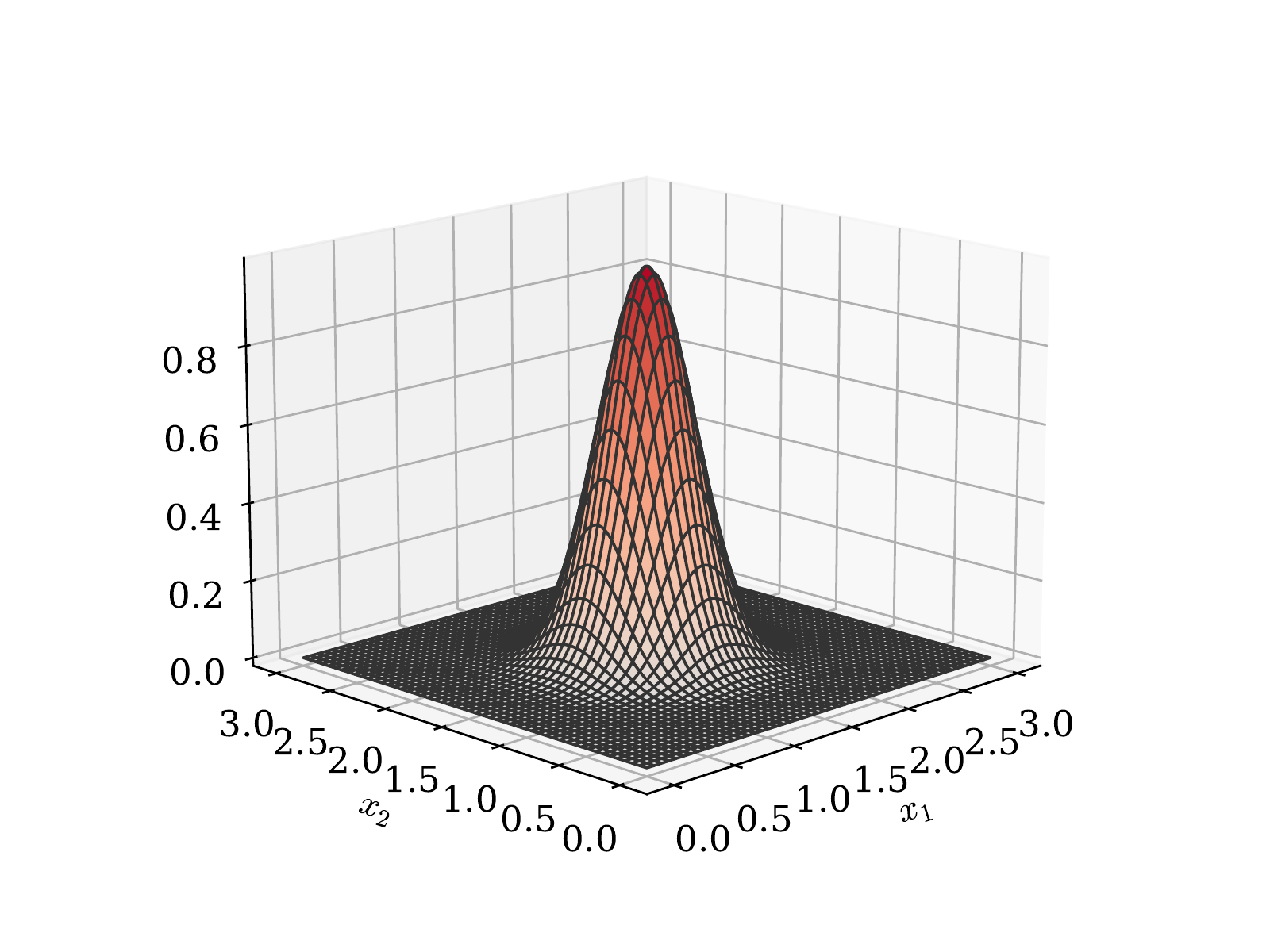}
    \end{subfigure}
    \begin{subfigure}{0.4\textwidth}
    \centering
    \includegraphics[trim={2cm 1cm 1cm 1cm}, width=\textwidth]{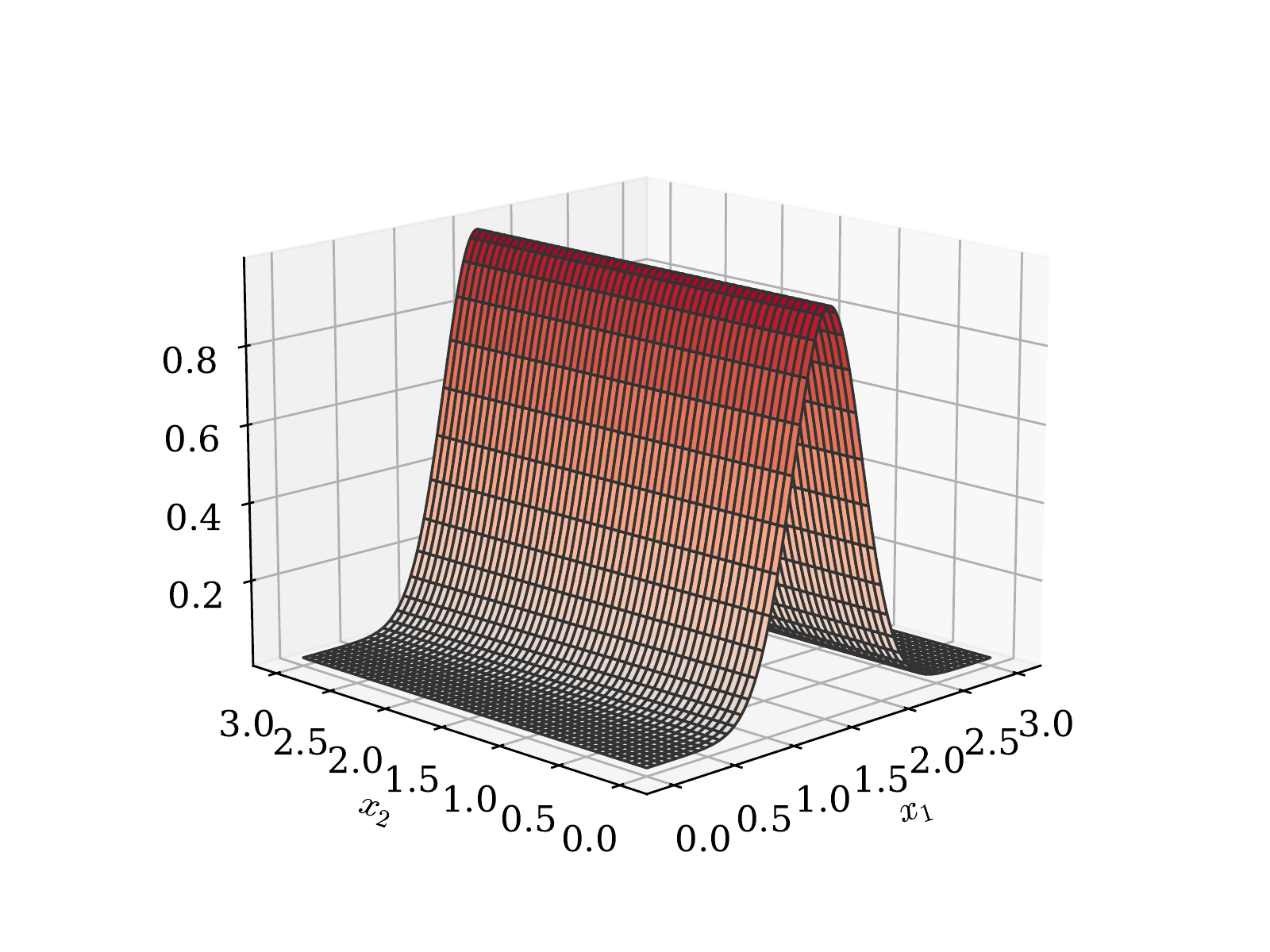}
    \end{subfigure}
    % \vspace{-20pt}
    \caption{Guidepost shapes on a square substrate domain: circular (left) and strip (right).}
    \label{fig:guideposts}
\end{figure}

We formulate the design problem based on guideposts, as is commonly done in the literature \cite{QinKhairaSuEtAl13,KhairaQinGarnerEtAl14,HannonDingBaiEtAl14}. In particular, we follow the approach of \cite{QinKhairaSuEtAl13}, in which the substrate interaction field is given by a linear combination of shape functions representing the guideposts. As mentioned previously, we define guideposts to be units of A-attracting regions ($f \leq 0$) prescribed in predefined shapes, such as circular spots or lengthwise strips, as illustrated in Figure \ref{fig:guideposts}. 

The design problem is then to optimize the locations of these guideposts. To do so, we first specify a fixed number of guideposts at given strengths, and take their locations within the substrate domain as the design variable to optimize. Two penalty terms are used to prevent the collapse of the guideposts and their departure from the domain. We note that this formulation ensures that the substrate pattern consists only of regular, manufacturable structures.

% \begin{figure}
% \centering
%     \begin{subfigure}{0.42\textwidth}
%     \centering
%     \includegraphics[width=\textwidth]{figures_pdf/spot.pdf}
%     \end{subfigure}
%     \begin{subfigure}{0.12\textwidth}
%     \centering
%     \includegraphics[width=0.8\textwidth]{figures_pdf/substrate_colorbar.pdf}
%     \end{subfigure}
%     \begin{subfigure}{0.42\textwidth}
%     \centering
%     \includegraphics[width=\textwidth]{figures_pdf/strip.pdf}
%     \end{subfigure}
%     % \vspace{-20pt}
%     \caption{Guidepost shapes on a square substrate domain: circular (left) and strip (right).}
%     \label{fig:guideposts}
% \end{figure}

\subsection{Circular guideposts}
\subsubsection{Definition and computation of derivatives}
We denote by $\mathbf{x} = (x_1, x_2)$ the coordinates on the two-dimensional substrate surface $\Omega_s$. Guideposts are defined based on a list of locations $\{\mathbf{r}_i\}_{i=1}^{N_p}$, where $\mathbf{r}_i = (r_{1,i}, r_{2,i}) \in \Omega_s \subset \mathbb{R}^2$ and $N_p$ is the number of guideposts.
% representing their locations. 
We define circular guideposts at a point $\mathbf{r}_i = (r_{1,i}, r_{2,i})$ by a shape function $\phi(\|\mathbf{x} - \mathbf{r}_i\|_2^2)$. In this work, we take $\phi(t) = w\exp(-t/2b^2)$, where $w$ controls the strength of the substrate chemical and $b$ controls the spread of the spot. The substrate interaction field is therefore
\beq
    f(x) = -\sum_{i=1}^{N_p} \tau(x_3) \phi(\|\mathbf{x} - \mathbf{r}_i\|_2^2),
\eeq
where the negative sign is introduced to make it A-attracting. Recall that the substrate design is defined only on the bottom boundary, and its decaying effect is carried in the $x_3$ direction by $\tau(x_3)$. If a thin film approximation is taken, we take $d=2$ and $\tau = 1$. 

The optimization variable is taken as $z = \{\mathbf{r}_i\}_{i=1}^{N_p}$, or
% now $z = \{\mathbf{r}_i\}_{i=1}^{N_p}$, alternatively expressed 
$z = (r_{1,i}, r_{2,i},...,r_{1,N_p}, r_{2,N_p})$ for convenience. We can compute the derivatives of $Q$ with respect to $z$ based on the functional gradient of $Q$ with respect to $f$ 
% and the chain rule with $\partial_z f$, 
i.e., 
\beq
    \frac{\partial Q}{\partial \mathbf{r}_i} 
        = \inner{D_f Q}{\partial_{\mathbf{r}_i}f} 
        = - \inner{D_f Q}{\tau \frac{\partial \phi(\|\mathbf{x} - \mathbf{r}_i\|_2^2)}{\partial \mathbf{r}_i}}.
\eeq
Note that
$$\frac{\partial Q}{\partial \mathbf{r}_i} = 
\left(\frac{\partial Q}{\partial r_{1,i}},
\frac{\partial Q}{\partial r_{2,i}}\right)^T
=
\left(
\frac{\partial Q}{\partial z_{2i-1}},
\frac{\partial Q}{\partial z_{2i}}
\right)^T.
$$
In the case of the Hessian action, we consider a direction $\hat{z} = \{\hat{\mathbf{r}}_j\}_{j=1}^{N_p}$ and compute the Hessian action
% \begin{linenomath}
% \begin{align*}
%     [\nabla^2 Q \hat{z}]_i &= \sum_{j=1}^{N_p} \frac{\partial^2 Q}{\partial \mathbf{r}_i\partial \mathbf{r}_j}\hat{\mathbf{r}}_j \\
%     &= \frac{\partial}{\partial \mathbf{r}_i} \left( \sum_{j=1}^{N_p} \frac{\partial Q}{\partial \mathbf{r}_j} \hat{\mathbf{r}}_j \right) \\
%     % &= \inner{D_f^2 Q \partial_{\mathbf{r}_i}f}{\sum_{j=1}^{N_p} \frac{\partial Q}{\partial \mathbf{r}_j} \hat{\mathbf{r}}_j} 
%     &= \inner{D_f^2 Q \partial_{\mathbf{r}_i}f}{\tau\frac{\partial \phi(\|\mathbf{x} - \mathbf{r}_i\|_2^2)}{\partial \mathbf{r}_i}} 
%             + \inner{D_f Q}{\sum_{j=1}^{N_p} \tau \frac{\partial^2 \phi(\|\mathbf{x}-\mathbf{r}_j\|_2^2)}{\partial \mathbf{r}_i \partial \mathbf{r}_j} \hat{\mathbf{r}}_j} \\
%     &= \inner{D_f^2 Q\left(\sum_{j=1}^{N_p} \tau \frac{\partial \phi(\|\mathbf{x} - \mathbf{r}_j\|_2^2)}{\partial \mathbf{r}_j} \cdot \hat{\mathbf{r}}_j \right)}{\tau\frac{\partial \phi(\|\mathbf{x} - \mathbf{r}_i\|_2^2)}{\partial \mathbf{r}_i}} 
%         \\
%         & + \inner{D_f Q}{\tau\frac{\partial^2 \phi(\|\mathbf{x}-\mathbf{r}_i\|_2^2)}{\partial \mathbf{r}_i^2} \hat{\mathbf{r}}_i}
% \end{align*}
% \end{linenomath}
\beq
    [\nabla^2 Q \hat{z}]_i =
    \inner{D_f^2 Q\left(\sum_{j=1}^{N_p} \tau \frac{\partial \phi(\|\mathbf{x} - \mathbf{r}_j\|_2^2)}{\partial \mathbf{r}_j} \cdot \hat{\mathbf{r}}_j \right)}{\tau\frac{\partial \phi(\|\mathbf{x} - \mathbf{r}_i\|_2^2)}{\partial \mathbf{r}_i}} 
    - \inner{D_f Q}{\tau\frac{\partial^2 \phi(\|\mathbf{x}-\mathbf{r}_i\|_2^2)}{\partial \mathbf{r}_i^2} \hat{\mathbf{r}}_i}
\eeq
by the chain and product rules. We note that 
\beq
\frac{\partial \phi(\|\mathbf{x} - \mathbf{r}_i)\|_2^2)}{\partial \mathbf{r}_i} = 2 \phi'(\|\mathbf{x} - \mathbf{r}_i)\|_2^2)(\mathbf{r}_i - \mathbf{x}),
\eeq
and
\beq
\frac{\partial^2 \phi(\|\mathbf{x} - \mathbf{r}_i)\|_2^2)}{\partial \mathbf{r}_i^2} = 4\phi''(\|\mathbf{x} - \mathbf{r}_i)\|_2^2)(\mathbf{r}_i - \mathbf{x})(\mathbf{r}_i - \mathbf{x})^T + 2 \phi'(\|\mathbf{x} - \mathbf{r}_i)\|_2^2) \mathbf{I},
\eeq
where $\mathbf{I} \in \mathbb{R}^{2\times2}$ is the identity matrix. The gradient of the full cost functional,
\beq
    \nabla_z \mathcal{J} = \nabla_z Q + \nabla_z \mathcal{P},
\eeq
and Hessian action for a given direction $\hat{z}$,
\beq
    \nabla_z^2 \mathcal{J} \hat{z} = \nabla_z^2 Q \hat{z} + \nabla_z^2 \mathcal{P} \hat{z},
\eeq
are now the sums of the design objective and penalty components.
We remark that the computation of the Hessian action requires one functional gradient calculation and one application of the functional Hessian in the term
\begin{linenomath}
$$ D_f^2 Q\left(\sum_{j=1}^{N_p} \tau \frac{\partial \phi(\|\mathbf{x} - \mathbf{r}_j\|_2^2)}{\partial \mathbf{r}_j} \cdot \hat{\mathbf{r}}_j \right). $$
\end{linenomath}

\subsubsection{Penalization} 
We introduce a repelling penalization term to prevent guideposts from collapsing onto the same location during optimization, and to control the spacing of the guideposts at their optimal configurations. The repelling penalization term uses a potential based on the distances between pairs of guideposts, and takes the form 
\beq \label{eq:repel_spot}
    \mathcal{P}_{\text{repel}} =\alpha \sum_{i \neq j} \psi(\|\mathbf{r}_i - \mathbf{r}_j\|_2^2),
\eeq
where we take $\psi(t) = t^{-1}$. In general, a weighted norm $\| \cdot \|_W$ of the relative positions of guideposts can be use in place of the Euclidean norm $\| \cdot \|_2$ to introduce anisotropy of the repelling force.
% The weighted norm $\|\mathbf{x}\|_W = \|W^{1/2}\mathbf{x}\|$ and the weighting matrix $W$ can be used to introduce anisotropy of the repelling force to encourage separation along one particular direction over the other. 
% \dc{I have not made use of W in the numerical results. Should I simply remove it from the equations and maybe leave a comment that it is possible to introduce anisotropy of the repulsion through W?} 
However, keeping with the Euclidean norm, we can compute the gradient of this penalization term as
\beq
    \frac{\partial \mathcal{P}_{\text{repel}}}{\partial \mathbf{r}_i}
     = \alpha \sum_{j=1, j\neq i}^{N} 2 \psi'(\|\mathbf{r}_i - \mathbf{r}_j\|_2^2) (\mathbf{r}_i - \mathbf{r}j),
\eeq
and components of the Hessian matrix as  
\beq
    \frac{\partial^2\mathcal{P}_{\text{repel}}}{\partial{\mathbf{r}_i}\partial{\mathbf{r}_j}} =
    \begin{cases}
        \alpha \left(
        4 \psi''(\|\mathbf{r}_i - \mathbf{r}_j\|_2^2)(\mathbf{r}_i - \mathbf{r}_j)(\mathbf{r}_i - \mathbf{r}_j)^T
                - 2 \psi'(\|\mathbf{r}_i - \mathbf{r}_j\|_2^2) \mathbf{I} \right) 
                & i \neq j, \\[10pt]
        \alpha \sum_{k=1, k\neq i}^{N} 
        	\left( 
            4\psi''(\|\mathbf{r}_i - \mathbf{r}_k\|_2^2)(\mathbf{r}_i - \mathbf{r}_k)(\mathbf{r}_i - \mathbf{r}_k)^T 
            + 2 \psi'(\|\mathbf{r}_i - \mathbf{r}_k\|_2^2) \mathbf{I} 
            \right) 
            & i = j. \\        
    \end{cases}
\eeq
% where $\mathbf{I}$ refers to the $2 \times 2$ identity matrix.
% \beq
%     \frac{\partial^2\mathcal{P}_{\text{repel}}}{\partial{\mathbf{r}_i}\partial{\mathbf{r}_j}} =
%     \begin{cases}
%         4 \psi''(\|\mathbf{r}_i - \mathbf{r}_j\|_W^2)W(\mathbf{r}_i - \mathbf{r}_j)(W(\mathbf{r}_i - \mathbf{r}_j))^T
%                 - 2 \psi'(\|\mathbf{r}_i - \mathbf{r}_j\|_W^2) W^{1/2} & i \neq j, \\[10pt]
%         \sum_{k=1, k\neq i}^{N} 
%             4\psi''(\|\mathbf{r}_i - \mathbf{r}_k\|_W^2)W^{1/2}(\mathbf{r}_i - \mathbf{r}_k)(W^{1/2}(\mathbf{r}_i - \mathbf{r}_k))^T 
%             + 2 \psi'(\|\mathbf{r}_i - \mathbf{r}_k\|_W^2) W^{1/2} & i = j. \\        
%     \end{cases}
% \eeq
%
% from which we can compute the full gradient and Hessian actions.
We also need to ensure that the guideposts remain inside the substrate domain $\Omega_s$. This is not enforced strongly, but instead encouraged through another penalization term 
\beq
    \mathcal{P}_{\text{wall}}(z) = \sum_{i=1}^{N} \left( e^{-a_1 r_{1,i}} + e^{-a_1(\lx - r_{1,i})} 
    + e^{-a_2 r_{2,i}} + e^{-a_2(\ly - r_{2,i})} \right),
\eeq
where $(0, \lx) \times (0, \ly)$ defines the rectangular domain $\Omega_s$, while $a_1$ and $a_2$ control the steepness of the potential in the $x_1$ and $x_2$ directions respectively. 
% An example of this potential is shown in Figure \ref{fig:wall} for a square domain. 
We can compute the gradient and Hessian as 
\beq
    \frac{\partial \mathcal{P}_{\text{wall}}}{\partial \mathbf{r}_i} = 
    \begin{bmatrix}
        -a_1(e^{-a_1 r_{1,i}} - e^{-a_1(\lx - r_{1,i})})\\[10pt]
        -a_2(e^{-a_2 r_{2,i}} - e^{-a_2(\ly - r_{2,i})}) 
    \end{bmatrix},
\eeq
and 
\beq
    \frac{\partial^2 \mathcal{P}_{\text{wall}}}{\partial_{\mathbf{r}_i}\partial_{\mathbf{r}_j}} = \delta_{ij}
    \begin{bmatrix}
        a_1^2 (e^{-a_1 r_{1,i}} + e^{-a_1(\lx - r_{1,i})}) & 0\\[10pt]
        0 & a_2^2 (e^{-a_2 r_{2,i}} + e^{-a_2(\ly - r_{2,i})})
    \end{bmatrix}.
\eeq
% from which we can compute the full gradient and Hessian actions.

% \begin{figure}[htpb!]
%     \centering
%     \includegraphics[width=0.75\textwidth]{figures_pdf/wall_2d.pdf}
%     \caption{Penalization to prevent the guideposts from leaving the domain.
%     % enforce the boundedness of guidepost locations
%     }
%     \label{fig:wall}
% \end{figure}
 
\subsection{Strip guideposts}
Strip guideposts are patterns that are uniform in the $x_2$ direction (length) of the substrate domain. Thus, they are parameterized by only their location along $x_1$. A simple approach to represent the strip guidepost is to consider it as 
% the 1D case of 
the circular guidepost 
% formulation, reduced to 
restricted to the one-dimensional domain $\tilde{\Omega}_s := (0,\lx)$. 
% representing the width of the domain. 
The optimization variable is now $z = (r_1, r_2, ..., r_{N_p})$, describing strip locations $r_i \in \tilde{\Omega}_s$. This leads to the substrate interaction field
\beq
    f(x) = -\sum_{i=1}^{N_p} \tau(x_3) \phi(|x_1 - r_i|^2).
\eeq
Derivatives are also computed in the same manner. We have gradients
\beq
    \frac{\partial Q}{\partial r_i} 
        = \inner{D_f Q}{\partial_{r_i}f} 
        = -\inner{D_f Q}{\tau \frac{\partial \phi(|x_1 - r_i|^2)}{\partial r_i}},
\eeq
and Hessian actions in direction $\hat{z} = (\hat{r}_1, \hat{r}_2, ..., \hat{r}_{N_p})$ as 
\begin{linenomath}
\begin{align}
    [\nabla^2 Q \hat{z}]_i
    &= \inner{D_f^2 Q\left(\sum_{j=1}^{N_p} \tau 
    \frac{\partial \phi(|x_1 - r_j|^2)}{\partial r_j} \hat{r}_j \right)}
    {\tau\frac{\partial \phi(|x_1 - r_i|^2)}{\partial r_i}} 
    - \inner{D_f Q}{\tau\frac{\partial^2 \phi(|x_1-r_i|^2)}{\partial r_i^2} \hat{r}_i}.
\end{align}
\end{linenomath}
% Again, we have made use of functional Gradient and Hessian actions computed via the adjoint method. 
The derivatives of the shape functions are now
\beq
\frac{\partial \phi(|x_1 - r_i|^2)}{\partial r_i} = 2 \phi'(|x_1 - r_i|^2)(r_i - x_1),
\eeq
and
\beq
\frac{\partial^2 \phi(|x_1 - r_i|^2)}{\partial r_i^2} = 4\phi''(|x_1 - r_i|^2)(r_i - x_1)^2 + 2 \phi'(|x_1 - r_i|^2).
\eeq

For the repelling penalization function, we also reduce the dimension to obtain
\beq \label{eq:repel_strip}
    \mathcal{P}_{\text{repel}} =\alpha \sum_{i \neq j} \psi(|r_i - r_j|^2).
\eeq
%
% noting that the matrix $W$ is no longer used. 
This leads to the derivatives
\beq
     \frac{\partial \mathcal{P}_{\text{repel}}} {\partial r_i} = \alpha \sum_{j=1, j\neq i}^{N} 2 \psi'(|r_i - r_j|^2) (r_i - r_j),
\eeq
and 
\beq
    \frac{\partial^2\mathcal{P}_{\text{repel}}}{\partial{r_i}\partial{r_j}} =
    \begin{cases}
        \alpha \left(
        4 \psi''(|r_i - r_j|^2)(r_i - r_j)^2
                - 2 \psi'(|r_i - r_j|^2) 
        \right) & i \neq j, \\[10pt]
        \alpha 
        \sum_{k=1, k\neq i}^{N} \left(
            4\psi''(|r_i - r_k|^2)(r_i - r_k)^2
            + 2 \psi'(|r_i - r_k|^2) 
            \right) & i = j .\\        
    \end{cases}
\eeq
In the case of the strip guideposts, the bounds only need to be defined on $r_i \in \tilde{\Omega}_s = (0, \lx)$. We therefore take
\beq
    \mathcal{P}_{\text{wall}} = \sum_{i=1}^{N} e^{-a_1 r_i} + e^{-a_1(\lx - r_j)},
\eeq
with its derivatives
\beq
    \frac{\mathcal{P}_{\text{wall}}}{\partial r_i} = 
        -a_1(e^{-a_1 r_i} - e^{-a_1(\lx - r_i)})
\eeq
and 
\beq
    \frac{\partial^2 \mathcal{P}_{\text{wall}}}{\partial_{r_i}\partial_{r_j}} = \delta_{ij}
        a_1^2 (e^{-a_1 r_i} + e^{-a_1(\lx - r_i)}).
\eeq

\section{Solving the optimal design problem}\label{sec:optimizer}
\subsection{Inexact Newton-CG with bounded step length}
% \dc{Added some details.} 
We propose a variant of the inexact Newton-CG optimization algorithm \cite{NocedalWright06} to solve the optimal design problem \eqref{eq:optimal_design}--\eqref{eq:pde_constraint} for the guidepost locations. Starting at some initial guess of guidepost locations $z^{(0)}$, at iteration $k = 0, 1, \dots $, we update the design variable as 
$$ 
z^{(k+1)} = z^{(k)} + \beta_z^{(k)}\delta z^{(k)}. 
$$
% \pc{$\alpha^{(k)}$ is replaced by $\beta^{(k)}$ everywhere since $\alpha$ has been used as a regularization parameter} 
Here $\beta_z^{(k)}$ is a step size and $\delta z^{(k)}$ is a descent direction, computed by solving the Newton system 
$$ 
\nabla^2 \mathcal{J}^{(k)} \delta z^{(k)} = - \nabla \mathcal{J}^{(k)}
$$
% where we use the notation $\nabla^2 \mathcal{J}^{(k)} := \nabla^2 \mathcal{J}(z^{(k)})$ and $\nabla \mathcal{J}^{(k)} := \nabla \mathcal{J}(z^{(k)})$. 
where $\nabla \mathcal{J}^{(k)}$ and $\nabla^2 \mathcal{J}^{(k)}$ are the gradient and Hessians of the cost function evaluated at $z^{(k)}$. 
The system is solved using the CG algorithm, terminating when the residual satisfies
$$ \|\nabla^2 \mathcal{J}^{(k)} \delta z^{(k)} + \nabla \mathcal{J}^{(k)}\| \leq r_{tol} \|\nabla \mathcal{J}^{(k)}\| $$
with tolerance
$r_{tol} = \min(0.5, \sqrt{\|\nabla \mathcal{J}^{(k)}\|/\|\nabla \mathcal{J}^{(0)}\|}),$
or upon the detection of negative curvature (see Algorithm 7.2, \cite{NocedalWright06}). This choice ensures that the Newton system is not solved to an excessively fine tolerance at early optimization iterations, when the design is far from the optimum. At the same time, the CG termination condition ensures asymptotic superlinear convergence of optimization iterations.
% We use a slightly modified form of this algorithm in which we do not compute a 

We cannot insist on monotonic descent in the cost, as changes to the solution branch of $u$ during optimization are often unpredictable and lead to discontinuities in the cost with respect to the design variable $z$. Thus, instead of a backtracking line search, we propose to clip the Newton step and select $\beta_z^{(k)}$ by
\beq\label{eq:step-k}
\beta_z^{(k)} = \min(1, \stepbound /\|\delta z^{(k)}\|_{\infty}),
\eeq
so that $\|\beta_z^{(k)} \delta{z}^{(k)}\|_{\infty} < \stepbound$, where $\stepbound > 0$ is a bound on the length of the optimization step. Near a local minimum, the proposed Newton steps $\delta z^{(k)}$ are small and automatically satisfy this bound. Therefore, we recover the natural step size of $\beta_z^{(k)} = 1$ and preserve the asymptotic convergence rates of the Newton method. 

The bound on the step size is used to improve the performance of the optimizer for our particular problem. The repelling potential can be extremely steep when guideposts start close to each other, and can result in large steps that move the guideposts out of the domain. The bound limits the distance moved by the guideposts to prevent overshooting. 
Moreover, having a small bound on the step also helps to achieve continuity of the state at successive iterations in the optimization scheme when a continuation scheme is used, as discussed below. 
% In this work, we select the value of $\stepbound$ to be on the order of the feature sizes of typical equilibrium morphologies encourage smooth transitions between the equilibrium states. 
In this work, we select the value of $\stepbound$ to be on the order of the width of the guideposts to encourage smooth transitions between the equilibrium states.

\subsection{Choice of initial guesses for the state equation solver}\label{sec:initial_guess_opt}
% \dc{Rewrote later parts}
Due to the non-uniqueness of the state problem, the choice of the initial guess for the state equation solver during optimization is important and nuanced. % The most fundamental 
In our formulation, we can use the solution operator $\mathcal{S}$ to write $u = \mathcal{S}(f(z), u^{(0)})$ for the solution to the state problem given some design variable $z \in \mathcal{Z}$ and an initial guess $u^{(0)} \in H^1(\Omega)$. One approach is to fix a sample from $\mathcal{D}$ as the initial guess for every iteration of the optimizer. That is, we use the same initial guess every time a state equation solved. However, the solutions $u = \mathcal{S}(f(z), u^{(0)})$ are sensitive to the substrate pattern, and small perturbations to the design variable $z$ can lead to convergence of $u$ to a different solution branch (as shown in the numerical example in Section \ref{sec:differentiability}). Thus, there is frequently a discontinuity of the mapping between $z$ and $u$ when a fixed initial guess is taken. We recall that $z$-derivatives in $\nabla_z Q$, $\nabla^2_z Q$ are local to the state $u$ used to compute them. This means they often do not accurately reflect the behavior of $Q(z)$ when changes in $z$ push the equilibrium state to other solution branches.
%, with morphologies changing significantly. 
This poses significant challenges to our derivative-based optimization methods.

Instead, one can try to preserve continuity (and smoothness) between $z$ and $u$ by a continuation scheme. We adopt the simple approach introduced in Section \ref{sec:differentiability}, which uses the current equilibrium state as the initial guess to the subsequent state equation solve. That is, using $u(z^{(k)})$ to denote the state corresponding to the $k$-th optimization iteration for the design variable, we compute the state in the next iteration, $u(z^{(k+1)})$, by
\beq
    u(z^{(k+1)}) = \mathcal{S}(f(z^{(k+1)}), u(z^{(k)})).
\eeq
This encourages the solution at $z^{(k+1)}$ to converge within the same branch as $u(z^{(k)})$, thereby achieving continuity of the state between successive iterations of the design variable, provided the step size $z^{(k+1)} - z^{(k)}$ is not too large. This can be enforced by using a small bound on the step length in the optimization algorithm. 

\subsection{Assessment of robustness to initial guesses}
% \dc{Rewrote this section.} 
Since the same design variable $z$ can yield different equilibrium states $u$ and, hence, different gradients $D_z \mathcal{J}$ corresponding to different initial guesses $u^{(0)}$, a design that is optimal for a particular state solution may be suboptimal when the state equation is solved starting at a different initial guess. It is therefore necessary to evaluate the robustness of the candidate optimal designs with respect to various initial guesses.  

In this work, we propose to first solve the optimization problem for an optimal design $z^*$, during which the state equations are solved using initial guesses given by the continuation scheme. Upon convergence to an optimal design, we then assess its robustness to the initial guess for the optimal design. To do so, we solve the state equation using $N$ samples of the initial guess $u_i^{(0)} \sim \mathcal{D}$ and inspect the various equilibrium states $u_i = \mathcal{S}(f(z^*), u_i^{(0)})$. We also compute the associated free energies $\mathcal{F}(u_i)$, and use the state with the minimal free energy as a proxy for the most likely equilibrium solution
\begin{linenomath}
$$u^* = \argmin_{i=1,...,N} \mathcal{F}(u_i),
$$
\end{linenomath}
and assess its optimality. Furthermore, we consider sample statistics of the design objective over the samples $Q_i = Q(\mathcal{S}(f(z^*), u_i^{(0)}))$ as an indicator for the robustness of the design.

We remark that this scheme of warm-starting the forward solver builds in continuity and smoothness to the optimization problem and aids convergence to the minimizer. However, by construction, it biases the equilibrium states computed in successive iterations to lie in a particular solution branch. These branches may lead to more optimal equilibrium states than otherwise achievable when using unbiased initial guesses such as those sampled from $u^{(0)} \sim \mathcal{D}$. An assessment of the optimal designs is therefore all the more important.

\subsection{Computational cost}
Most of computational cost in the optimization algorithm is in solving the PDEs to evaluate the objective function and its derivatives. Due to the bilevel nature of the optimal design problem, we will refer to solving the energy-stable Newton equation for state problem as the \textit{inner} Newton iterations to distinguish them from the \textit{outer} Newton iterations of the guidepost optimization algorithm. In our algorithm, we solve the state equation once for every optimization iteration. This itself requires a number of inner Newton iterations until convergence, with each inner Newton iteration being a linear PDE solve. The algorithm then requires one adjoint equation solve to compute the gradient, and a pair of incremental forward and adjoint equation solves for every application of the Hessian action in the CG algorithm. 

Suppose a total of $M$ optimizations steps are taken, and $N_{s}^{(k)}$ inner Newton iterations and $N_{c}^{(k)}$ Hessian actions are required at step $k$ for $k = 0,...,M$, then the total number of linear PDE solves required, $N_t$, is approximately
\beq
N_t = \underbrace{\sum_{k=0}^{M} N_{s}^{(k)}}_{\text{state}} + \underbrace{M}_{\text{adjoint}} + \underbrace{\sum_{k=0}^{M-1} 2 N_{c}^{(k)}}_{\text{incremental}}.
\eeq
We note that within one optimization step, the operators associated with the incremental forward and adjoint equations are constant across CG iterations. Thus, the cost of repeatedly applying the Hessian action can be ameliorated by reusing the factorization of the PDE operators if using a direct solver, or by building a good preconditioner if using an iterative solver. Therefore, the linear PDE solves of the inner Newton iterations remains the dominant cost of the algorithm. The continuation scheme along with the bound on the step size is also used to reduce the number of inner Newton solves $N_{s}^{(k)}$ required converge in the state problem. As demonstrated in \cite{CaoGhattasOden22}, the number of inner Newton iterations remains small and independent of the state dimension, as is characteristic of Newton methods for many operator equations. Moreover, as demonstrated below, the number of outer iterations depends only weakly on the design variable dimension. Thus, overall, we arrive at a method for optimal design that is fast and scalable with respect to state and design dimensions.

\section{Numerical Results}\label{sec:results}
\subsection{General setting}{}
Here, we present results from numerical experiments in which we solve the optimal design problem for various target morphologies and guidepost configurations. 
% The model parameters adopted in the subsequent numerical experiments are shown in Table \ref{tab:material_params}. 
Starting from an initial guidepost configuration, we optimize guidepost locations using the bounded step-length Newton CG optimizer with continuation of the state solver, as described in Section \ref{sec:optimizer}. We then assess the converged optimal designs by solving the state problem with initial guesses $u^{(0)}_i$ drawn independently from the transformed Gaussian random field. As suggested previously, we implement the finite element discretization using FEniCS, and use LU solvers from PETSc \cite{BalayAbhyankarAdamsEtAl15} to solve the resulting linear systems. We also draw from hIPPYlib \cite{VillaPetraGhattas20} to implement the optimization algorithm and sampling from the initial guess distribution $\mathcal{D}$. 

Across our experiments, we assume fixed model parameters of $m = 0$, $\varepsilon = 0.08$, $\sigma = 12.8$, corresponding to a lamella producing material. Based on our model parameters, we use $(\delta_G, \gamma_G) = (0.8, 0.02)$ for the covariance operator of the Gaussian random field and $s=1$ as the scaling parameter.  This selection ensures fluctuations in the random initial guess are over much finer scales than the lamellae structures in typical equilibrium morphologies. We adopt a strength of $w = 0.5$ for the guideposts. 

% Based on our model parameters, we select a step size bound of $\stepbound = 0.2$ for the Newton algorithm, which is on the same order as the size of typical features exhibited by the equilibrium morphologies. 

% \begin{table}[htbp!]
%     \centering
%     \begin{tabular}{l r}
%     \toprule
%     Model parameters  \\
%     \midrule
%     $m$ & 0 \\
%     $\varepsilon$ & 0.08 \\
%     $\sigma$ & 12.8 \\
%     $w$ & 0.5 \\
%     \bottomrule
%     \end{tabular}
%     \caption{Material parameters used in the numerical experiments.}
%     \label{tab:material_params}
% \end{table}

\subsection{Design with strip guideposts in 2D}
\subsubsection{Problem set-up}
\begin{figure}[htbp!]
\centering
\begin{subfigure}{0.32\textwidth}
    \centering
    \includegraphics[width=0.95\textwidth]{figures_pdf/state_colors.pdf}
\end{subfigure}

\begin{subfigure}{0.32\textwidth}
    \centering
    \includegraphics[width=0.95\textwidth]{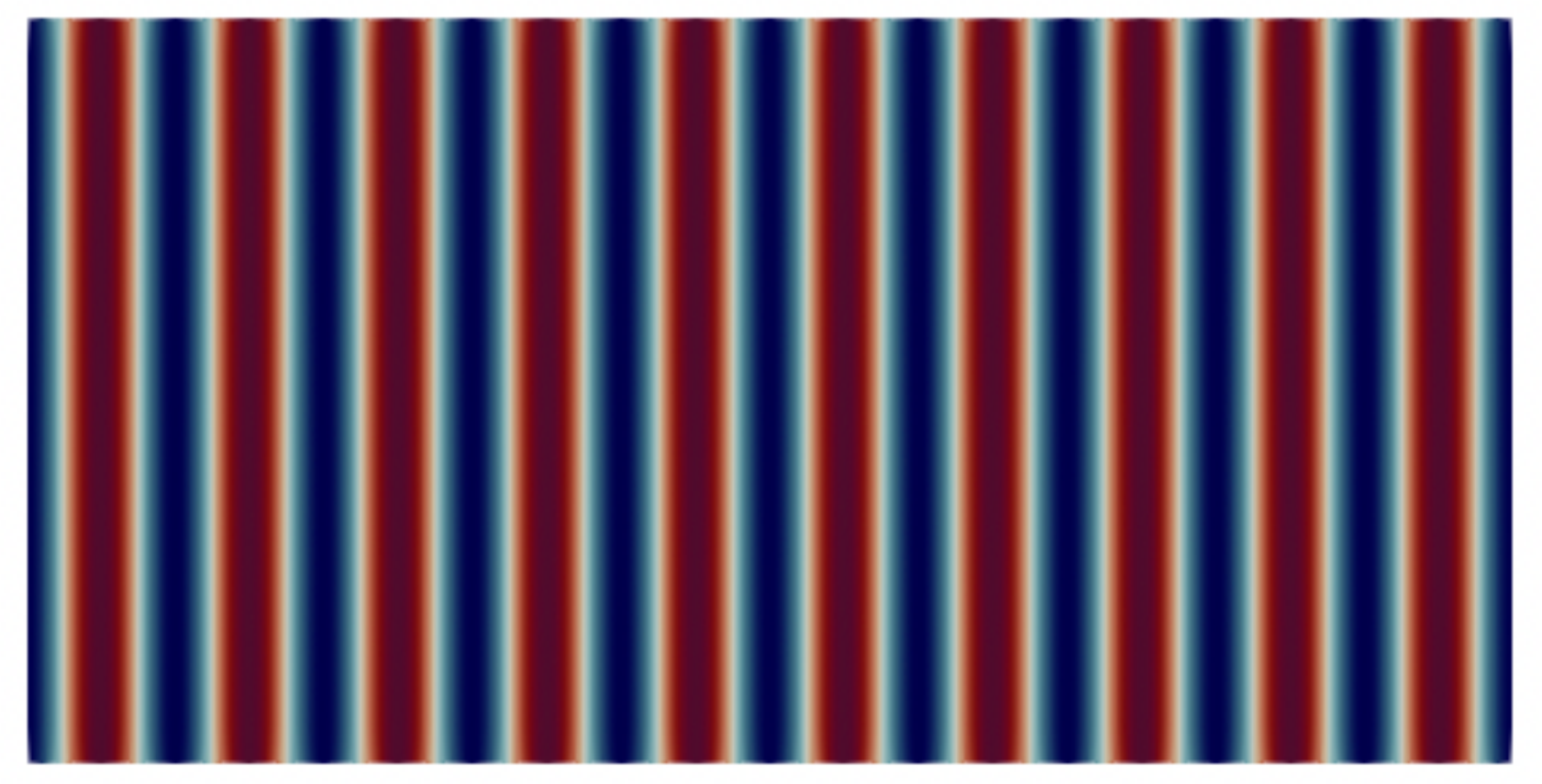}
    \caption{Target morphology}
    \label{fig:2d_strip_target_new}
\end{subfigure}
\begin{subfigure}{0.32\textwidth}
    \centering
    \includegraphics[width=0.95\textwidth]{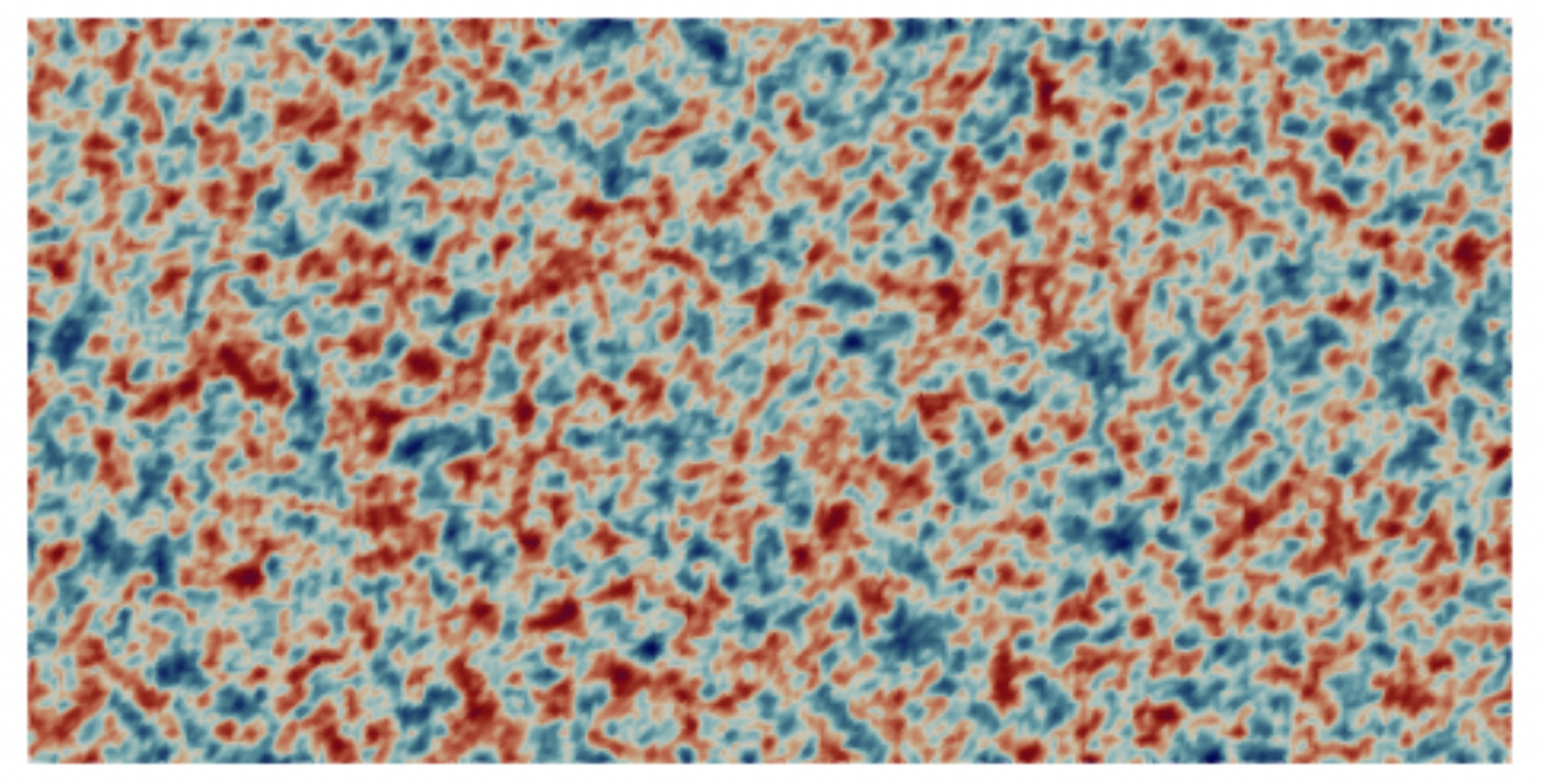}
    \caption{Initial guess}
    \label{fig:2d_strip_initial}
\end{subfigure}
\begin{subfigure}{0.32\textwidth}
    \centering
    \includegraphics[width=0.95\textwidth]{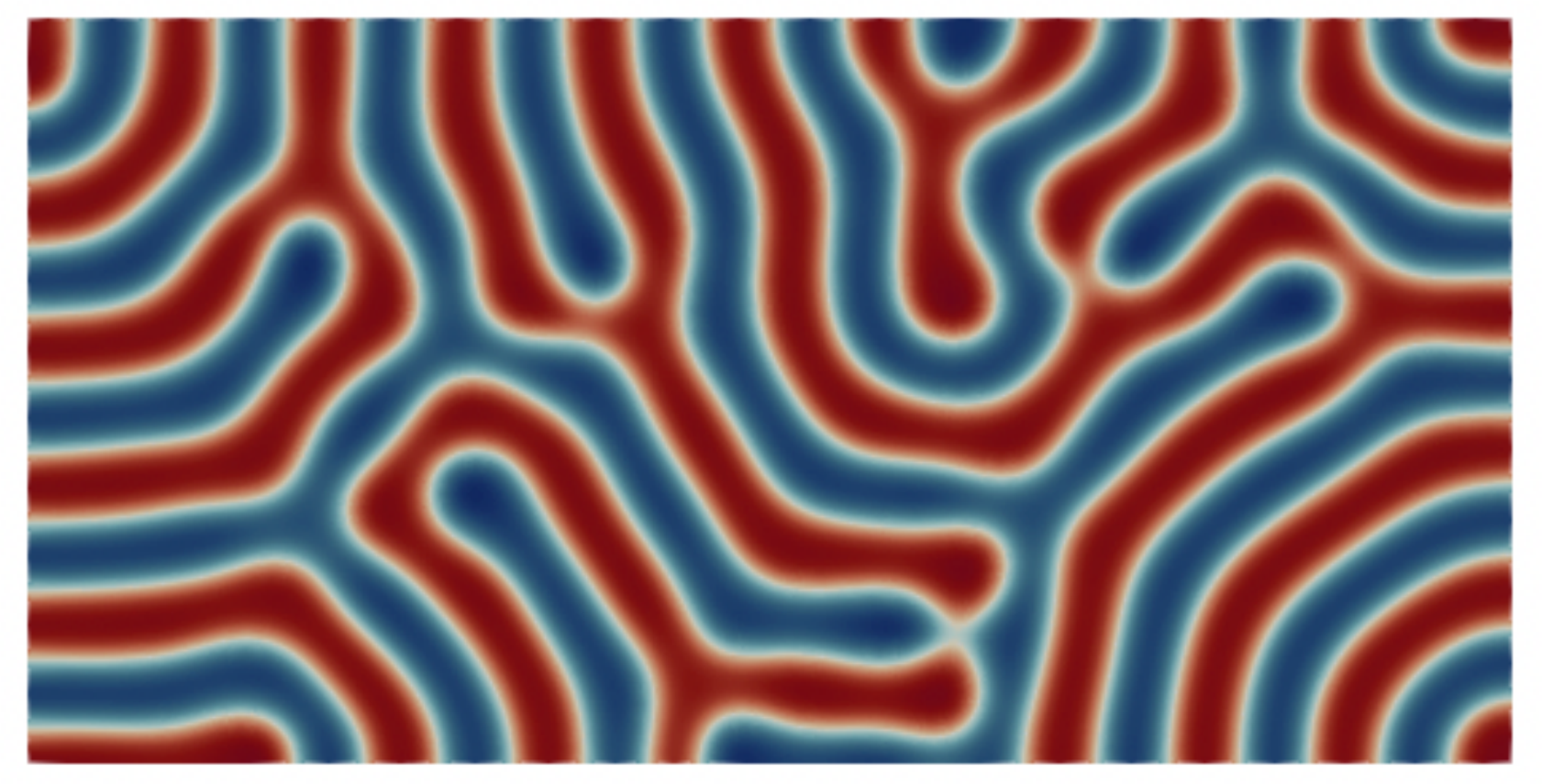}
    \caption{Equilibrium state}
    \label{fig:2d_strip_equilibrium}
\end{subfigure}
\caption{Strips target morphology in 2D (left), shown alongside an initial guess drawn from $\mathcal{D}$ (center) and the equilibrium state
computed using the initial guess over a fully neutral substrate for a comparison of the feature sizes.}
\label{fig:2d_strip_features}
\end{figure}

We first consider a problem in which the target morphology consists of equally spaced vertical strips spaced at a distance $l_{target}= 1$ in a domain $\Omega = [0, 10]\times[0,5]$, as shown in Figure \ref{fig:2d_strip_features}. This is selected to approximately match the natural size of the lamella features corresponding to our chosen model parameters. For comparison, we also plot an initial guess drawn from the transformed Gaussian distribution $\mathcal{D}$ along with the equilibrium state produced with a neutral substrate, showing the relative sizes of the features. We then optimize the locations of four strip shaped guideposts with width $b=0.2$ for the shape function, using a maximum step size of $\stepbound = 0.2$. 

% \begin{figure}[htbp!]
% \centering
% \includegraphics[width=0.3\textwidth]{figures_pdf/strip_results/strip_target.pdf}
% \caption{Target morphology consisting of vertical strips in 2D.}
% \label{fig:2d_strip_target}
% \end{figure}

This is a common test problem in both experimental and computational investigations (e.g. \cite{LiuRamirezHernandezHanEtAl13, KhairaQinGarnerEtAl14}), since the target morphology is simple and of practical importance. In this case, one expects the optimal guideposts to be regularly spaced, aligning with the strips in the target morphology. However, as an optimal design problem, it still exhibits the typical challenges we face, namely the non-convexity of the design objective and the difficulties arising from the non-uniqueness of the state problem. 

\begin{figure}[htbp!]
\centering
\begin{subfigure}{0.45\textwidth}
\includegraphics[width=0.95\textwidth]{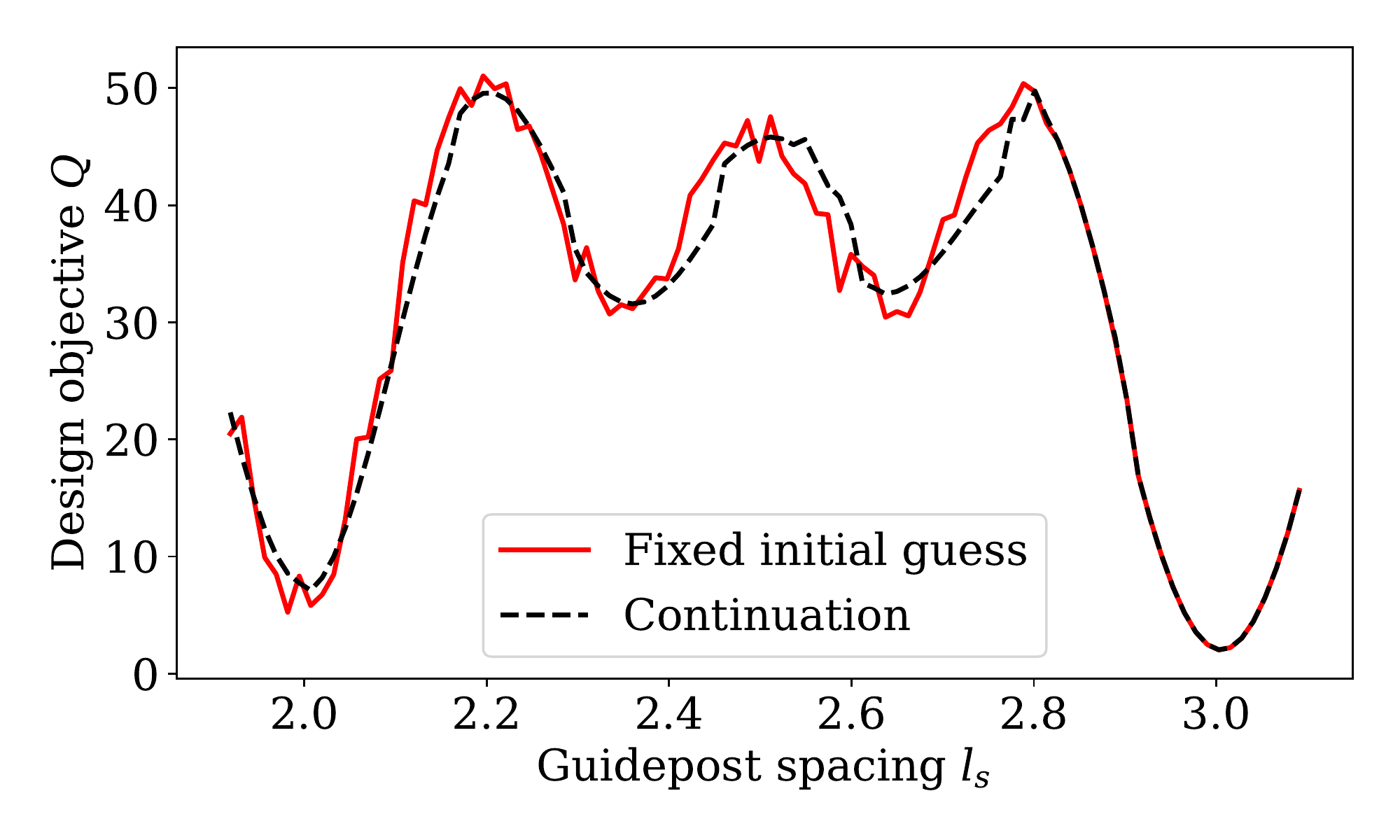}
\end{subfigure}
\begin{subfigure}{0.45\textwidth}
\includegraphics[width=0.95\textwidth]{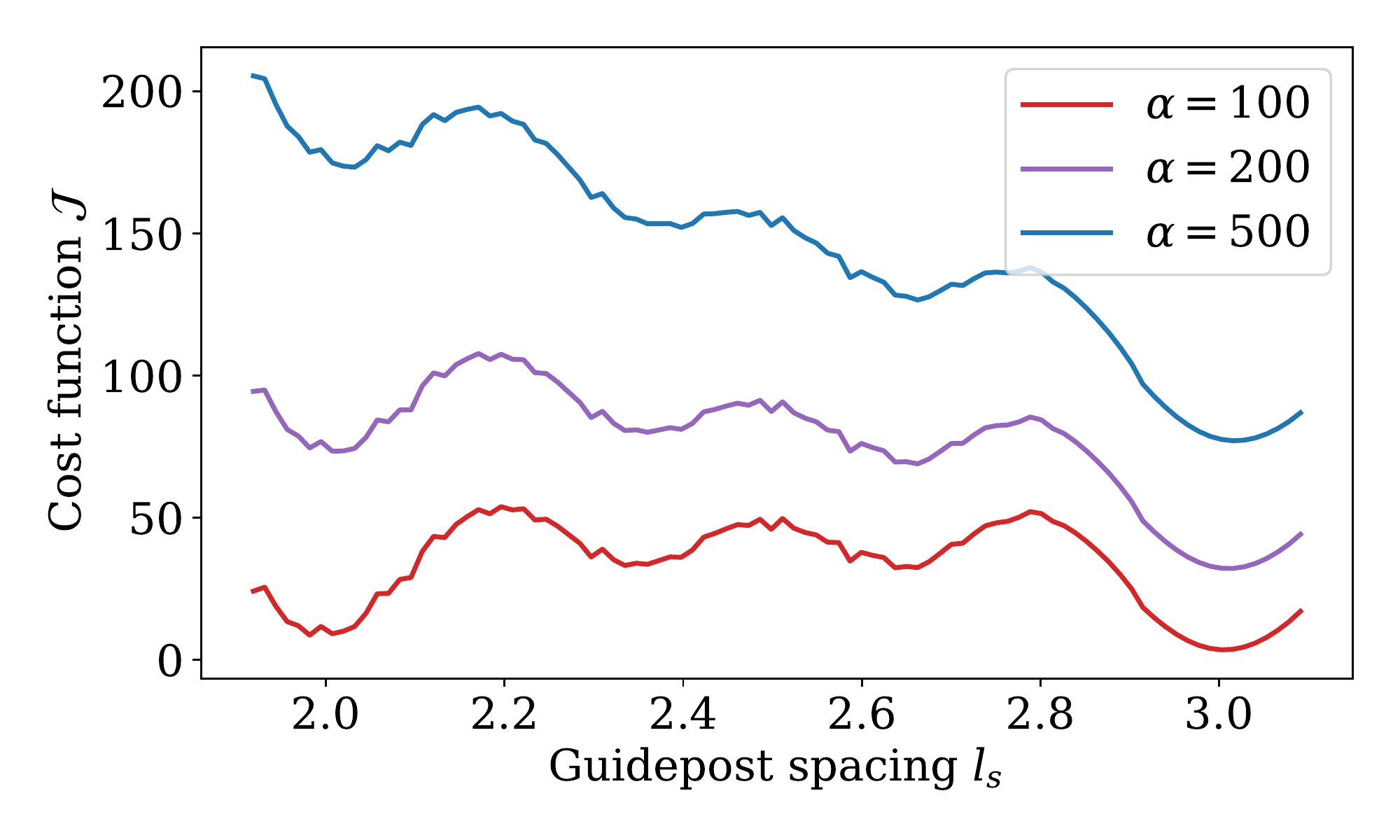}
\end{subfigure}

\caption{Left: Design objective $Q(u;u_d) = \|u - u_d\|_{L^2(\Omega)}$ computed for four equally spaced strip guideposts of various spacing $l_s$. For the red line, the state problems are solved using a fixed initial guess drawn from $\mathcal{D}$. For the black dashed line, the state problems are solved with the proposed continuation scheme. Right: Cost function with repulsion penalization term added, $Q + \mathcal{P}_{\text{repel}}$, with $\alpha = 100$, 200, and 500, using objective values from the fixed initial guess.}

\label{fig:2d_strip_obj}
\end{figure}

To visualize the aforementioned challenges, we first consider a simple representation in which the four substrate guideposts are equally spaced at some distance $l_s$, with the leftmost guidepost fixed at $x_1 = 1/2.$ We plot the design objective in terms of the guidepost spacing. To do so, we fix an initial guess drawn from $\mathcal{D}$ and use it to solve for the state and evaluate the objective function for increasing values of guidepost spacing. For comparison, we also evaluate the objective function using the continuation scheme. That is, as we increase the spacing $l_s$, we use the solution from the previous design as the initial guess for the next. We plot resulting objective values in Figure \ref{fig:2d_strip_obj}, along with the cost function computed from adding the repulsion penalization with various strengths $\alpha$. 

In this setting, the design objective exhibits several local minima corresponding to configurations where the guideposts coincide with the strips in the target morphology ($l_s = 2$ and $l_s = 3$). This is regularized by the addition of the repulsion penalization, making the minimum near $l_s = 3$ more pronounced.
% Furthermore, for each guidepost spacing, the different initial guesses can produce different values of the design objective, demonstrating the necessity to evaluate sensitivity of candidate designs to initial guesses of the state problem. 
We also observe that with the fixed initial guess, the design objective tends to be non-smooth, and jumps in the objective function correspond to instances where the state problem converges to an alternative solution branch as the spacing is perturbed. These jumps disappear near $l_s = 3$, where the guidepost configuration constrains the solution to a single branch and instead gives rise to a smooth curve. In contrast, the continuation scheme leads to a much smoother curve throughout, as the state problem is encouraged to remain on a particular solution branch. This highlights the role played by the continuation scheme in producing meaningful derivative computations for the optimizer.

% This suggests that the continuation scheme is important to achieving convergence by our optimizer.

% We make several observations. First, the design objective exhibits several local minima corresponding to configurations where the guideposts coincide with the strips in the target morphology. This behaviour is largely independent of the initial guess used. Furthermore, for each guidepost spacing, the different initial guesses can produce different values of the design objective, demonstrating the necessity to evaluate sensitivity of candidate designs to initial guesses of the state problem. We can also observe that with a fixed initial guess, the value of the design objective tends to be non-smooth, with the jumps corresponding to instances where the state problem converges to an alternative solution branch as the spacing is perturbed. In contrast, the continuation scheme leads to a much smoother curve, particularly near the local minima, as the state problem is encouraged to remain on a particular solution branch. This suggests that the continuation scheme is crucial in achieving convergence by our optimizer.

% Furthermore, jumps in the design objective arise from the state converging to different solution branches as the substrate design changes. 
% In particular, we wish to obtain optimal designs where the spacing is larger than $l_c$ to achieve density multiplication. 
% For our numerical experiments, we use a domain $\Omega = [0, 10]\times[0,5]$ and define the target morphology as in Figure \ref{fig:2d_strip_target}.

\subsubsection{Optimization results}
\begin{figure}[h]
\centering
\begin{subfigure}{0.32\textwidth}
    \centering
    \includegraphics[width=0.95\textwidth]{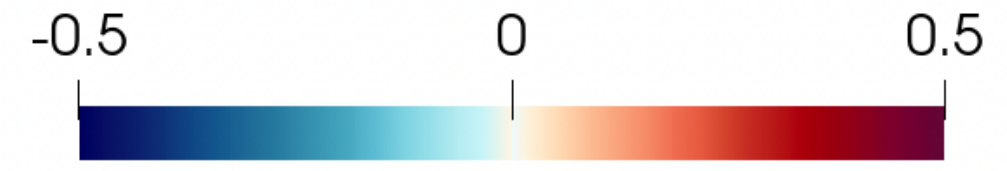}
\end{subfigure}

\begin{subfigure}{0.32\textwidth}
    \centering
    \includegraphics[width=0.95\textwidth]{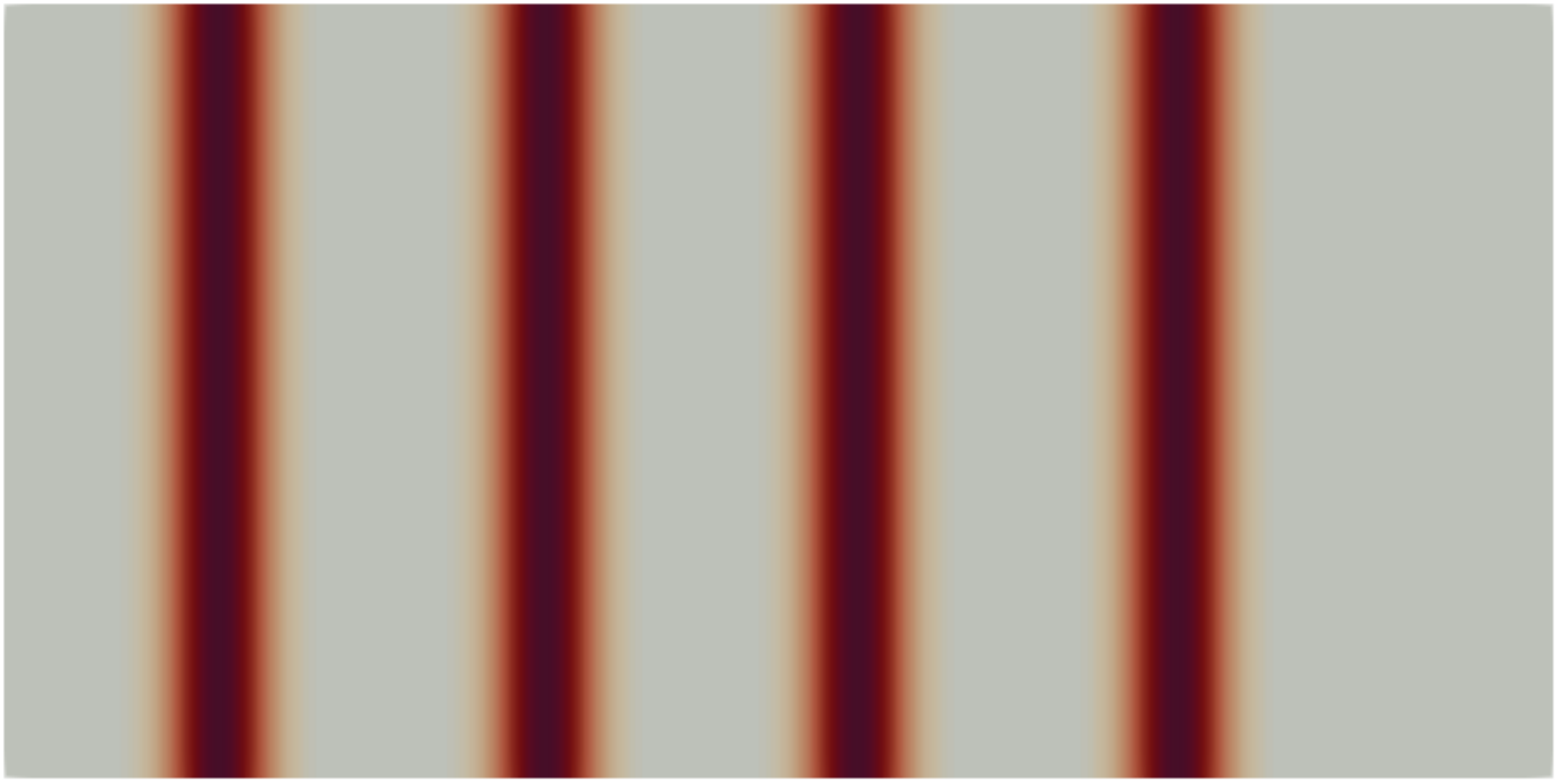}
    \caption{Optimal design $\alpha = 10^{2}$}
\end{subfigure}
\begin{subfigure}{0.32\textwidth}
    \centering
    \includegraphics[width=0.95\textwidth]{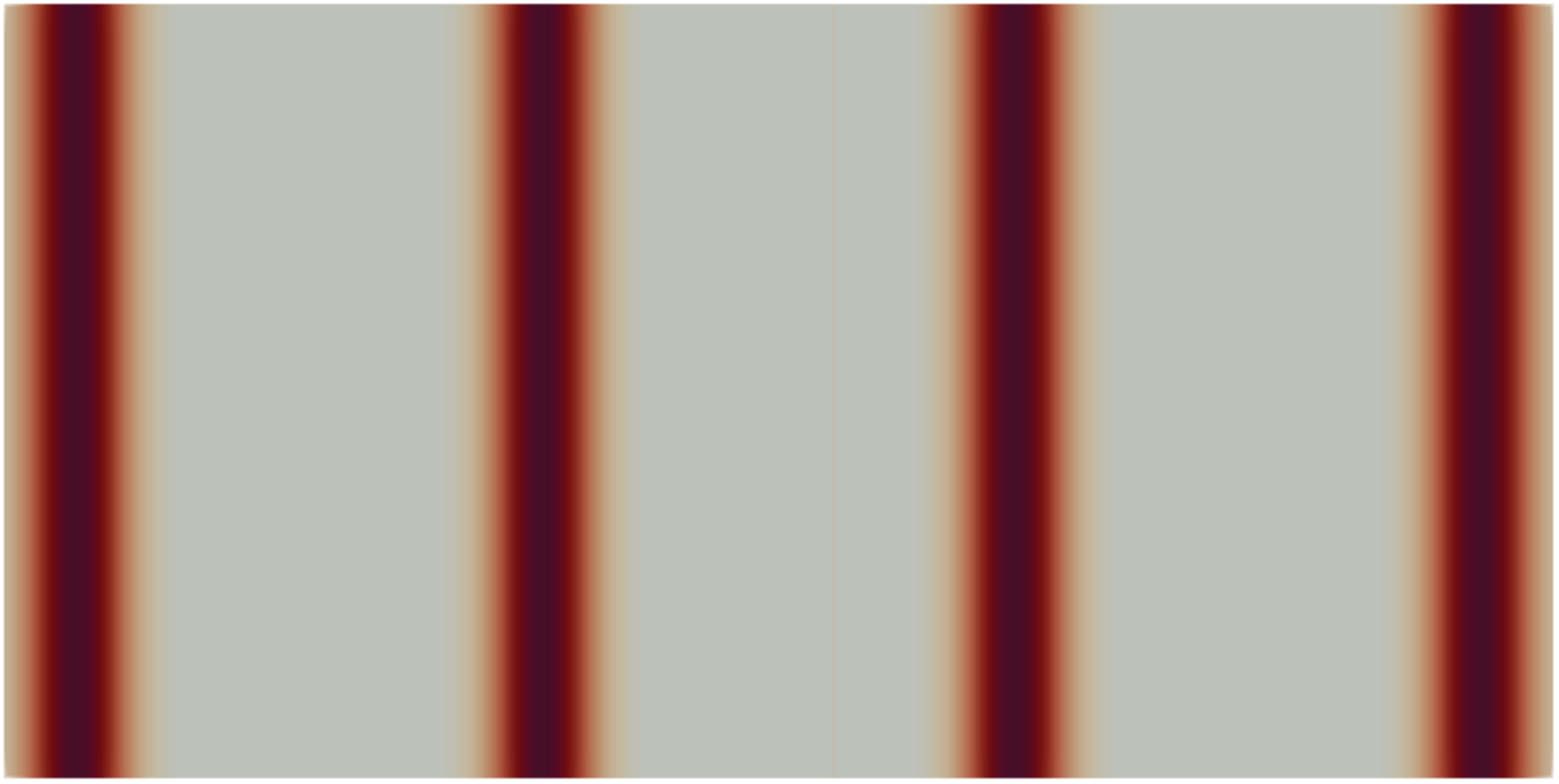}
    \caption{Optimal design $\alpha = 10^{3}$}
\end{subfigure}
\begin{subfigure}{0.32\textwidth}
    \centering
    \includegraphics[width=0.95\textwidth]{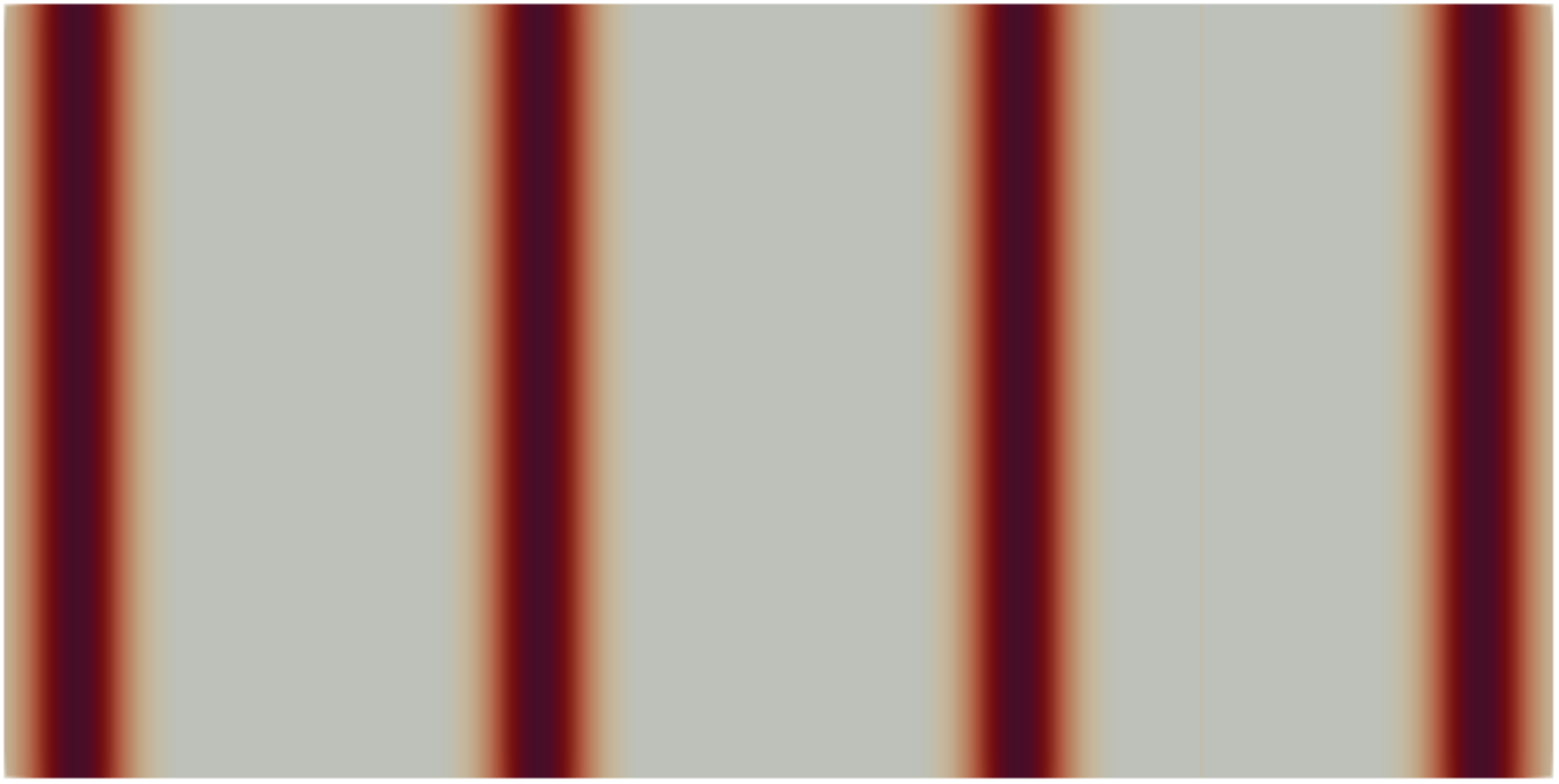}
    \caption{Optimal design $\alpha = 10^{4}$}
\end{subfigure}

\begin{subfigure}{0.32\textwidth}
    \centering
    \includegraphics[width=0.95\textwidth]{figures_pdf/state_colors.pdf}
\end{subfigure}

\begin{subfigure}{0.32\textwidth}
    \centering
    \includegraphics[width=0.95\textwidth]{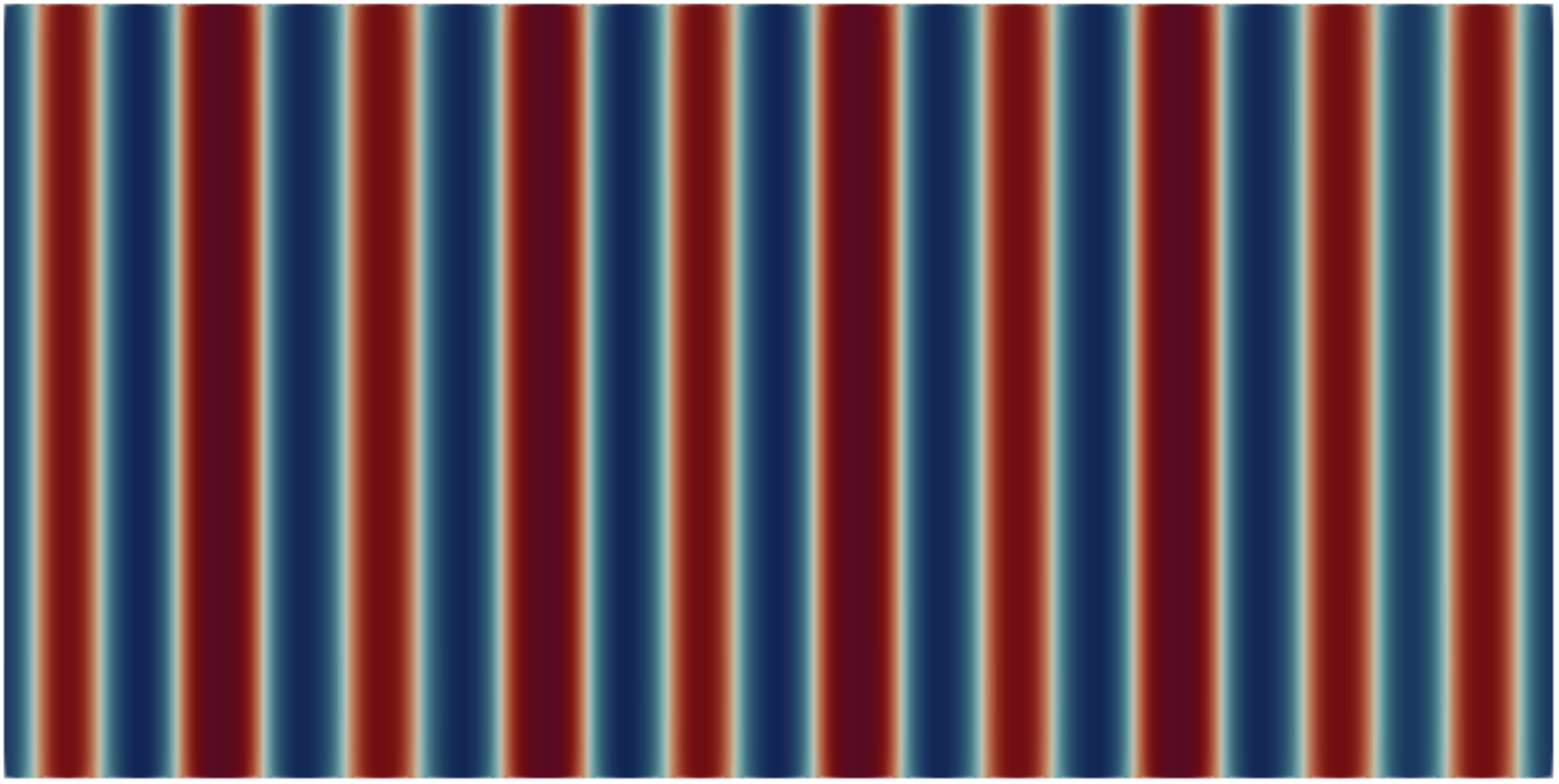}
    \caption{Optimal state $\alpha = 10^{2}$}
\end{subfigure}
\begin{subfigure}{0.32\textwidth}
    \centering
    \includegraphics[width=0.95\textwidth]{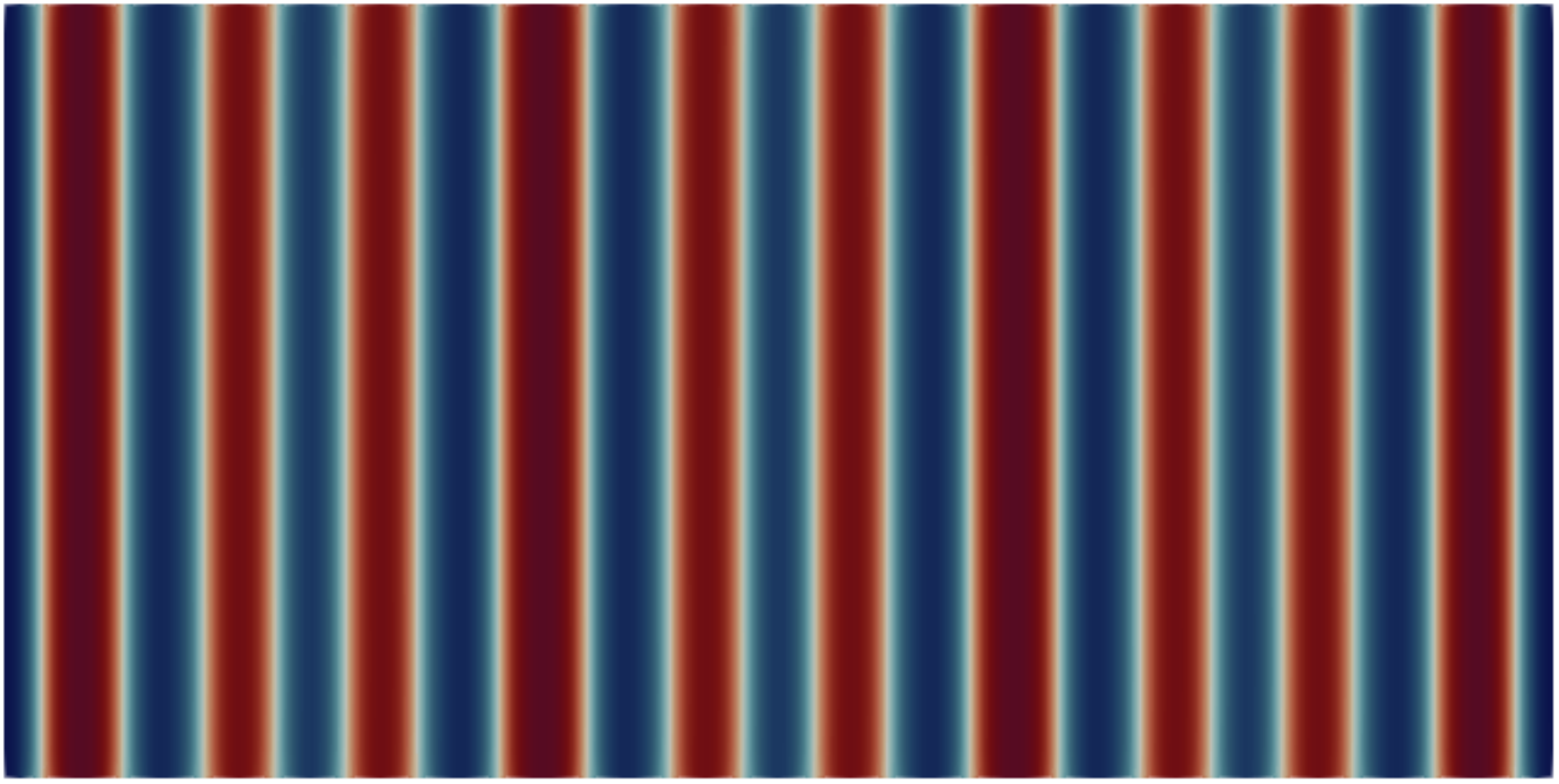}
    \caption{Optimal state $\alpha = 10^{3}$}
\end{subfigure}
\begin{subfigure}{0.32\textwidth}
    \centering
    \includegraphics[width=0.95\textwidth]{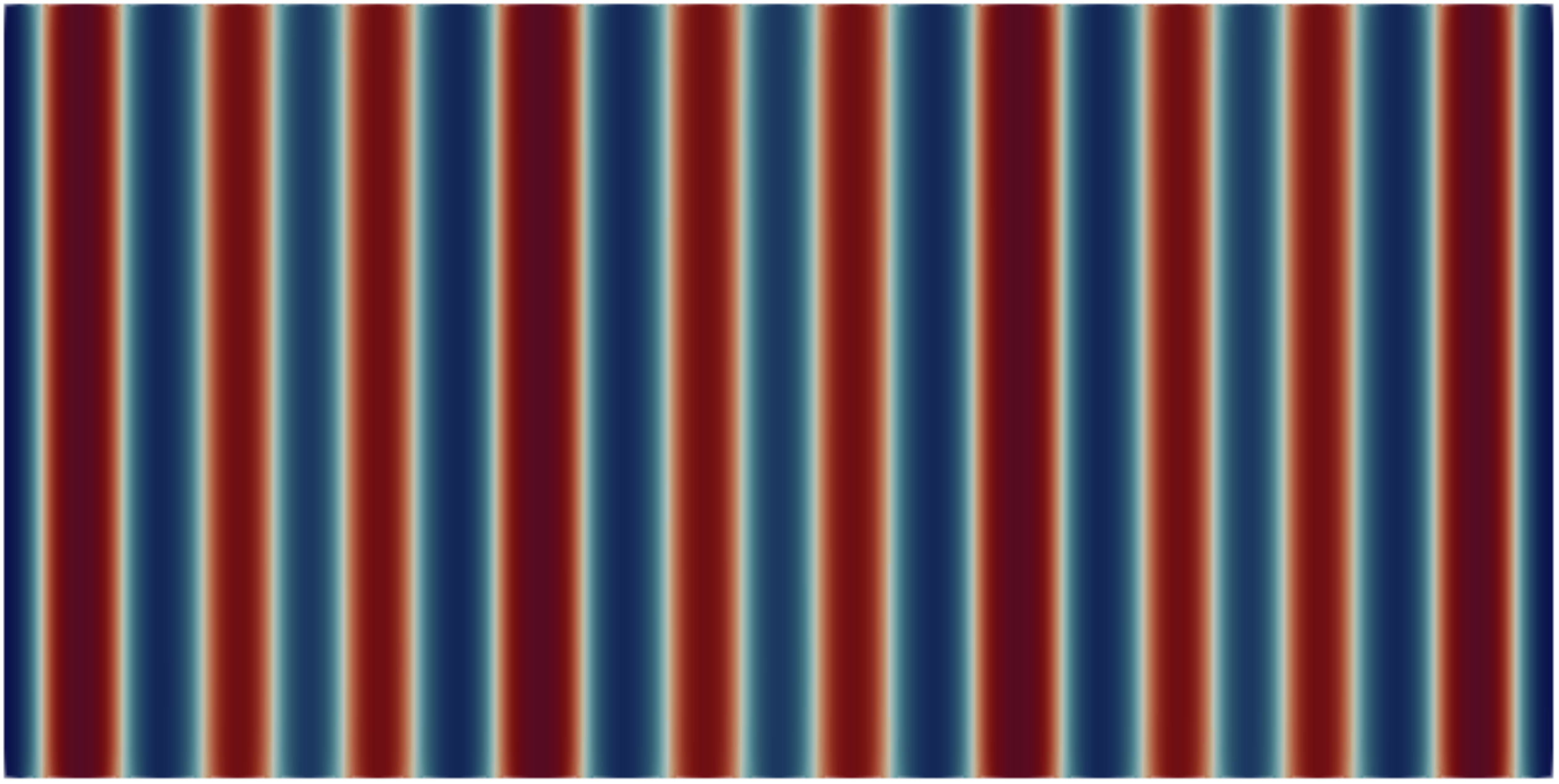}
    \caption{Optimal state $\alpha = 10^{4}$}
\end{subfigure}

\caption{Optimal designs (top) and the corresponding states (bottom) computed at the final optimization iteration for the strip target morphology. Results are shown for different values of the distance penalty parameter $\alpha = 10^{2}$ (left), $\alpha = 10^{3}$ (middle) and $\alpha = 10^{4}$ (right).}
\label{fig:strip_results}
\end{figure}

\begin{figure}[h]
\centering
\captionsetup[subfigure]{justification=centering}

\begin{subfigure}{0.32\textwidth}
    \centering
    \includegraphics[width=0.95\textwidth]{figures_pdf/state_colors.pdf}
\end{subfigure}

\begin{subfigure}{0.32\textwidth}
    \centering
    \includegraphics[width=0.95\textwidth]{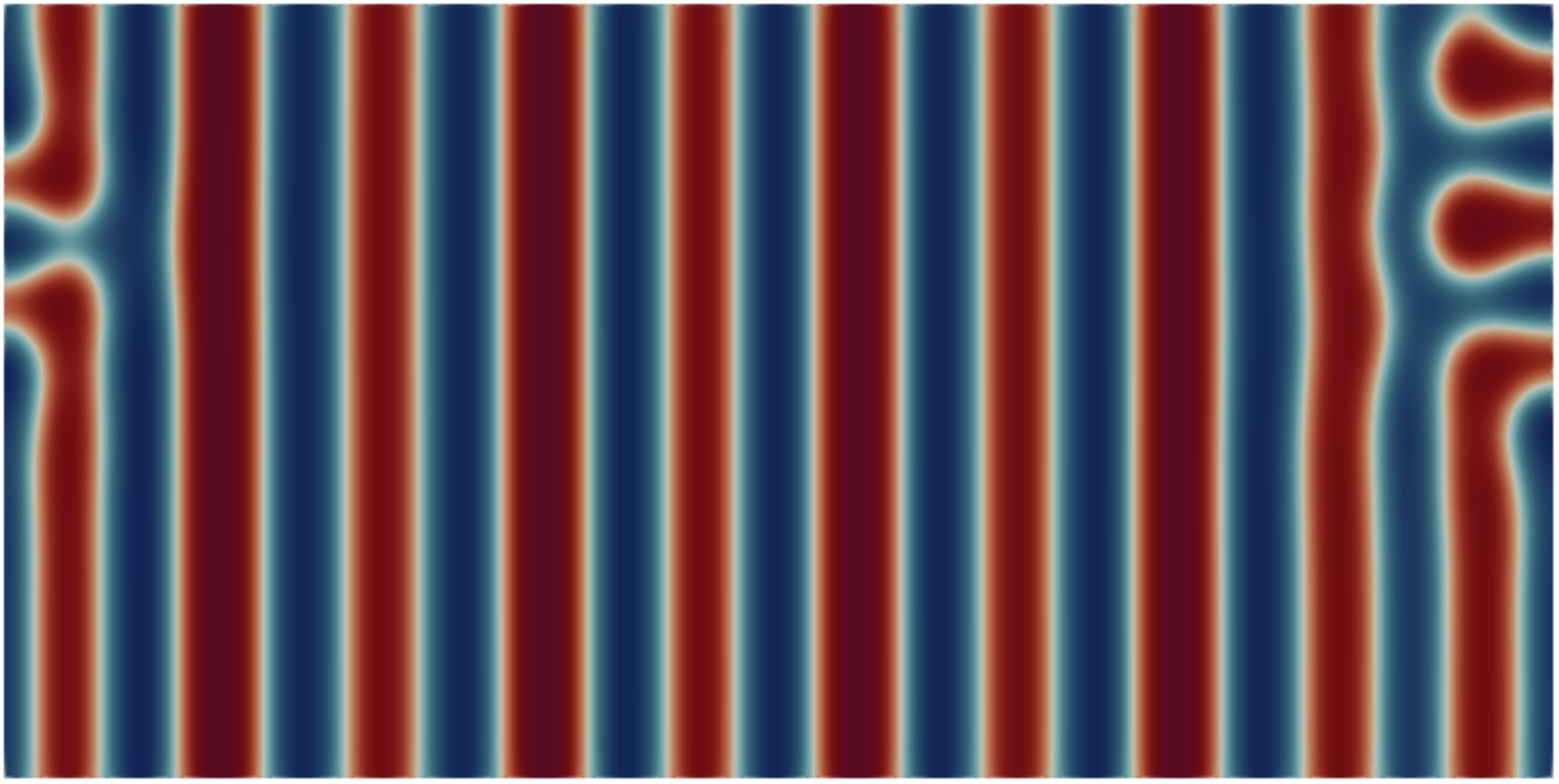}
    \caption{Minimum energy state $\alpha=10^{2}$}
\end{subfigure}
\begin{subfigure}{0.32\textwidth}
    \centering
    \includegraphics[width=0.95\textwidth]{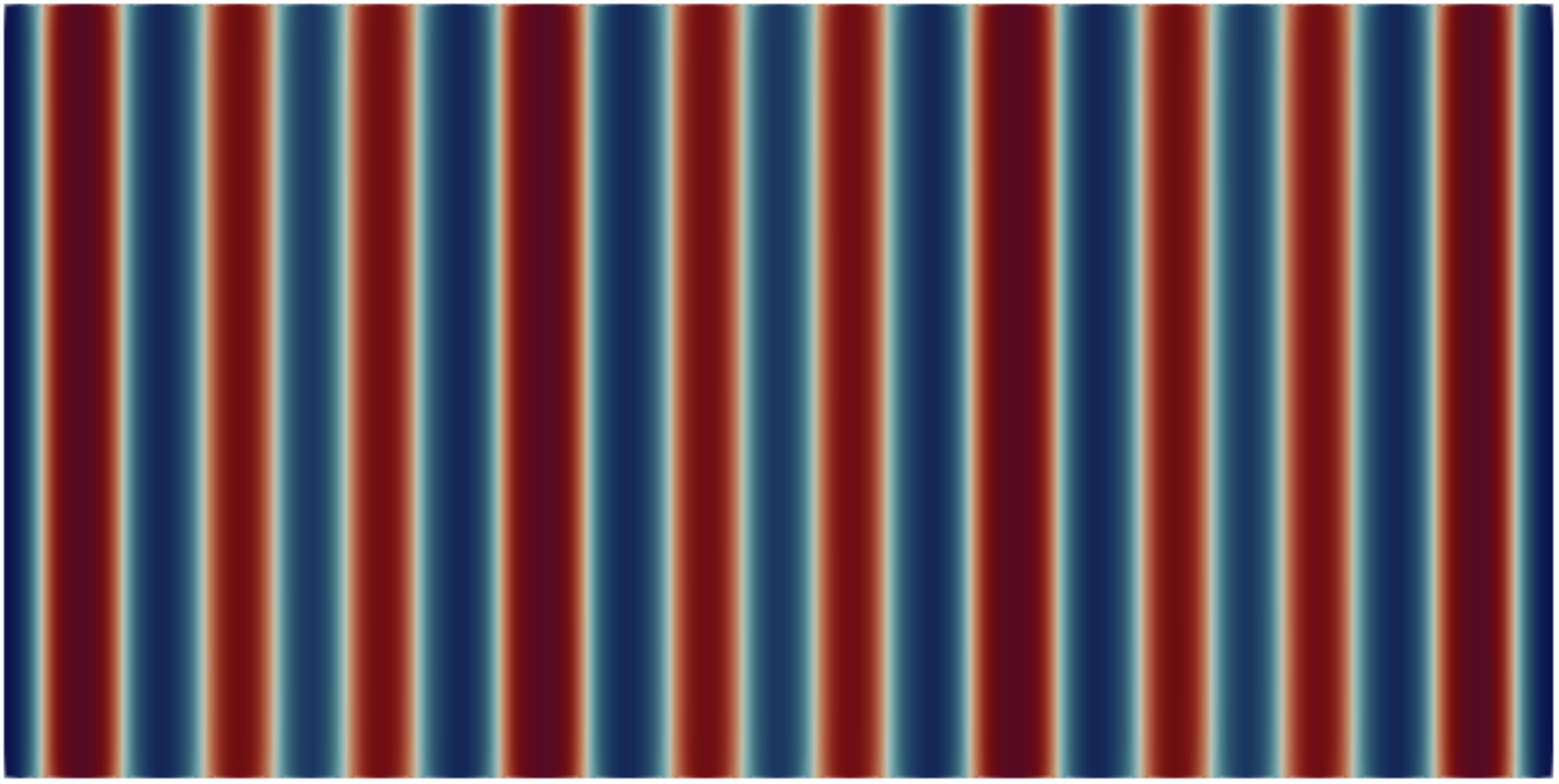}
    \caption{Minimum energy state $\alpha=10^{3}$}
\end{subfigure}
\begin{subfigure}{0.32\textwidth}
    \centering
    \includegraphics[width=0.95\textwidth]{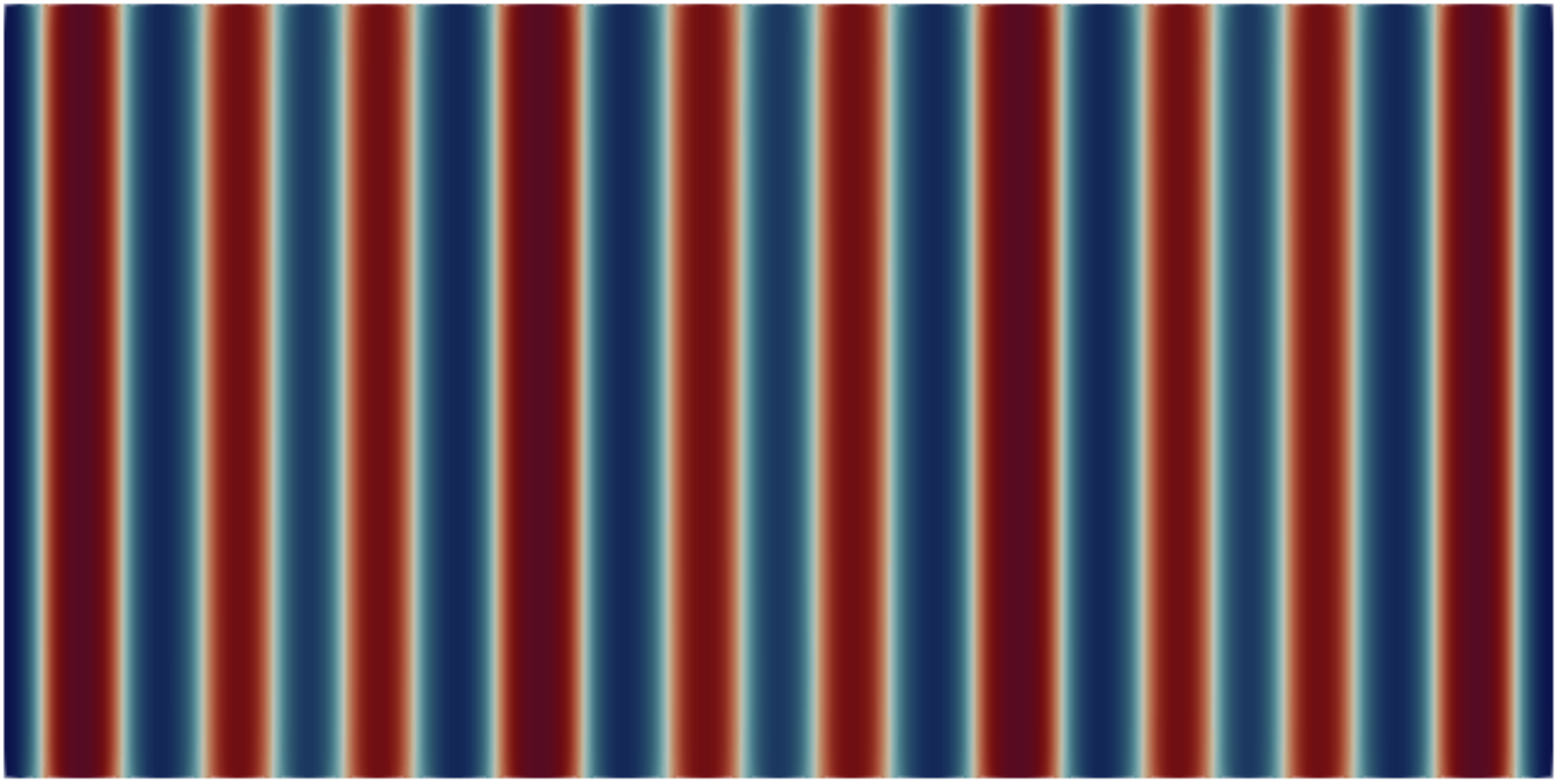}
    \caption{Minimum energy state $\alpha=10^{4}$}
\end{subfigure}

\begin{subfigure}{0.32\textwidth}
    \centering
    \includegraphics[width=0.95\textwidth]{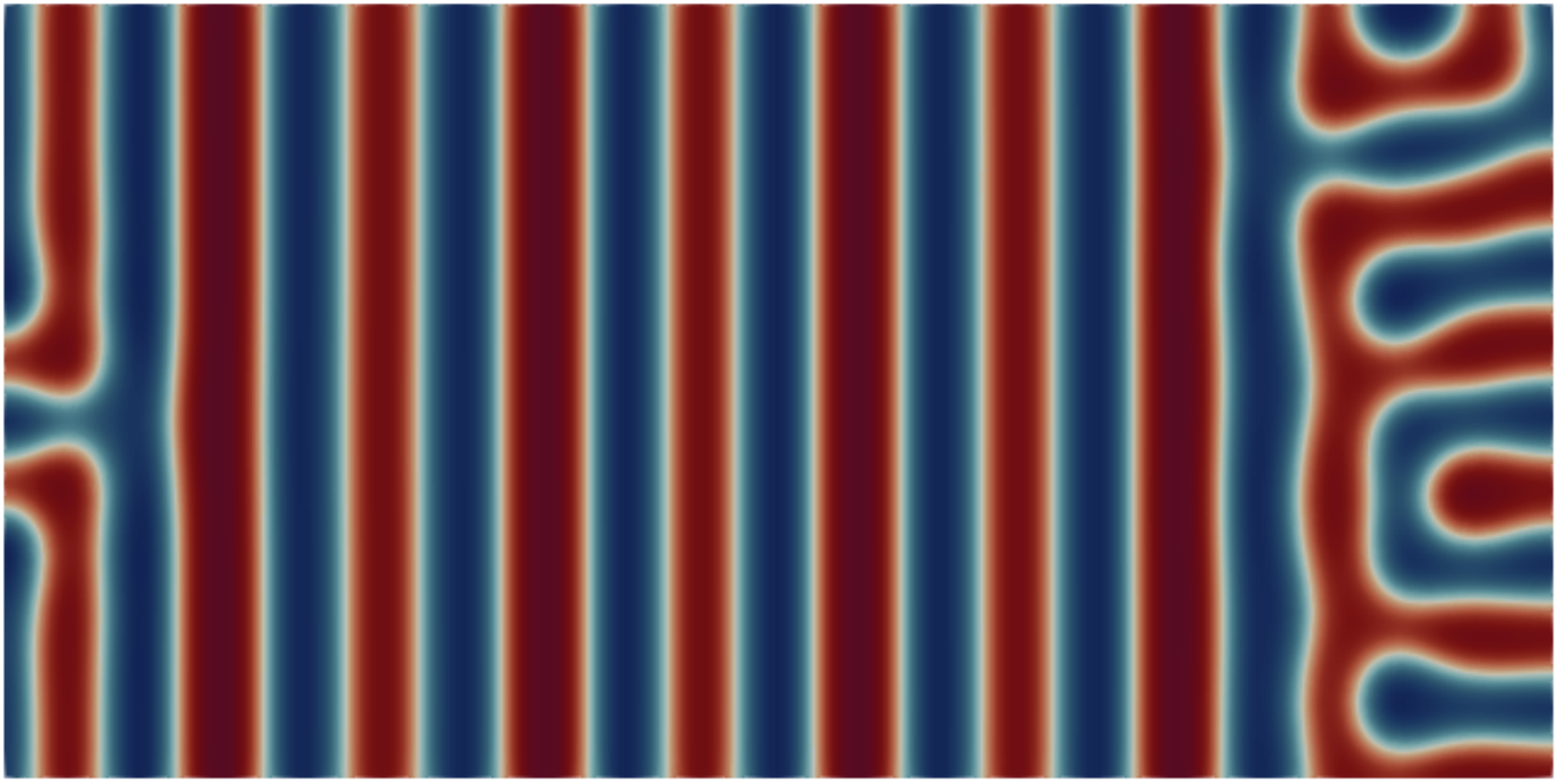}
    \caption{Sample state $\alpha=10^{2}$}
\end{subfigure}
\begin{subfigure}{0.32\textwidth}
    \centering
    \includegraphics[width=0.95\textwidth]{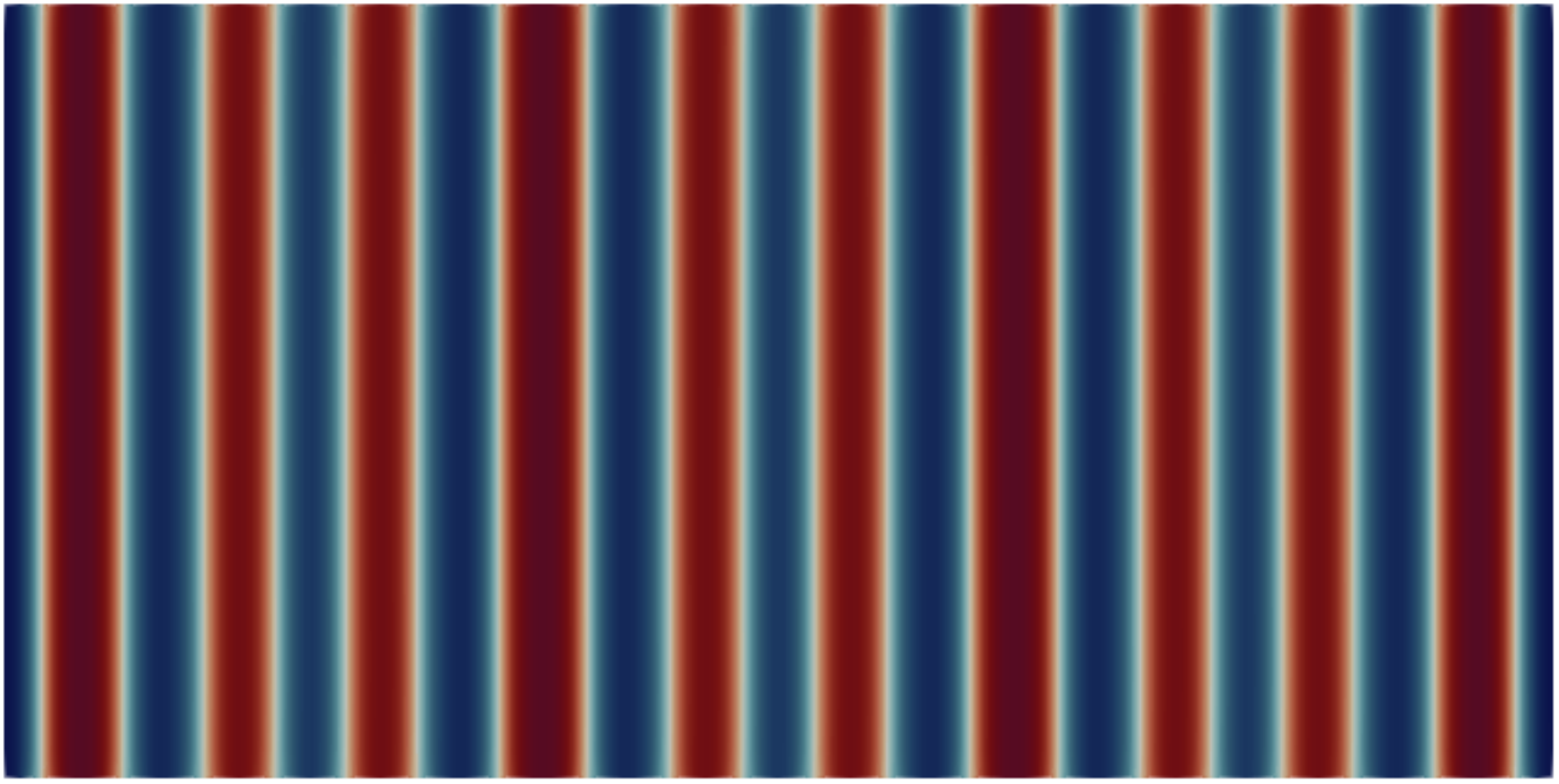}
    \caption{Sample state $\alpha=10^{3}$}
\end{subfigure}
\begin{subfigure}{0.32\textwidth}
    \centering
    \includegraphics[width=0.95\textwidth]{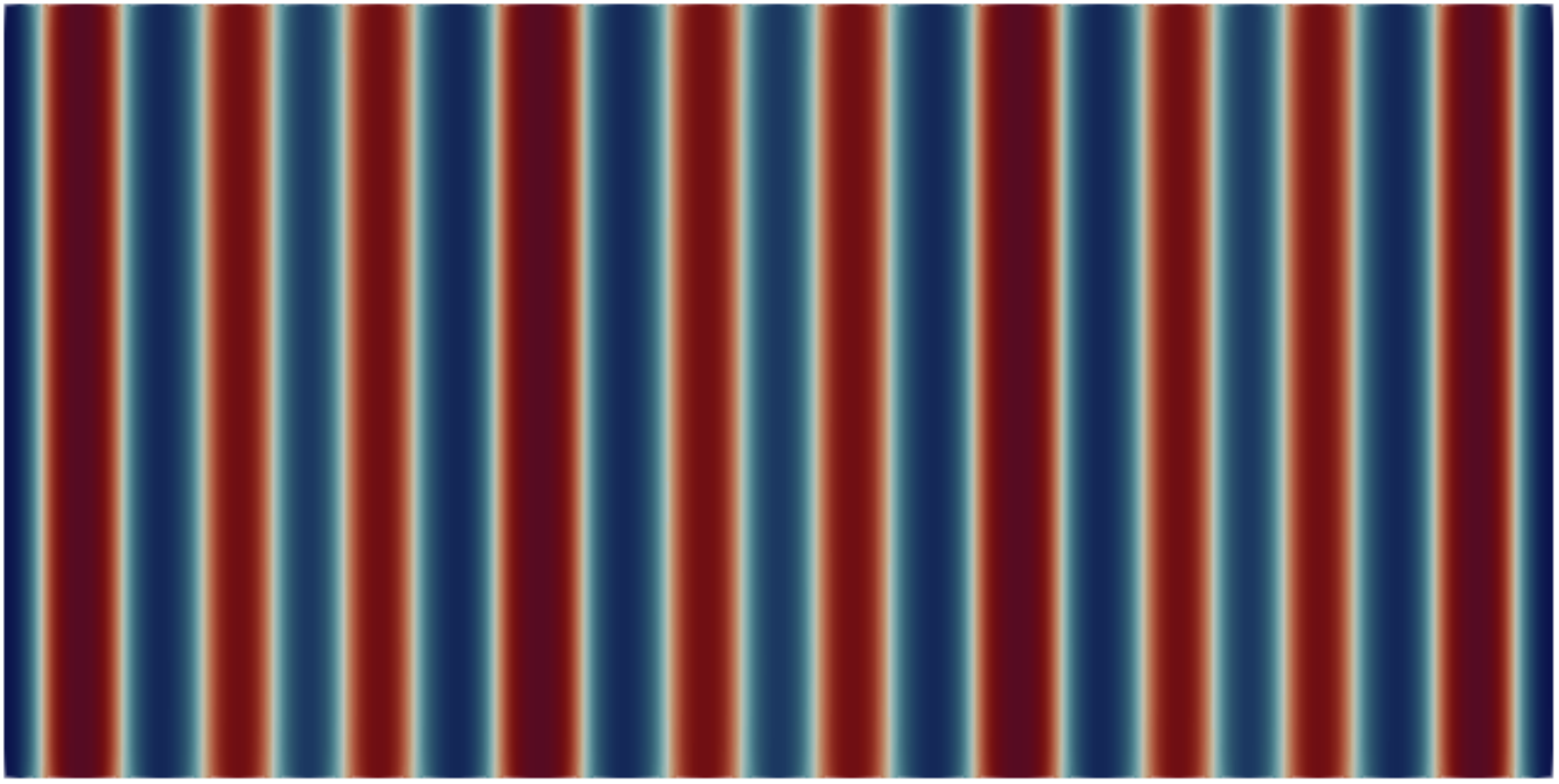}
    \caption{Sample state $\alpha=10^{4}$}
\end{subfigure}

\begin{subfigure}{0.32\textwidth}
    \centering
    \includegraphics[width=0.95\textwidth]{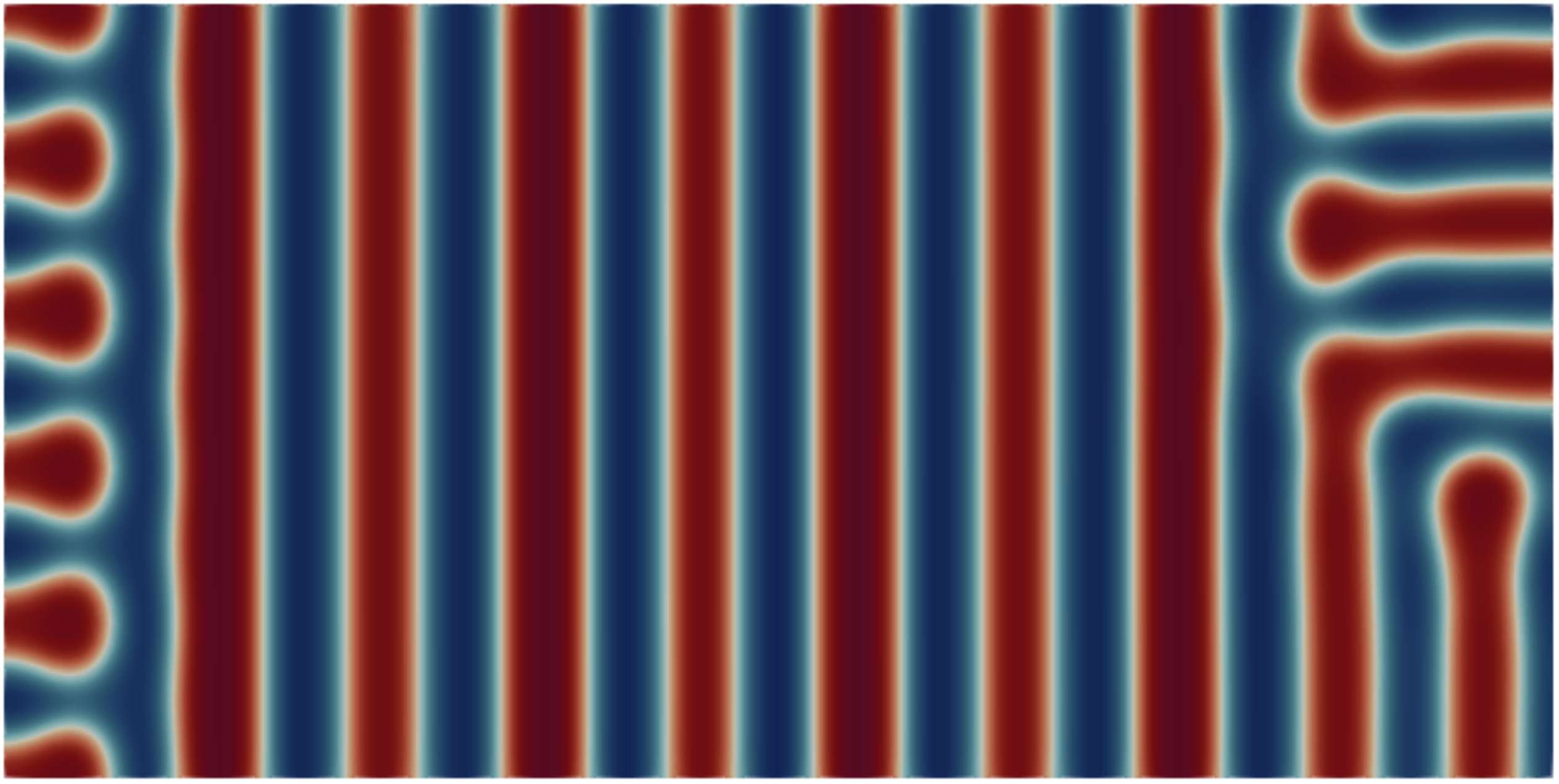}
    \caption{Sample state $\alpha=10^{2}$}
\end{subfigure}
\begin{subfigure}{0.32\textwidth}
    \centering
    \includegraphics[width=0.95\textwidth]{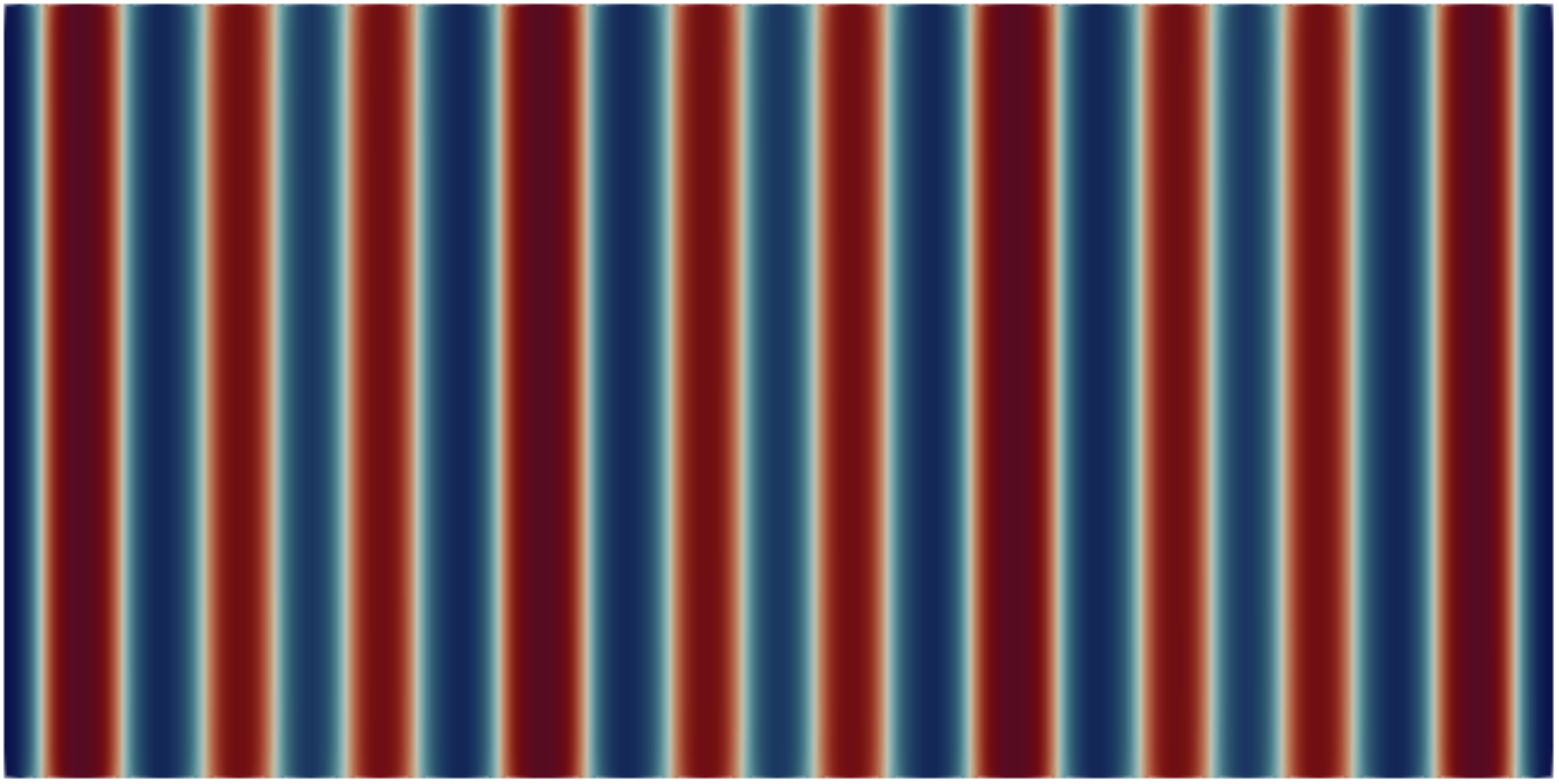}
    \caption{Sample state $\alpha=10^{3}$}
\end{subfigure}
\begin{subfigure}{0.32\textwidth}
    \centering
    \includegraphics[width=0.95\textwidth]{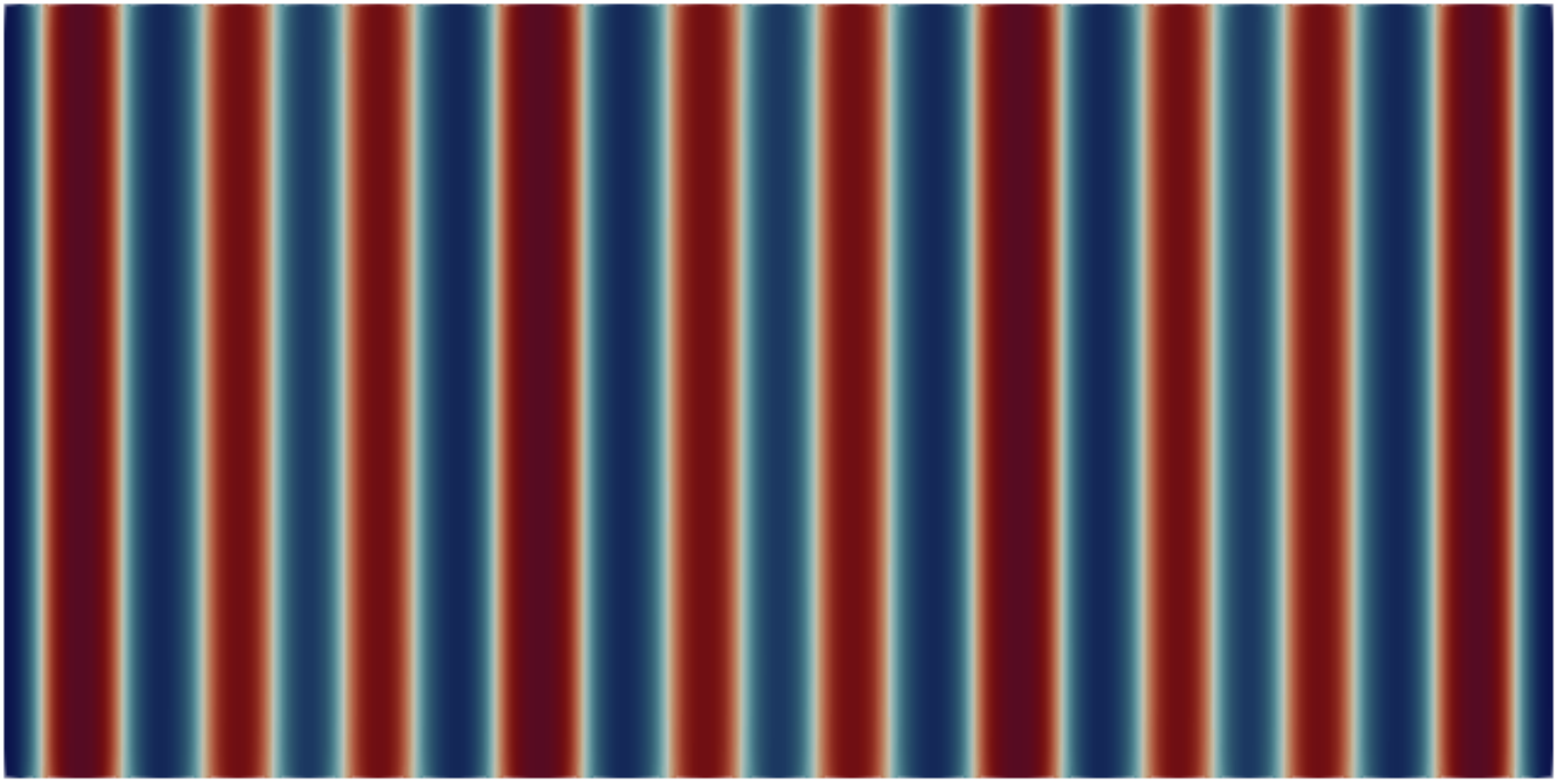}
    \caption{Sample state $\alpha=10^{4}$}
\end{subfigure}
% \caption{Optimal designs (top) for a strip target morphology using different values of the distance penalty parameter $\alpha = 10^{2}$ (left), $\alpha = 10^{3}$ (middle) and $\alpha = 10^{4}$ (right). The bottom three rows show three sample equilibrium states computed from different random initial guesses at the optimal designs, with the minimum energy states shown first. 
% % The target morphology is plotted for comparison.
% }
\caption{Sample sample equilibrium states computed from random initial guesses each using the optimal designs for the strip target morphology. Results are shown for different values of the distance penalty parameter $\alpha = 10^{2}$ (left), $\alpha = 10^{3}$ (middle) and $\alpha = 10^{4}$ (right). The minimum energy states are shown in the top row.
% The target morphology is plotted for comparison.
}
\label{fig:strip_samples}
\end{figure}

We set the initial locations of the four guideposts as randomly distributed with uniform distribution in the domain and carry out the optimization. We consider three values of the penalty parameter in the repelling term, $\alpha = 10^{2}$, $10^{3}$, and $10^{4}$, to investigate its effect. In Figure \ref{fig:strip_results}, we show the resulting optimal designs (with substrate strength in absolute values) and the optimal states---states computed from the optimal designs at the final iteration of the optimization algorithm using the continuation scheme. As a comparison, we show in Figure \ref{fig:strip_samples} the states with minimum free energy among the sample states computed for the optimal design, using randomly drawn initial guesses. Along with that, we also present a few of the other sample states in Figure \ref{fig:strip_samples} to give an indication of the possible defects.

% The resulting optimal designs (with substrate strength shown in absolute values) and sample states are shown in Figure \ref{fig:strip_results}. 
For $\alpha = 10^2$, the optimal guideposts are located with $2\times$ density multiplication; i.e., skipping one target strip, while for $\alpha = 10^{3}$ and $10^{4}$, the optimal guideposts are located with $3\times$ density multiplication; i.e., skipping two target strips. 
However, the sample states for the design with $\alpha = 10^{2}$ are defective, despite the optimal state (computed by the continuation scheme during optimization) being defect free. On the other hand, the optimal designs in the latter two cases are robust to different initial guesses, for which the sample equilibrium states are all defect-free.
% For $\alpha=10^{2}$, the sample equilibrium states at the optimal design expose defects. 
% \pc{Do the defects exist at the converged state in the optimization?}{}

As previously shown in Figure \ref{fig:2d_strip_obj}, the objective function exhibits several multiple local minima as the guideposts move in and out of alignment with the target morphology. Therefore, the penalization term is important to producing well separated guideposts and achieving higher levels of density multiplication. This helps the optimizer to navigate through these local minima and converge to more robust designs. 

% In this problem, there are several local minima of the design objective corresponding to instances when the guideposts are located at some of the target strips. \pc{can you show the local minima as you showed in one presentation?} The results demonstrate that stronger penalization helps the optimizer to navigate through these local minima and converge to more robust designs. 

% To investigate the influence of the penalty parameter $\alpha$ in the repelling term \eqref{eq:repel_strip} for the robustness of the optimal design, we take $\alpha = 10^{2}$, $10^{3}$, and $10^{4}$. The robustness is assessed by solving the forward problem at the converged optimal designs starting from random initial guesses drawn independently from that used for the optimization.
% Beginning with a random set of initial guidepost locations drawn from a uniform distribution on the substrate domain, 

% The optimal designs (with substrate strength shown in absolute values) and sample states are shown in Figure \ref{fig:strip_results}. 
% For simplicity, we plot the magnitude of the substrate interaction field for the optimal designs instead of its value, since we have used $f < 0$ to attract species A.

\subsection{Designs with circular guideposts in 2D}
\begin{figure}[h]
\centering
\begin{subfigure}{0.19\textwidth}
    \centering
    \includegraphics[width=0.9\textwidth]{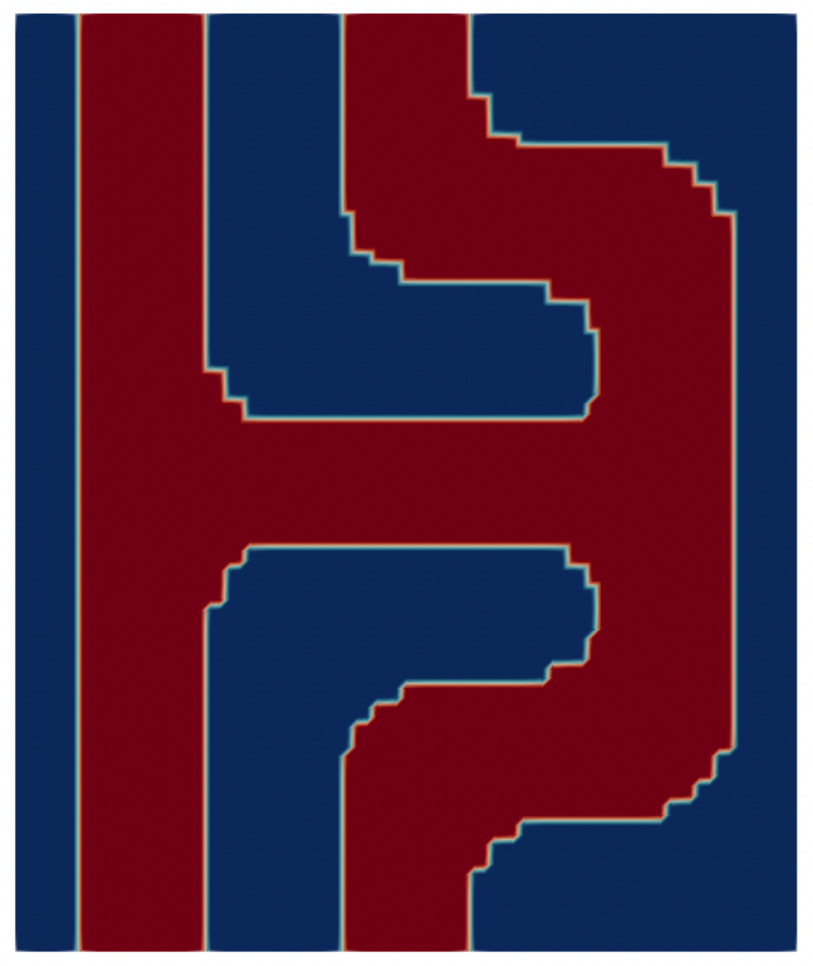}
    \caption{Junction}
    \label{fig:junction2d_target}
\end{subfigure}
\begin{subfigure}{0.19\textwidth}
    \centering
    \includegraphics[width=0.9\textwidth]{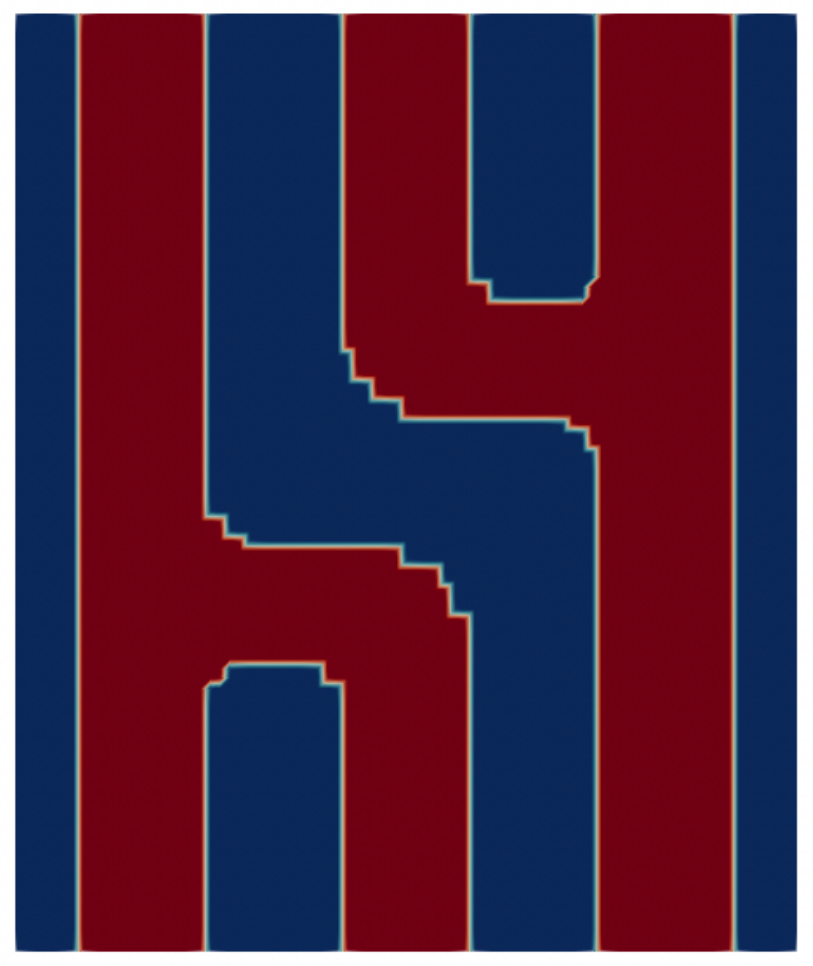}
    \caption{Jog}
    \label{fig:jog2d_target}
\end{subfigure}
\caption{Target morphologies in 2D, defined over $\Omega = [0, L]\times[0, 7L/6]$ with $L=3$.}
\label{fig:2d_targets}
\end{figure}

\begin{figure}
    \centering
    \captionsetup[subfigure]{justification=centering}
    \begin{subfigure}{0.32\textwidth}
        \centering
        \includegraphics[width=1.0\textwidth]{figures_pdf/design_colors.pdf}
    \end{subfigure}

    \begin{subfigure}{0.19\textwidth}
        \centering
        \includegraphics[width=0.9\textwidth]{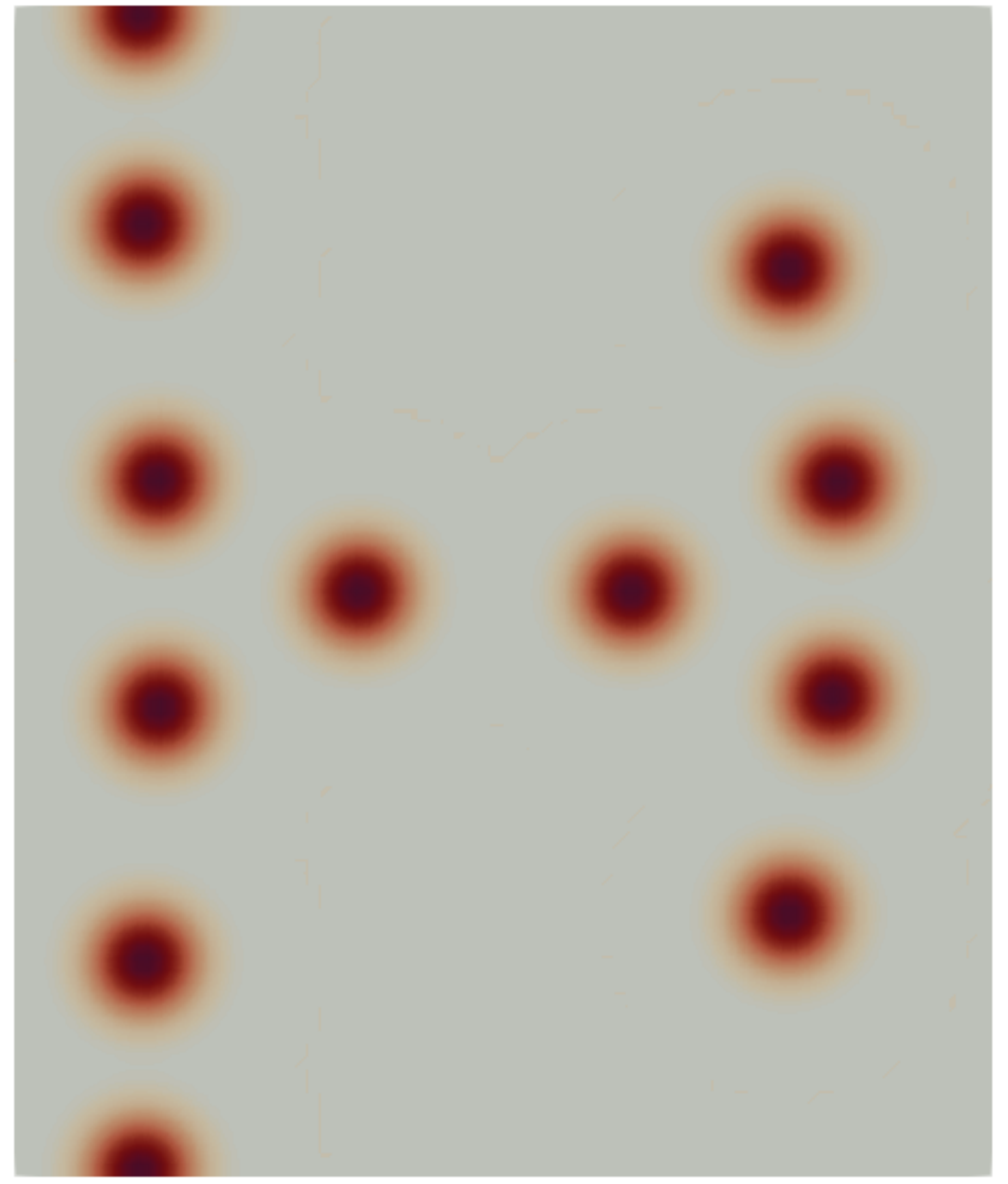}
        \caption{Optimal design $N_p = 12$}
    \end{subfigure}
    \begin{subfigure}{0.19\textwidth}
        \centering
        \includegraphics[width=0.9\textwidth]{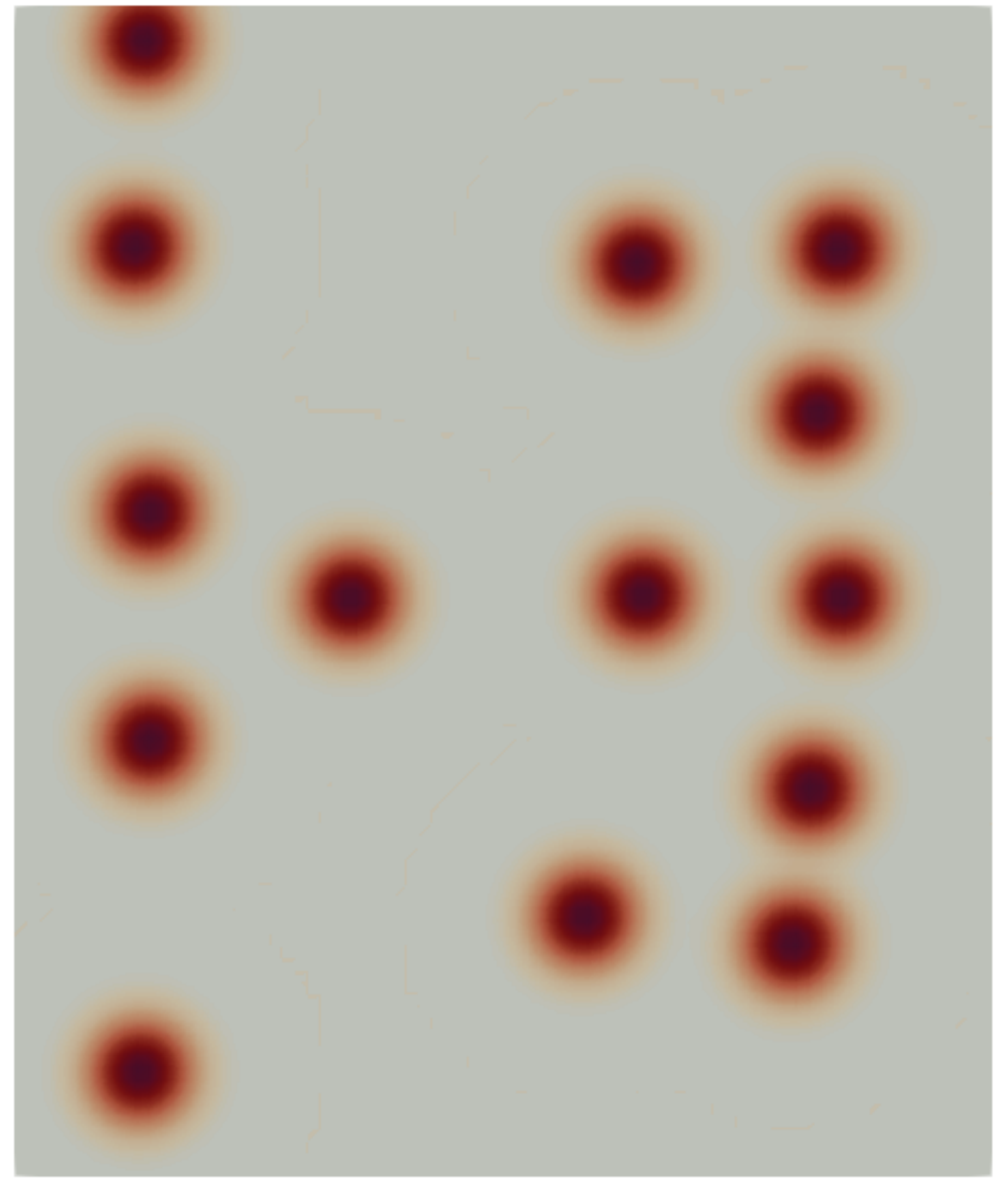}
        \caption{Optimal design $N_p = 14$}
    \end{subfigure}
    \begin{subfigure}{0.19\textwidth}
        \centering
        \includegraphics[width=0.9\textwidth]{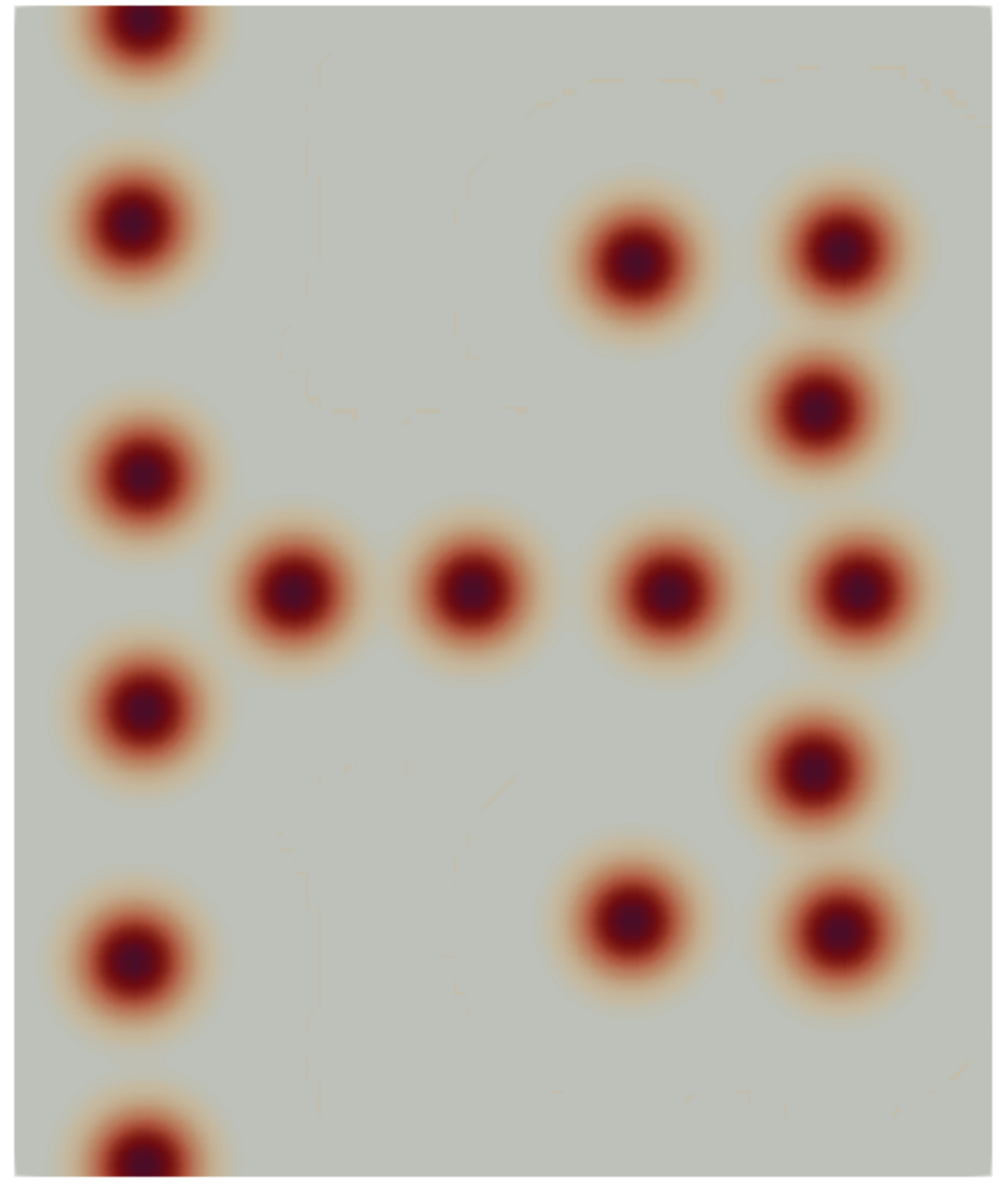}
        \caption{Optimal design $N_p = 16$}
    \end{subfigure}
    \begin{subfigure}{0.19\textwidth}
        \centering
        \includegraphics[width=0.9\textwidth]{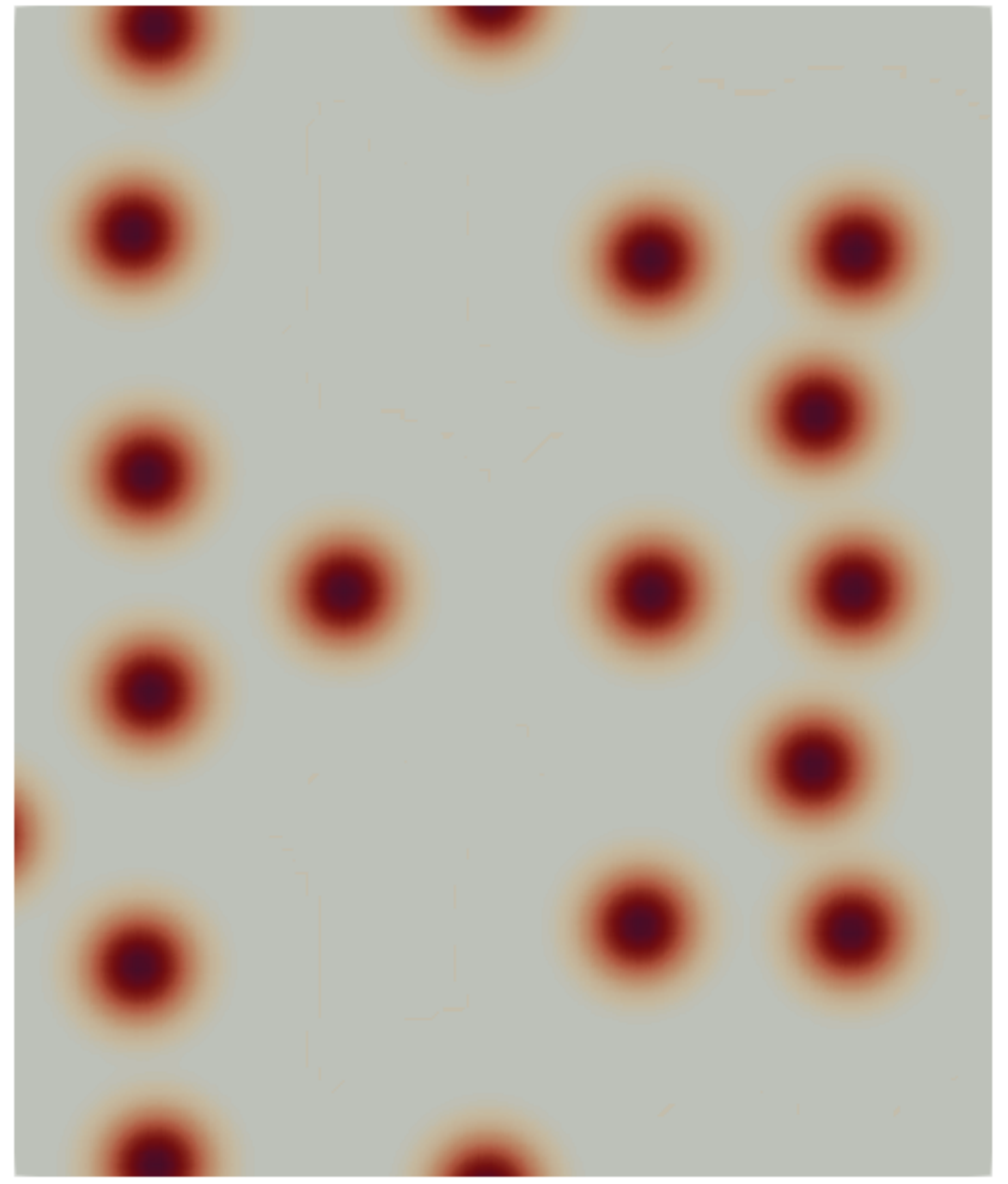}
        \caption{Optimal design $N_p = 16$}
    \end{subfigure}
    \begin{subfigure}{0.19\textwidth}
        \centering
        \includegraphics[width=0.9\textwidth]{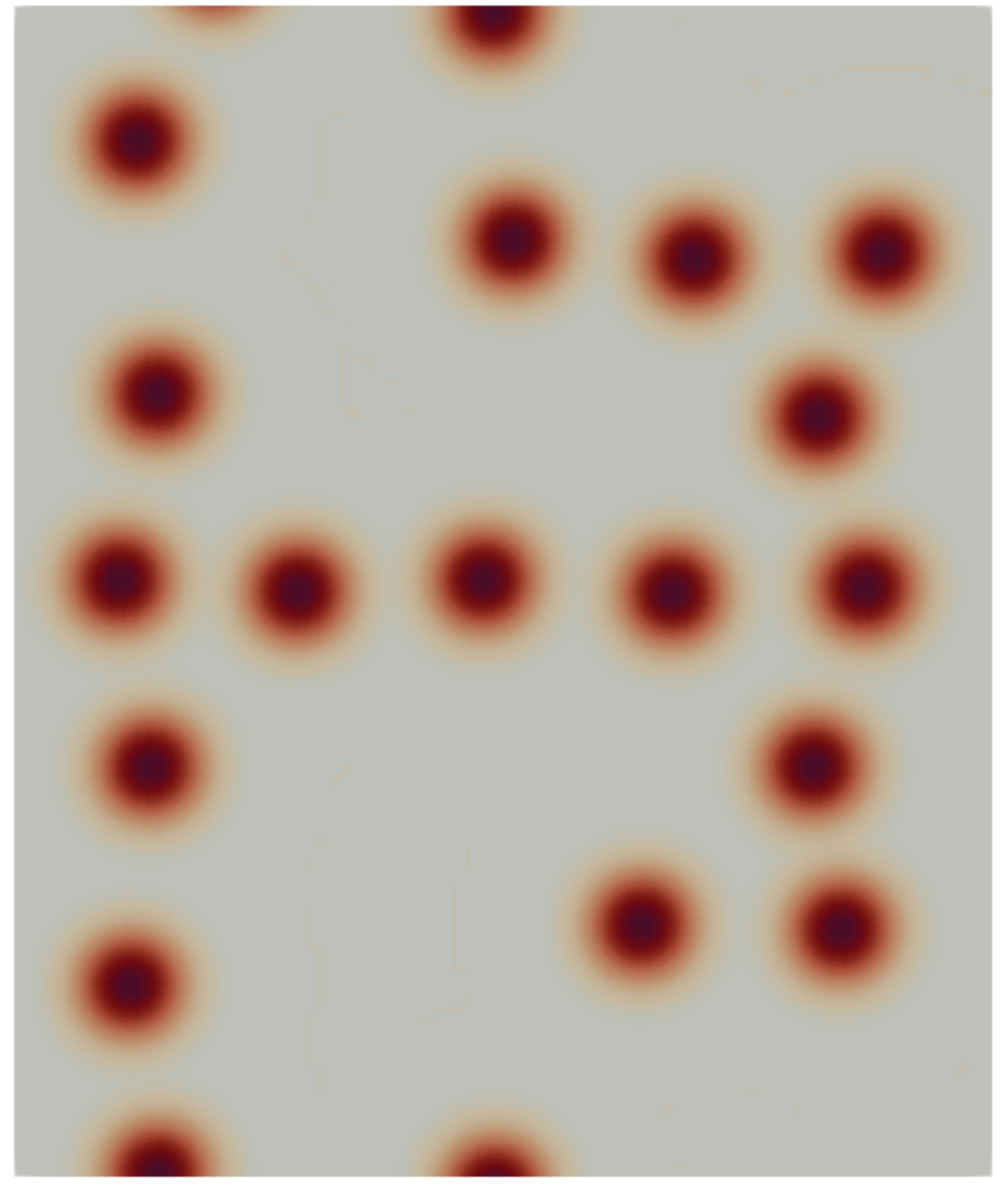}
        \caption{Optimal design $N_p = 20$}
    \end{subfigure}

    \begin{subfigure}{0.32\textwidth}
        \centering
        \includegraphics[width=1.0\textwidth]{figures_pdf/state_colors.pdf}
    \end{subfigure}

    \begin{subfigure}{0.19\textwidth}
        \centering
        \includegraphics[width=0.9\textwidth]{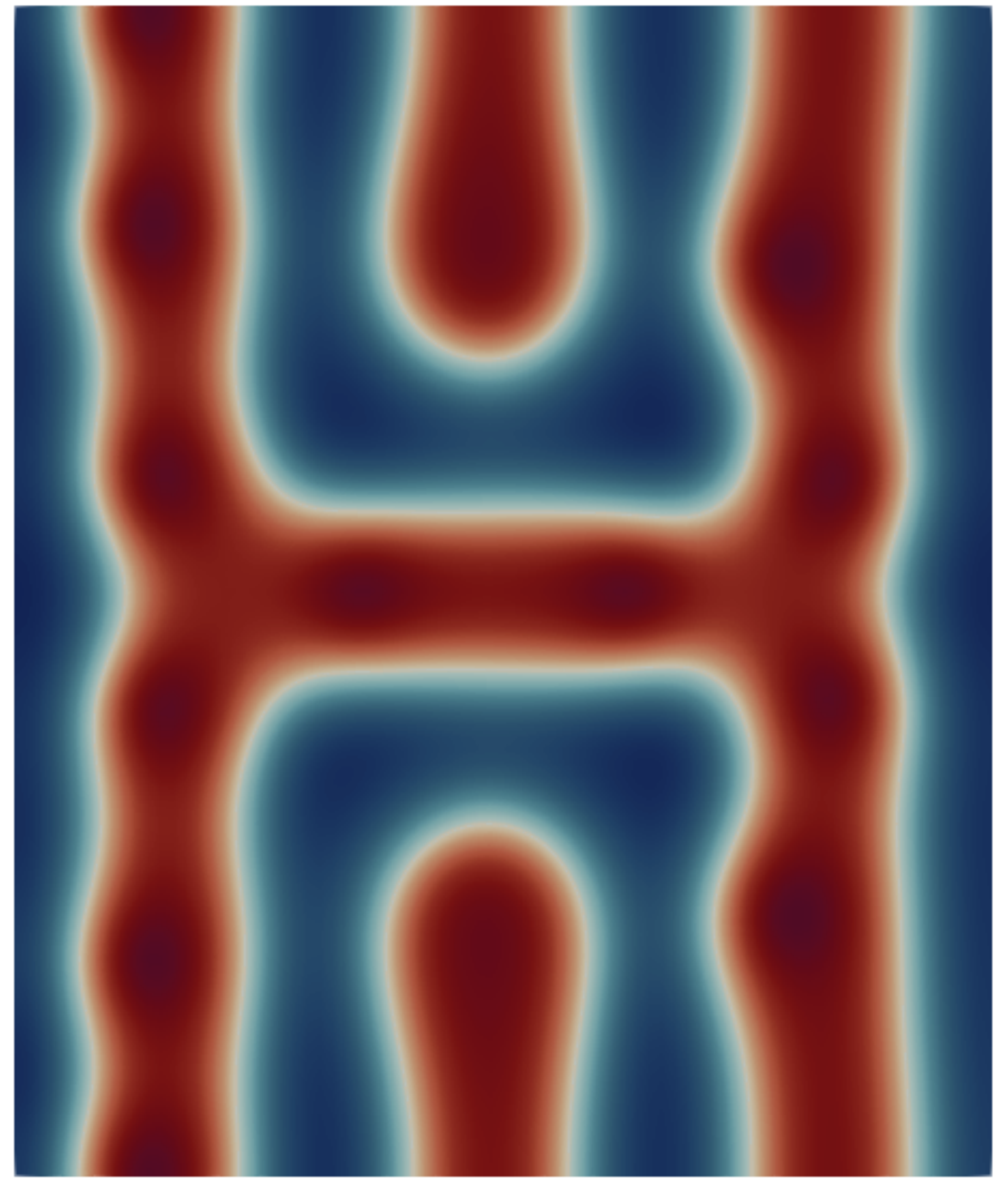}
        \caption{Minimum $\mathcal{F}(u)$ sample $N_p = 12$}
    \end{subfigure}
    \begin{subfigure}{0.19\textwidth}
        \centering
        \includegraphics[width=0.9\textwidth]{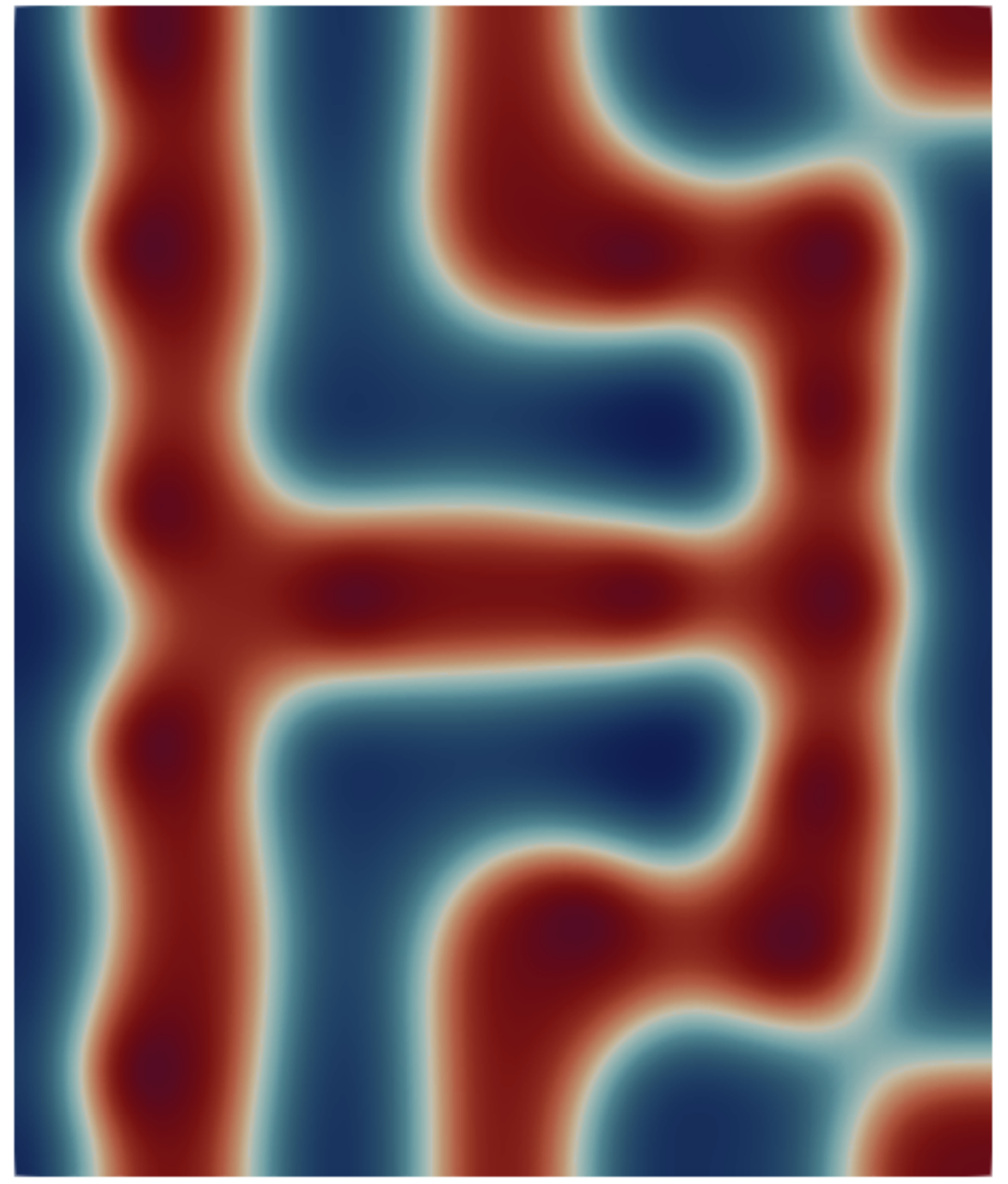}
        \caption{Minimum $\mathcal{F}(u)$ sample $N_p = 14$}
    \end{subfigure}
    \begin{subfigure}{0.19\textwidth}
        \centering
        \includegraphics[width=0.9\textwidth]{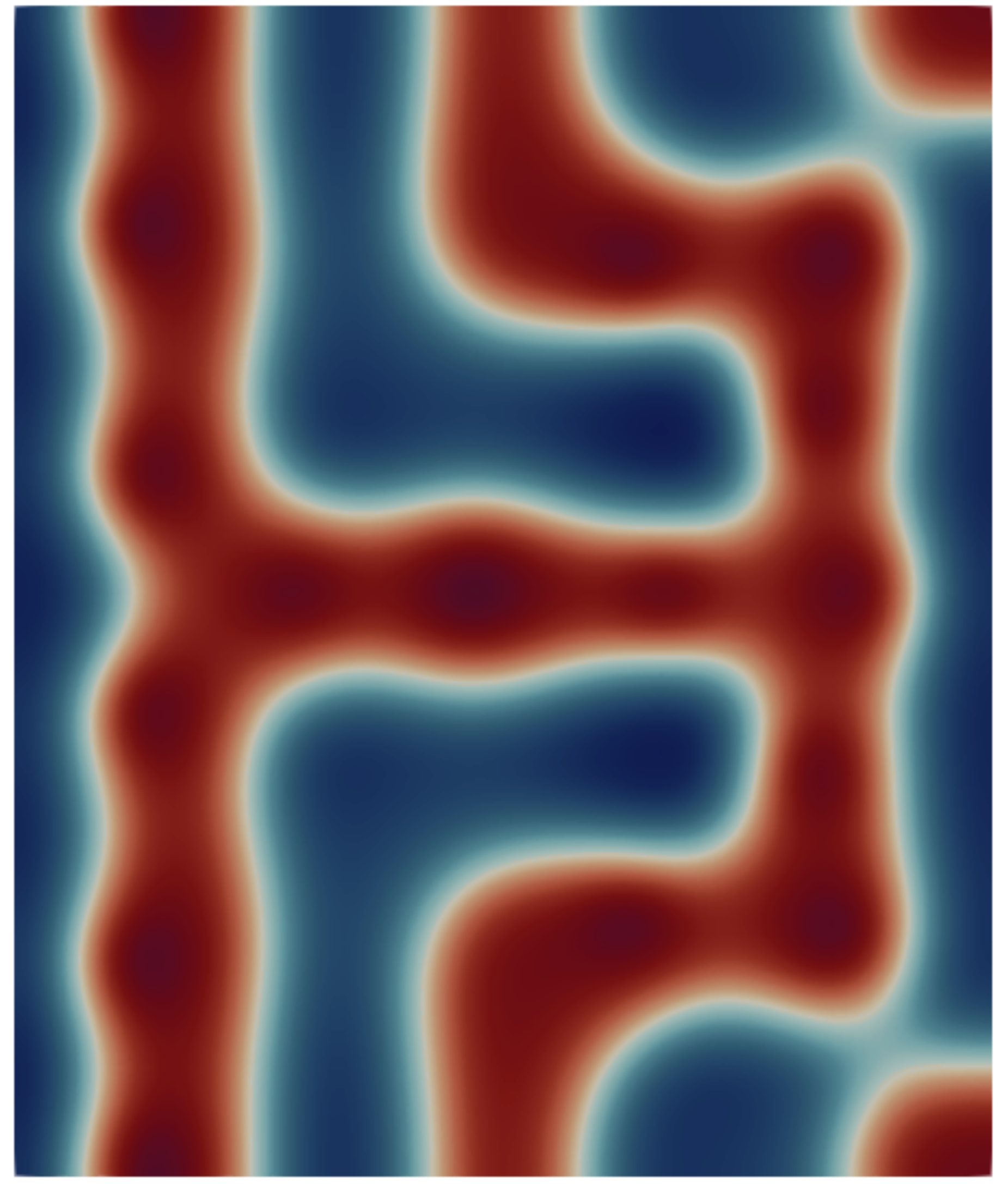}
        \caption{Minimum $\mathcal{F}(u)$ sample $N_p = 16$}
    \end{subfigure}
    \begin{subfigure}{0.19\textwidth}
        \centering
        \includegraphics[width=0.9\textwidth]{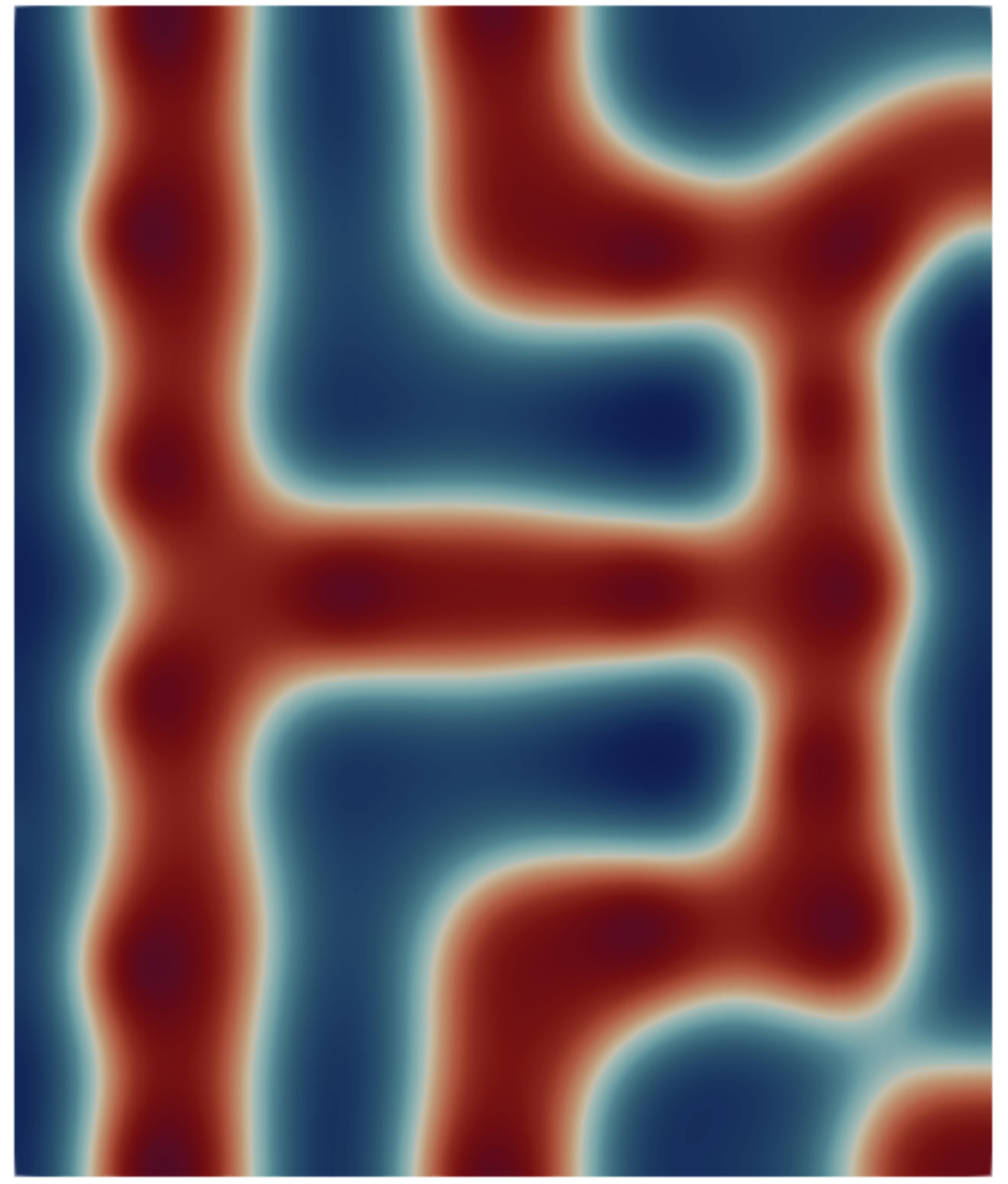}      
        \caption{Minimum $\mathcal{F}(u)$ sample $N_p = 16$}
    \end{subfigure}
    \begin{subfigure}{0.19\textwidth}
        \centering
        \includegraphics[width=0.9\textwidth]{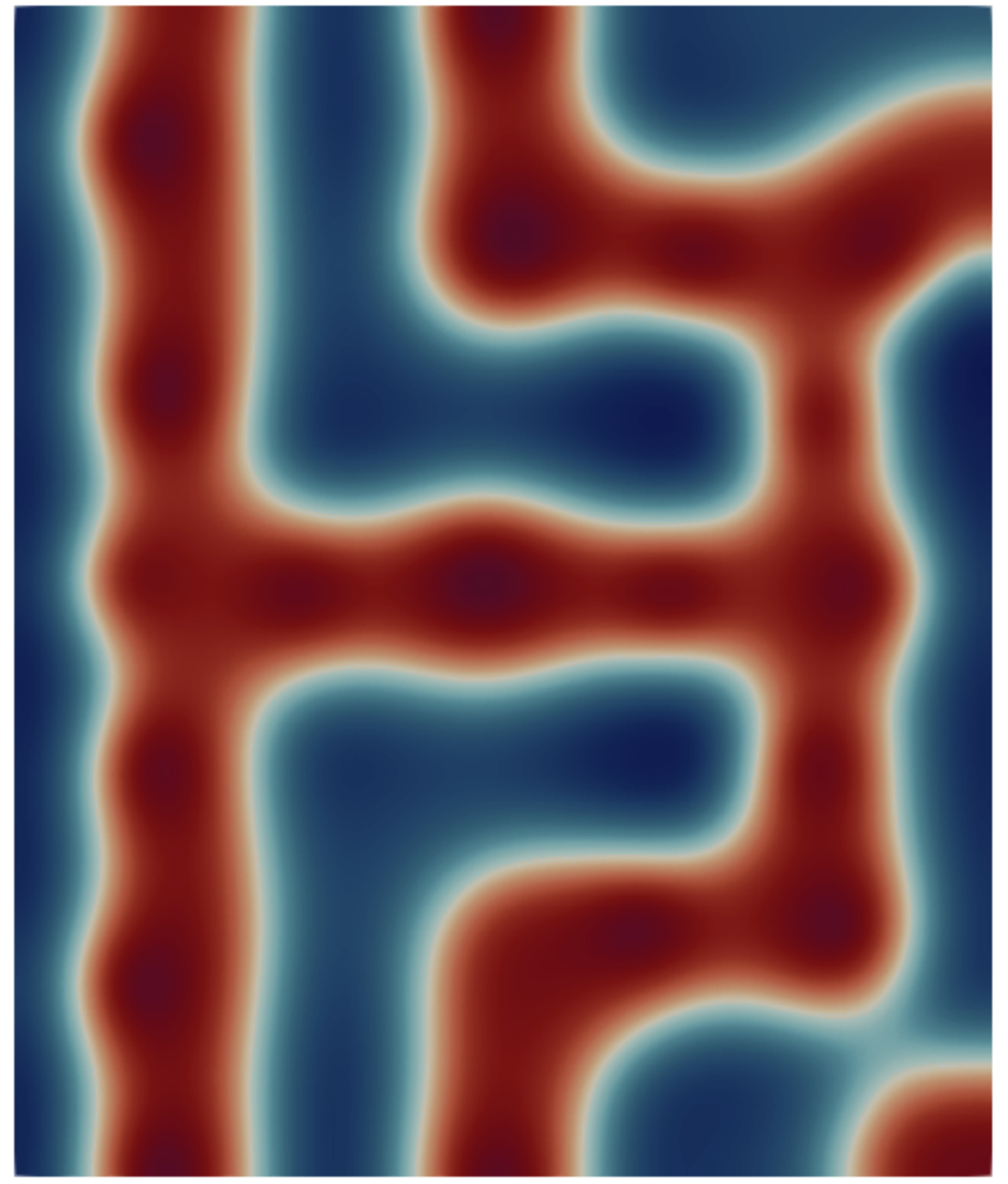}
        \caption{Minimum $\mathcal{F}(u)$ sample $N_p = 20$}
    \end{subfigure}

    \begin{subfigure}{0.19\textwidth}
        \centering
        \includegraphics[width=0.9\textwidth]{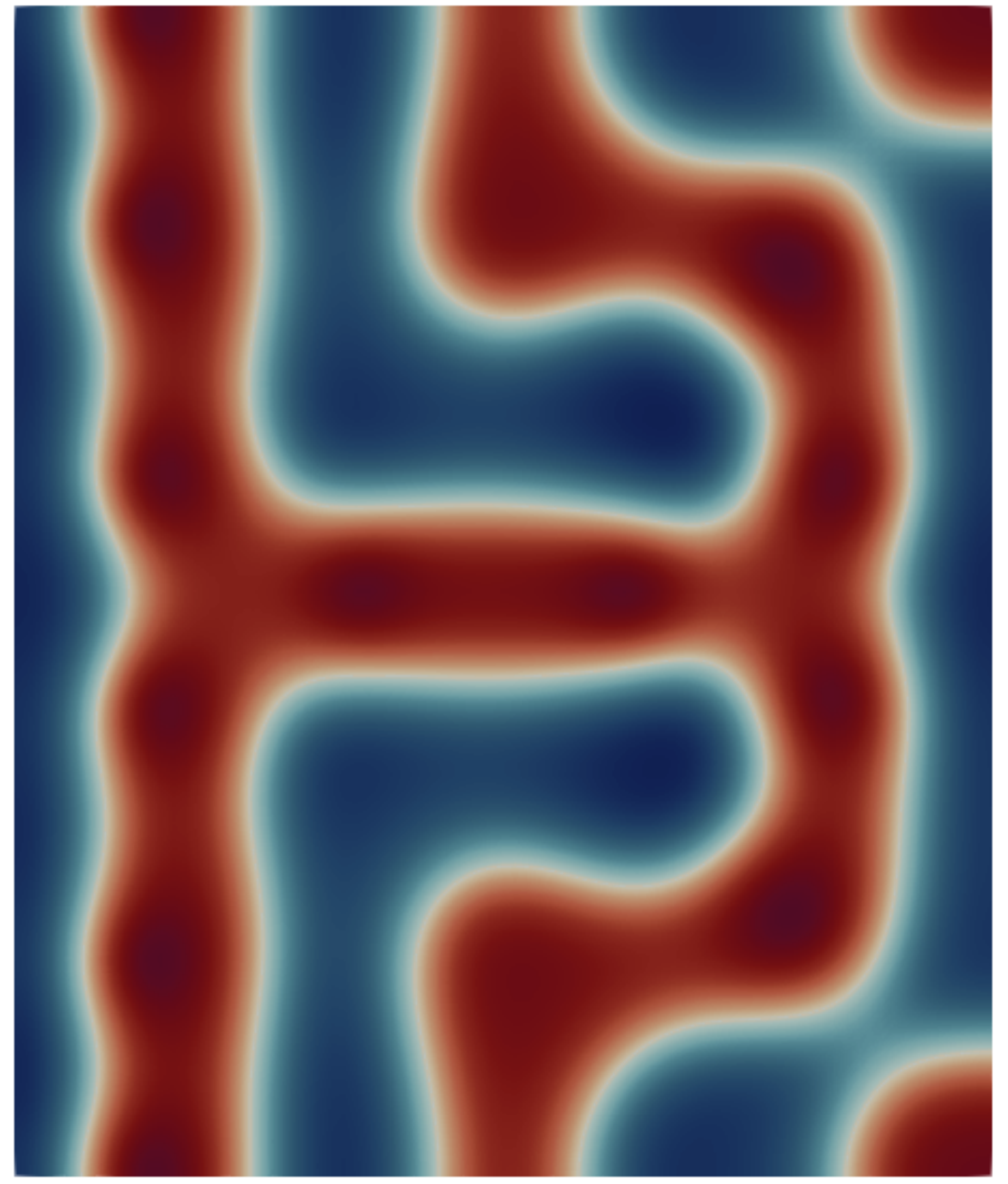}
        \caption{Sample 2 $N_p = 12$}
    \end{subfigure}
    \begin{subfigure}{0.19\textwidth}
        \centering
        \includegraphics[width=0.9\textwidth]{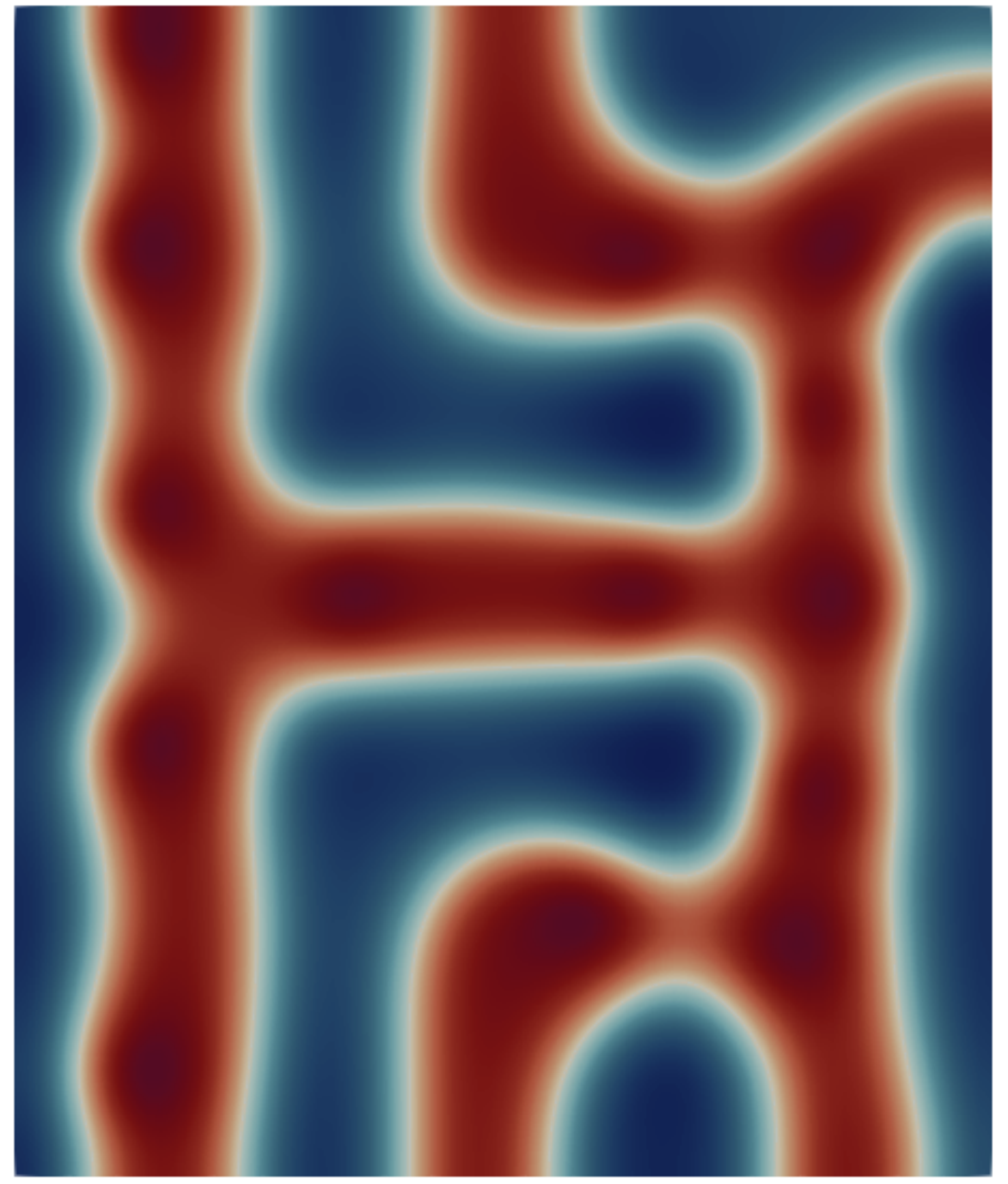}
        \caption{Sample 2 $N_p = 14$}
    \end{subfigure}
    \begin{subfigure}{0.19\textwidth}
        \centering
        \includegraphics[width=0.9\textwidth]{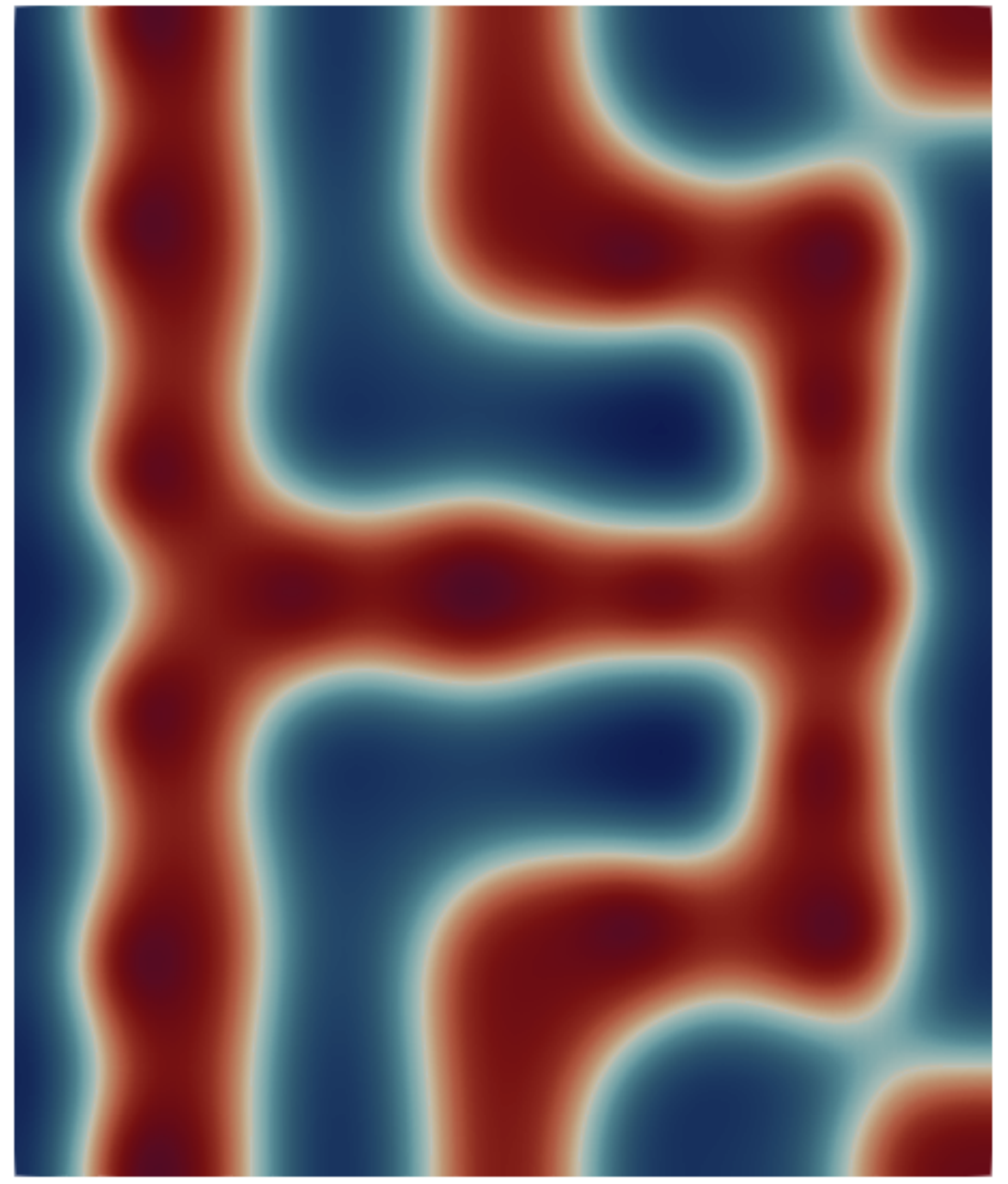}
        \caption{Sample 2 $N_p = 16$}
    \end{subfigure}
    \begin{subfigure}{0.19\textwidth}
        \centering
        \includegraphics[width=0.9\textwidth]{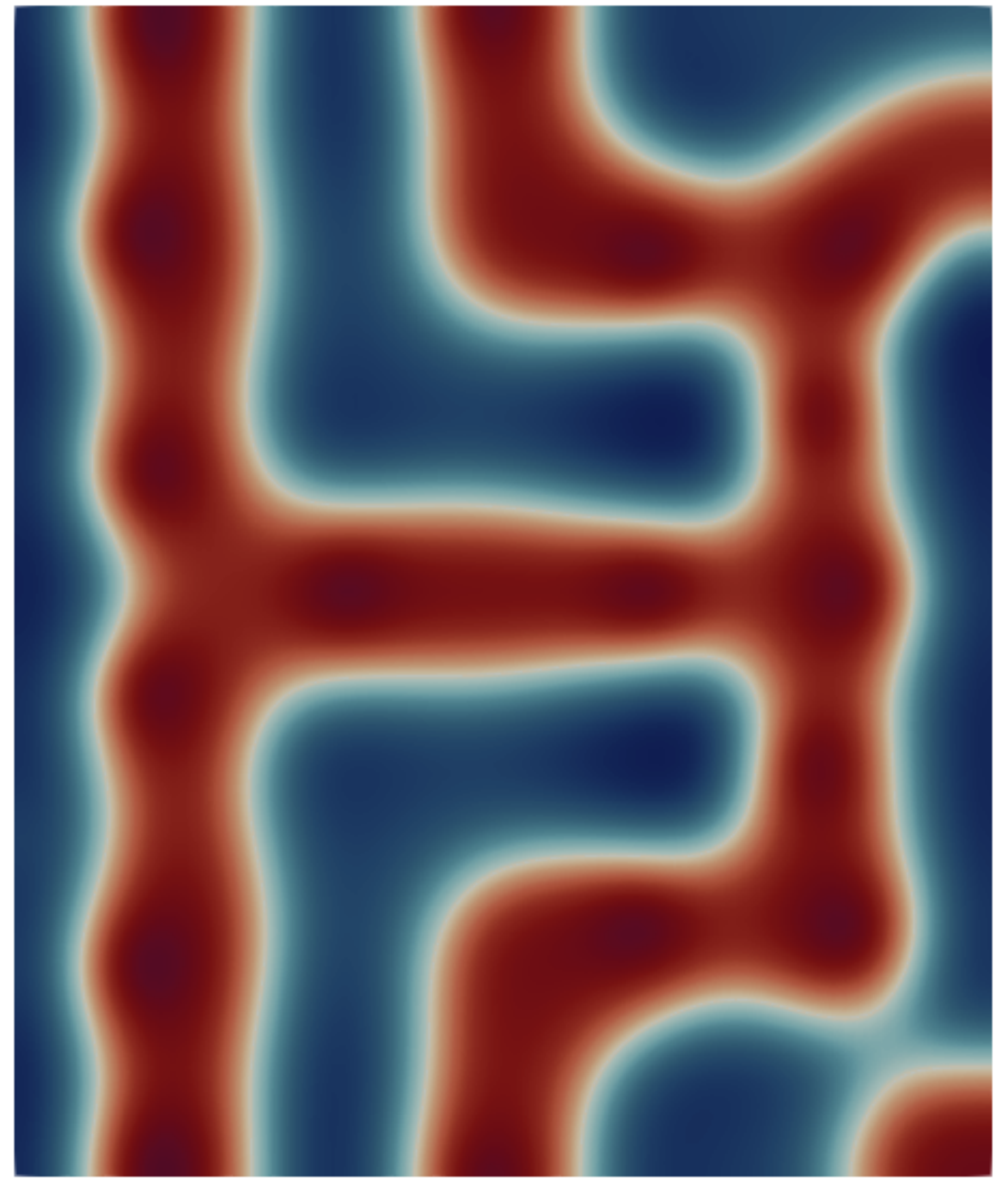}      
        \caption{Sample 2 $N_p = 16$}
    \end{subfigure}
    \begin{subfigure}{0.19\textwidth}
        \centering
        \includegraphics[width=0.9\textwidth]{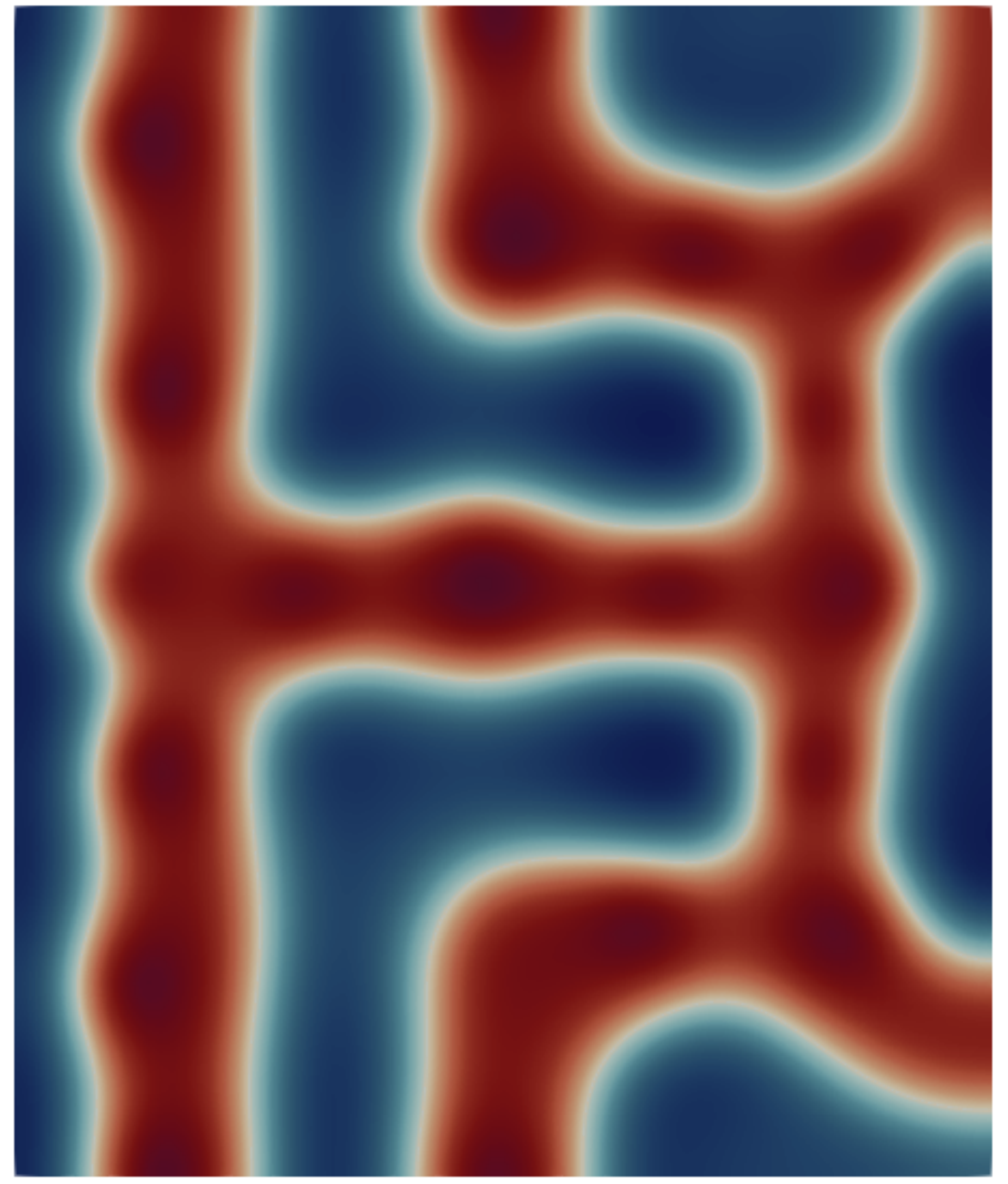}
        \caption{Sample 2 $N_p = 20$}
    \end{subfigure}

    \begin{subfigure}{0.19\textwidth}
        \centering
        \includegraphics[width=0.9\textwidth]{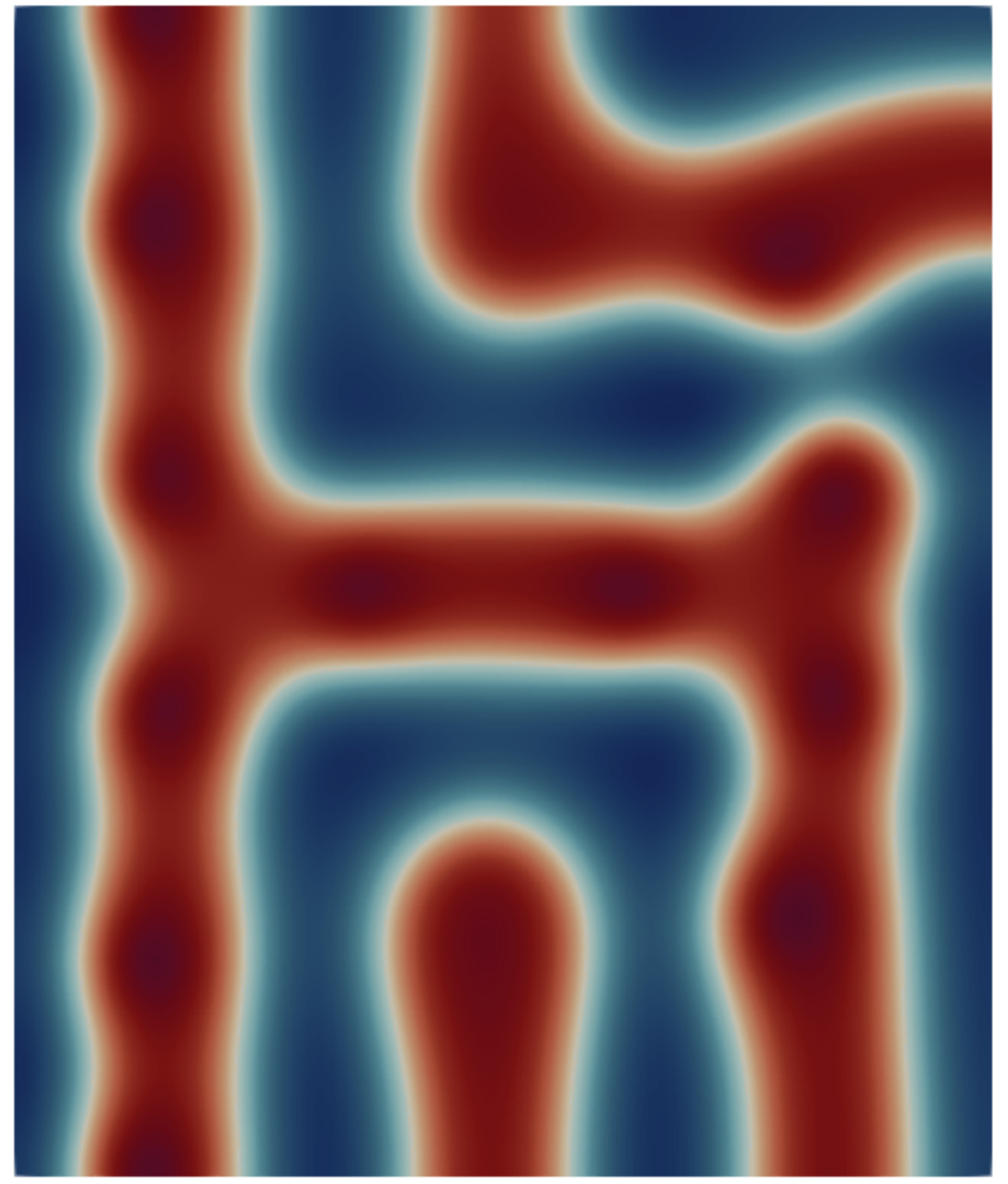}
        \caption{Sample 3 $N_p = 12$}
    \end{subfigure}
    \begin{subfigure}{0.19\textwidth}
        \centering
        \includegraphics[width=0.9\textwidth]{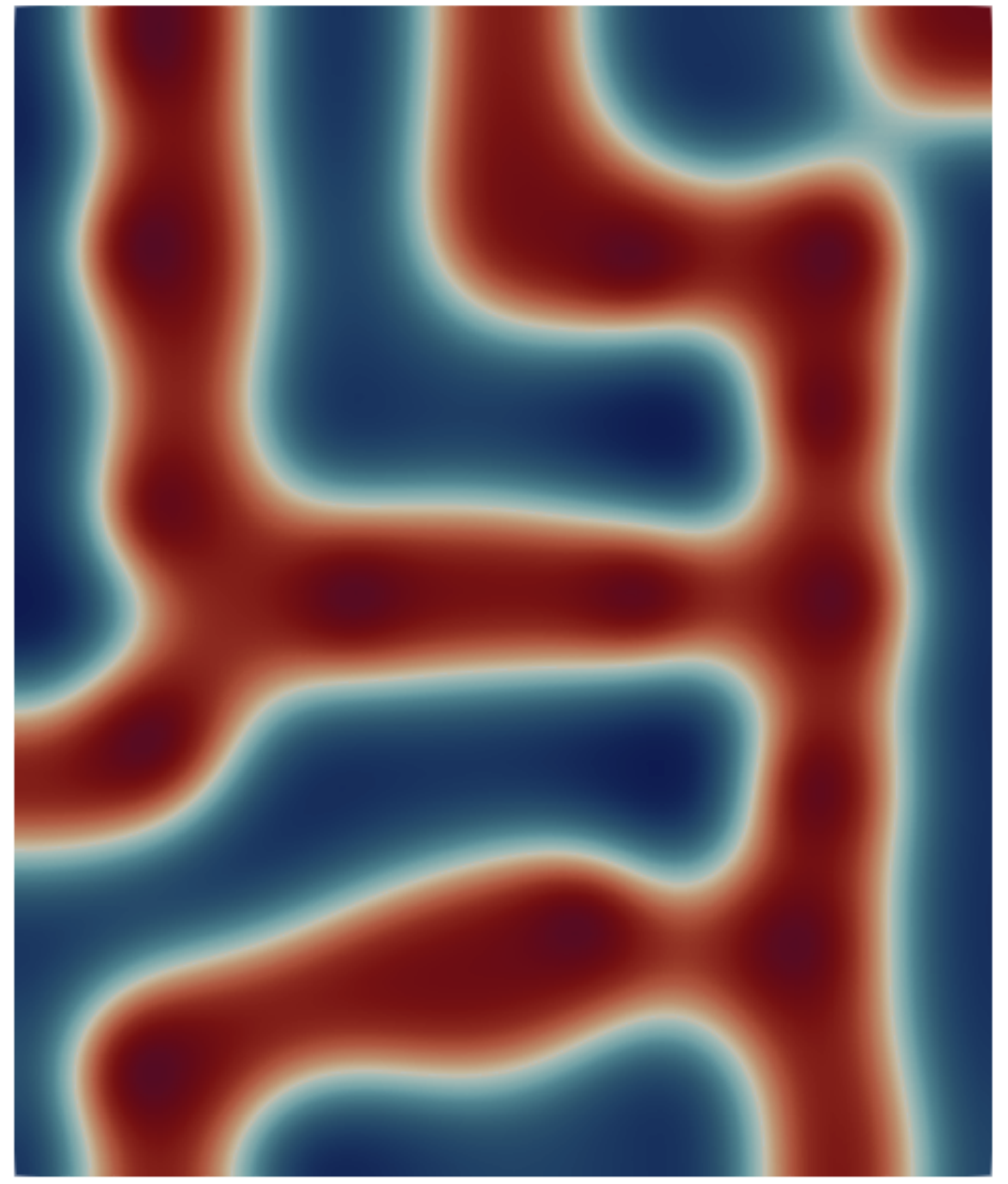}
        \caption{Sample 3 $N_p = 14$}
    \end{subfigure}
    \begin{subfigure}{0.19\textwidth}
        \centering
        \includegraphics[width=0.9\textwidth]{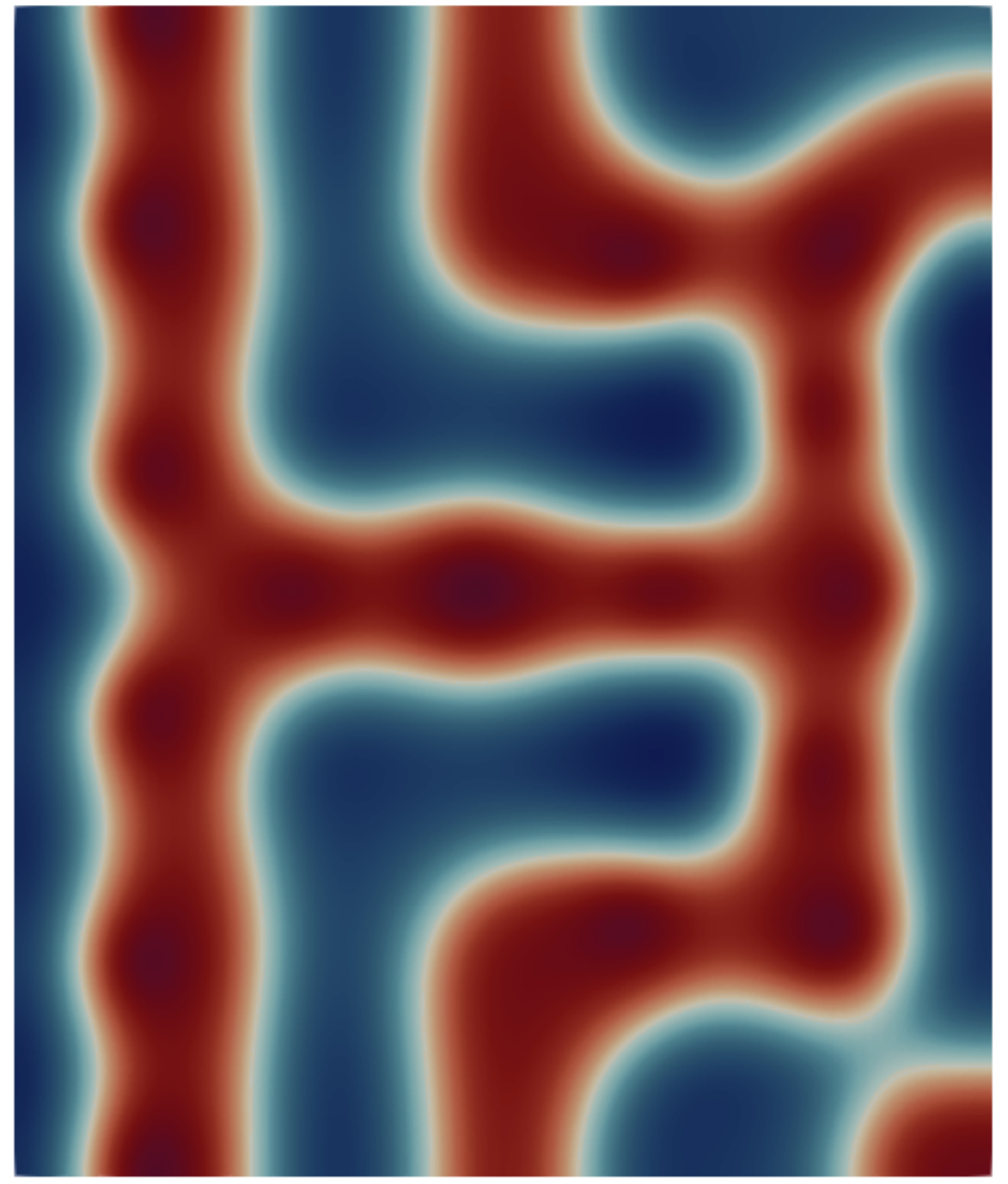}
        \caption{Sample 3 $N_p = 16$}
    \end{subfigure}
    \begin{subfigure}{0.19\textwidth}
        \centering
        \includegraphics[width=0.9\textwidth]{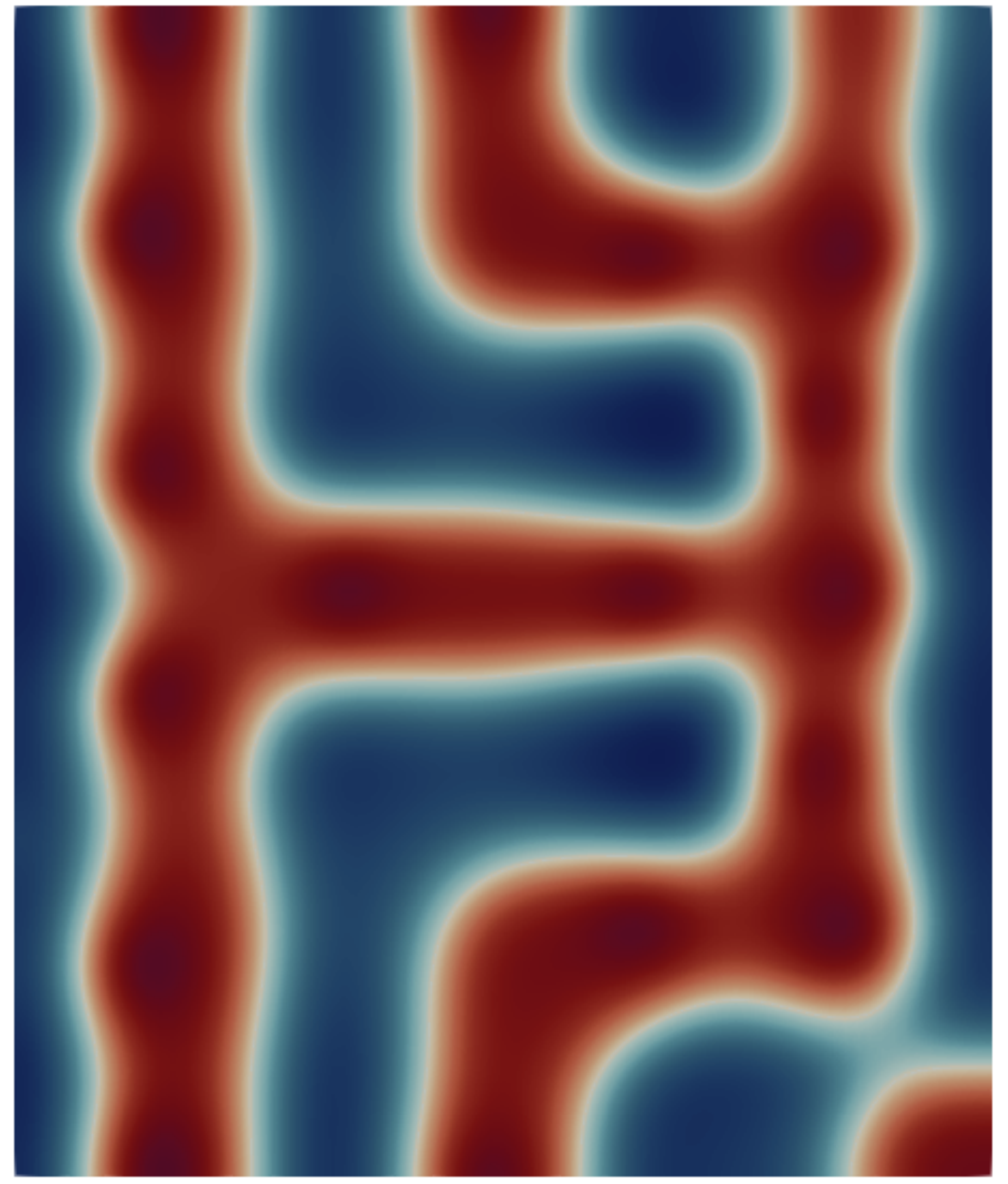}      
        \caption{Sample 3 $N_p = 16$}
    \end{subfigure}
    \begin{subfigure}{0.19\textwidth}
        \centering
        \includegraphics[width=0.9\textwidth]{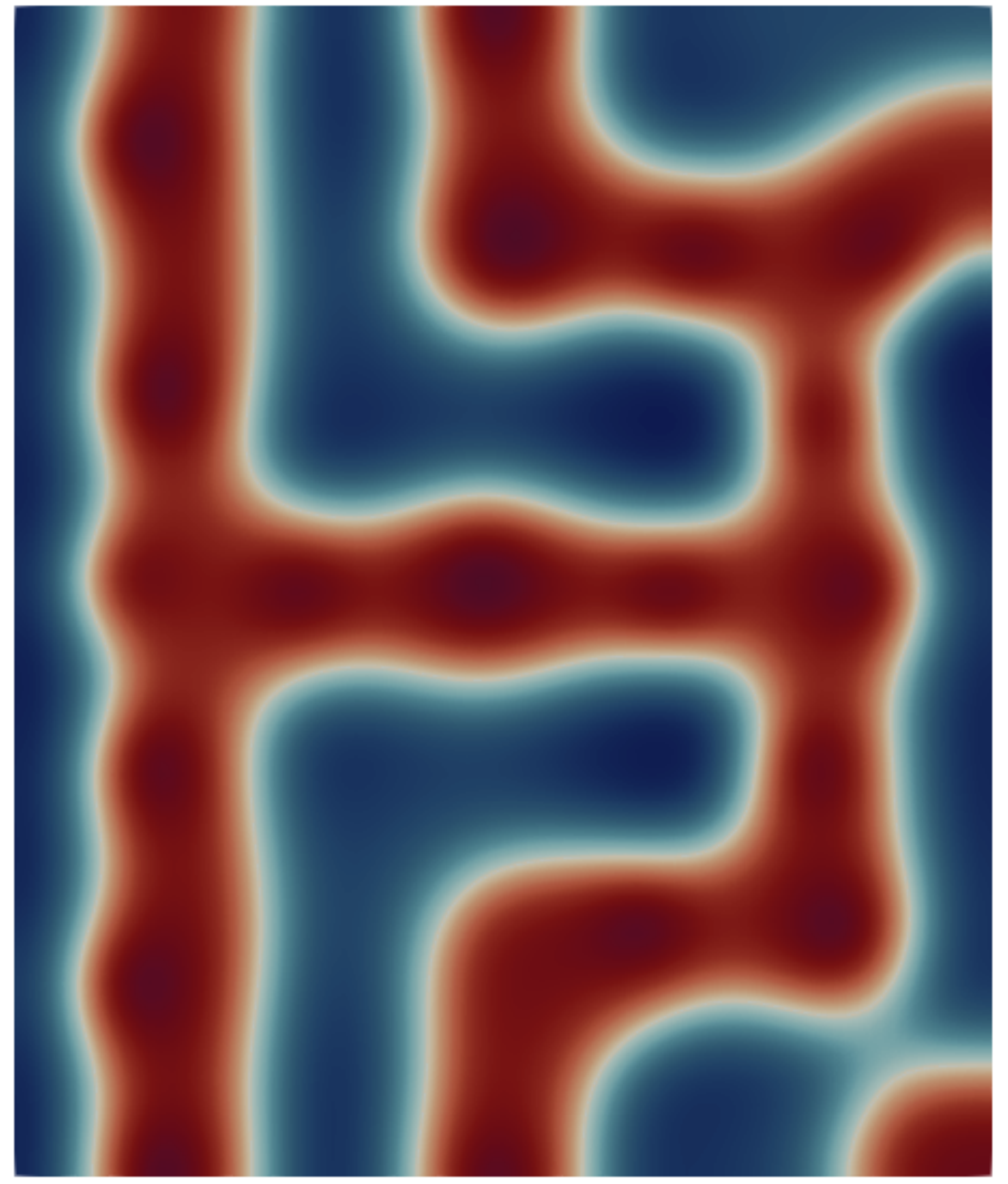}
        \caption{Sample 3 $N_p = 20$}
    \end{subfigure}
    \caption{Optimal design (top) for the junction target morphology using  $N_p = 12, 14, 16, 18, 20$ guideposts  and repulsion strength $\alpha = 5 \times 10^{-3}$. 
    % The bottom three rows show three sample equilibrium states corresponding to the optimal designs, computed from three different sample initial guesses. The minimum energy states are the samples with lowest free energy. 
    The bottom three rows show sample equilibrium states computed from different sample initial guesses $u_0$ at the optimal designs, with the minimum energy states shown first. 
    % The minimum energy states are the samples with lowest free energy.
    % In particular, at the top of this section are the samples states with the minimum free energy.
    }
    \label{fig:junction2d_Np}
\end{figure}

\begin{figure}
    \centering
    \captionsetup[subfigure]{justification=centering}
    \begin{subfigure}{0.32\textwidth}
        \centering
        \includegraphics[width=1.0\textwidth]{figures_pdf/design_colors.pdf}
    \end{subfigure}

    \begin{subfigure}{0.19\textwidth}
        \centering
        \includegraphics[width=0.9\textwidth]{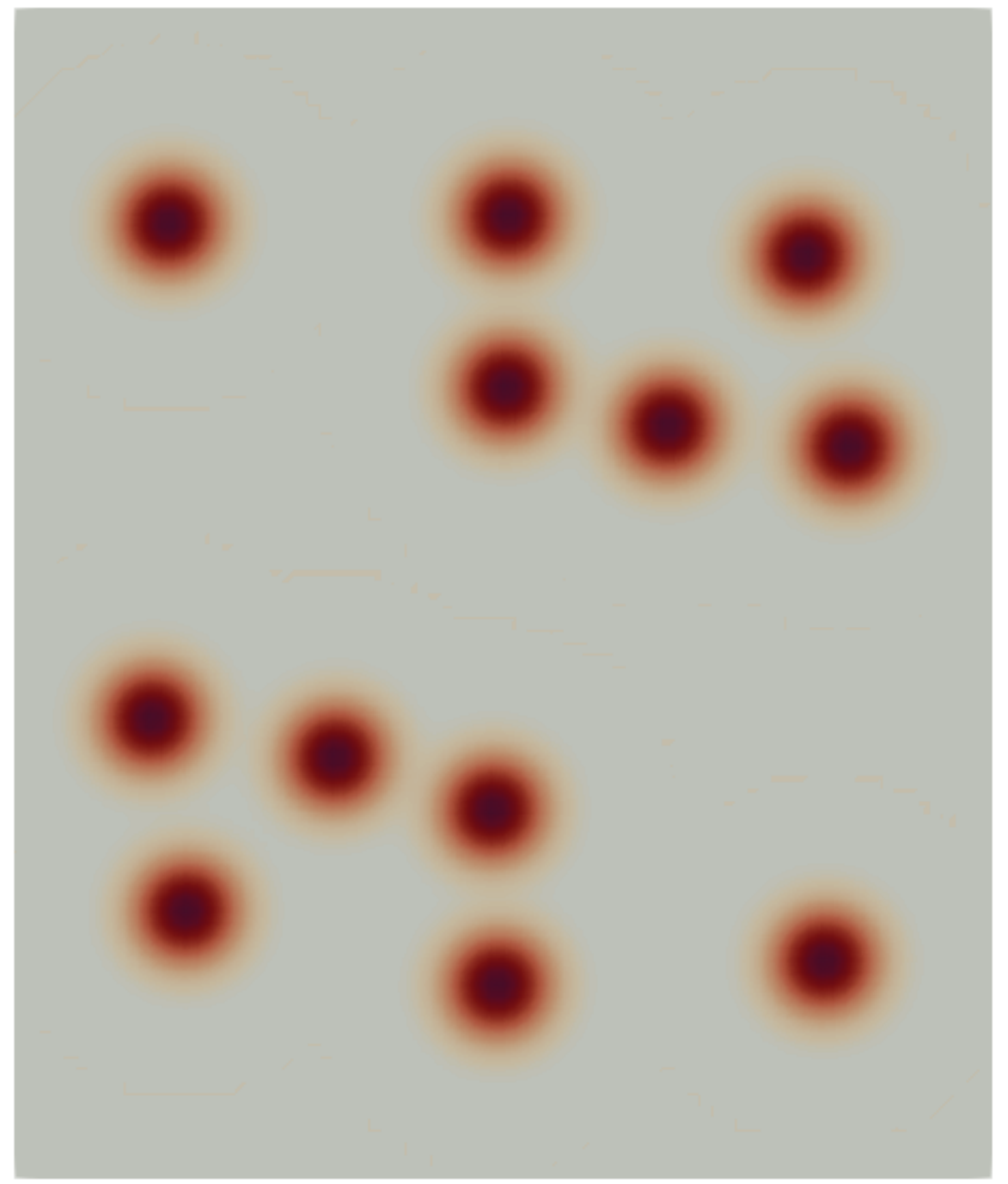}
        \caption{Optimal design $N_p = 12$}
    \end{subfigure}
    \begin{subfigure}{0.19\textwidth}
        \centering
        \includegraphics[width=0.9\textwidth]{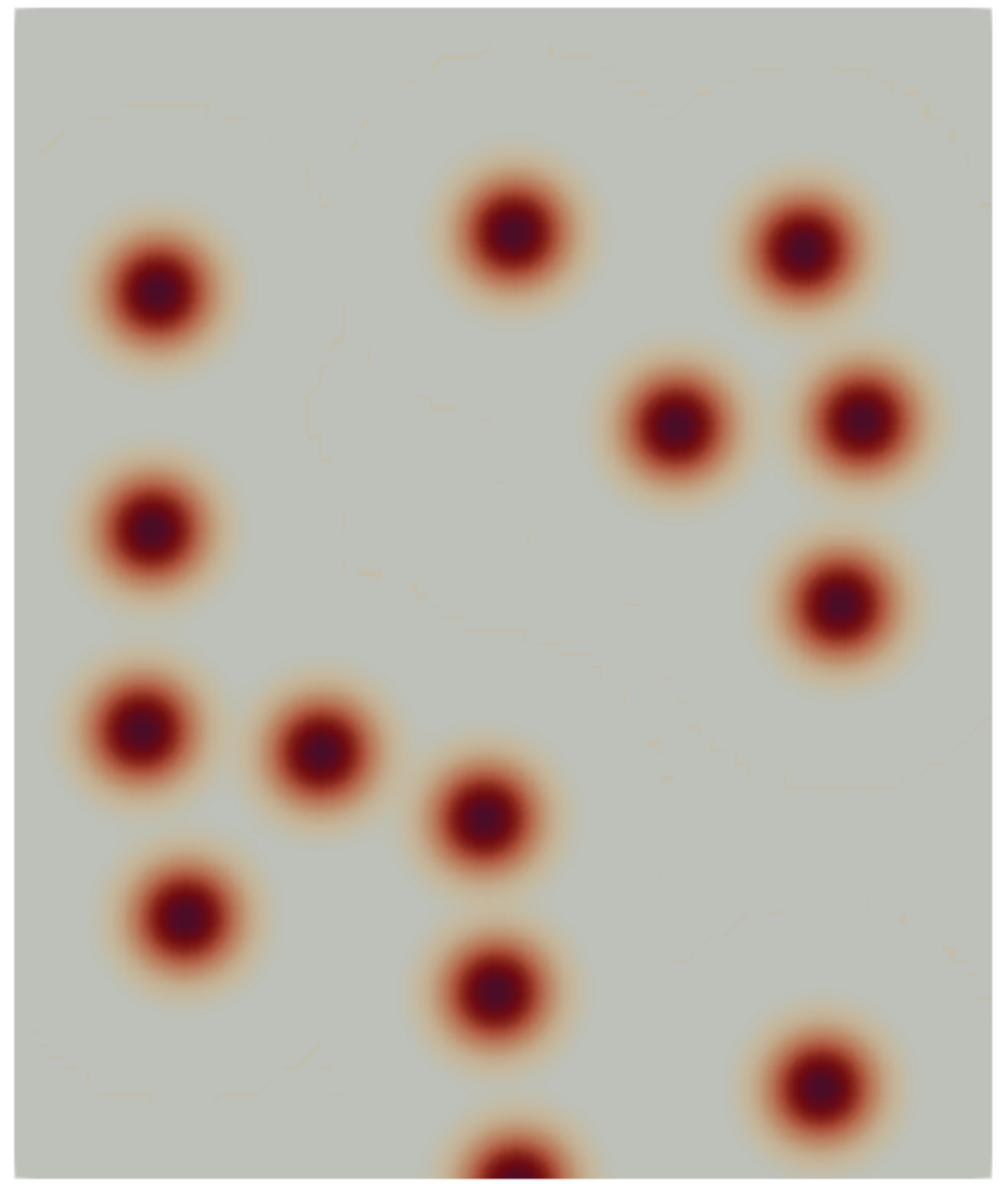}
        \caption{Optimal design $N_p = 14$}
    \end{subfigure}
    \begin{subfigure}{0.19\textwidth}
        \centering
        \includegraphics[width=0.9\textwidth]{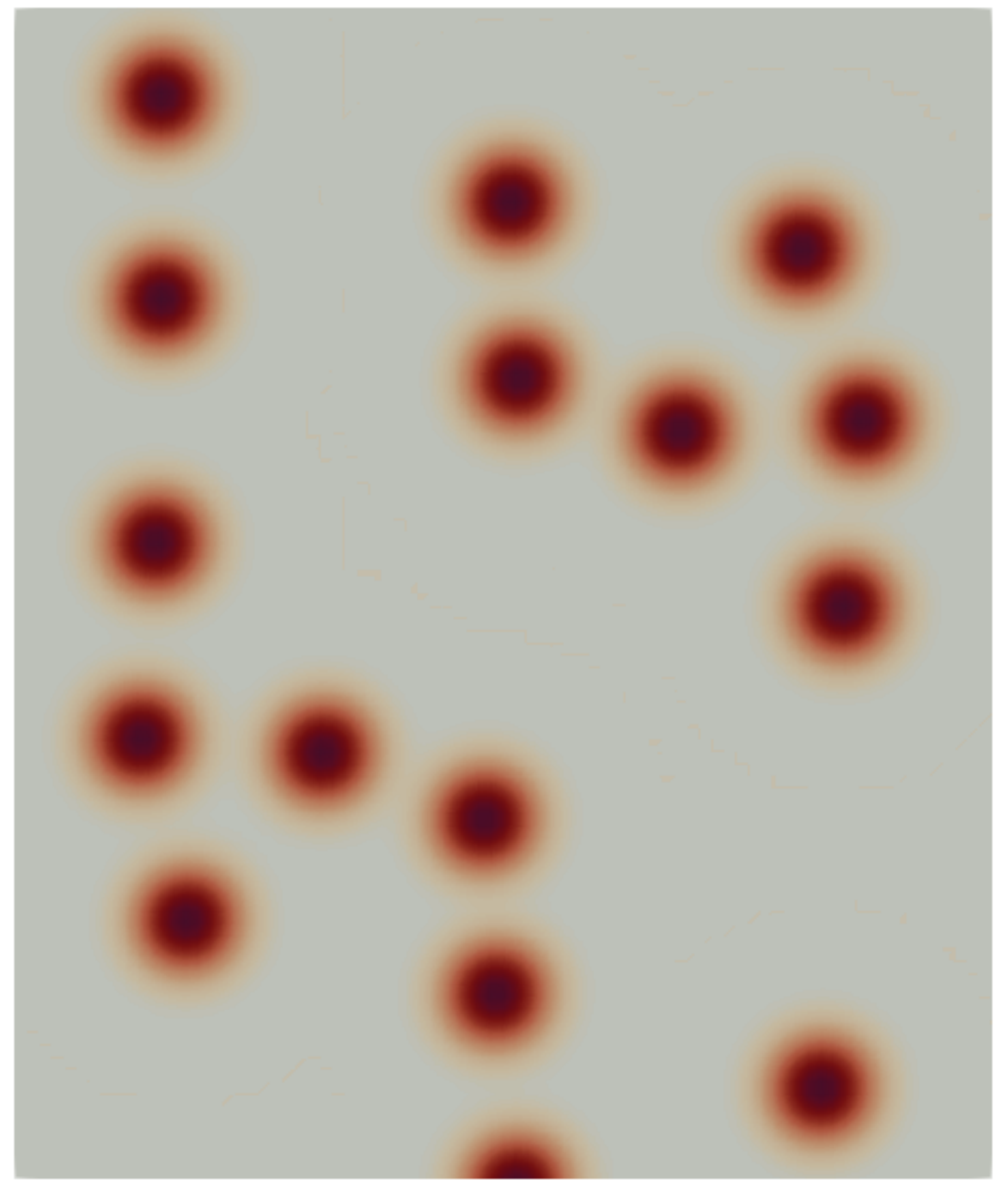}
        \caption{Optimal design $N_p = 16$}
    \end{subfigure}
    \begin{subfigure}{0.19\textwidth}
        \centering
        \includegraphics[width=0.9\textwidth]{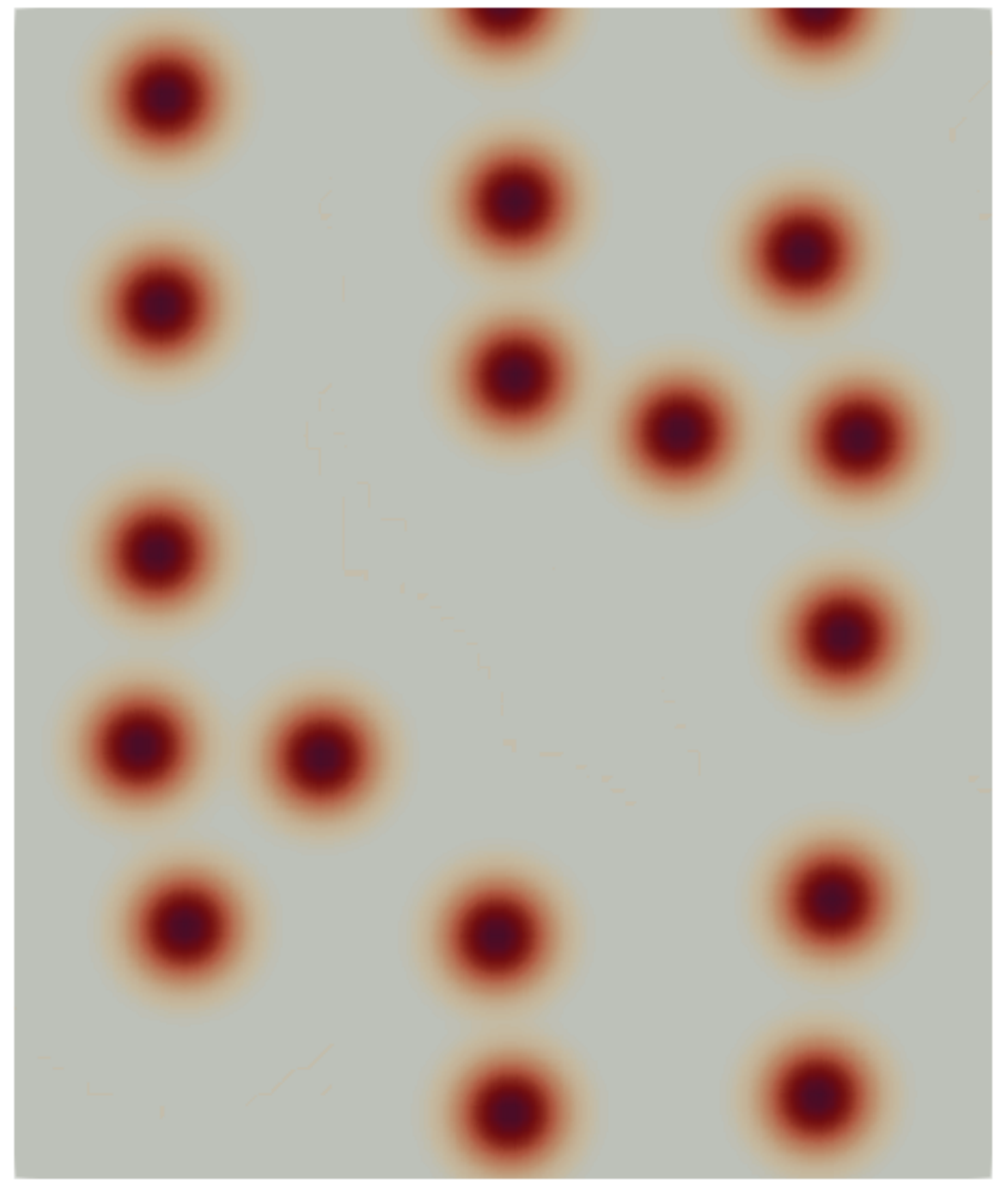}
        \caption{Optimal design $N_p = 16$}
    \end{subfigure}
    \begin{subfigure}{0.19\textwidth}
        \centering
        \includegraphics[width=0.9\textwidth]{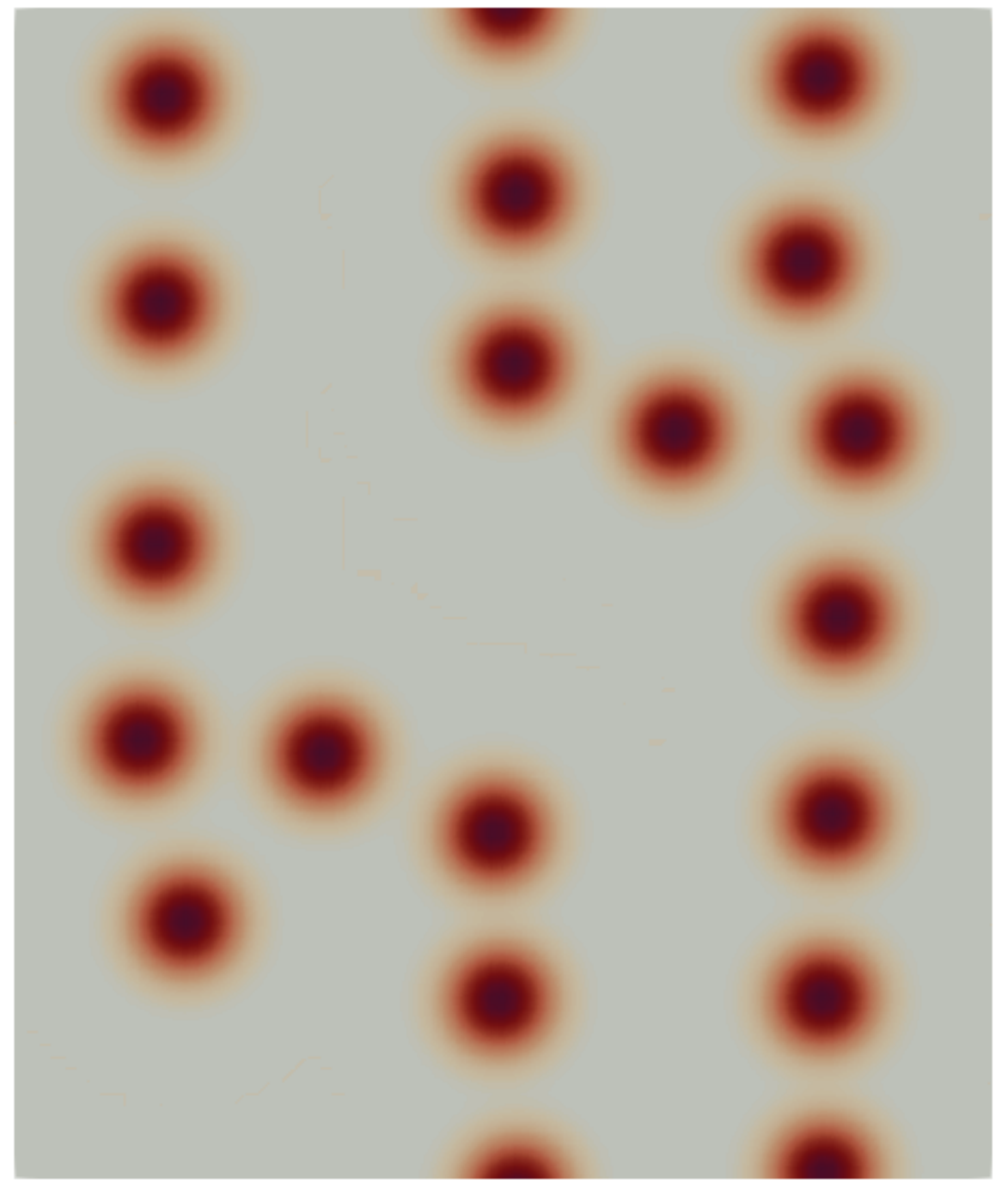}
        \caption{Optimal design $N_p = 20$}
    \end{subfigure}

    \begin{subfigure}{0.32\textwidth}
        \centering
        \includegraphics[width=1.0\textwidth]{figures_pdf/state_colors.pdf}
    \end{subfigure}

    \begin{subfigure}{0.19\textwidth}
        \centering
        \includegraphics[width=0.9\textwidth]{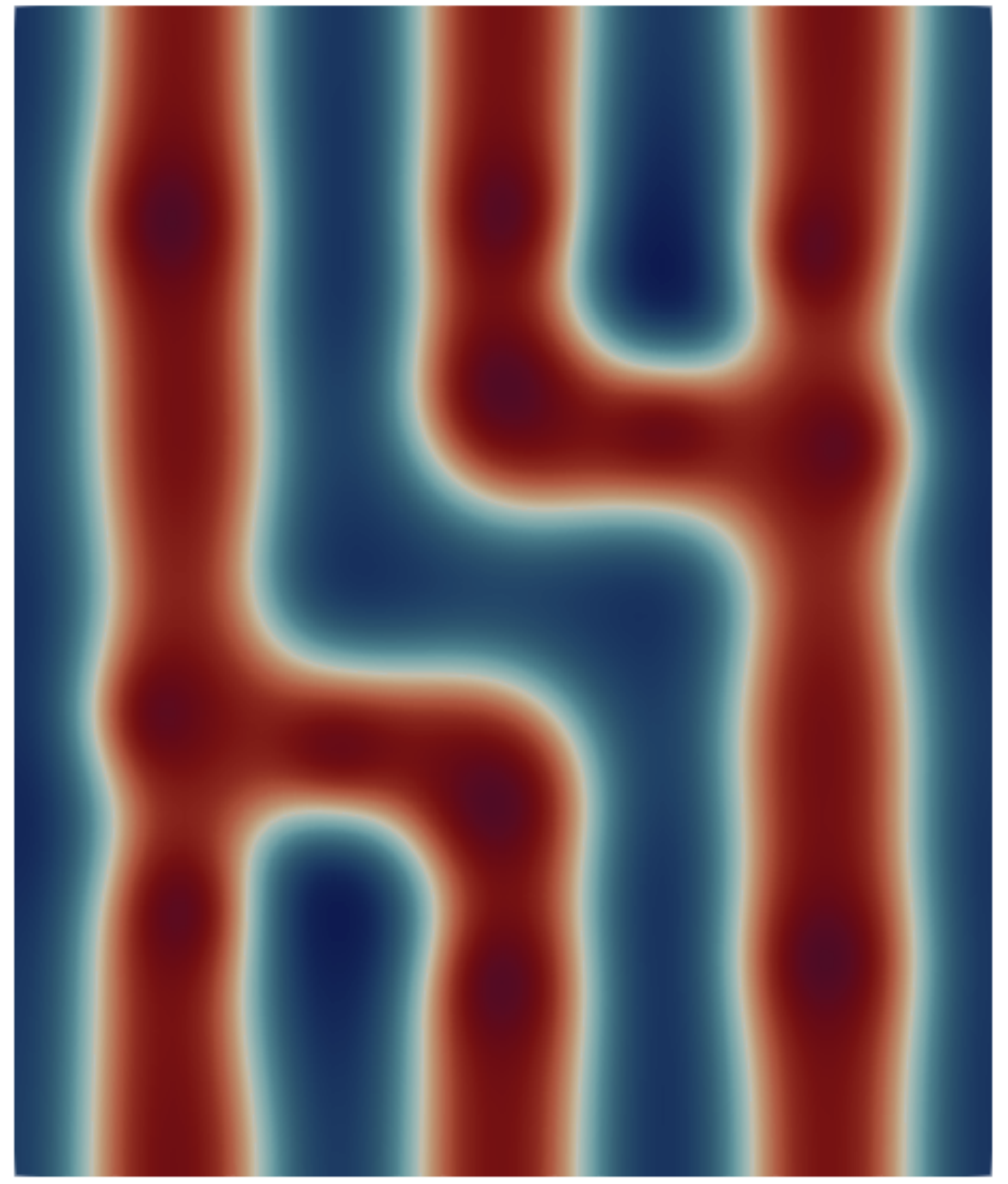}
        \caption{Minimum $\mathcal{F}(u)$ sample $N_p = 12$}
    \end{subfigure}
    \begin{subfigure}{0.19\textwidth}
        \centering
        \includegraphics[width=0.9\textwidth]{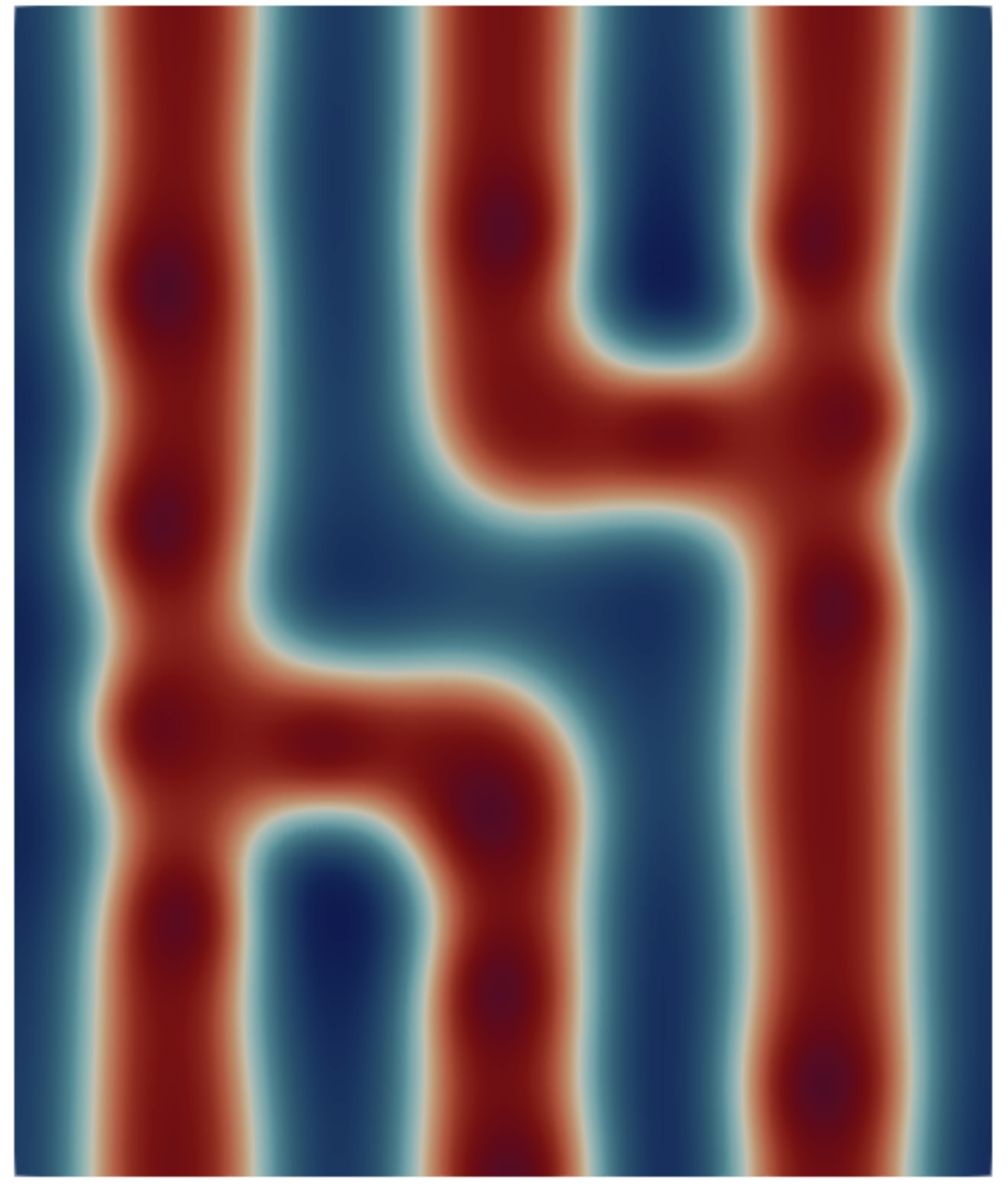}
        \caption{Minimum $\mathcal{F}(u)$ sample $N_p = 14$}
    \end{subfigure}
    \begin{subfigure}{0.19\textwidth}
        \centering
        \includegraphics[width=0.9\textwidth]{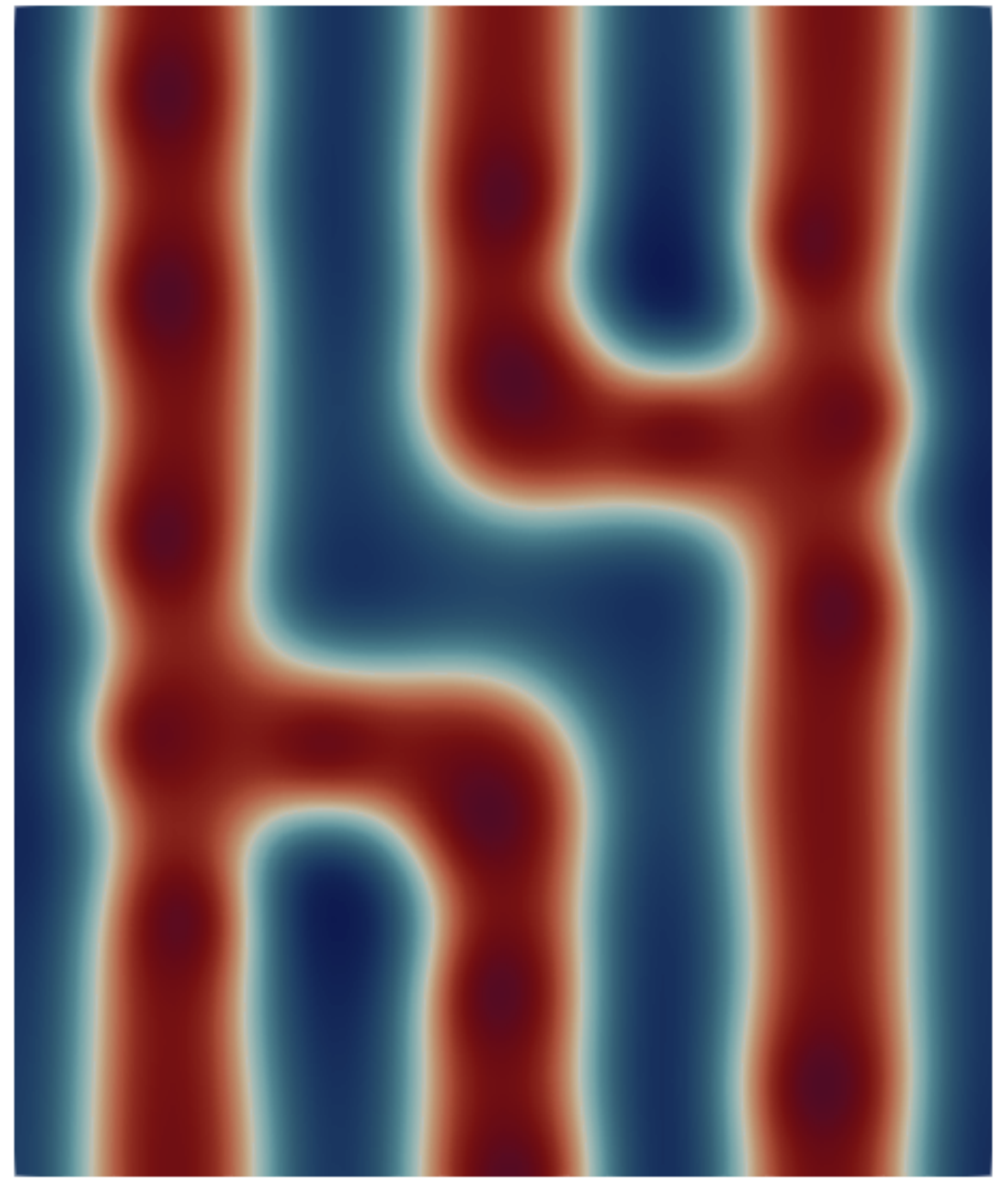}
        \caption{Minimum $\mathcal{F}(u)$ sample $N_p = 16$}
    \end{subfigure}
    \begin{subfigure}{0.19\textwidth}
        \centering
        \includegraphics[width=0.9\textwidth]{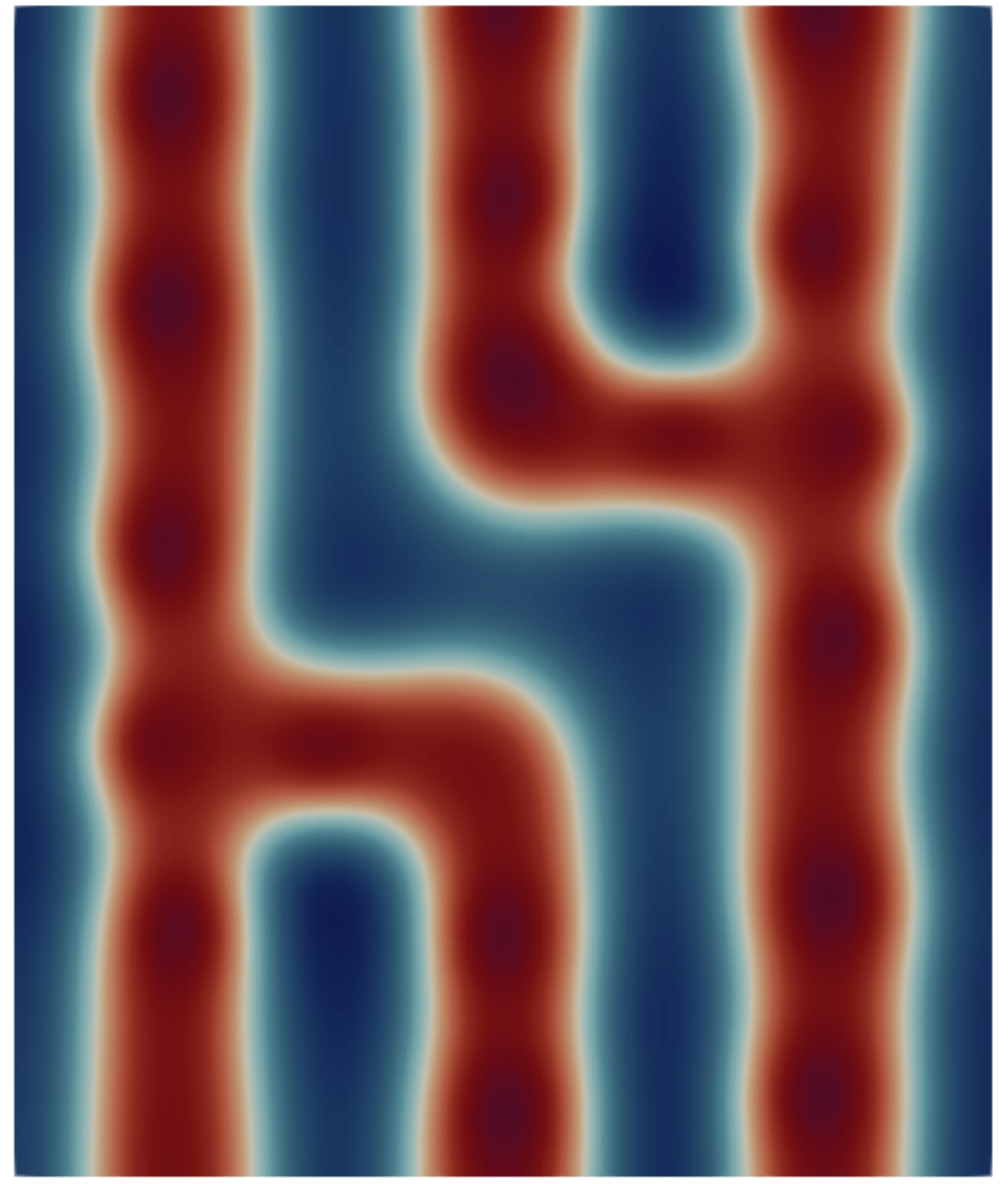}      
        \caption{Minimum $\mathcal{F}(u)$ sample $N_p = 16$}
    \end{subfigure}
    \begin{subfigure}{0.19\textwidth}
        \centering
        \includegraphics[width=0.9\textwidth]{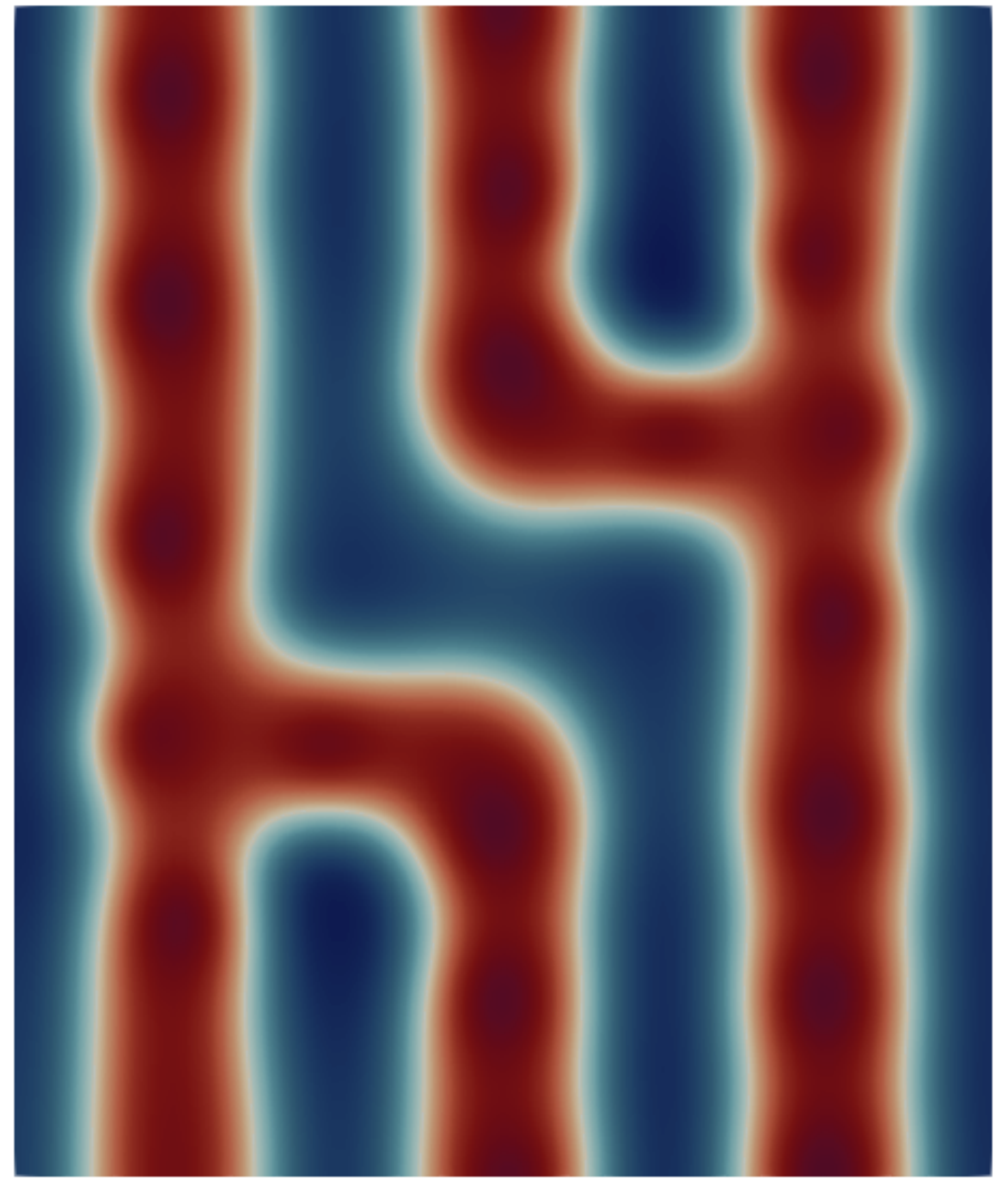}
        \caption{Minimum $\mathcal{F}(u)$ sample $N_p = 20$}
    \end{subfigure}

    \begin{subfigure}{0.19\textwidth}
        \centering
        \includegraphics[width=0.9\textwidth]{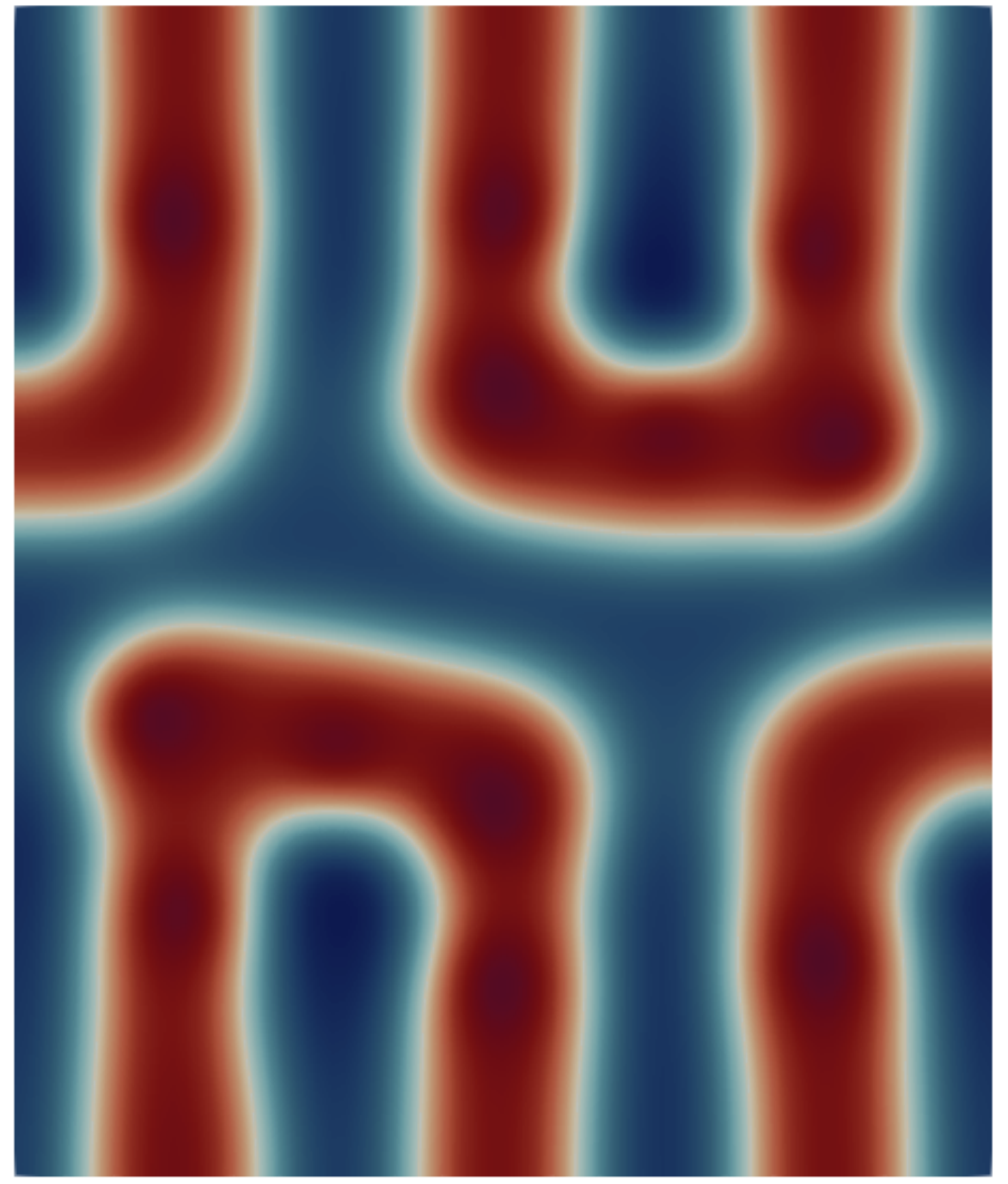}
        \caption{Sample 2 $N_p = 12$}
    \end{subfigure}
    \begin{subfigure}{0.19\textwidth}
        \centering
        \includegraphics[width=0.9\textwidth]{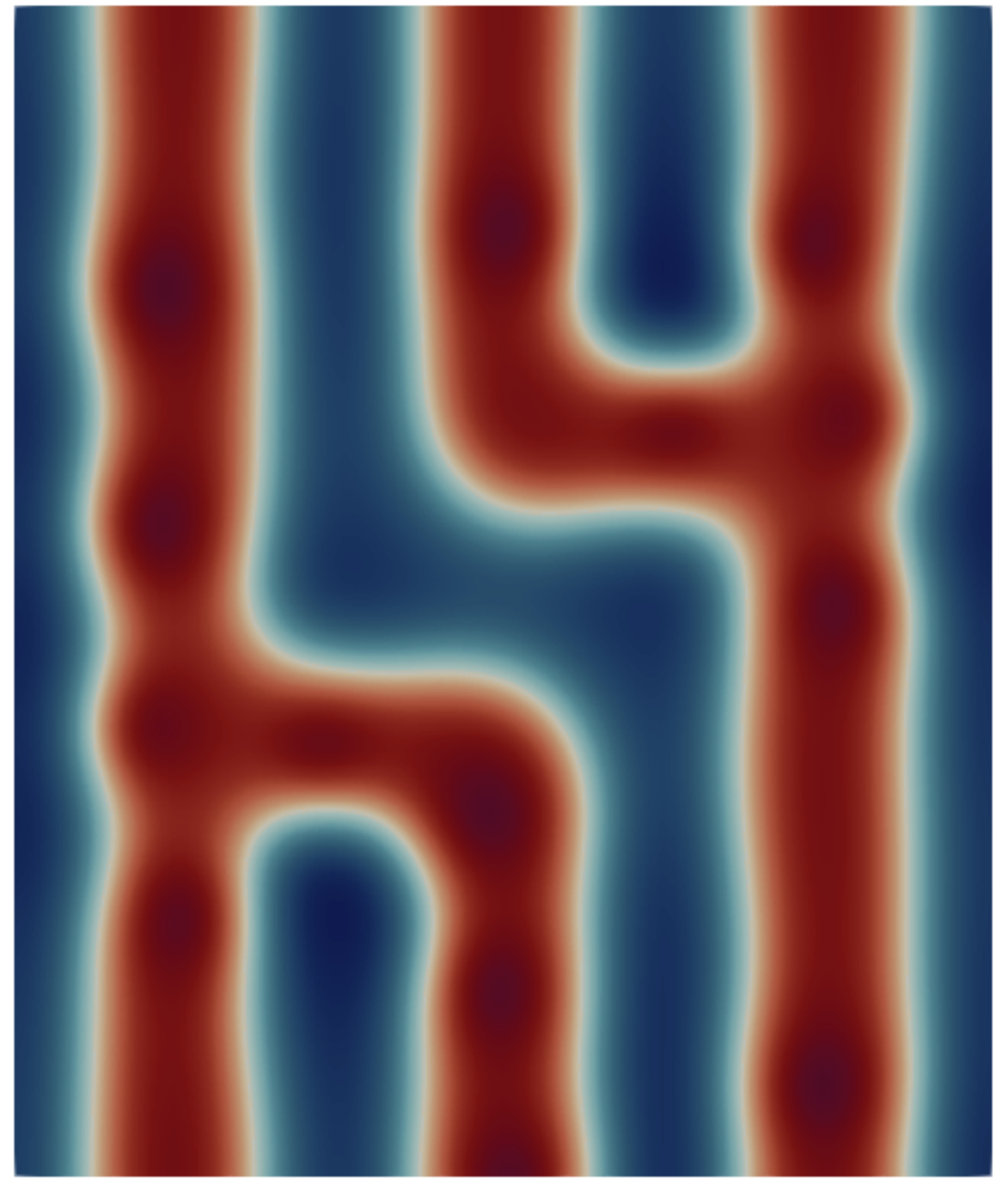}
        \caption{Sample 2 $N_p = 14$}
    \end{subfigure}
    \begin{subfigure}{0.19\textwidth}
        \centering
        \includegraphics[width=0.9\textwidth]{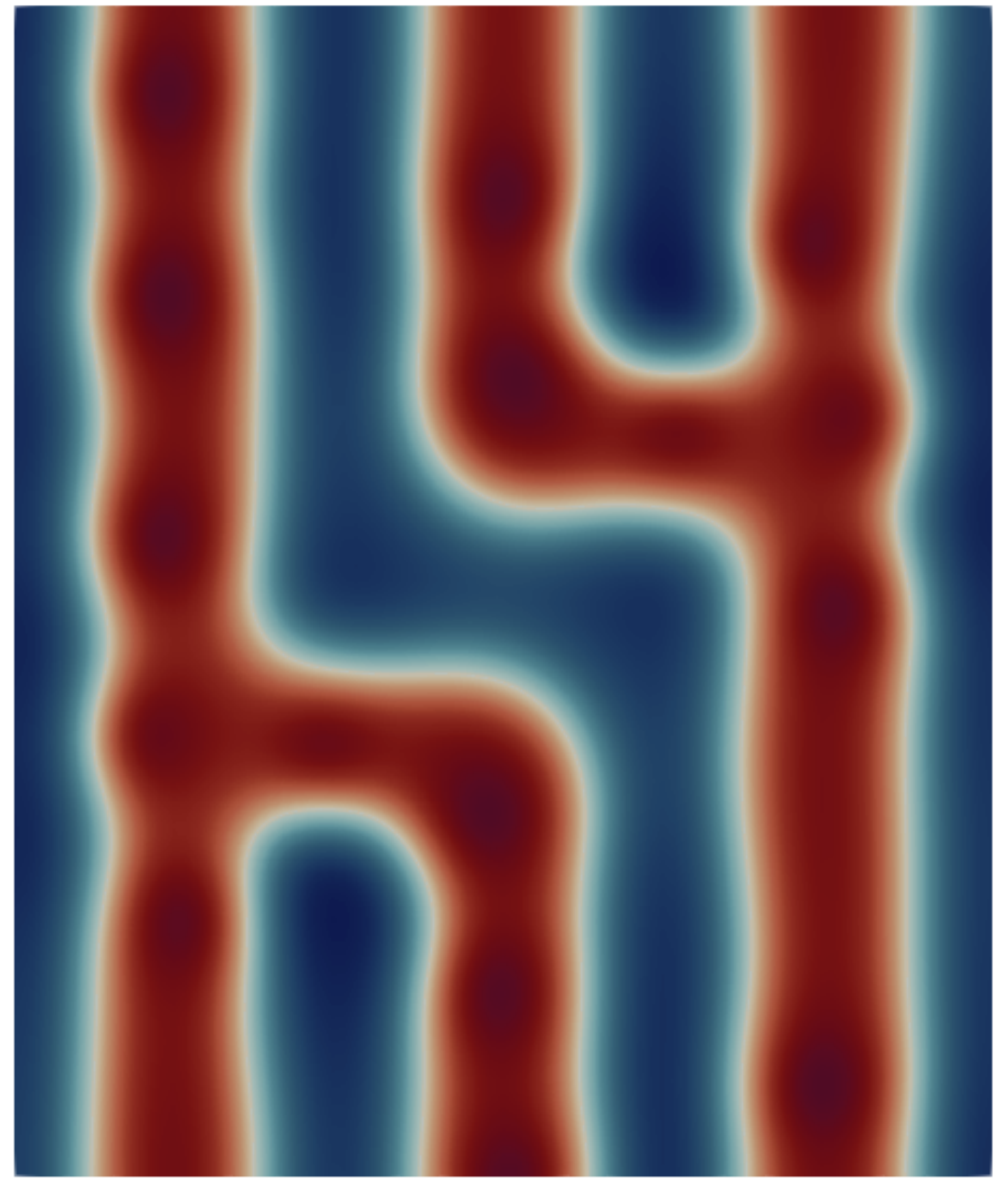}
        \caption{Sample 2 $N_p = 16$}
    \end{subfigure}
    \begin{subfigure}{0.19\textwidth}
        \centering
        \includegraphics[width=0.9\textwidth]{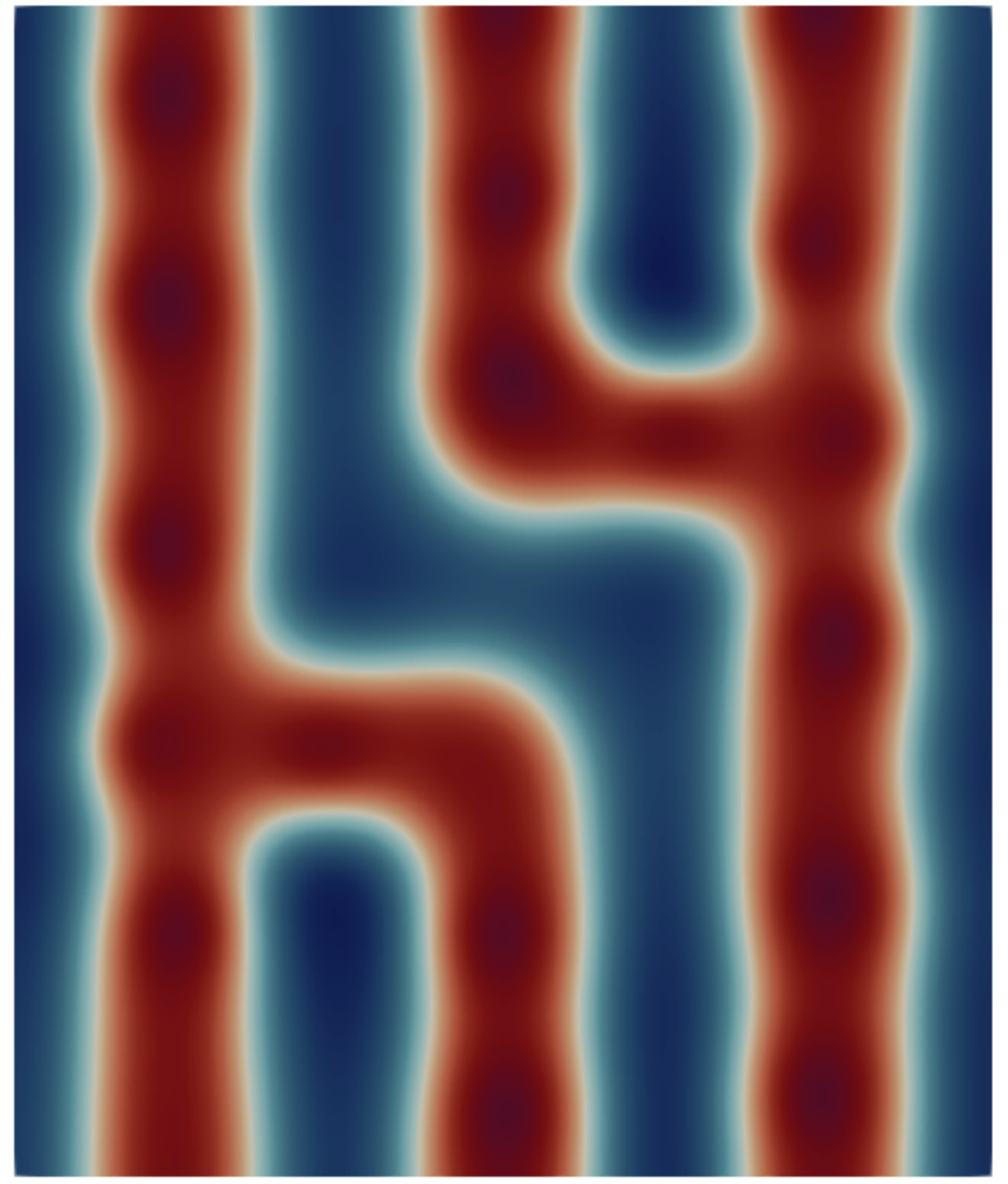}      
        \caption{Sample 2 $N_p = 16$}
    \end{subfigure}
    \begin{subfigure}{0.19\textwidth}
        \centering
        \includegraphics[width=0.9\textwidth]{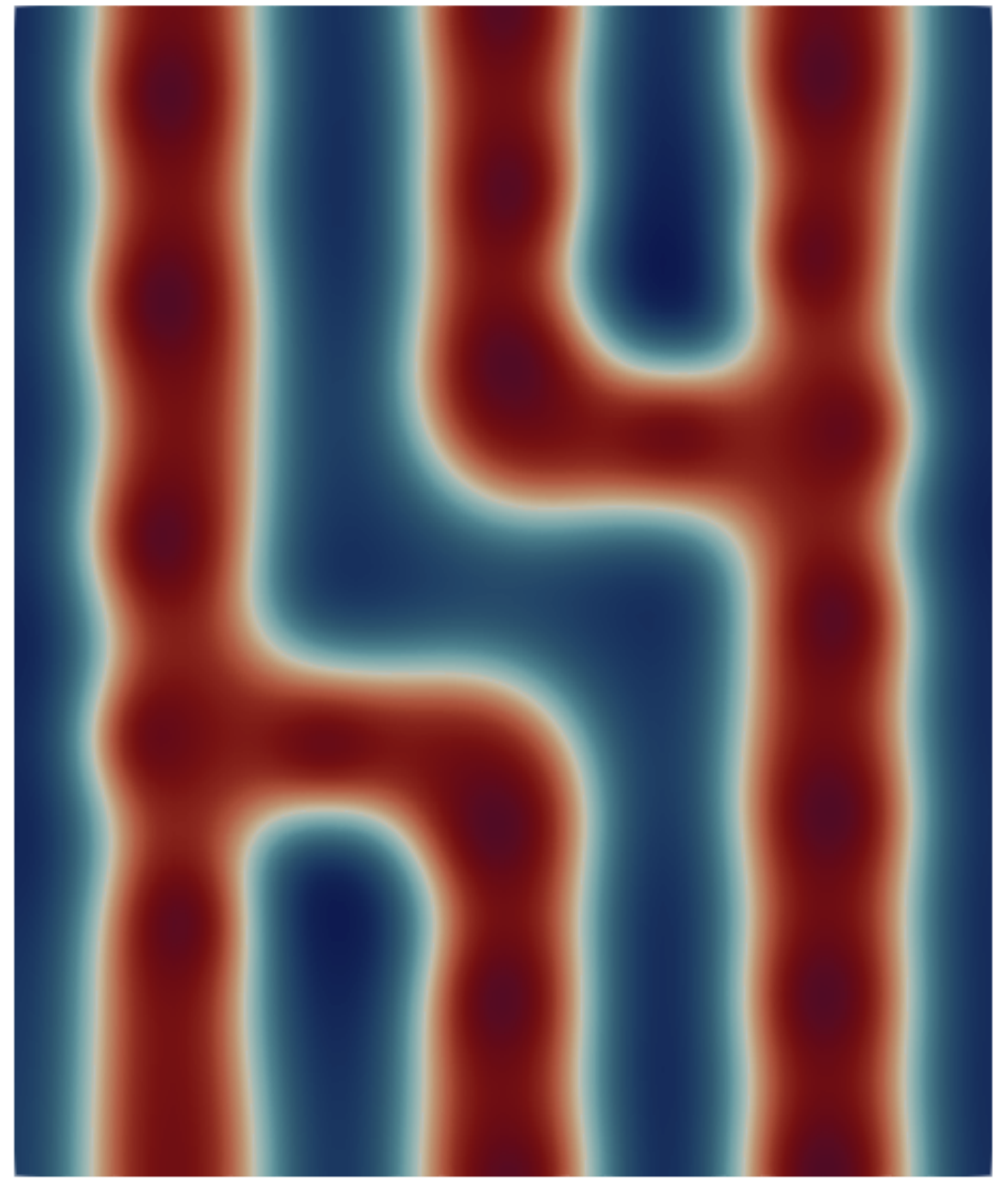}
        \caption{Sample 2 $N_p = 20$}
    \end{subfigure}

    \begin{subfigure}{0.19\textwidth}
        \centering
        \includegraphics[width=0.9\textwidth]{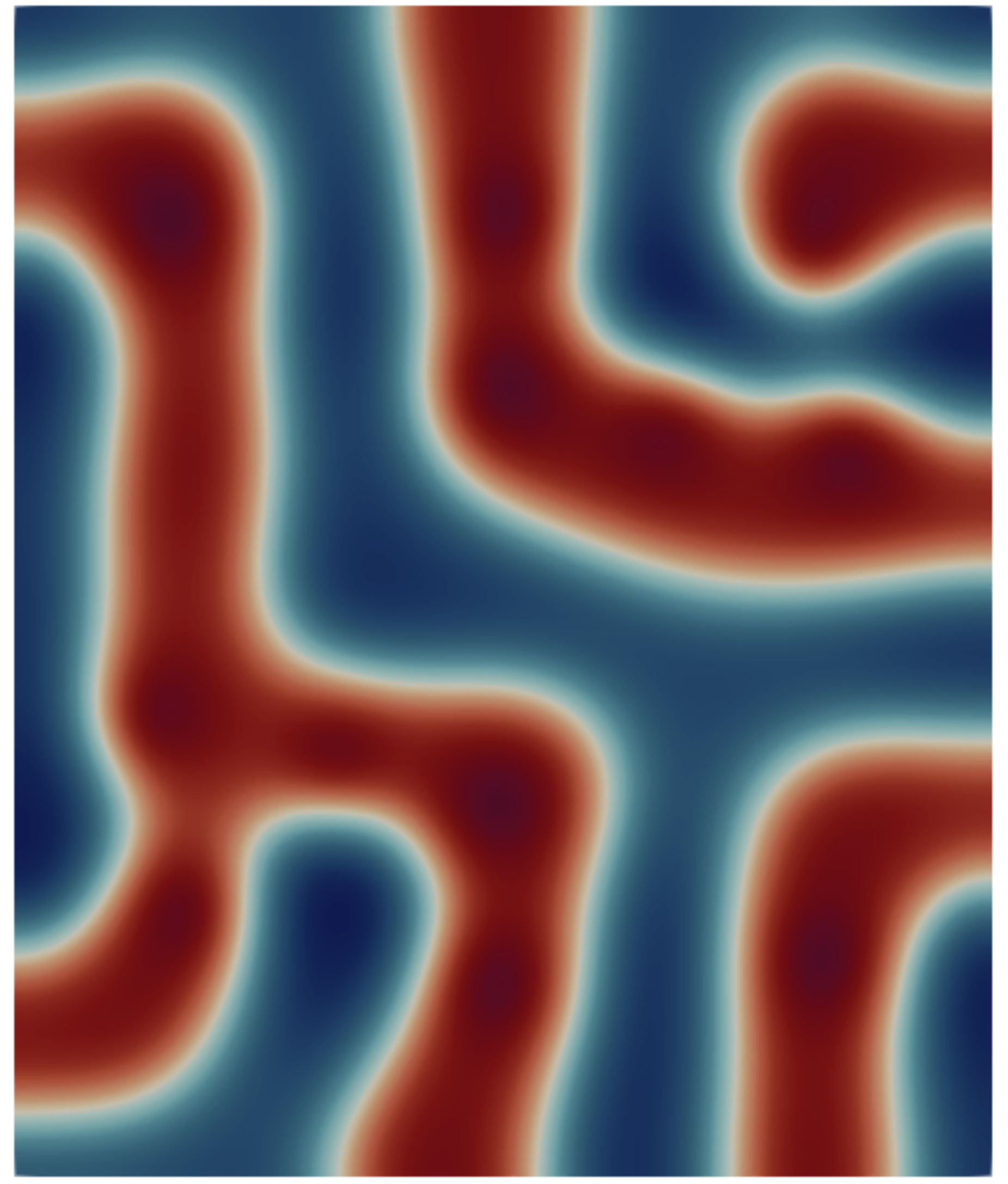}
        \caption{Sample 3 $N_p = 12$}
    \end{subfigure}
    \begin{subfigure}{0.19\textwidth}
        \centering
        \includegraphics[width=0.9\textwidth]{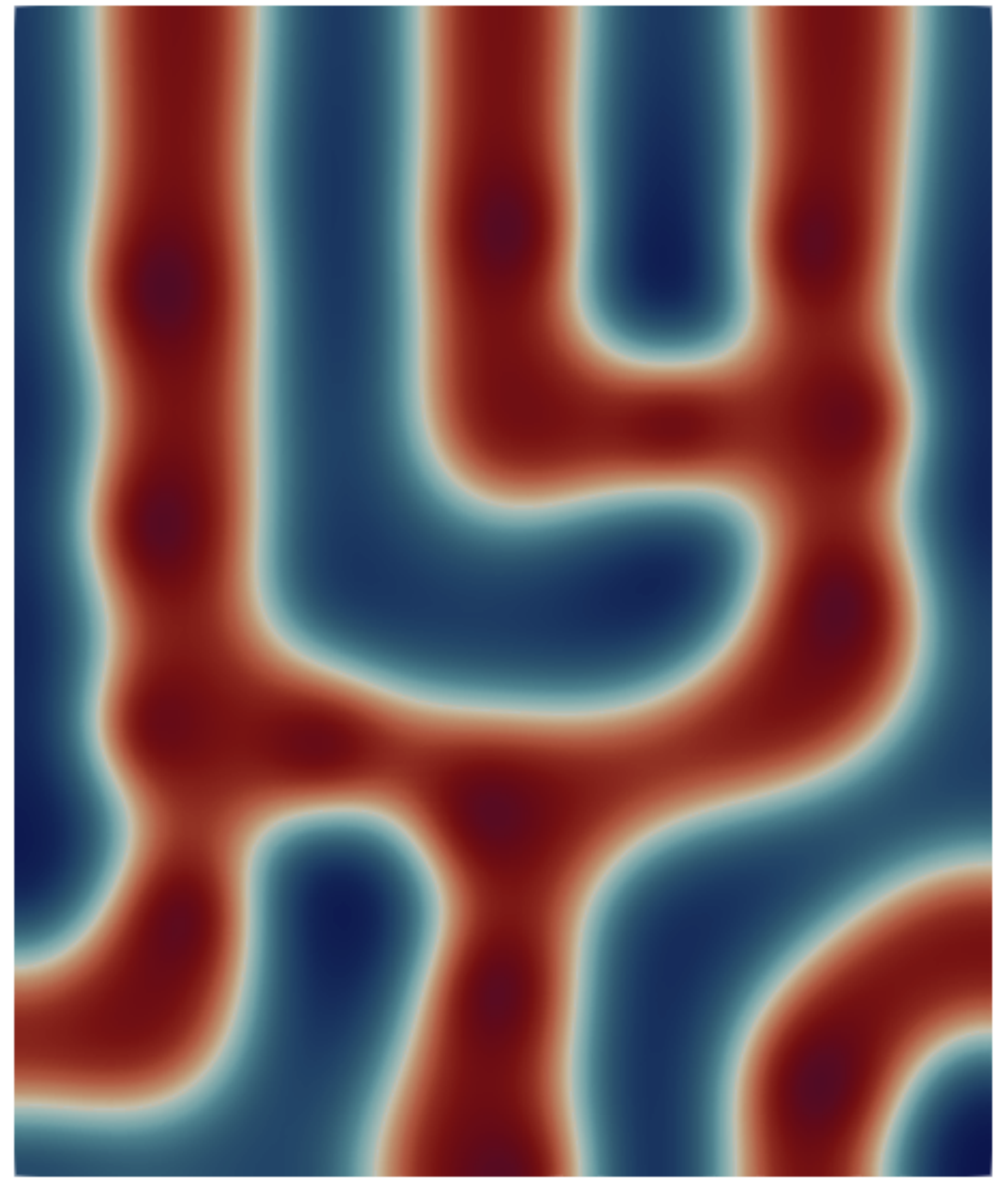}
        \caption{Sample 3 $N_p = 14$}
    \end{subfigure}
    \begin{subfigure}{0.19\textwidth}
        \centering
        \includegraphics[width=0.9\textwidth]{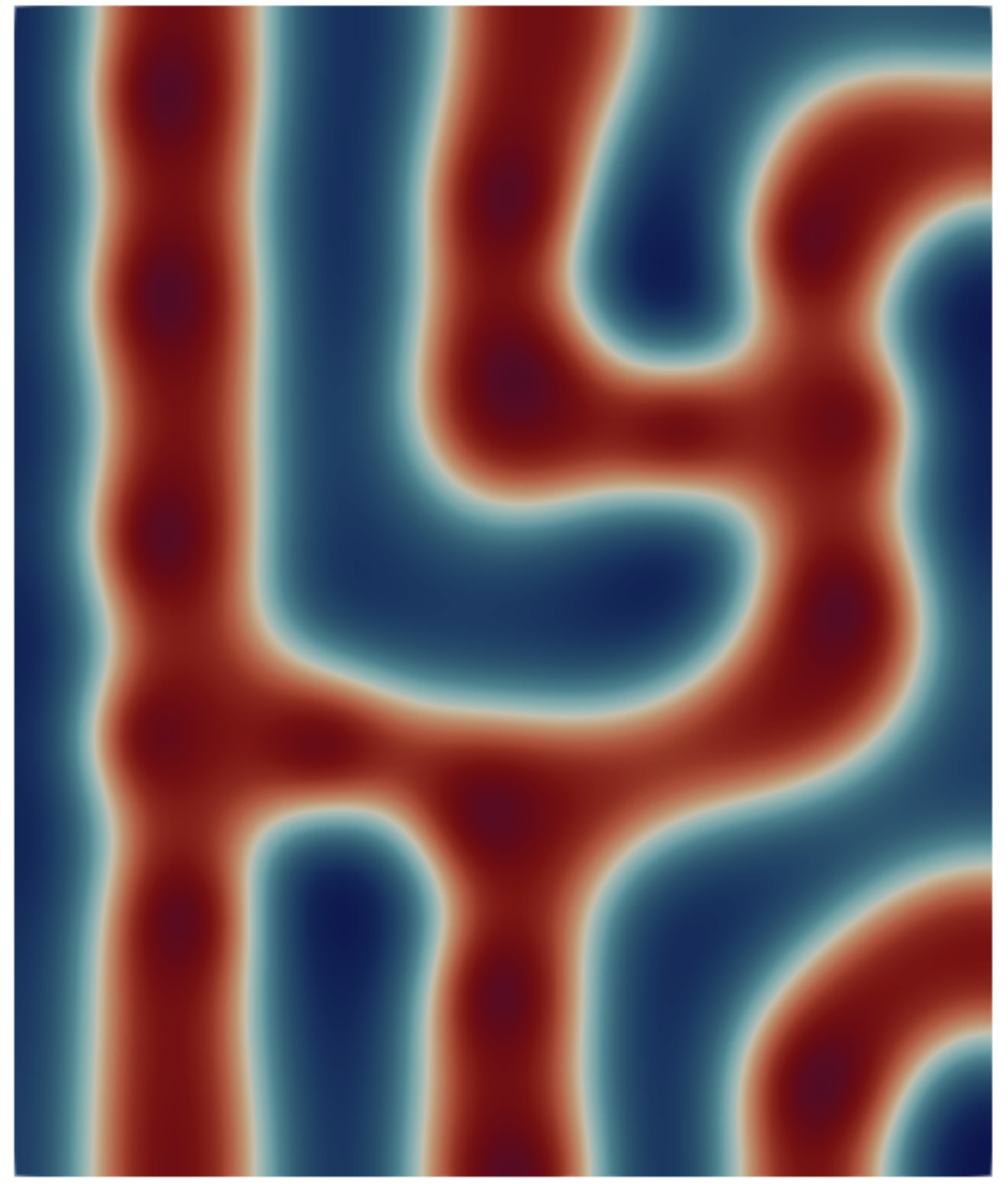}
        \caption{Sample 3 $N_p = 16$}
    \end{subfigure}
    \begin{subfigure}{0.19\textwidth}
        \centering
        \includegraphics[width=0.9\textwidth]{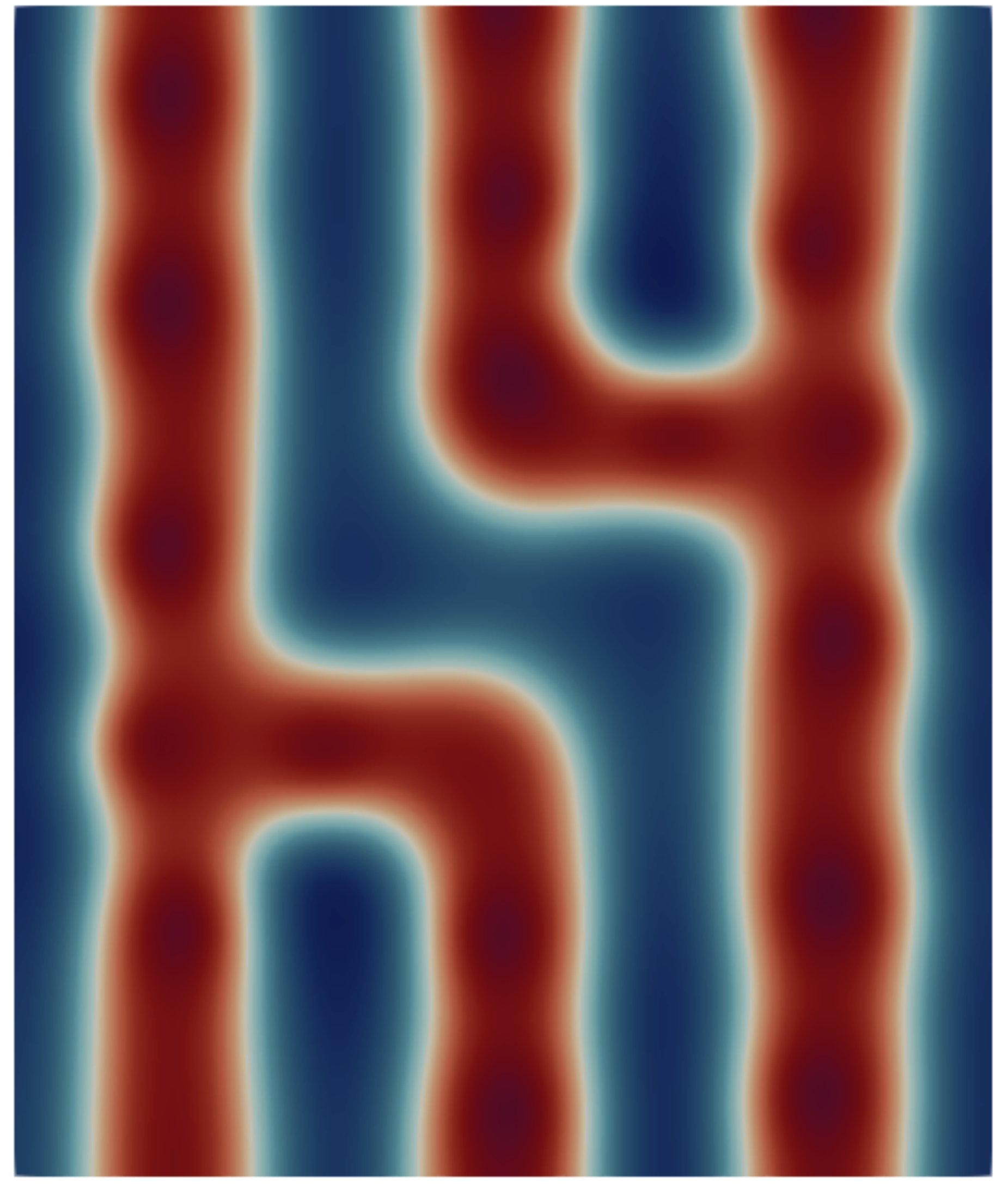}      
        \caption{Sample 3 $N_p = 16$}
    \end{subfigure}
    \begin{subfigure}{0.19\textwidth}
        \centering
        \includegraphics[width=0.9\textwidth]{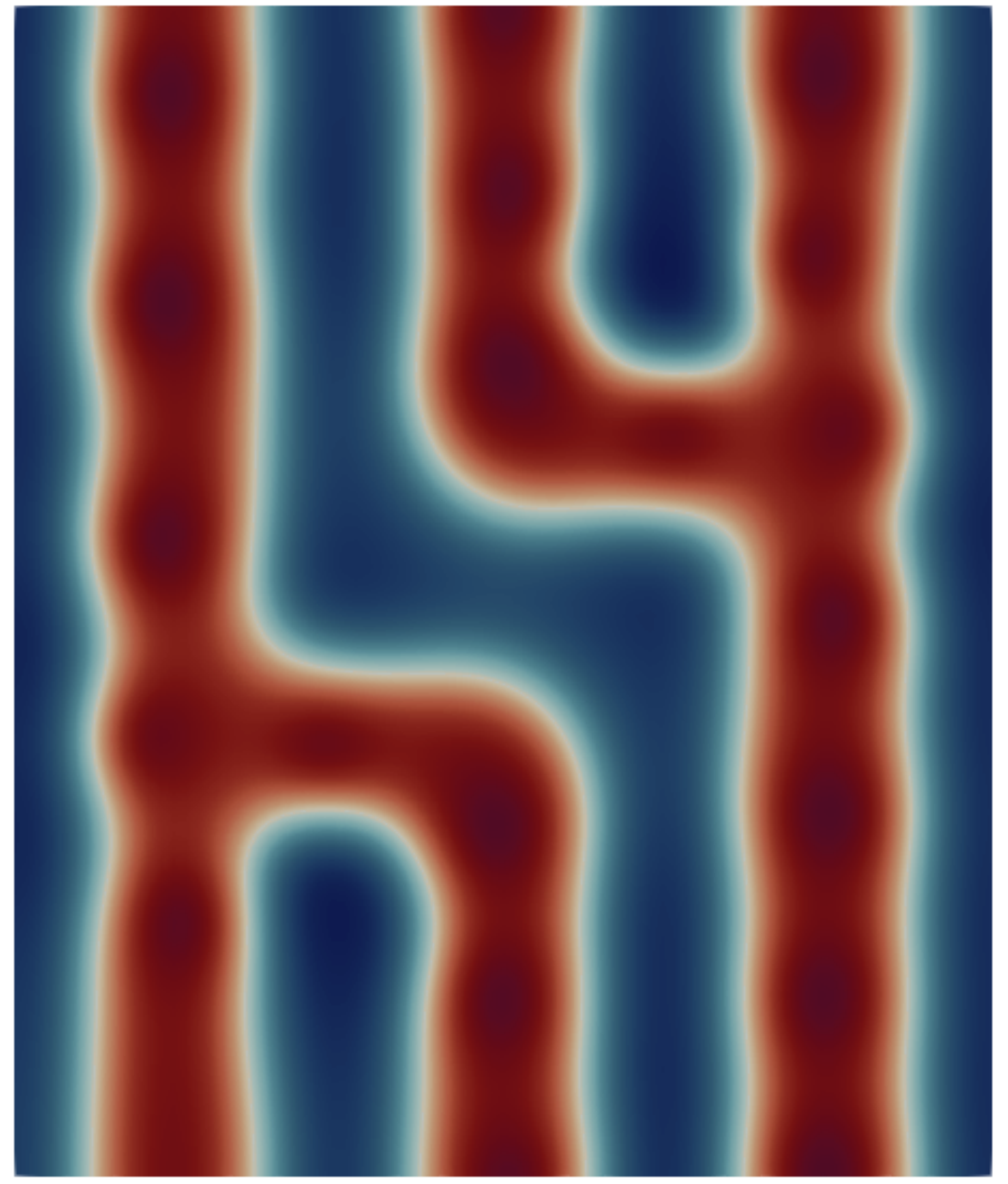}
        \caption{Sample 3 $N_p = 20$}
    \end{subfigure}
    \caption{Optimal design results for the jog target morphology using guidepost numbers $N_p = 12, 14, 16, 18, 20$. The final three rows are solutions of the state problem computed using samples of the initial guess $u_0$. In particular, at the top of this section are the samples states with the minimum free energy.}
    \label{fig:jog_Np}
\end{figure}

We showcase the flexibility of our formulation by considering two additional non-trivial target morphologies, the junction and jog morphologies, as shown in Figure \ref{fig:2d_targets}. In particular, the junction morphology (Figure \ref{fig:2d_targets}, left) is based on a target considered in \cite{HannonGotrikRossEtAl13}, while the jog morphology (Figure \ref{fig:2d_targets}, right) is presented in \cite{StoykovichKangDaoulasEtAl07} as an essential circuit geometry. We use $N_p$ circular guideposts of width $b=0.08$, and solve the optimal design problem for both target morphologies with $N_p = $ 12, 14, 16, 18, and 20. For the penalization parameter, we select $\alpha = 5\times 10^{-3}$ for the junction and $\alpha = 10^{-3}$ for the jog. We optimize with a step bound of $\stepbound=0.05$. The optimal designs and sample states are presented in Figure \ref{fig:junction2d_Np} for the junction target morphology and in Figure \ref{fig:jog_Np} for the jog target morphology.

\begin{figure}
\centering
    \begin{subfigure}{0.45\textwidth}
        \centering
        \includegraphics[width=1.0\textwidth]{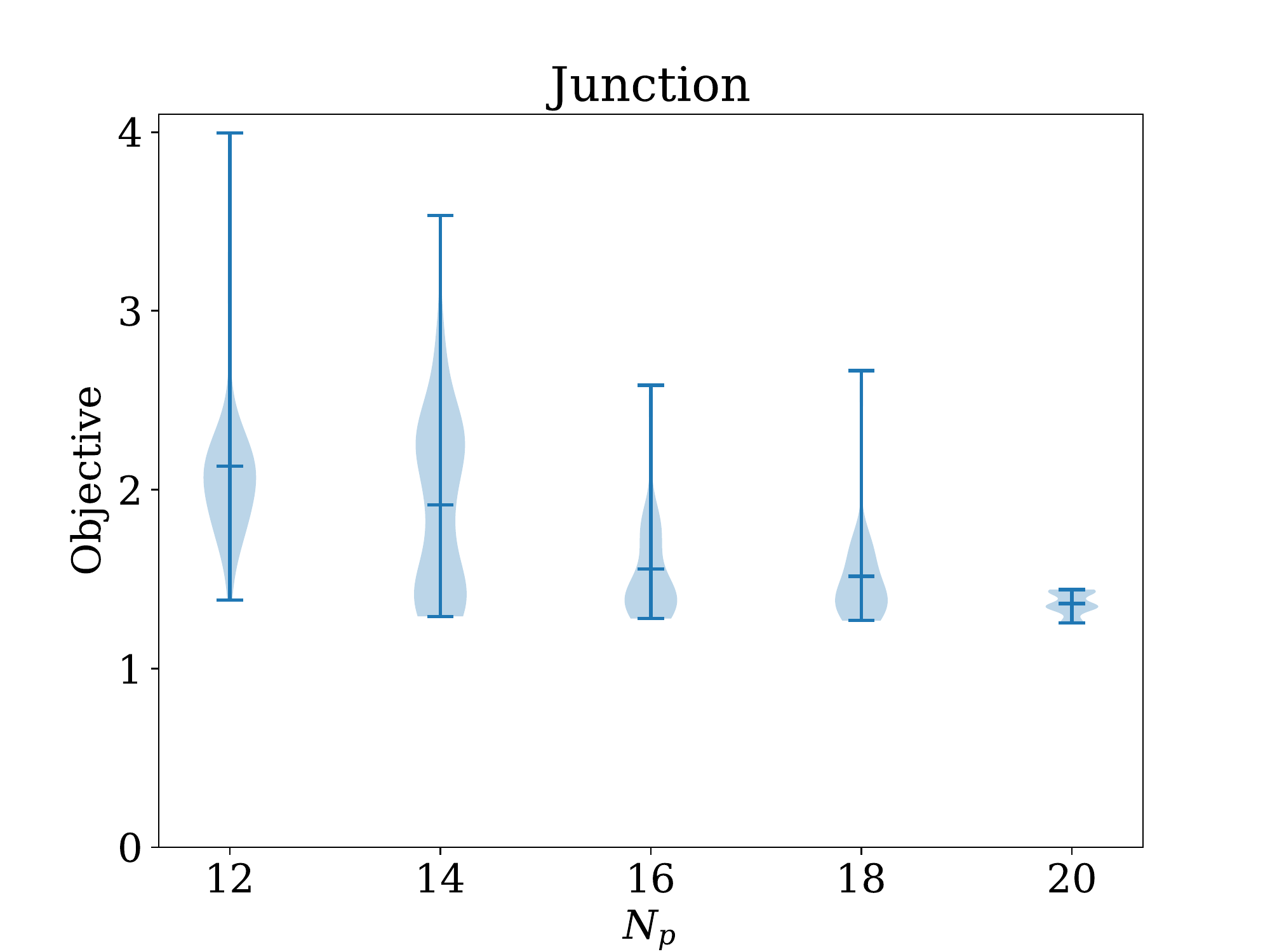}
    \end{subfigure}
    \begin{subfigure}{0.45\textwidth}
        \centering
        \includegraphics[width=1.0\textwidth]{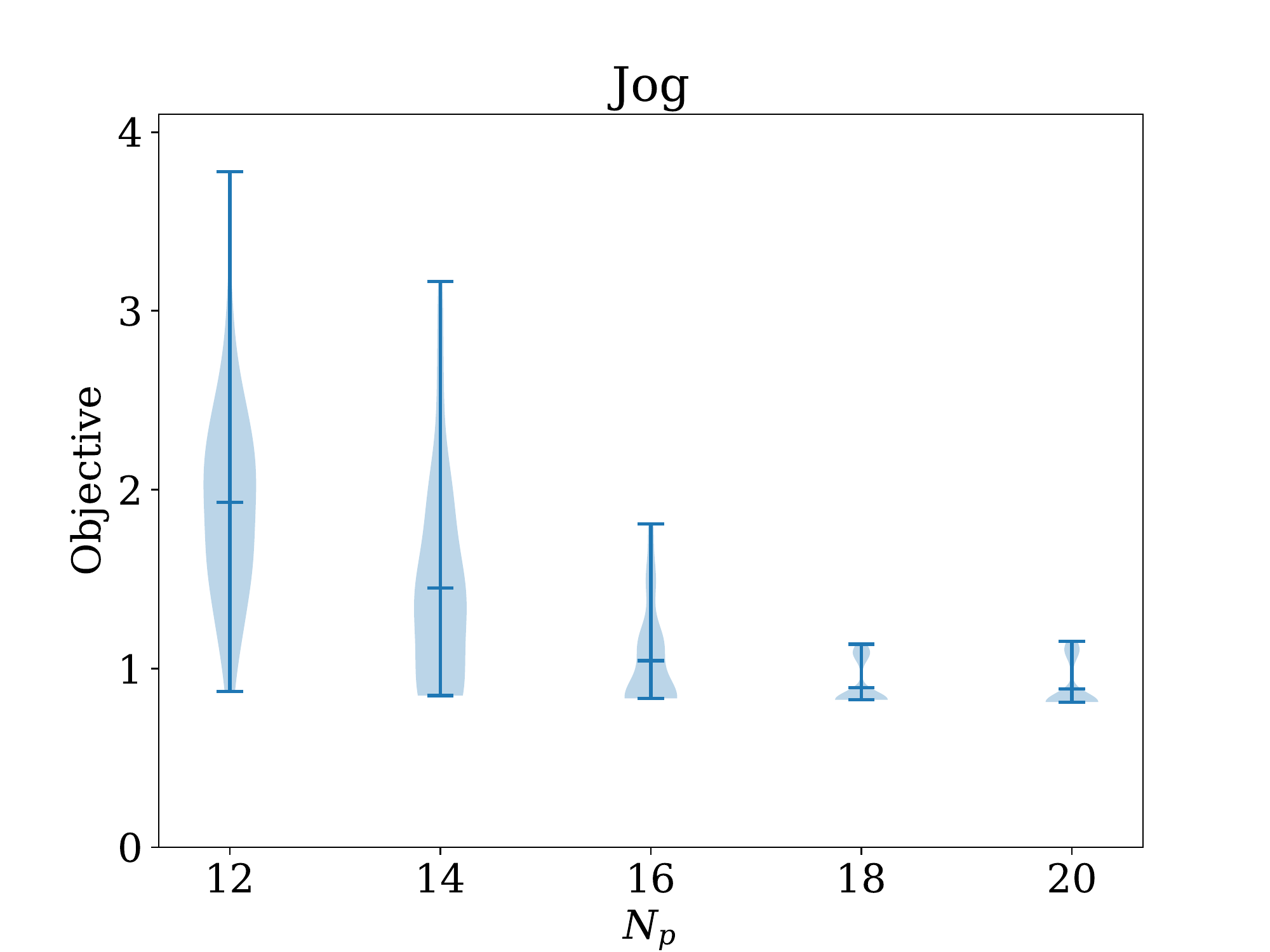}
    \end{subfigure}

    \caption{Distribution of objective function evaluations using 100 initial guess samples for the optimal design of the junction (left) and jog (right) target morphologies.}
    \label{fig:spot_obj_dist}
\end{figure}

In almost all of the cases, the sample states with minimum free energies highly resemble the target morphology. However, we observe more defective sample states for designs with fewer guideposts, such as that for $N_p = 12$ (visually, we may consider a state to be defective if it has missing connections between the red regions). As more guideposts are added to the system, the optimal designs tend to become less sensitive to different initial guesses. This is quantitatively shown in Figure \ref{fig:spot_obj_dist}, where we plot distributions of the objective function evaluated using 100 sample states. Smaller means and variances of the objective function are associated with designs using more guideposts, indicating that such designs are more robust to initial guesses.

% However, the optimal designs also represent a trade-off between the sparsity of the substrate pattern with different number of guideposts and the robustness of the optimal design. For both the junction and jog target morphologies, the number of defective sample states at the optimal designs is largest for $N_p = 12$, and tends to decrease as more guideposts are added to the system. 
% \pc{clarify: how many samples did you compute? how many are defective? how to quantify/determine if a state is defective, visually?}

\begin{figure}
    \centering
    \captionsetup[subfigure]{justification=centering}
    \begin{subfigure}{0.32\textwidth}
        \centering
        \includegraphics[width=1.0\textwidth]{figures_pdf/design_colors.pdf}
    \end{subfigure}

    \begin{subfigure}{0.19\textwidth}
        \centering
        \includegraphics[width=0.9\textwidth]{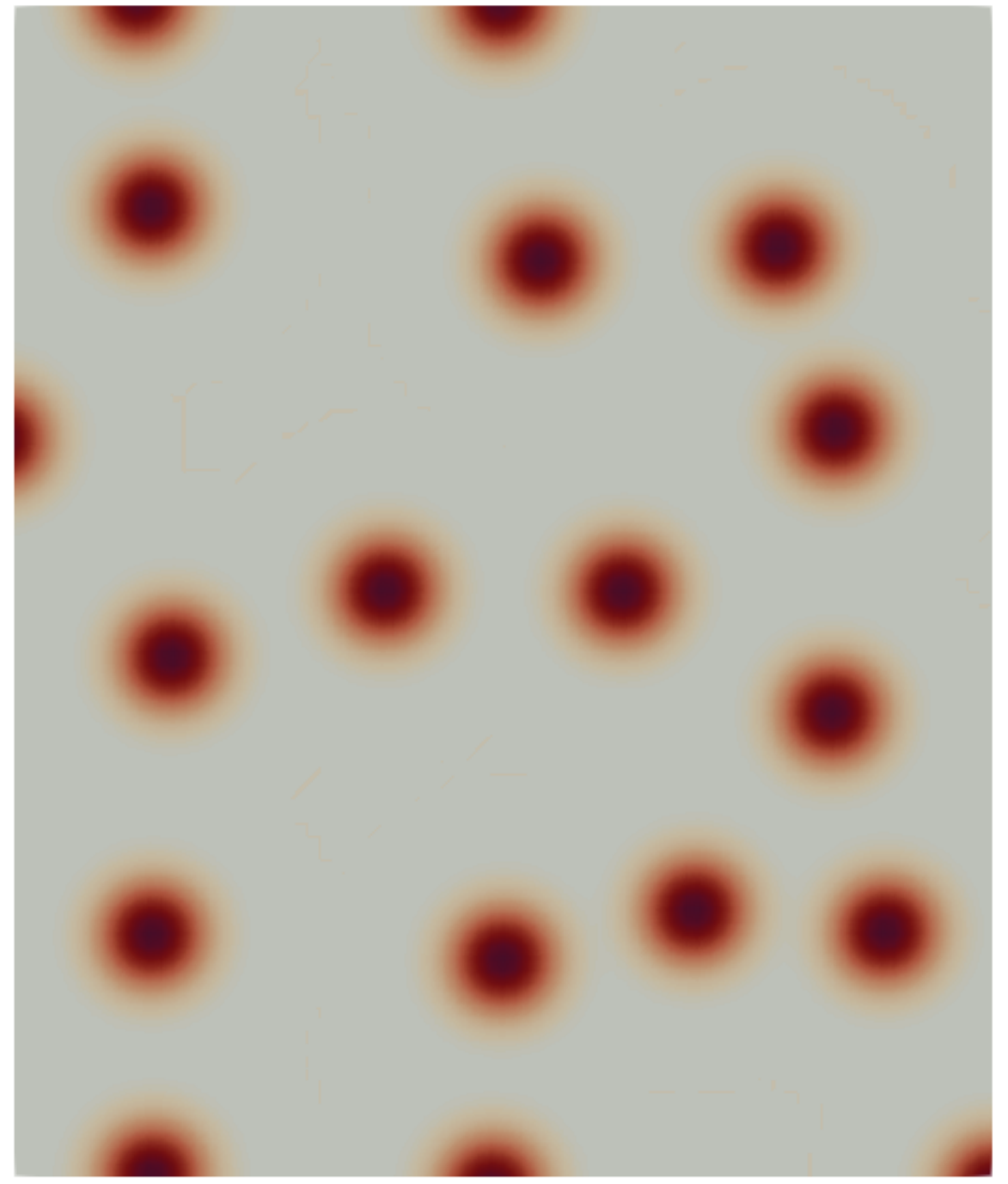}
        \caption{Optimal design $\alpha = 1\times 10^{-2}$}
    \end{subfigure}
    \begin{subfigure}{0.19\textwidth}
        \centering
        \includegraphics[width=0.9\textwidth]{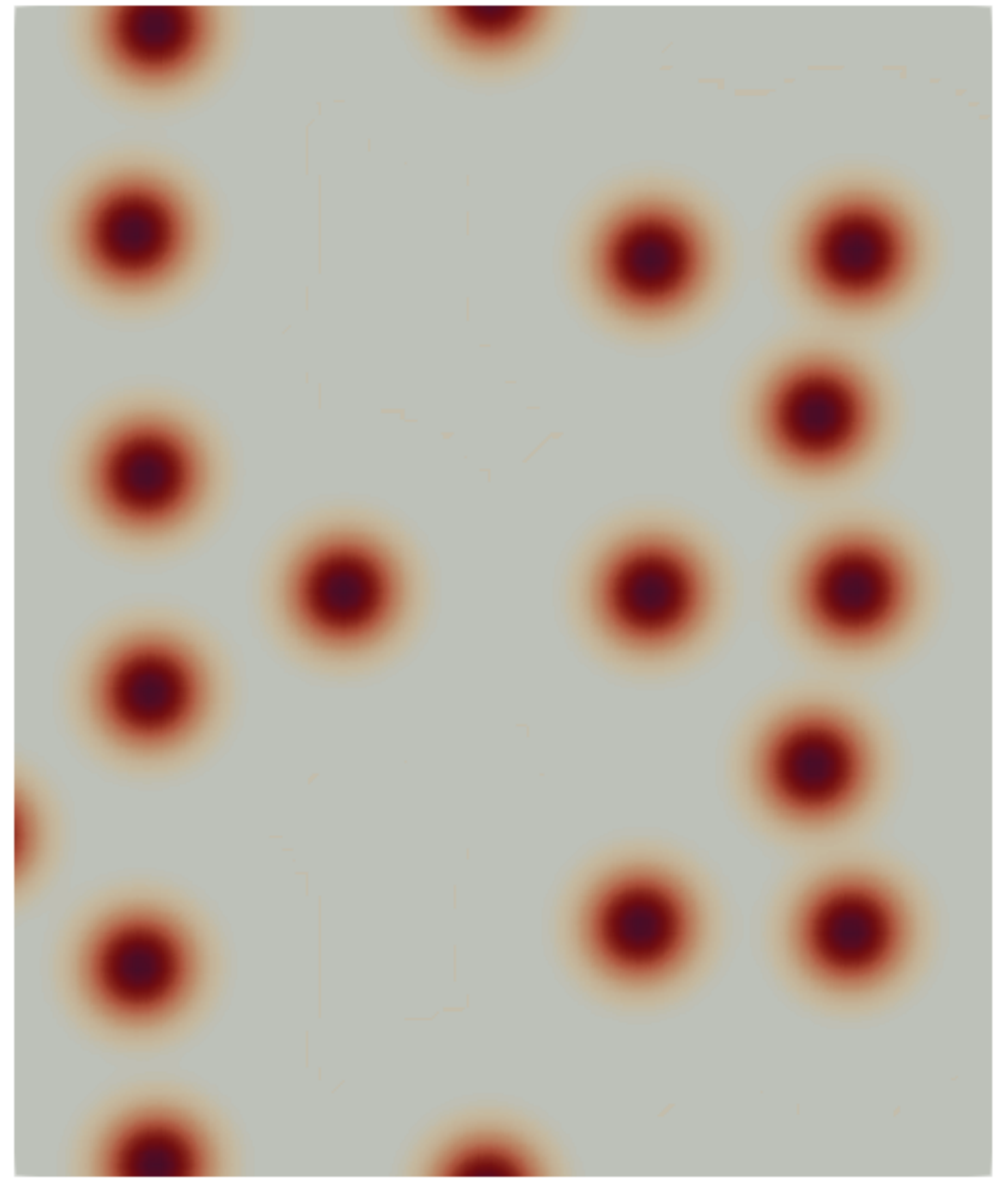}
        \caption{Optimal design $\alpha = 5\times 10^{-3}$}
    \end{subfigure}
    \begin{subfigure}{0.19\textwidth}
        \centering
        \includegraphics[width=0.9\textwidth]{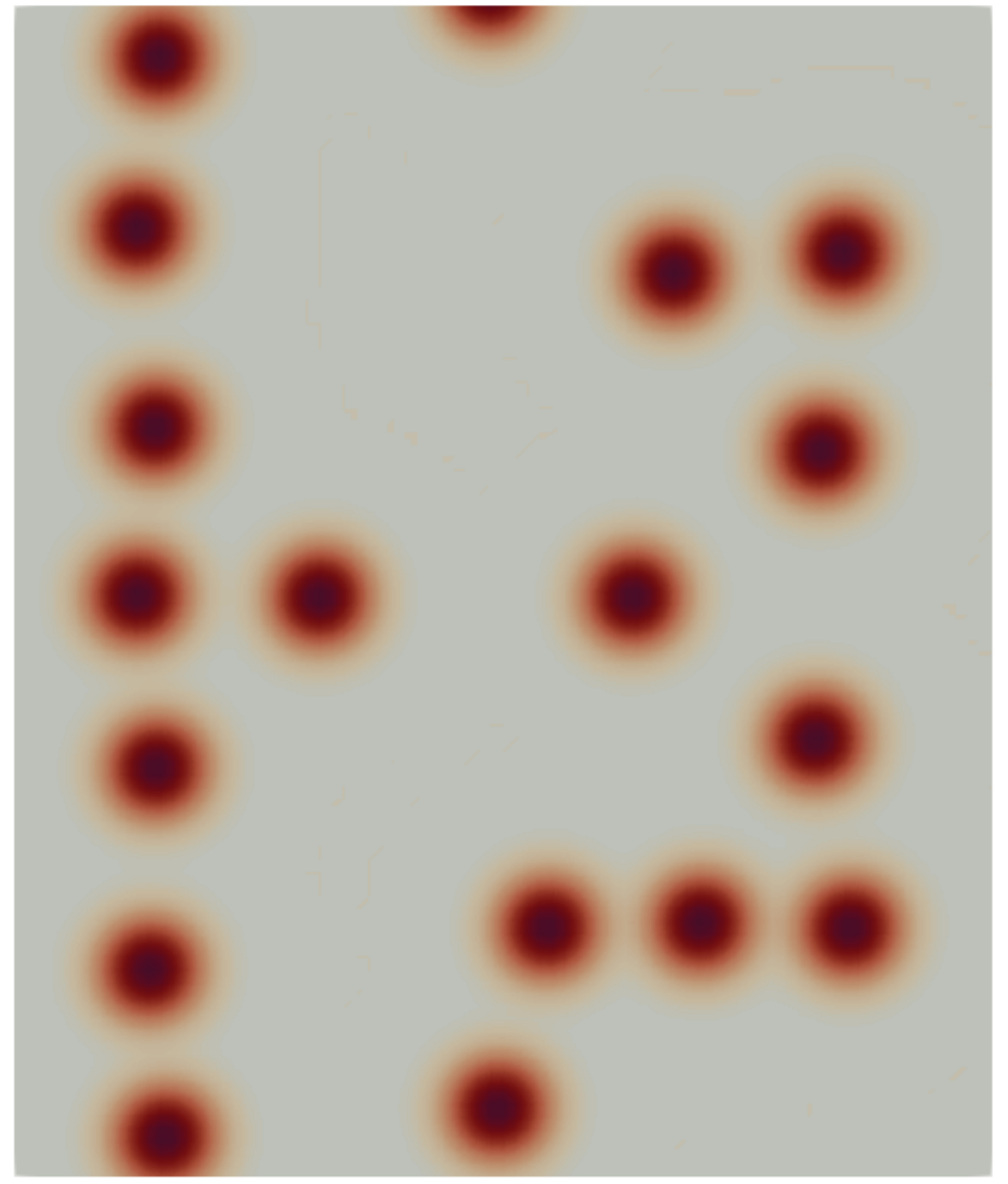}
        \caption{Optimal design $\alpha = 1\times 10^{-3}$}
    \end{subfigure}
    \begin{subfigure}{0.19\textwidth}
        \centering
        \includegraphics[width=0.9\textwidth]{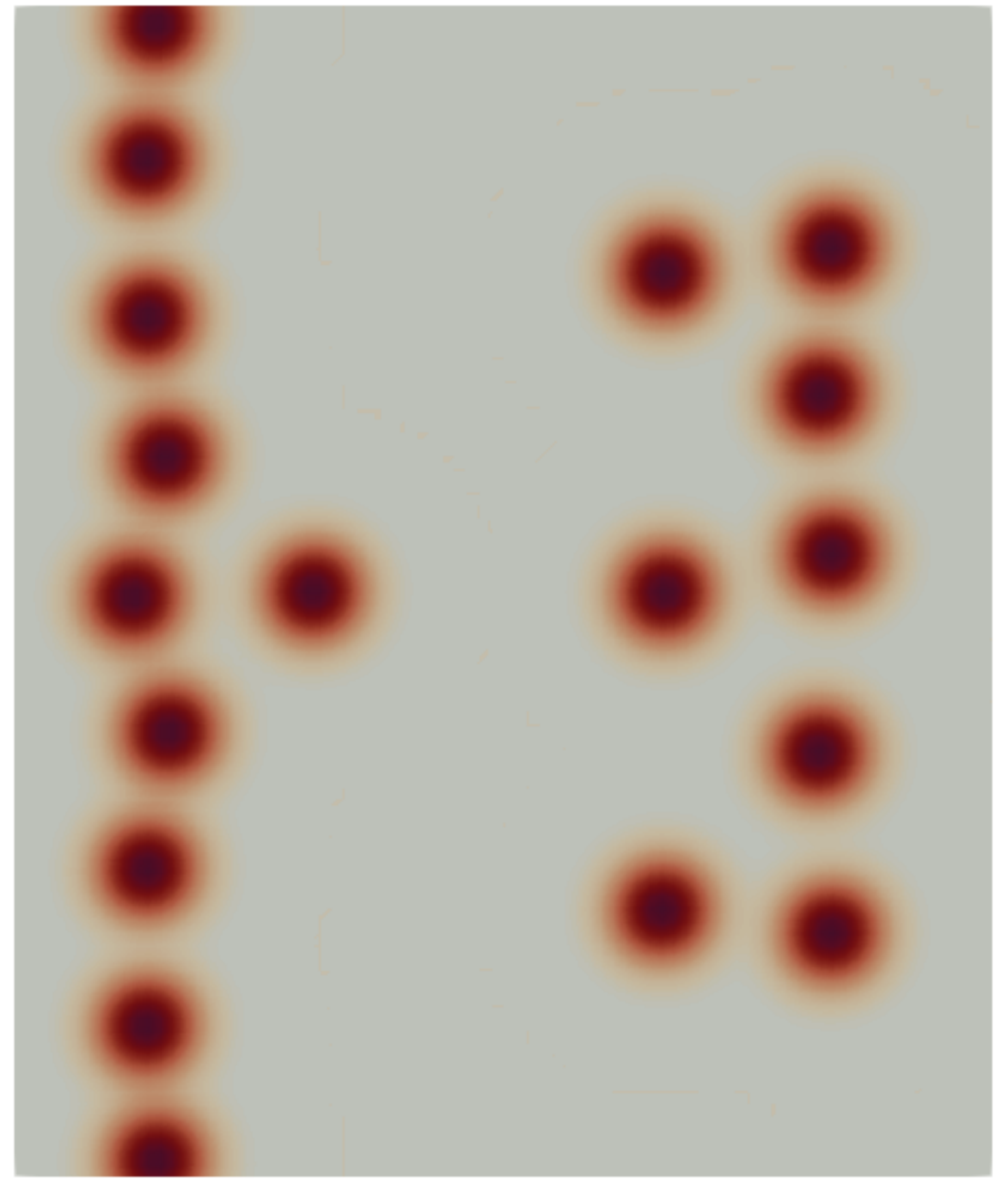}
        \caption{Optimal design $\alpha = 5\times 10^{-4}$}
    \end{subfigure}

    \begin{subfigure}{0.32\textwidth}
        \centering
        \includegraphics[width=1.0\textwidth]{figures_pdf/state_colors.pdf}
    \end{subfigure}

    \begin{subfigure}{0.19\textwidth}
        \centering
        \includegraphics[width=0.9\textwidth]{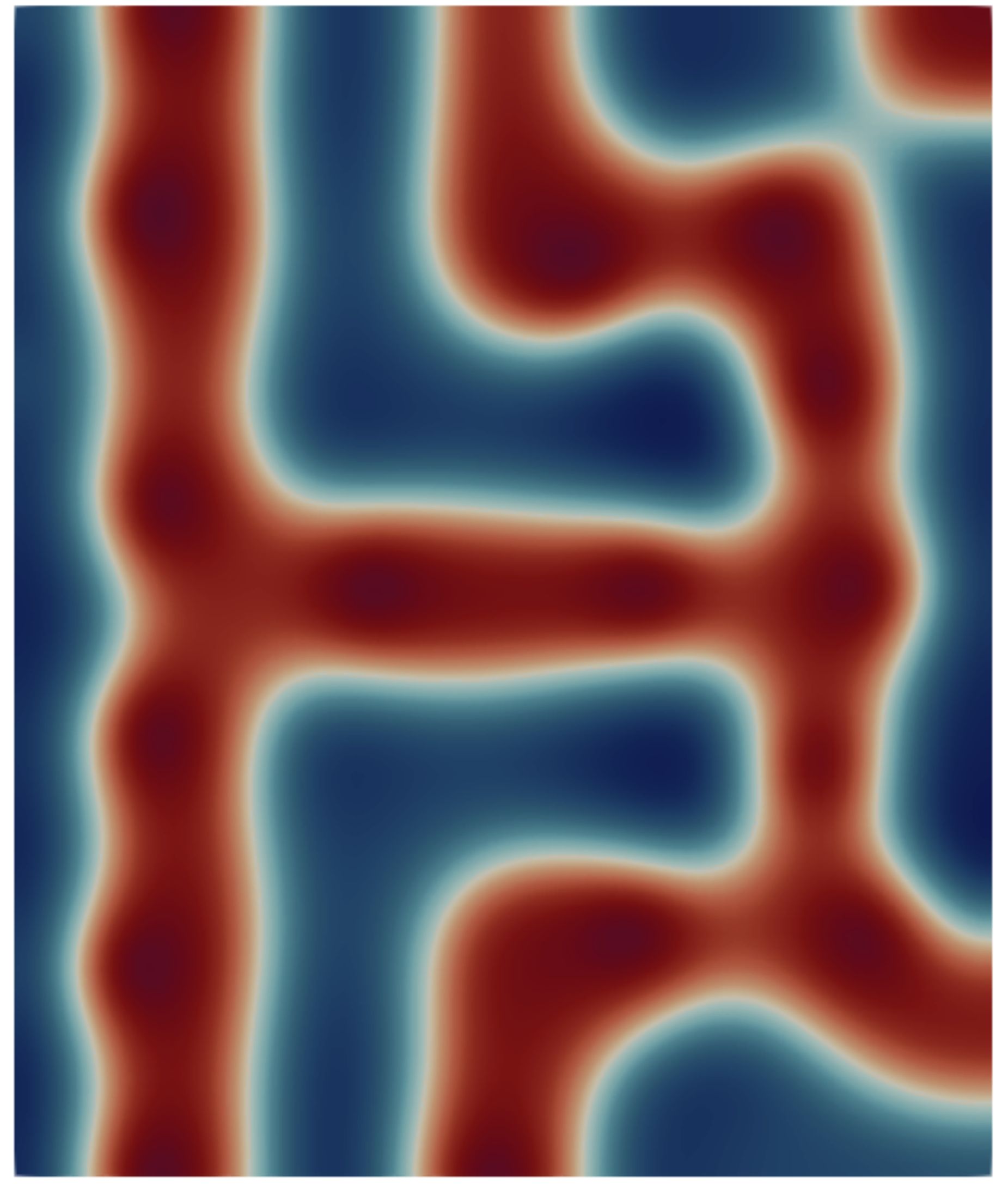}
        \caption{Minimum $\mathcal{F}(u)$ sample $\alpha = 1\times 10^{-2}$}
    \end{subfigure}
    \begin{subfigure}{0.19\textwidth}
        \centering
        \includegraphics[width=0.9\textwidth]{figures_pdf/junction2d/elbow_state_18spots_5e-3_min.pdf}
        \caption{Minimum $\mathcal{F}(u)$ sample $\alpha = 5\times 10^{-3}$}
    \end{subfigure}
    \begin{subfigure}{0.19\textwidth}
        \centering
        \includegraphics[width=0.9\textwidth]{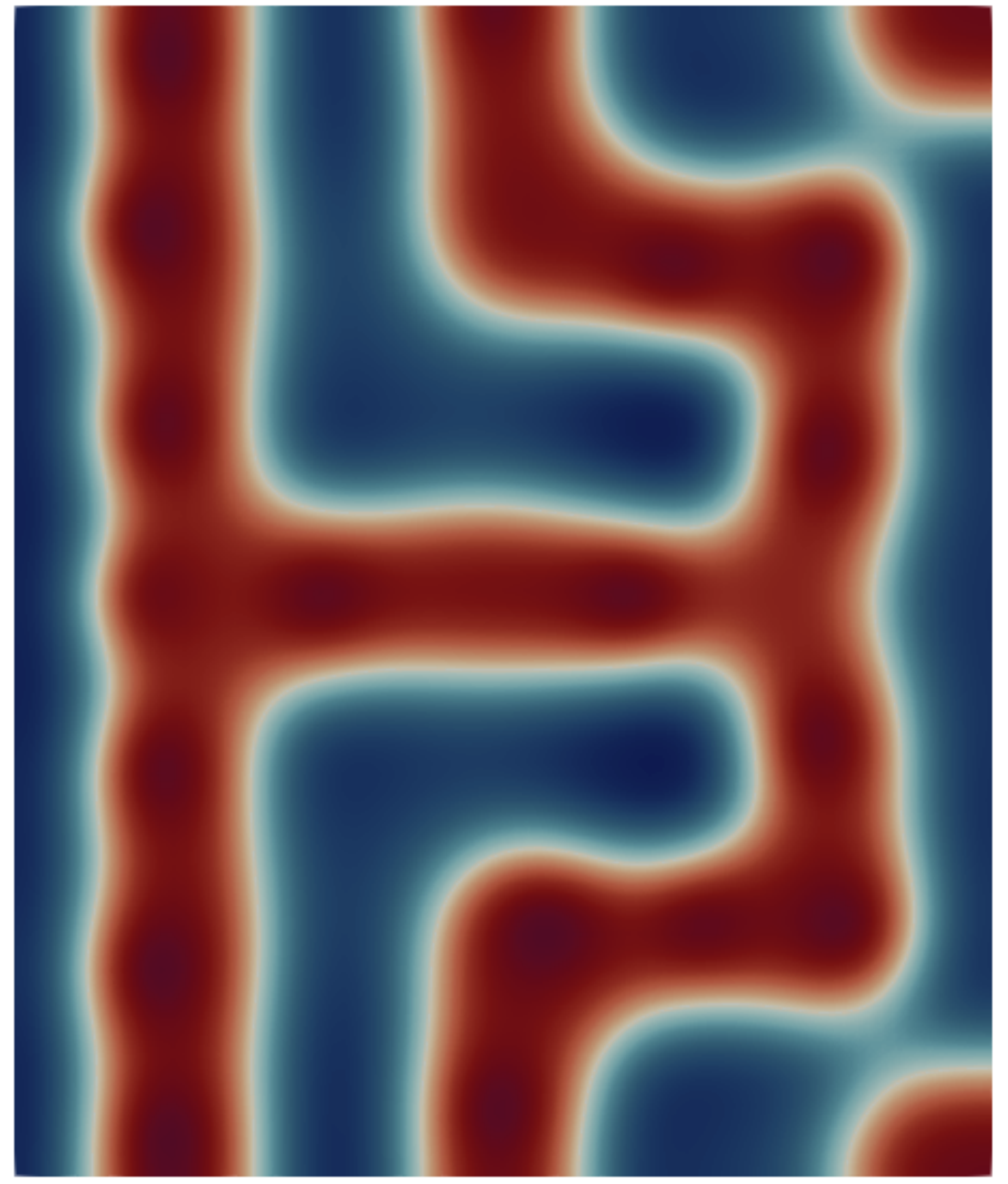}
        \caption{Minimum $\mathcal{F}(u)$ sample $\alpha = 1\times 10^{-3}$}
    \end{subfigure}
    \begin{subfigure}{0.19\textwidth}
        \centering
        \includegraphics[width=0.9\textwidth]{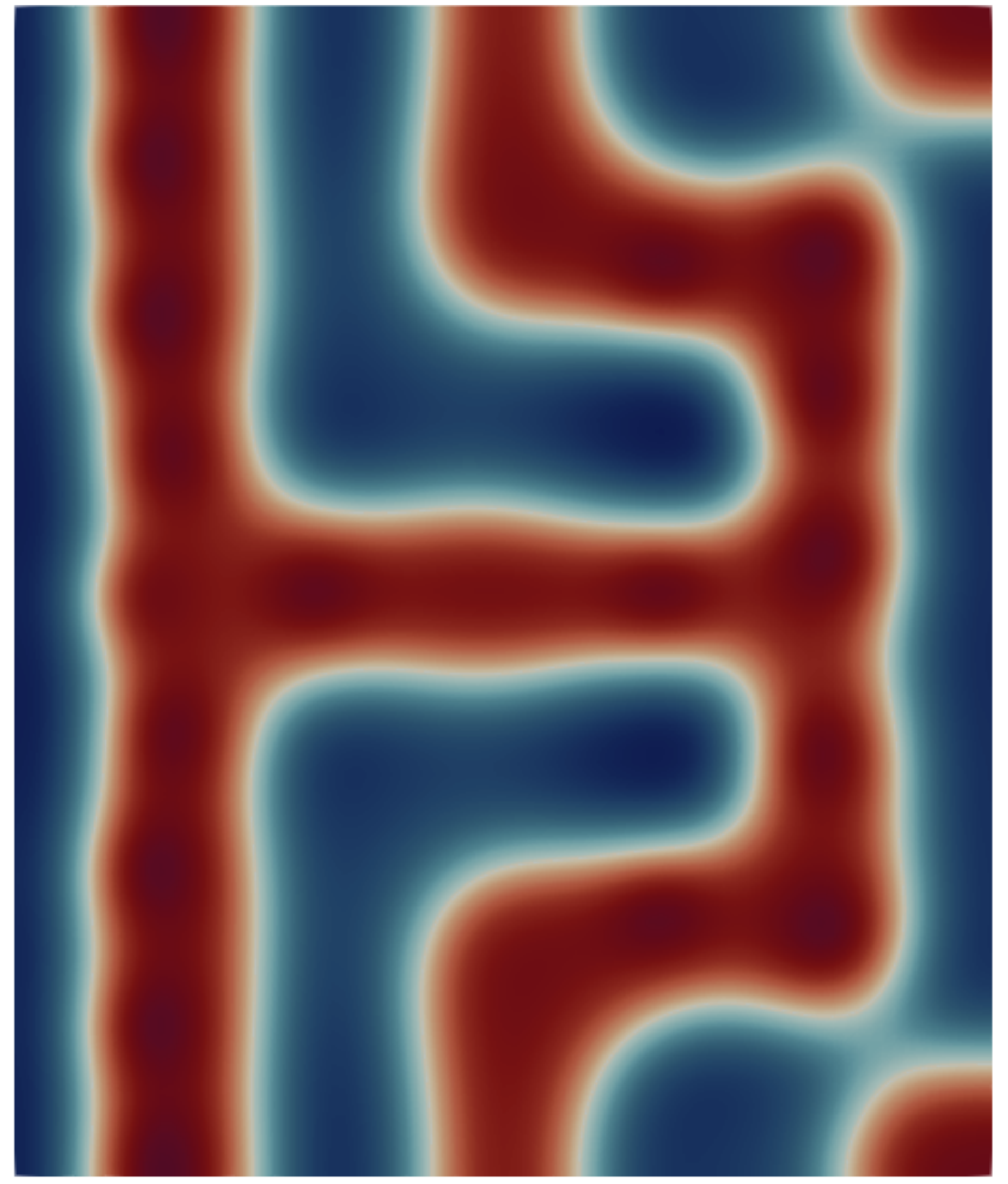}
        \caption{Minimum $\mathcal{F}(u)$ sample $\alpha = 5\times 10^{-4}$}
    \end{subfigure}

    \begin{subfigure}{0.19\textwidth}
        \centering
        \includegraphics[width=0.9\textwidth]{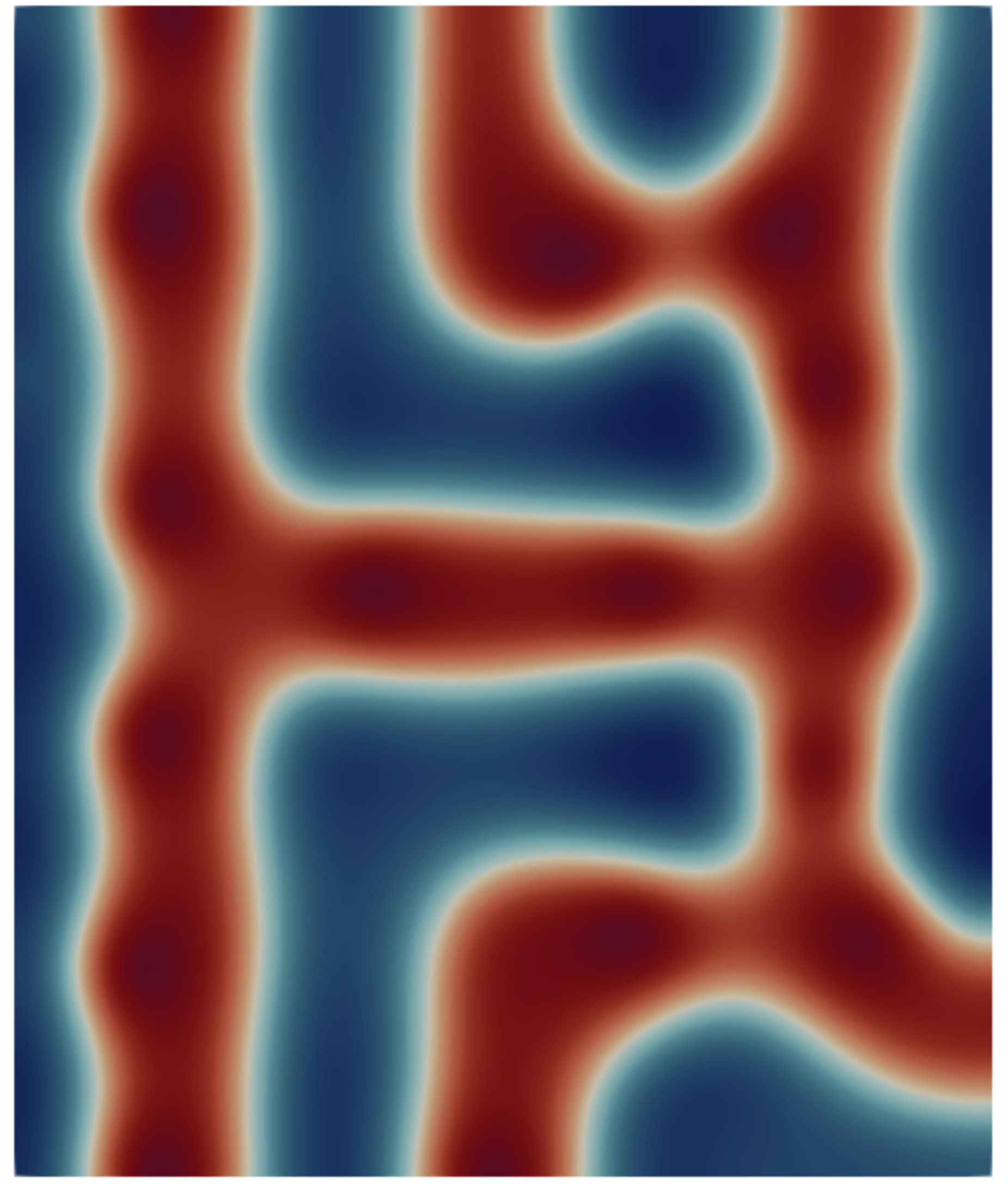}
        \caption{Sample 2 $\alpha = 1 \times 10^{-2}$}
    \end{subfigure}
    \begin{subfigure}{0.19\textwidth}
        \centering
        \includegraphics[width=0.9\textwidth]{figures_pdf/junction2d/elbow_state_18spots_5e-3_0.pdf}
        \caption{Sample 2 $\alpha = 5 \times 10^{-3}$}
    \end{subfigure}
    \begin{subfigure}{0.19\textwidth}
        \centering
        \includegraphics[width=0.9\textwidth]{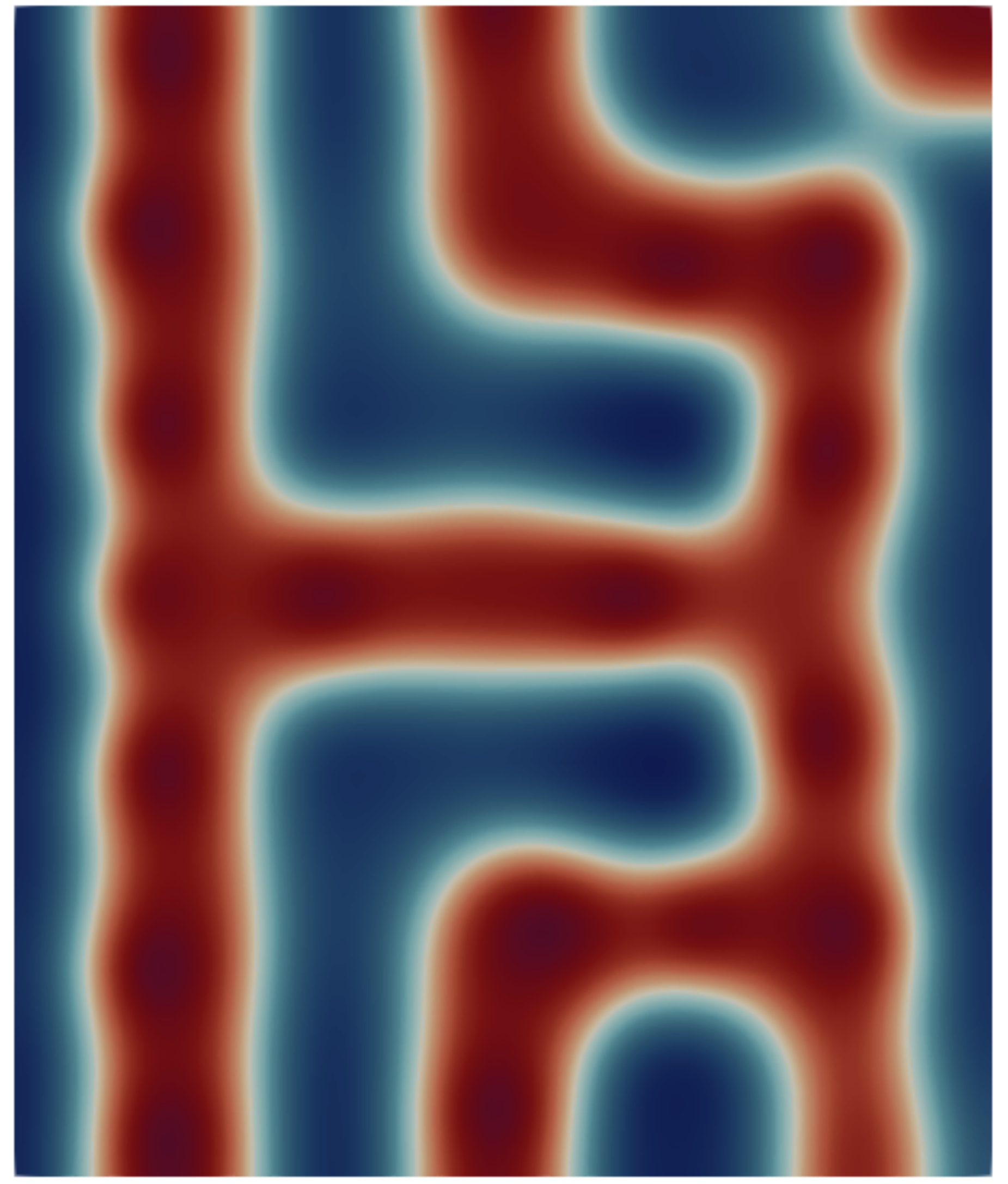}
        \caption{Sample 2 $\alpha = 1 \times 10^{-3}$}
    \end{subfigure}
    \begin{subfigure}{0.19\textwidth}
        \centering
        \includegraphics[width=0.9\textwidth]{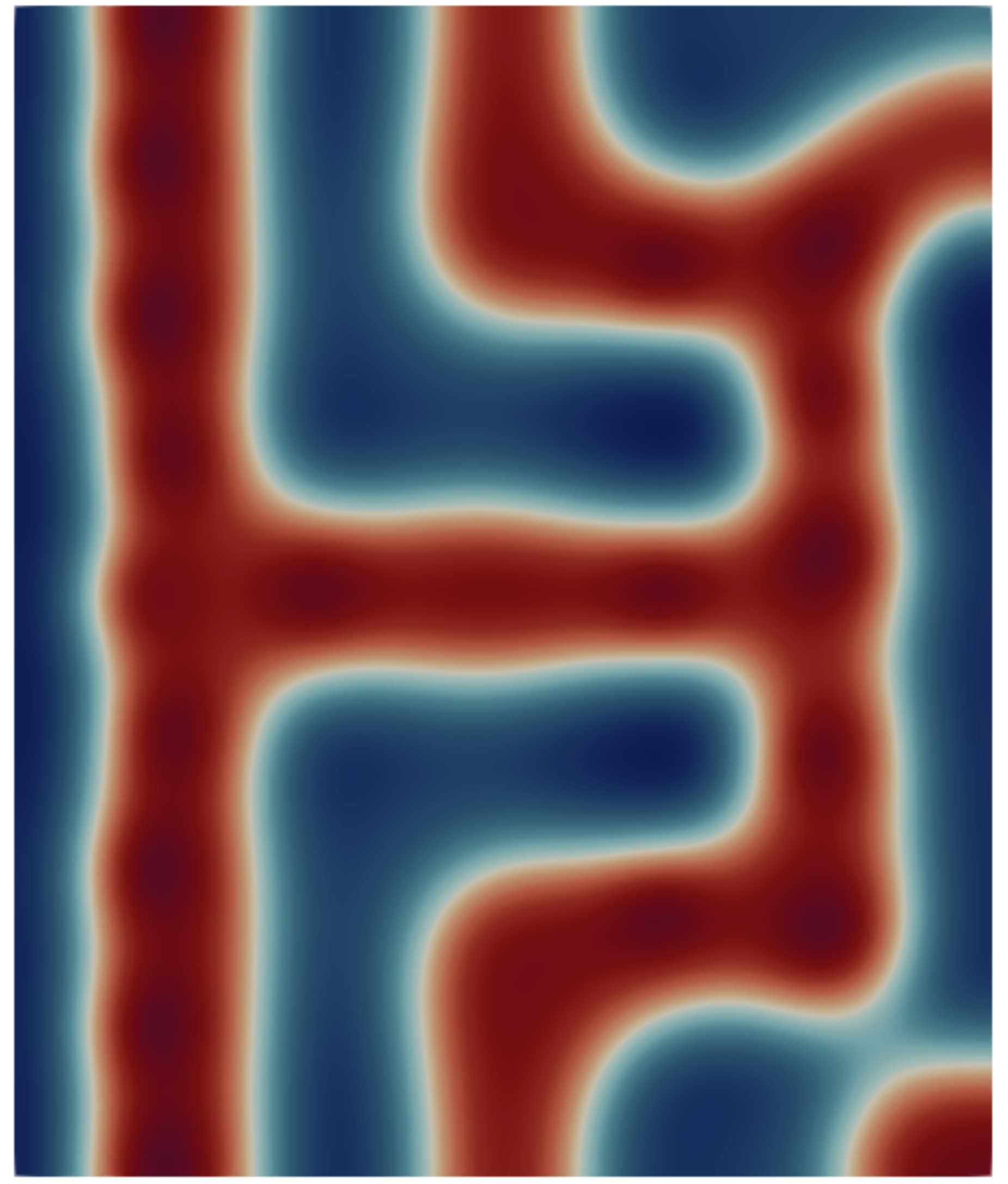}
        \caption{Sample 2 $\alpha = 5 \times 10^{-4}$}
    \end{subfigure}

    \begin{subfigure}{0.19\textwidth}
        \centering
        \includegraphics[width=0.9\textwidth]{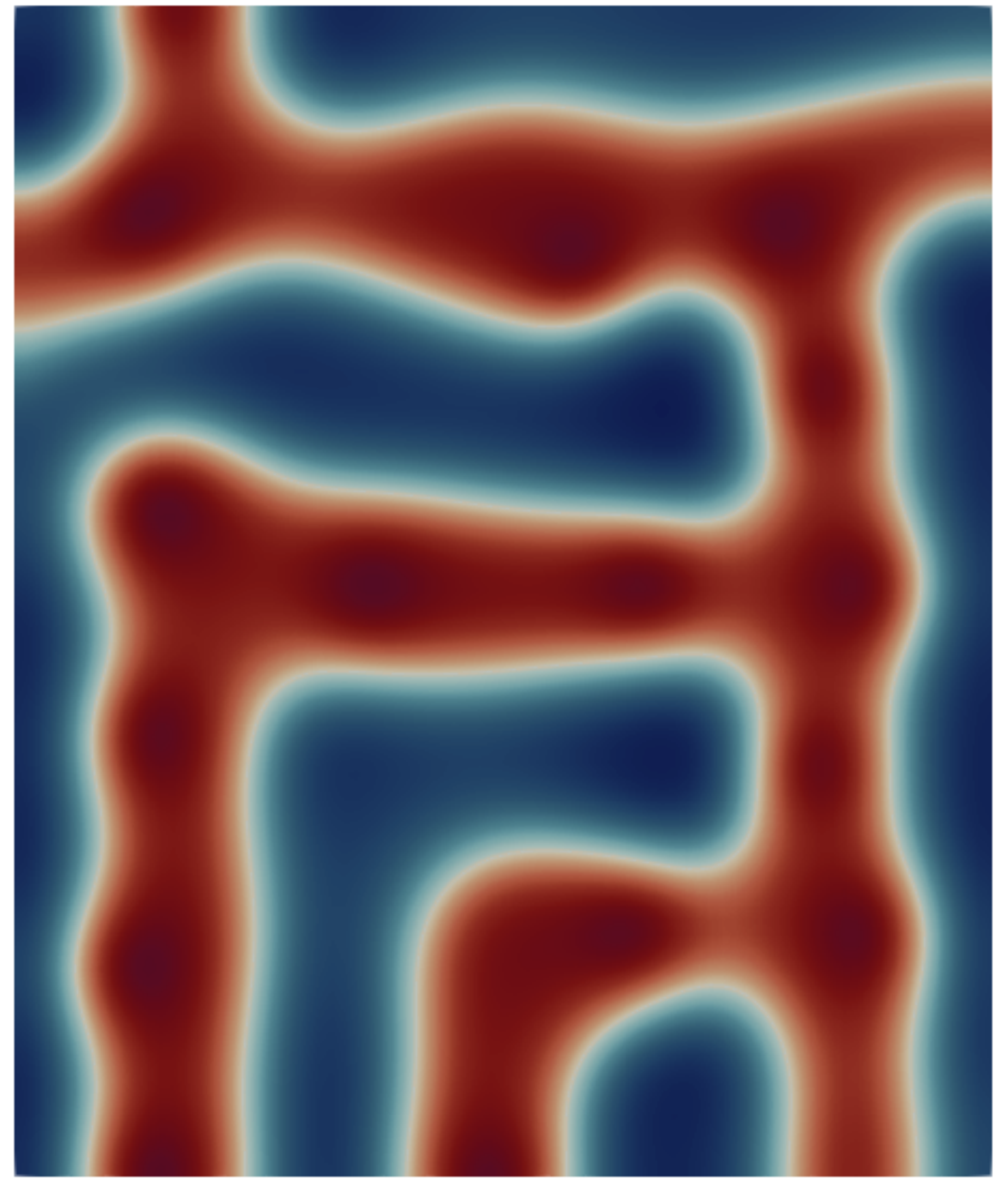}
        \caption{Sample 3 $\alpha = 1 \times 10^{-2}$}
    \end{subfigure}
    \begin{subfigure}{0.19\textwidth}
        \centering
        \includegraphics[width=0.9\textwidth]{figures_pdf/junction2d/elbow_state_18spots_5e-3_1.pdf}
        \caption{Sample 3 $\alpha = 5 \times 10^{-3}$}
    \end{subfigure}
    \begin{subfigure}{0.19\textwidth}
        \centering
        \includegraphics[width=0.9\textwidth]{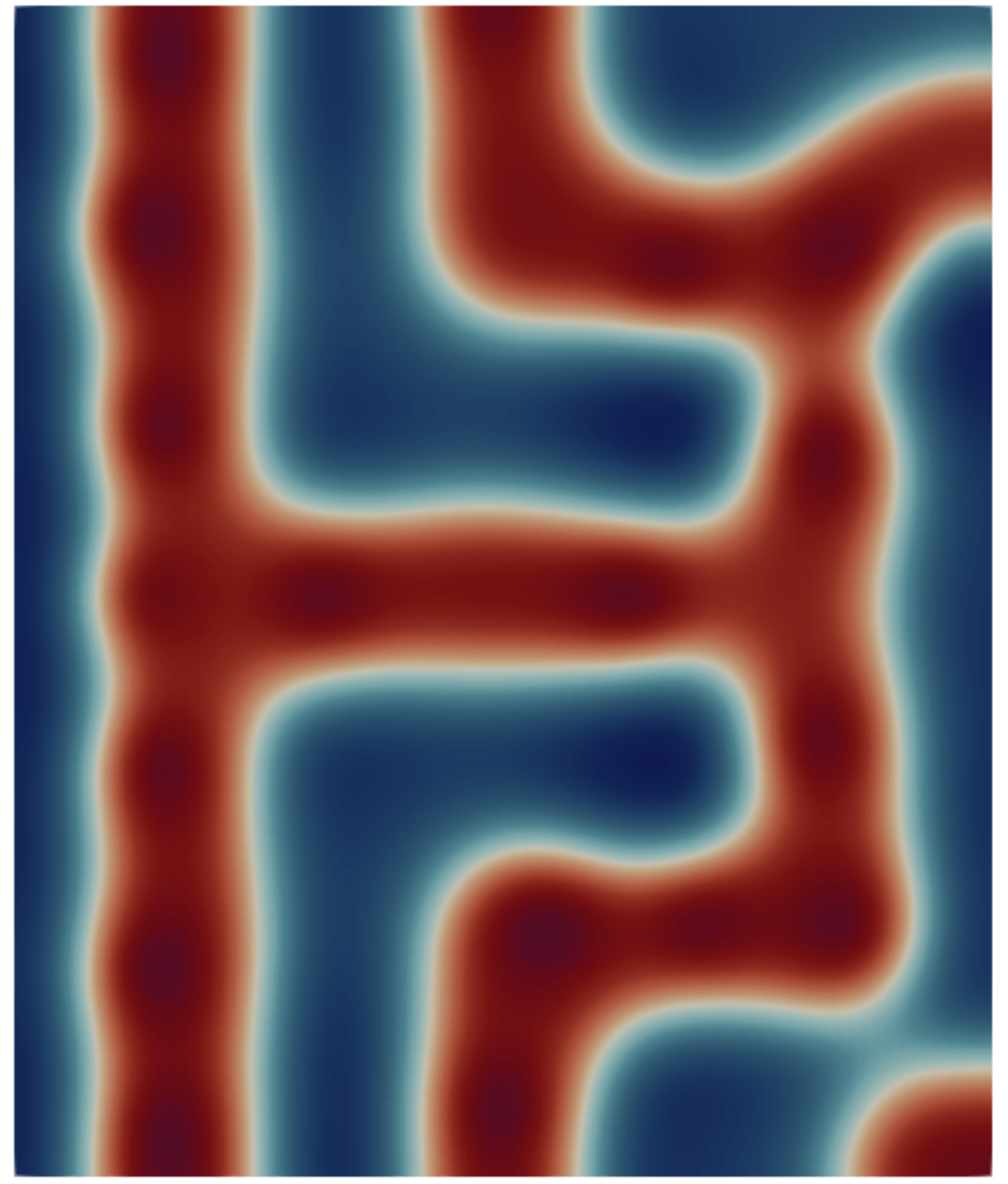}
        \caption{Sample 3 $\alpha = 1 \times 10^{-3}$}
    \end{subfigure}
    \begin{subfigure}{0.19\textwidth}
        \centering
        \includegraphics[width=0.9\textwidth]{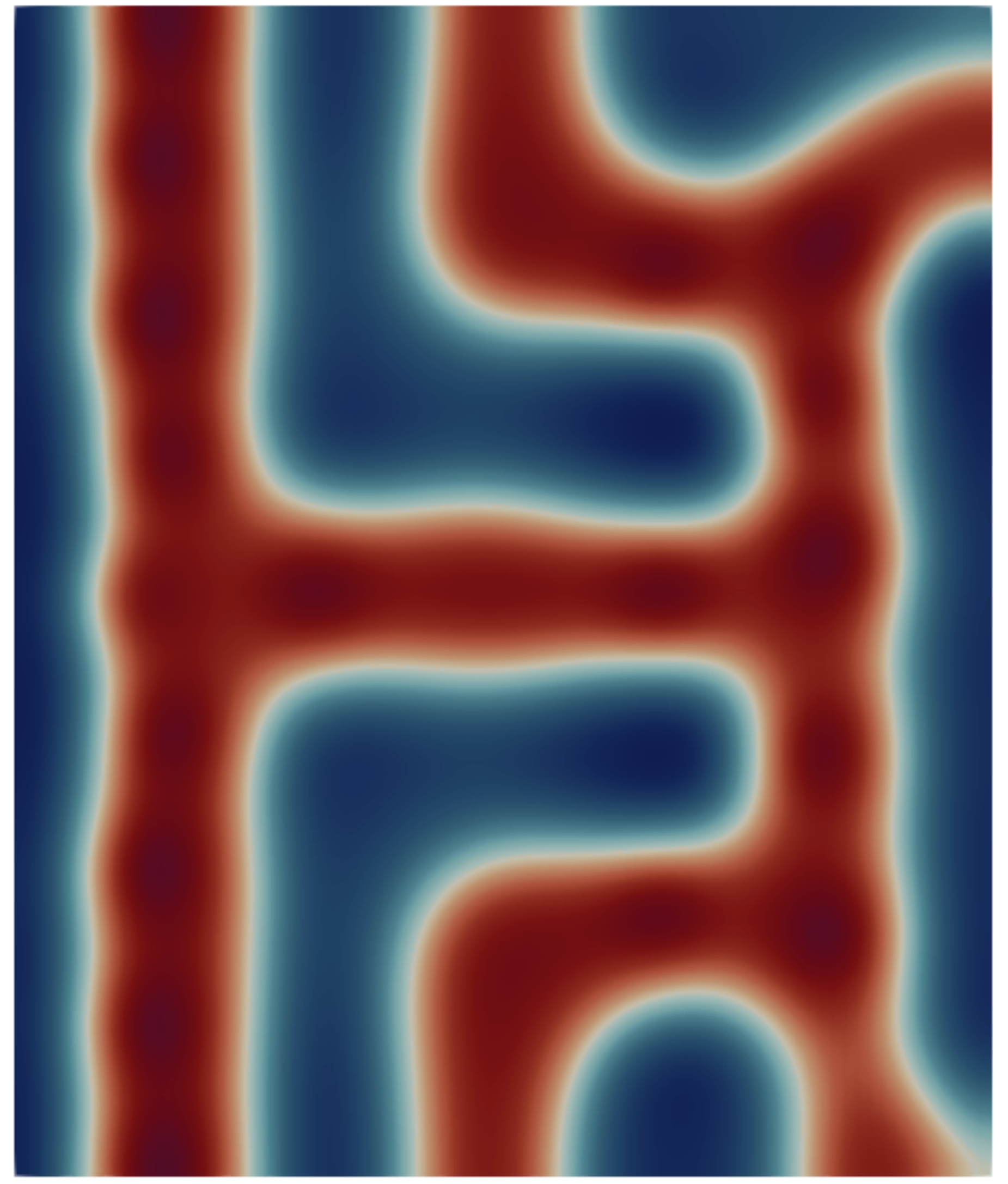}
        \caption{Sample 3 $\alpha = 5 \times 10^{-4}$}
    \end{subfigure}

    \caption{Optimal designs and sample states for the junction target morphology using repulsion strengths $\alpha = 10^{-2}, 5\times 10^{-3}, 10^{-3}, 5\times 10^{-4}$. The bottom three rows show sample equilibrium states computed from different sample initial guesses $u_0$ at the optimal designs, with the minimum energy states shown first.}
    \label{fig:junction2d_alpha}
\end{figure}

We also consider separately the influence of the repelling potential by solving the optimal design problem over a range of $\alpha = $ $10^{-2}$, $5\times 10^{-3}$, $10^{-3}$, $5 \times 10^{-4}$ using 18 guideposts. We present the results for just the junction target morphology, showing the optimal designs and sample states in Figure \ref{fig:junction2d_alpha}. 

We observe that as $\alpha$ reduces, guideposts are allowed to be closer to each other in the optimal designs. As there are numerous guideposts in this setup, the optimal designs using small $\alpha$ (small repulsion) can reliably achieve sample states close to the target morphology. On the other hand, with large $\alpha$ (repulsion), the design is dominated by the effect of the repulsion, pushing guideposts out of alignment with the target morphology in order to maintain separation. 
This is contrasted with the previous example for the strip target morphology, where large values of $\alpha$ are required to separate the guideposts and navigate the numerous local minima. 
Therefore, in practice, the $\alpha$ parameter should be tuned based on the target morphology, number of guideposts, as well as the desired level of separation between the guideposts.

% To demonstrate the flexibility of our method for a variety of target morphologies, we consider a different target morphology defined over the domain $\Omega = [0, L] \times [0, 7L/6]$ with $L=3$, as shown in Figure \ref{fig:jog_target}. We solve the optimal design problem with $N_p = $ 12, 14, 16, 18, and 20 guideposts, and present the optimal designs and sample states in Figure \ref{fig:jog_Np}. We can draw similar conclusion as before, i.e., the sample states with minimum energies have less defects \pc{is this explained in the last example?} for all choices of $N_p$, and the robustness of the optimal designs improve as more guideposts are used.

\subsubsection{Assessment of computational cost}

We assess the computational cost of the optimization algorithm by repeating the design optimization for the junction target morphology over 20 different initial guidepost configurations. We consider the cases of $N_p = 12$, 16, 20, 30, and 40. Since we have a large number of guideposts within a relatively small domain for the latter two cases, we adopt a reduced repulsion strength of $\alpha = 5\times 10^{-4}$ such that the problem is not masked by the repulsion. 

\begin{figure}
    \centering
    \begin{subfigure}{0.45\textwidth}
        \centering
        \includegraphics[width=0.9\textwidth]{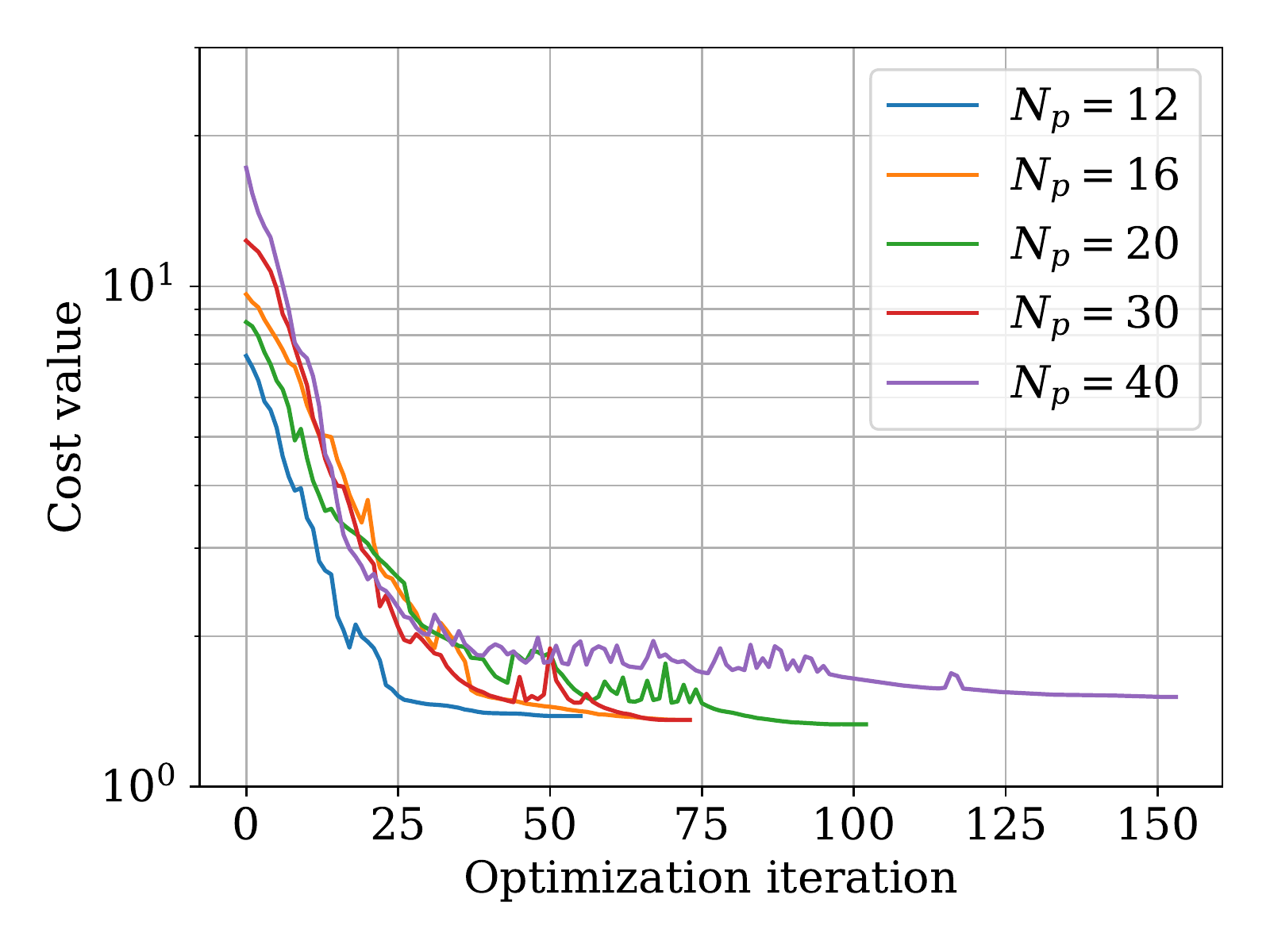}
        % \caption{Cost function}
    \end{subfigure}
    \begin{subfigure}{0.45\textwidth}
        \centering
        \includegraphics[width=0.9\textwidth]{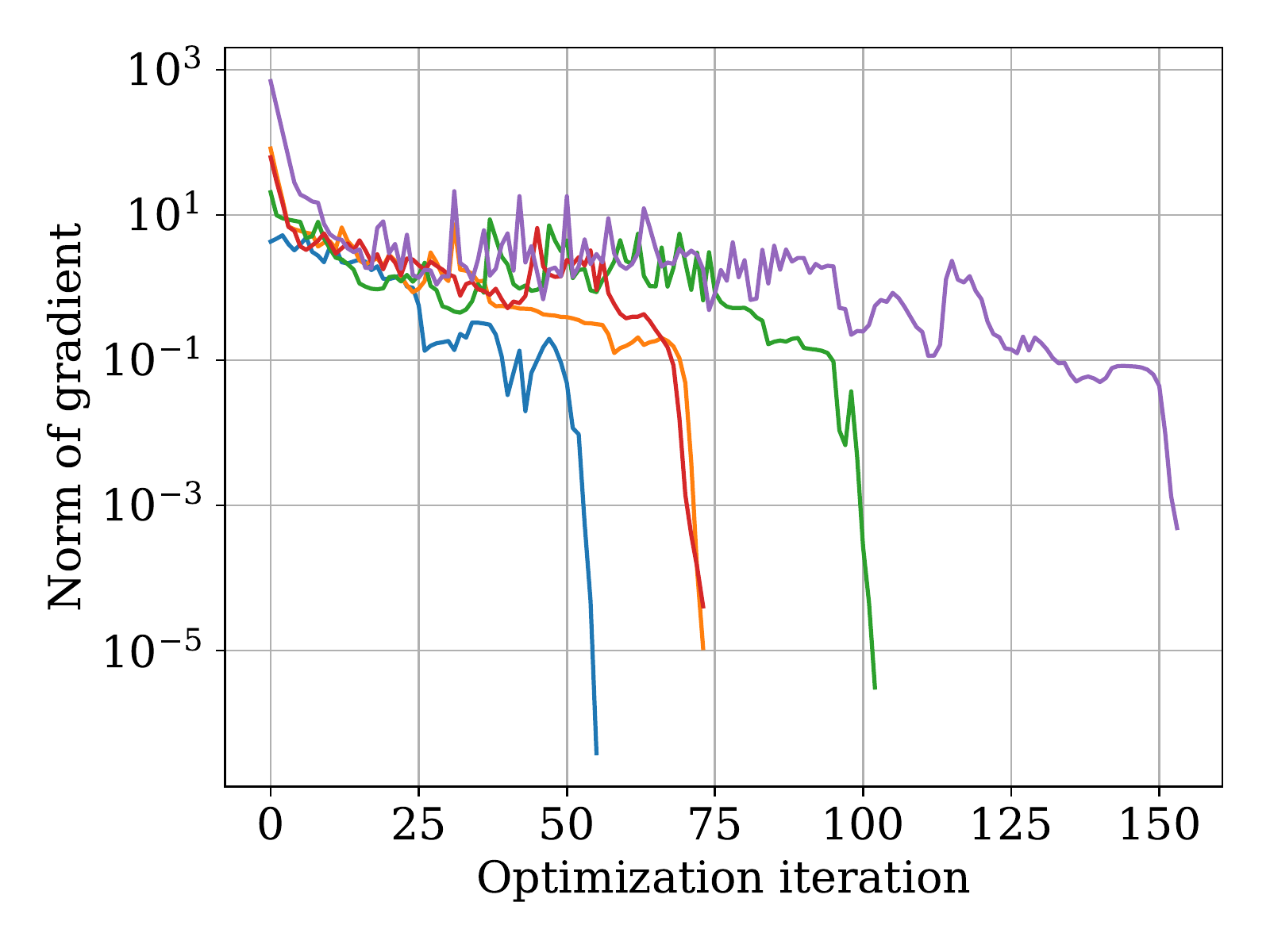}
        % \caption{Gradient norm of cost function}
    \end{subfigure}
    \caption{Values of the cost function (left) and norm of its gradients (right) at each optimization iteration. Plots are shown for three optimization trials with different initial guesses, with the number of guideposts $N_p$ = 12, 16, and 20.}
    \label{fig:opt_iters}
\end{figure}

In Figure \ref{fig:opt_iters} we plot the values of the cost function and the norm of its gradients with respect to the design variables for a single run for each of $N_p = 12$, 16, 20, 30, and 40. We observe that significant reduction of both the cost and gradient norm is typically achieved within 100 optimization iterations. The gradient norm decreases rapidly in the first few iterations where the guideposts separate quickly due to the repulsion penalty, and in the last few iterations where we observe the rapid (quadratic) convergence of the outer Newton algorithm close to the local minimum.

\begin{figure}
    \centering
    \begin{subfigure}{0.45\textwidth}
        \centering
        \includegraphics[width=0.9\textwidth]{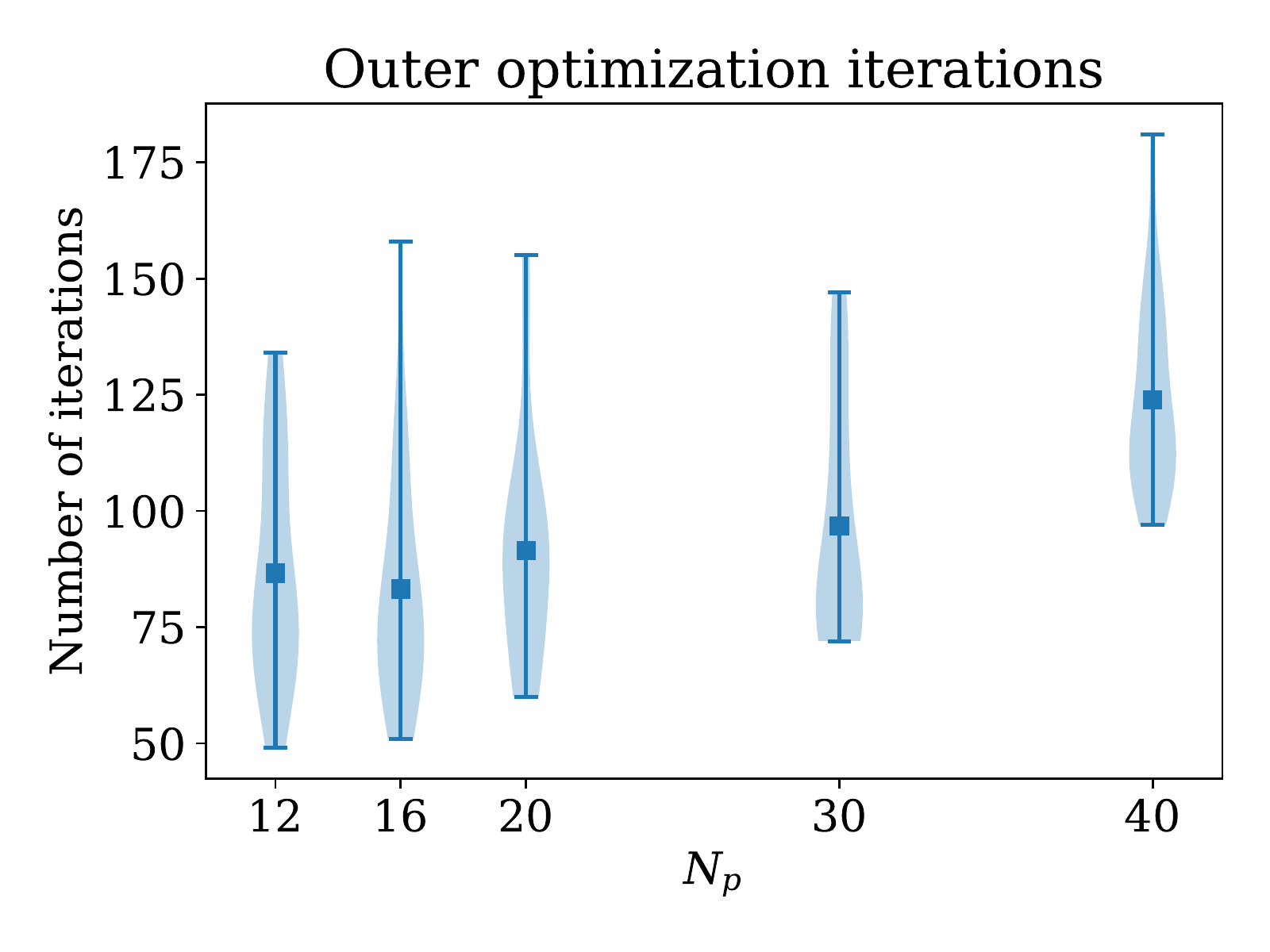}
        % \caption{Optimization}
    \end{subfigure}
    \begin{subfigure}{0.45\textwidth}
        \centering
        \includegraphics[width=0.9\textwidth]{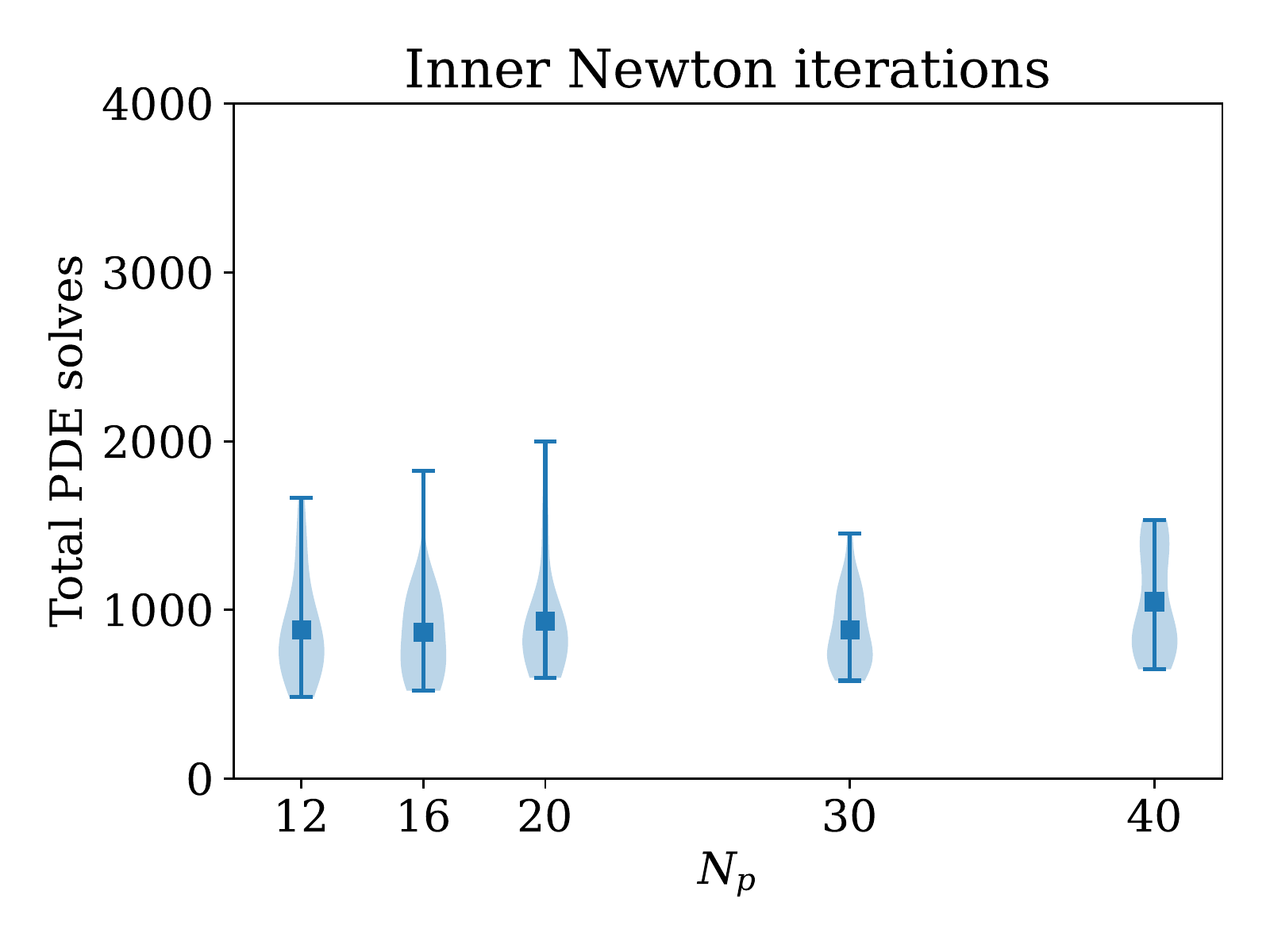}
        % \caption{Gradient norm of cost function}
    \end{subfigure}
    \caption{Distribution of the number of outer optimization iterations (left) and total number linear PDE solves associated with the inner Newton iterations for the state equation over the entire optimization run (right), starting at 20 different initial guidepost configurations. Results are shown for $N_p = 12$, 16, 20, 30, and 40.}
    \label{fig:opt_distributions}
\end{figure}

Furthermore, we plot in Figure \ref{fig:opt_distributions} the distributions from the 20 sample runs for both the number of optimization iterations and the total number of inner (energy-stable) Newton iterations accumulated over the entire optimization run. From these results, we again see that the optimization problem tends to converge in $\mathcal{O}(10^2)$ iterations, requiring on average $\mathcal{O}(10^3)$ linearized state PDE solves in total. The total number of inner Newton iterations remain relatively consistent across the different numbers of guideposts, despite requiring more optimization iterations on average. This suggests that the overall optimization cost does not increase significantly with the number of guideposts and hence dimension of the optimization problem, within the range considered. 

\begin{figure}
    \centering
    \includegraphics[width=0.5\textwidth]{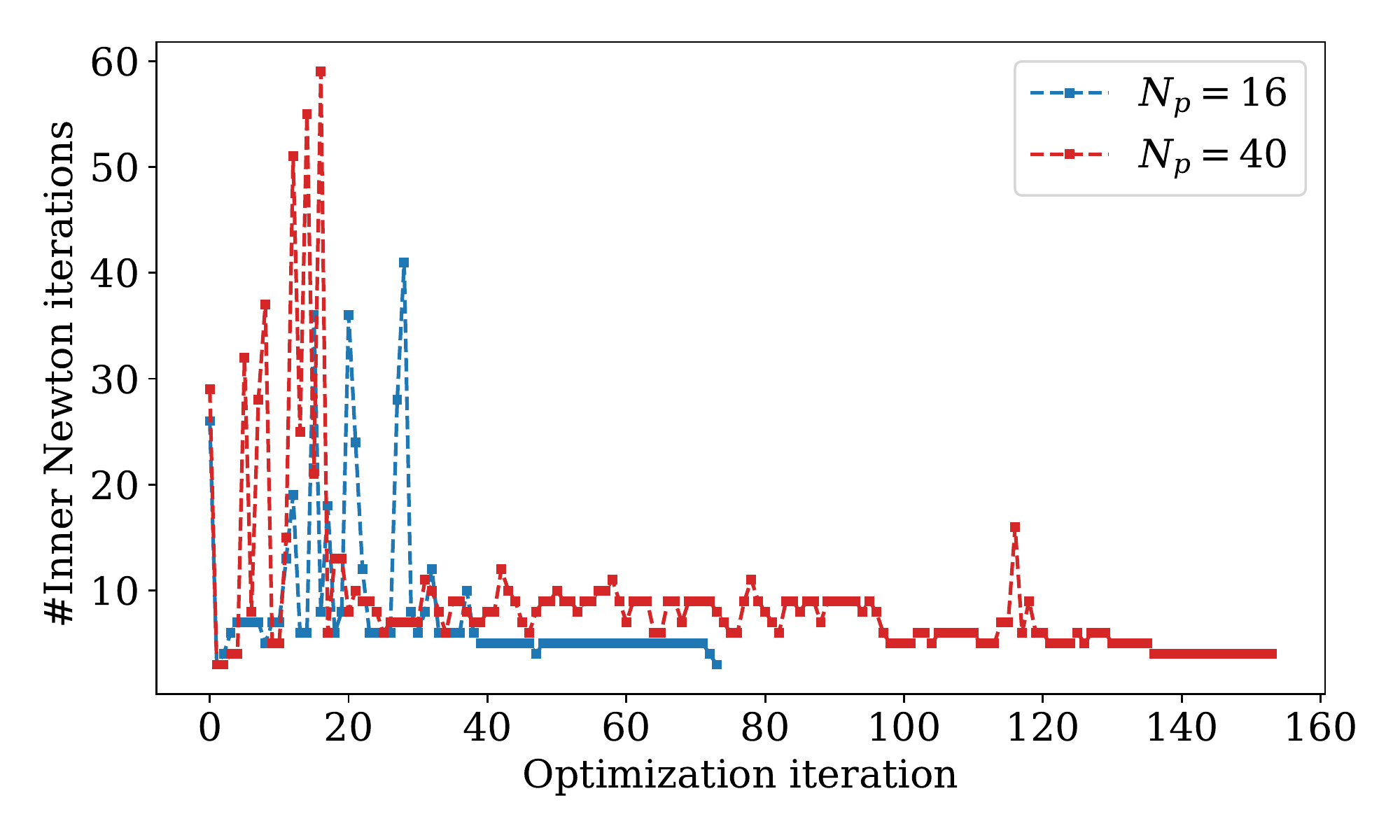}
    % \caption{Optimization}
    \caption{Number of inner Newton iterations at each outer optimization iteration. Results shown for one $N_p = 16$ run and one $N_p=40$ run.}
    \label{fig:opt_solves}
\end{figure}

We compare in detail one optimization run using $N_p = 16$ with another using $N_p = 40$. We plot the number of inner Newton iterations used for solving the state problem at each optimization iteration for the two cases in Figure \ref{fig:opt_solves}. In the figure, we observe peaks in the number of inner Newton iterations followed by a significant reduction for the subsequent steps. This is a result of the continuation scheme used in the optimization. The peaks in the iteration numbers reflect instances when the state moves to a different solution branch, thus requiring a large number of state iterations to converge. Subsequent state PDE solves are then initialized at the converged solution of the previous steps. When the change of the substrate design is small, the subsequent solutions remain close to the previous solutions and hence require relatively few of state iterations. In particular, the peaks tend to be concentrated near the beginning of the optimization where the guidepost configuration is largely suboptimal and does not tend to constrain the solution. On the other hand, the later optimization iterations tend to remain on the same solution branch and require very few state iterations. This explains why the $N_p = 40$ case does not require more state iterations in total, despite requiring more optimization iterations compared to the $N_p = 16$ case. 

We compare our results with existing studies that also adopt the
Cahn--Hilliard formulation, such as \cite{QinKhairaSuEtAl13}, which
optimizes the guidepost locations using the derivative free CMA-ES
algorithm. For a modest-sized problem, this study uses the CMA-ES
algorithm with a population size of 28 ($\mathcal{O}(10^1)$), and
reports convergence to optimal designs in $\mathcal{O}(10^2)$
generations. As a conservative estimate, this requires a total of
$\mathcal{O}(10^3)$ evaluations of the non-linear, time-evolution PDE,
for which the study uses $\mathcal{O}(10^5)$ time steps to solve. This
amounts to $\mathcal{O}(10^8)$ time-stepping iterations in total, each
of which is comparable to a single energy-stable Newton minimization step. 

In contrast, our formulation attains five orders of magnitude
reduction in cost ($\mathcal{O}(10^3)$ state Newton iterations versus
$\mathcal{O}(10^8)$ state time steps). Computational efficiency in our
approach is achieved through a combination of the fast energy-stable
Newton algorithm for solving the state problem, the continuation
scheme during optimization, and the inexact Newton design optimization algorithm that is tailored to capitalize on the energy-minimization formulation and the continuation scheme for the state problem. We also remark that in our approach, the optimization problem is solved to local optimality due to the rapid local convergence of the Newton algorithm, while heuristic algorithms like the CMA-ES only provide approximate solutions.

% At every optimization iteration, we need to solve once the state and adjoint PDE to obtain the value and gradient of the cost functional, as well 

\subsection{Design in 3D}
Finally, we leverage the efficiency of the proposed computational procedure and solve the optimization problem using the full three-dimensional model. As far as we know, the problem of optimizing for the substrate design using a 3D model has not been studied beyond simple target morphologies with low dimensional optimization variables, such as that in \cite{KhairaQinGarnerEtAl14}, because of the prohibitive computational cost of conventional methods. We define 3D target morphologies by the extrusion of the 2D junction and jog morphologies in the domain $\Omega = [0, L] \times [0, 7L/6] \times [0, L/8]$. This is discretized by linear tetrahedral elements, leading to 204,282 degrees of freedom for the combined state variables $(u_h, \mu_h)$. We again use circular guideposts of width $b=0.08$ and set $c_s = 0.75$ in \eqref{eq:ds} as the parameter defining the decay length scale of the substrate attraction. For the repulsion strength, we adopt $\alpha = 5\times 10^{-3}$ for the junction design and $\alpha = 10^{-3}$ for the jog design. 

% We consider a 3D target morphology by extrusion of the 2D junction morphology in the domain $\Omega = [0, L] \times [0, 7L/6] \times [0, L/8]$. We set $d_s = 0.8$ in \eqref{eq:tau-3} as the parameter characterizing the decay of the attractive strength of the substrate. 

\begin{figure}[htbp!]
    \centering
    \begin{subfigure}{0.32\textwidth}
        \includegraphics[width=0.99\textwidth]{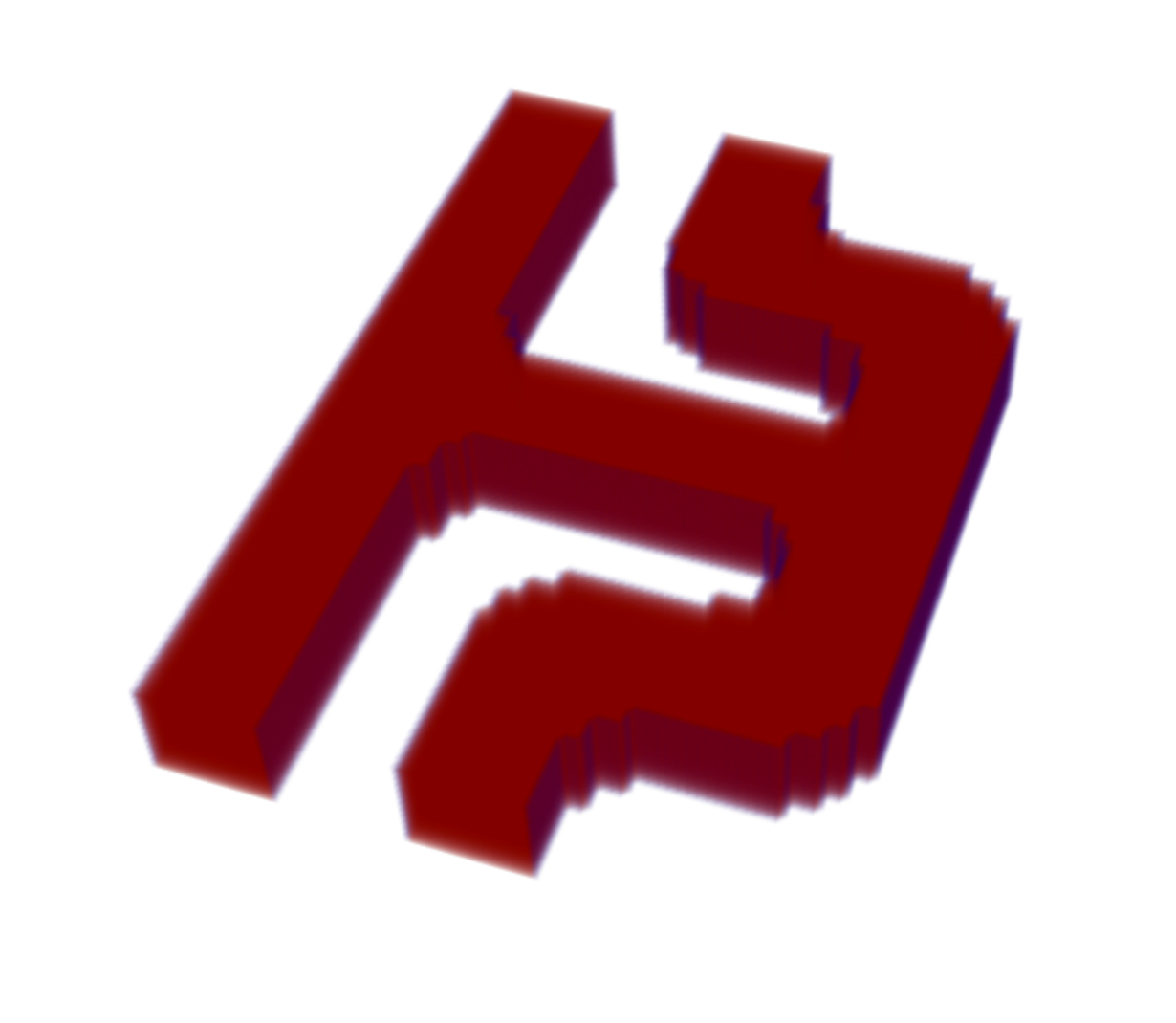}
        \caption{3D junction target morphology}
    \end{subfigure}
    \begin{subfigure}{0.32\textwidth}
        \centering
        \includegraphics[width=0.99\textwidth]{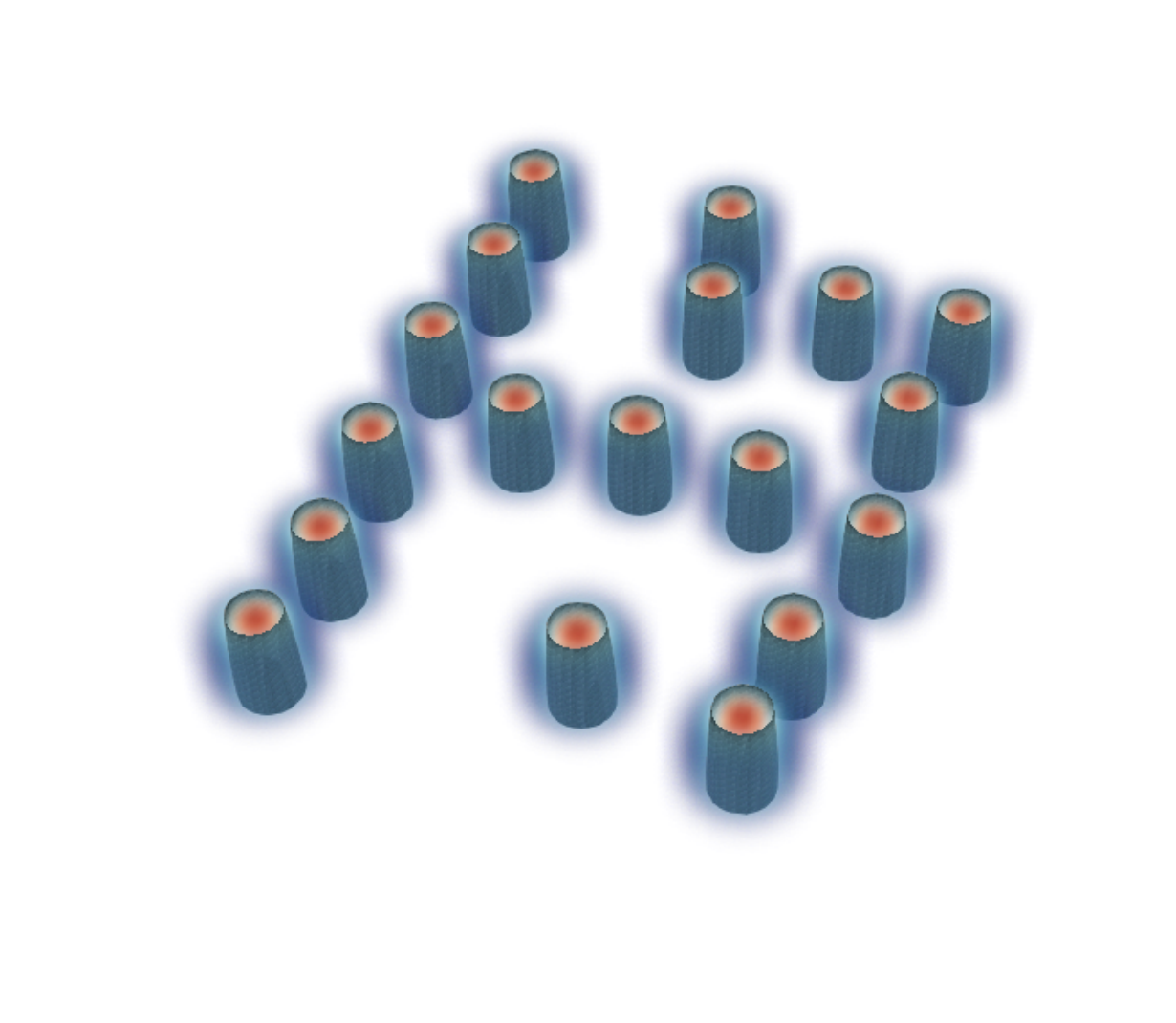}
        \caption{Substrate strength at optimal design}
    \end{subfigure}

    \begin{subfigure}{0.32\textwidth}
        \centering
        \includegraphics[width=0.99\textwidth]{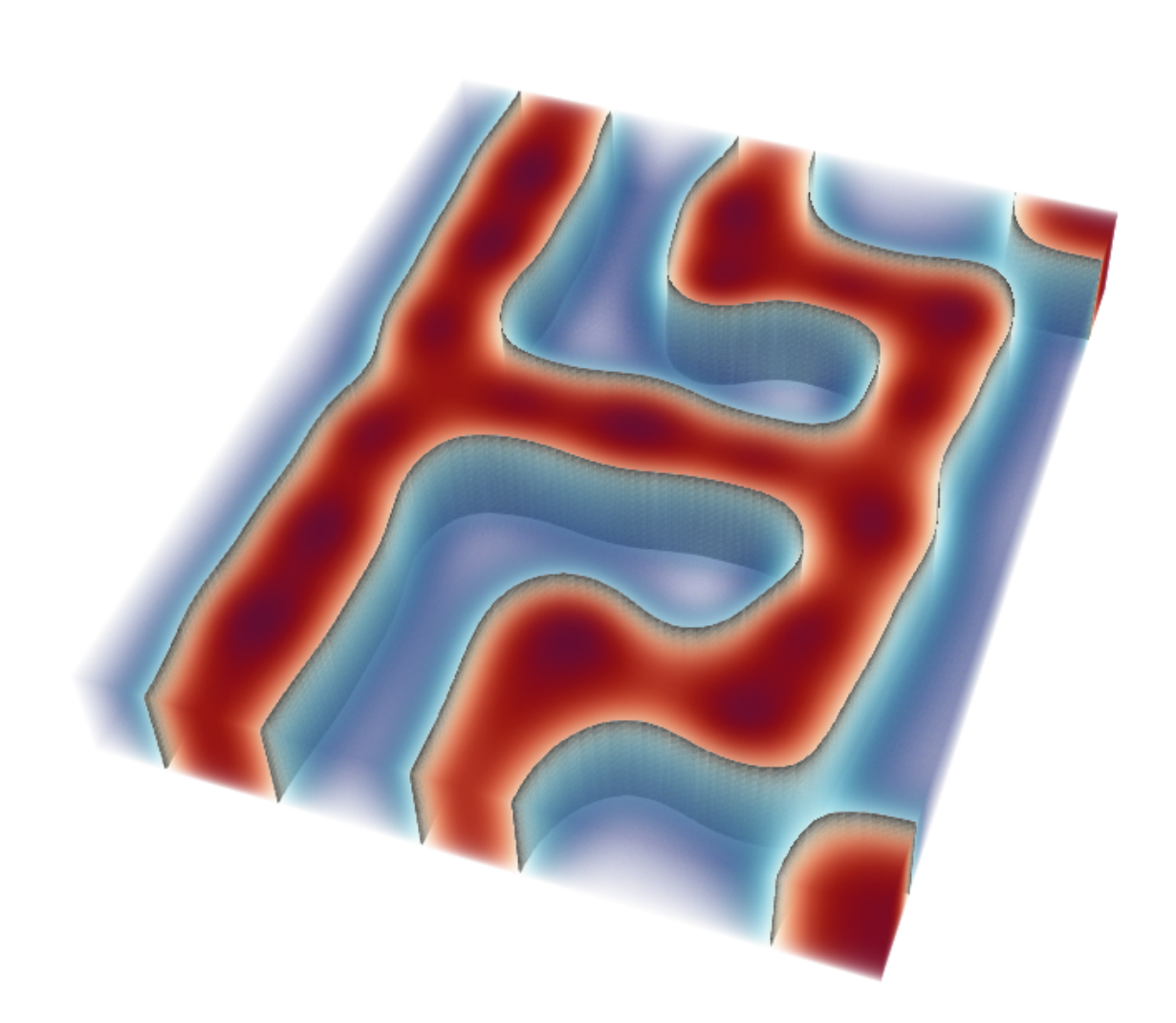}
        \caption{Minimum $\mathcal{F}(u)$ sample}
    \end{subfigure}
    \begin{subfigure}{0.32\textwidth}
        \centering
        \includegraphics[width=0.99\textwidth]{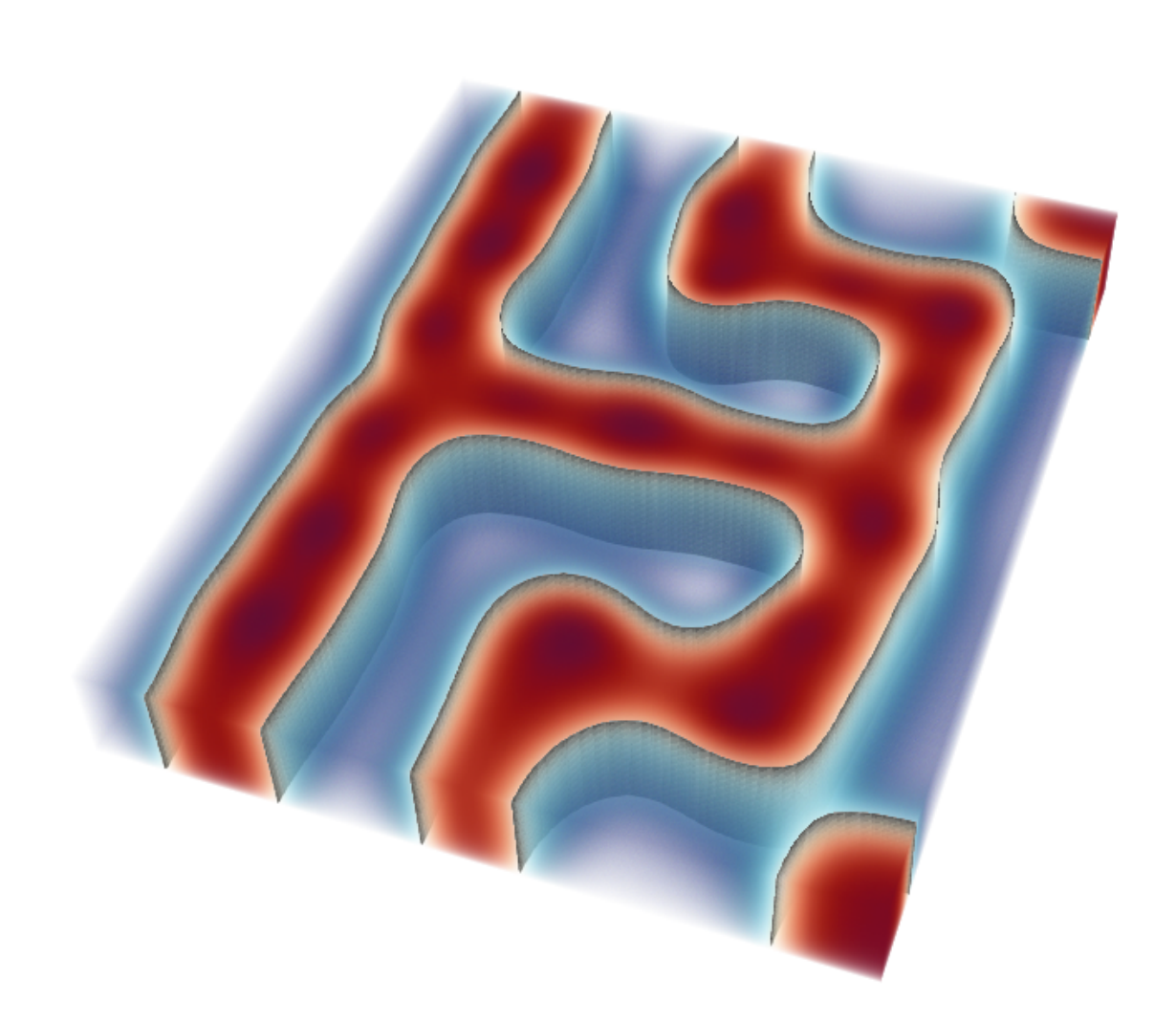}
        \caption{Sample state 2}
    \end{subfigure}
    \begin{subfigure}{0.32\textwidth}
        \centering
        \includegraphics[width=0.99\textwidth]{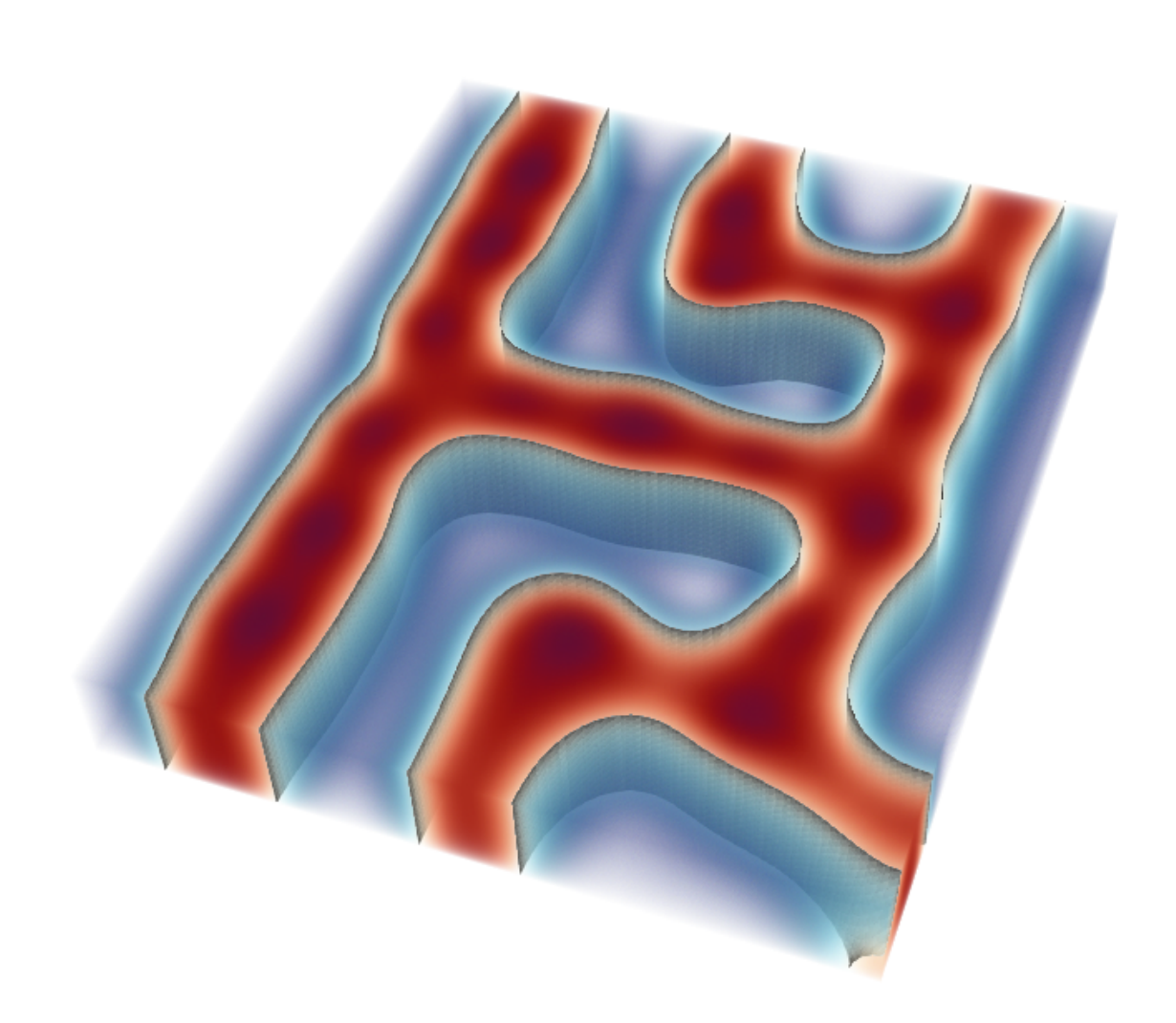}
        \caption{Sample state 3}
    \end{subfigure}
    \caption{3D junction target morphology, substrate interaction field at the optimal design, and three sample equilibrium states for different initial guesses at the optimal design. The one with minimum energy is shown in bottom left.  % Observe the failure modes that are now possible due to three dimensional effects.
    }
    \label{fig:junction3d_results}
\end{figure}

\begin{figure}[htbp!]
    \centering
    \begin{subfigure}{0.32\textwidth}
        \includegraphics[width=0.99\textwidth]{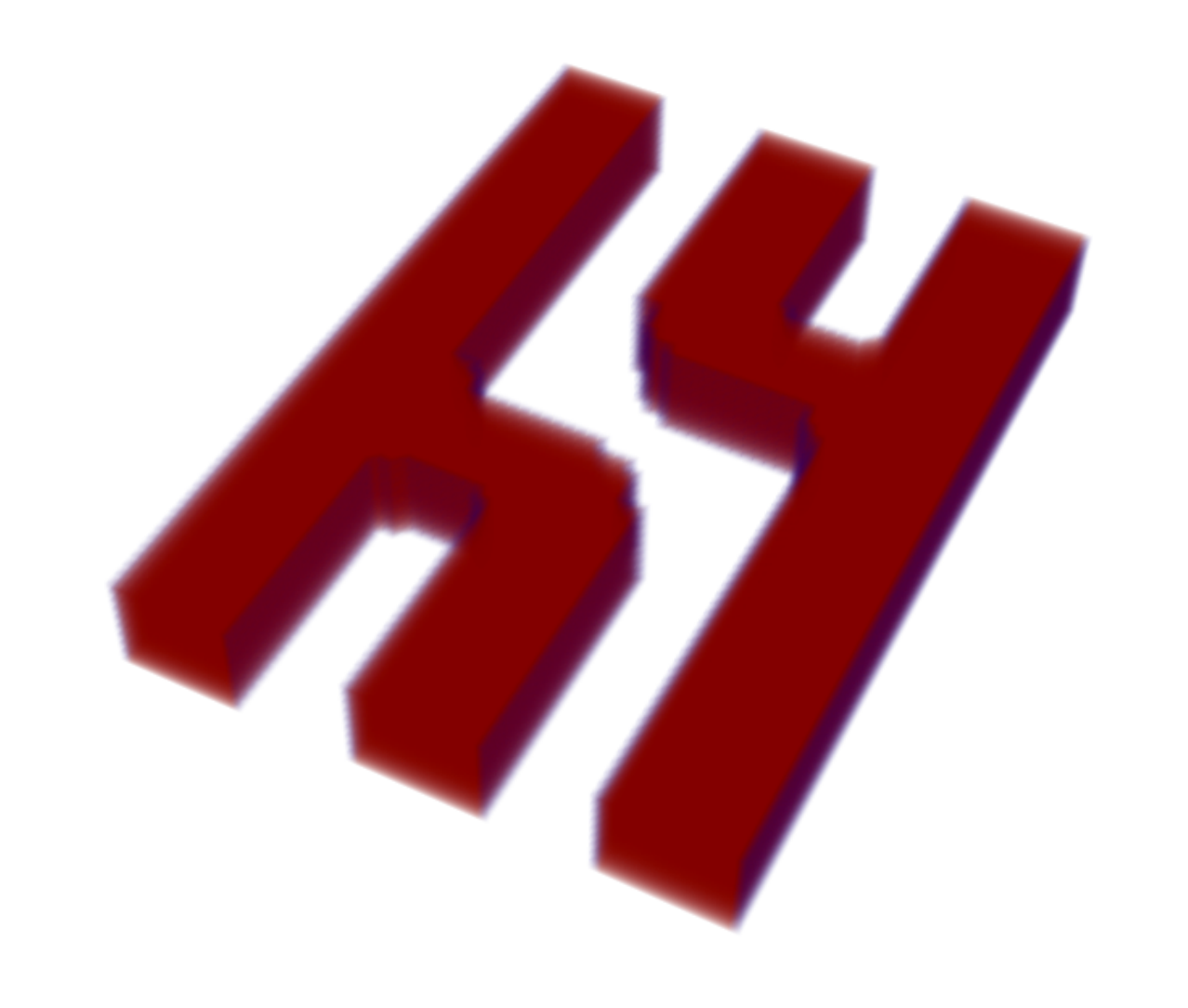}
        \caption{3D junction target morphology}
    \end{subfigure}
    \begin{subfigure}{0.32\textwidth}
        \centering
        \includegraphics[width=0.99\textwidth]{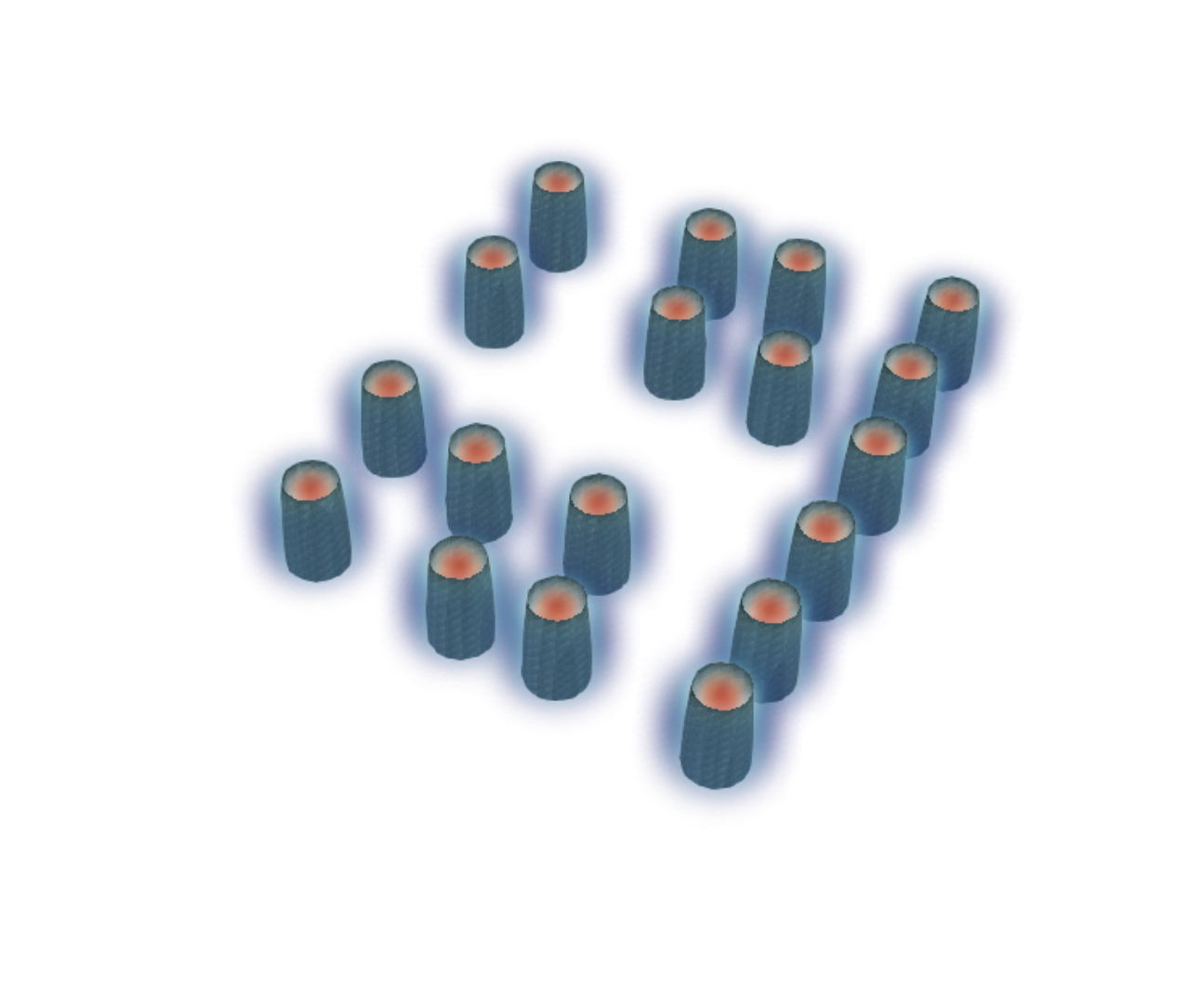}
        \caption{Substrate strength at optimal design}
    \end{subfigure}

    \begin{subfigure}{0.32\textwidth}
        \centering
        \includegraphics[width=0.99\textwidth]{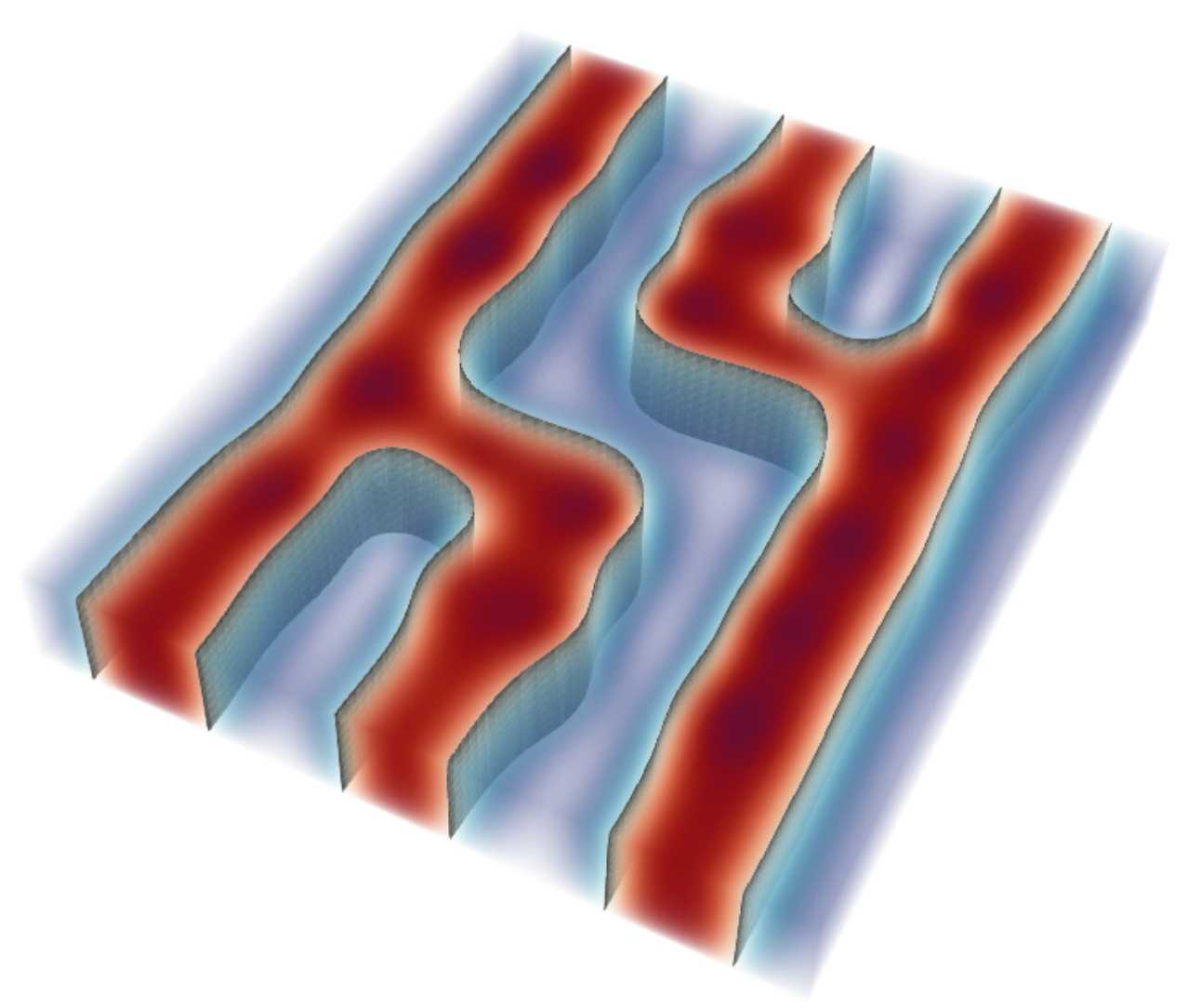}
        \caption{Minimum $\mathcal{F}(u)$ sample}
    \end{subfigure}
    \begin{subfigure}{0.32\textwidth}
        \centering
        \includegraphics[width=0.99\textwidth]{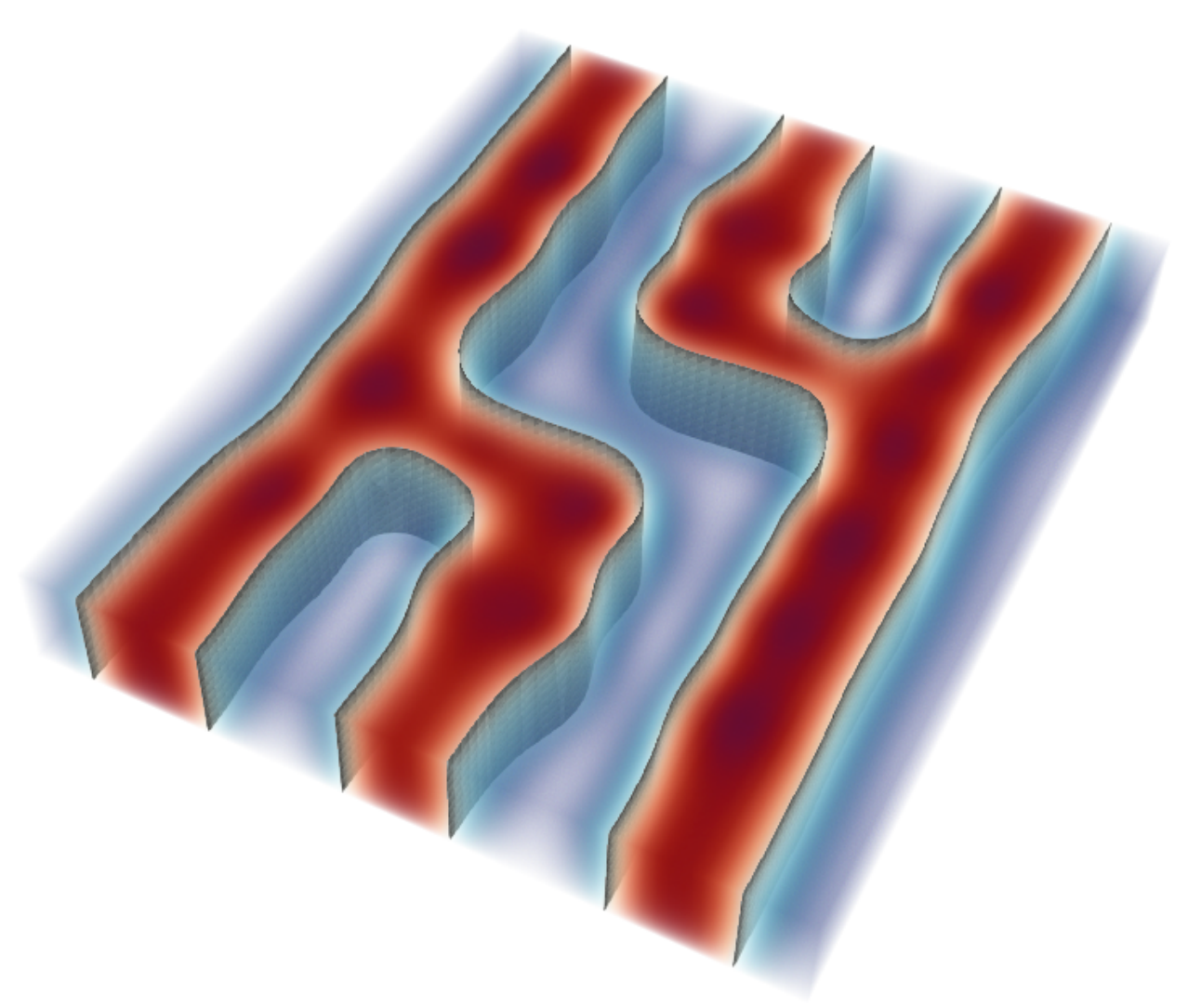}
        \caption{Sample state 2}
    \end{subfigure}
    \begin{subfigure}{0.32\textwidth}
        \centering
        \includegraphics[width=0.99\textwidth]{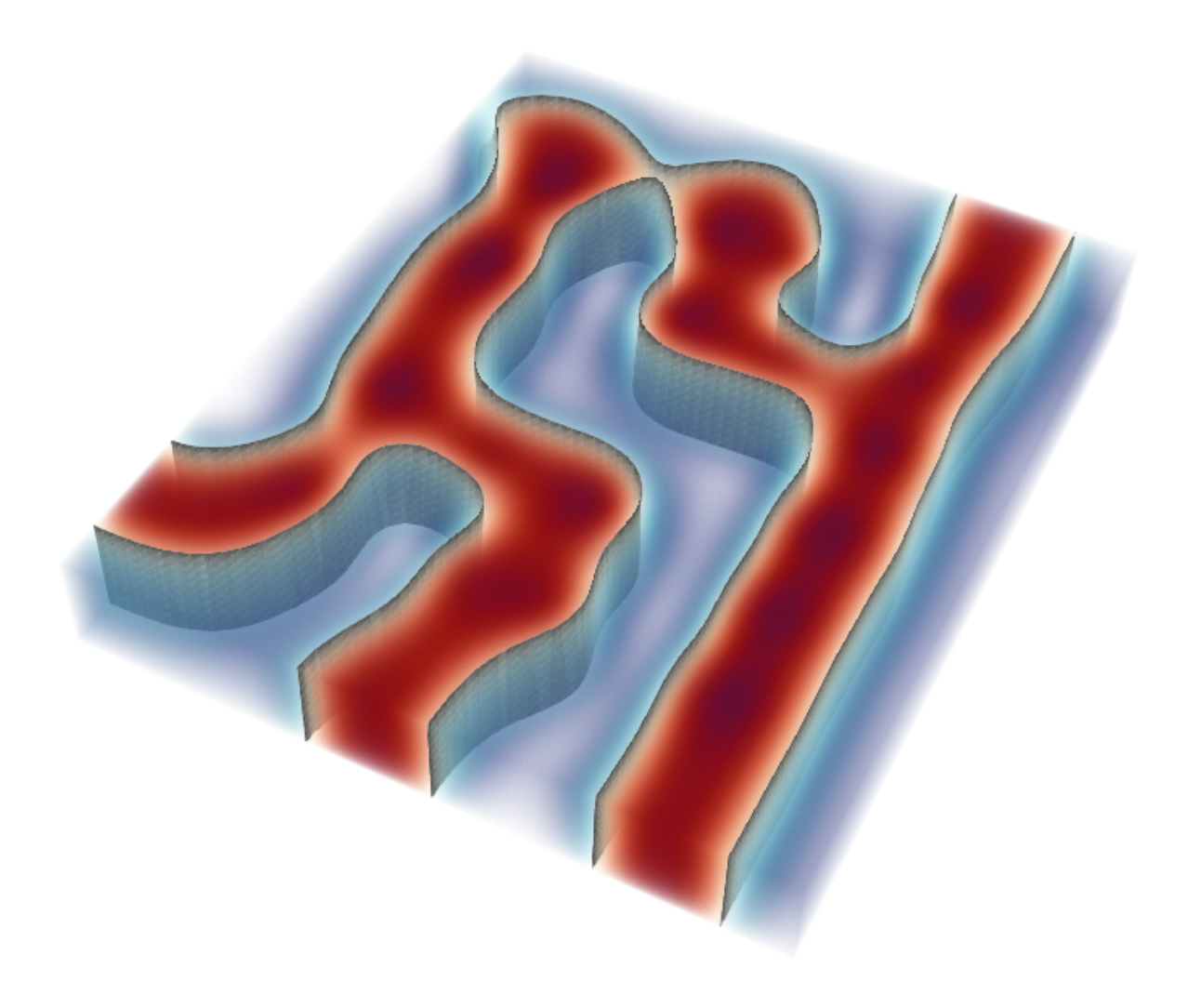}
        \caption{Sample state 3}
    \end{subfigure}
    \caption{3D jog morphology, substrate interaction field at the optimal design, and three sample equilibrium states for different initial guesses at the optimal design. The one with minimum energy is shown in bottom left.  % Observe the failure modes that are now possible due to three dimensional effects.
    }
    \label{fig:jog3d_results}
\end{figure}

We present the optimal design and sample states in Figure \ref{fig:junction3d_results} for the junction morphology and Figure \ref{fig:jog3d_results} for the jog morphology. In general, the 3D problem is more difficult than the 2D problem as the extra dimension introduces additional modes for defects. Nevertheless, in our results the sample state with the minimum energy closely resemble the target morphology. Moreover, these states can be consistently produced by the optimal designs across different initial guesses. 

% As seen in the junction design, we observe sample states that are defective due three-dimensional effects, which cannot manifest in 2D approximation of the domain. In general, we observe more defective sample states computed from random initial guesses at the optimal design for same number of guideposts in the 3D case compared to the 2D approximation. 

% In fact, we observe more defective sample states computed from random initial guesses at the optimal design due to the higher degrees of freedom in the three-dimensional domain than in the two-dimensional domain for the same number of guideposts. 
% In these results, we are able to produce accurately the 3D target morphology in the minimum energy sample. 

% \pc{Do you have examples for the 3D Bend design? Does that one look better? You can also report that one if you have.}
% As three-dimensional effects are now accounted for, we observe new modes of defects. 
% For example, in the third sample state, we notice the region on the right where the BCP remains connected on the bottom surface but not on the top.

\section{Conclusion}\label{sec:conclusion}
In this work, we developed a PDE-constrained optimization computational framework for the chemoepitaxy process for DSA of BCPs, utilizing a phase field model based on the Ohta--Kawasaki free energy. Specifically, we parameterized the substrate pattern using guideposts defined by shape functions and optimized for their locations. We proposed an adjoint-based optimization method---a modified inexact Newton-CG algorithm---to solve the optimization problem, which was demonstrated to converge rapidly, with the number of optimization iterations scaling almost independent of the number of guideposts.
% rapidly to the optimal designs. 
% This formulation allowed for the use of an adjoint based methods, and hence a modified inexact Newton-CG optimization algorithm which led to rapid convergence to optimal designs. 

% Significant computational reduction (compared to conventional methods) was achieved by our approach, due to 
As a result of the rapid convergence of both the outer Newton optimizer for the design optimization and the inner Newton solver for the energy minimization, along with the use of a continuation scheme, our method has been able to achieve significant computational efficiency, delivering five orders of magnitude speedup relative to previous work. The continuation scheme employed uses the state solution from the previous design iteration to initialize the energy minimization for the subsequent design iteration. This leads to a significant reduction in the number of iterations needed to arrive at the equilibrium state. The efficiency of this approach permits the design for non-trivial target morphologies in 3D, where the computational burden associated with repeatedly solving the state problem would otherwise render design optimization prohibitive.

%  In particular, 
% computational cost of the energy minimization was further reduced by the continuation scheme used in optimization, where the solution from previous optimization iteration was used as an initial guesse for the next step. 

A key challenge in this problem stemmed from the non-uniqueness of the equilibrium states, which in our formulation manifested as a dependence on the choice of initial guesses used in the energy minimization. This meant that robustness of any optimal design needed to be assessed by demonstrating its ability to produce the target morphologies over a range of initial guesses. The non-uniqueness also led to discontinuities in the objective function whenever changes to the design variable moved the equilibrium state from one solution branch to another.

%In this particular work, we adopted a largely deterministic approach 
In addressing these issues, we adopted a deterministic approach:
During optimization, we used the continuation scheme to remain in the same basin of attraction as the equilibrium solution corresponding to the previous design iteration, and thus preserve smoothness of the objective function. An assessment of the robustness was then made upon convergence to the optimal design by solving the state problem with random samples of initial guesses. In particular, the equilibrium state with minimum free energy among the sample equilibrium states was used as an indicator of the most likely state, while the remaining samples were used to give an indication of the robustness and possible defective states. We observed that while most optimal designs were able to produce the target morphology, the robustness to initial guesses largely depends on the number of guideposts used. The use of more guideposts tends to lead to a lower likelihood of defects in the equilibrium states. 

The deterministic approach presented in this study is effective while being computationally inexpensive. However, the robustness of the design was not accounted for during optimization. In future work, we look to explicitly account for the robustness 
% during optimization. This would lead to a 
by extending the deterministic formulation of the design optimization problem presented here to a stochastic one, in which optimization is carried out to minimize the occurrence defects based on the statistics of the design objective with respect to the distribution of random initial guesses. This opens up the new challenges of accurately and efficiently estimating these statistics during optimization.

\section{Acknowledgements}\label{sec:acknowledgements}
This work is supported in part by the U.S.\ Department of Energy,
Office of Science, Office of Advanced Scientific Computing Research
under award DE-SC0019303, and by the U.S.\ National Science Foundation
Division  of Mathematical Sciences under award 2012453. 

\bibliography{references}

\appendix
\section{The \texorpdfstring{$\hcirc$}{hcirc} and \texorpdfstring{$\hcircdual$}{hcirdual} spaces}\label{appendix:hcirc}

To rigorously formulate the free energy minimization, and in particular, interpret the inverse Laplacian operator $(-\Delta_N)^{-1}$, we utilize properties of the space $\hcirc$ and its dual. We present some of the fundamental properties of these spaces following \cite{CaoGhattasOden22}. Recall that $\hcirc$ defines the $H^1(\Omega)$ functions with zero average,
\beq
    \hcirc := \{ u \in H^1(\Omega): \int_{\Omega} u \dx = 0 \}.
\eeq 
This is a Hilbert space with inner product $(u, v)_{\hcirc} = \inner{\nabla u}{\nabla v}$, where $(\cdot, \cdot)$ denotes the $L^2(\Omega)$ inner product. We can then associate its dual, $\hcirc'$ with the following subspace of $H^1(\Omega)'$,
\beq
    \hcircdual := \{g \in H^1(\Omega)' : \langle g, c \rangle = 0 \quad \forall c \text{ const a.e.}\}
\eeq
We then define the inverse Laplacian by its weak form. For $w \in \hcirc$ and $g \in \hcircdual$, $w = (-\Delta_N)^{-1}g$ if 
\beq
    \inner{\nabla w}{\nabla v} = \langle g, v \rangle \quad \forall v \in H^1(\Omega).
\eeq
Thus, $(-\Delta_N)^{-1} : \hcircdual \rightarrow \hcirc$ is a continuous solution operator, and can be used to define an inner product on $\hcircdual$,
\begin{linenomath}
$$ (f,g)_{\hcircdual} :=((-\Delta_N)^{-1}f, (\Delta_N)^{-1}g)_{\hcirc} 
= \duality{f}{(-\Delta_N)^{-1}g} = \duality{g}{(-\Delta_N)^{-1}f}.$$

Finally, we have $R : \hcirc \rightarrow \hcircdual$ the $L^2(\Omega)$ Riesz map defined such that $R(u_0) = (u_0, \cdot) \; \forall u_0 \in \hcirc$. Using the definition of $(-\Delta_N)^{-1} : \hcircdual \rightarrow \hcirc$, the expression $(-\Delta_N)^{-1}(u-m)$ can also be written as $(-\Delta_N)^{-1}R(u-m)$. This leads to the definition of the nonlocal term as a norm,
\beq
\int_{\Omega} (u-m)(-\Delta_N)^{-1}(u-m)\dx = \duality{R(u-m)}{(-\Delta_N)^{-1}R(u-m)} = \|R(u-m)\|_{\hcircdual}^2.
\eeq
\end{linenomath}

\section{Continuity of the second derivative operator}\label{appendix:continuity}
We show the continuity of the second derivative operator assuming $\Omega \subset \mathbb{R}^{d}$ is a bounded Lipschitz domain. Recall that the second derivative operator is given by 
\beq
    \duality{D_u^2 \mathcal{F}(u)\hat{u}_0}{\tilde{u}_0} 
    = \inner{W''(u)\hat{u}_0}{\tilde{u}_0}
    + \inner{\varepsilon^2 \nabla\hat{u}_0}{\nabla\tilde{u}_0}
    + \inner{\sigma(-\Delta_N)^{-1}\hat{u}_0}{\tilde{u}_0} \quad \forall \hat{u}_0, \tilde{u}_0 \in \hcirc.
\eeq
where $W''(u) = 3u^2 - 1$. When using the formalism of $(-\Delta_N)^{-1} : \hcircdual \rightarrow \hcirc$, we write this more precisely as 
\beq
    \duality{D_u^2 \mathcal{F}(u)\hat{u}_0}{\tilde{u}_0} 
    = \inner{W''(u)\hat{u}_0}{\tilde{u}_0}
    + \inner{\varepsilon^2 \nabla\hat{u}_0}{\nabla\tilde{u}_0}
    + \inner{\sigma(-\Delta_N)^{-1}R(\hat{u}_0)}{\tilde{u}_0} \quad \forall \hat{u}_0, \tilde{u}_0 \in \hcirc.
\eeq
 We focus first on the double-well potential term, where $W''(u) = 3u^2 - 1$. By the generalized H\"{o}lder inequality, we have
\beq
    \inner{u^2\hat{u}_0}{\tilde{u}_0} 
    \leq \|u^2\|_{L^2(\Omega)}\|\hat{u}_0\|_{L^4(\Omega)}\|\tilde{u}_0\|_{L^4(\Omega)}
    = \|u\|_{L^4(\Omega)}^2\|\hat{u}_0\|_{L^4(\Omega)}\|\tilde{u}_0\|_{L^4(\Omega)}
\eeq
We note that for $d=2,3$, the Sobolev embedding theorem gives 
\begin{linenomath}
$$ \|v\|_{L^4(\Omega)} \leq C_s \|v\|_{H^1(\Omega)} \quad \forall v \in H^1(\Omega).$$
\end{linenomath}
Furthermore, for $u_0 \in \hcirc$, the Poincar\'{e} inequality implies the equivalence of the $\hcirc$ and $H^1(\Omega)$ norms so that
\begin{linenomath}
$$ C_1 \|u_0\|_{\hcirc} \leq \|u_0\|_{H^1(\Omega)} \leq C_2 \|u_0\|_{\hcirc}.$$
\end{linenomath}
Thus, 
\begin{linenomath}
$$ \inner{(3u^2-1)\hat{u}_0}{\tilde{u}_0} \leq C_w(u)\|\hat{u}_0\|_{\hcirc}\|\tilde{u}_0\|_{\hcirc}, $$
\end{linenomath}
where $C_w(u)$ is a constant depending on $u$. We now turn our attention to the nonlocal term. Using the Riesz map, we obtain
% Given $u_0 \in \hcirc$, we define the mapping $R: \hcirc \ni u_0 \rightarrow (u_0, \cdot) \in \hcircdual$.  Thus, we represent 
\begin{linenomath}
$$ \inner{\sigma(-\Delta_N)^{-1}R(\hat{u}_0)}{\tilde{u}_0} 
= \sigma\duality{R(\tilde{u}_0)}{(-\Delta_N)^{-1}R(\hat{u}_0)}
\leq \sigma \|R(\tilde{u}_0)\|_{\hcircdual} \|(-\Delta_N)^{-1}R(\hat{u}_0)\|_{\hcirc}.
% \leq \sigma \|(-\Delta_N)^{-1}\hat{u}_0\|_{\hcirc} \|(-\Delta_N)^{-1}\tilde{u}_0\|_{\hcirc}
$$
\end{linenomath}
By the continuity of the inverse Laplacian operator and the Riesz map, we have 
\begin{linenomath}
$$ \inner{\sigma(-\Delta_N)^{-1}\hat{u}_0}{\tilde{u}_0} \leq C_{\sigma} \|\hat{u}_0\|_{\hcirc} \|\tilde{u}_0\|_{\hcirc}.
$$
\end{linenomath}
Finally, continuity of the second term follows directly from the Cauchy-Schwarz inequality. This asserts the continuity of the bilinear form defined by the second derivative.
% of the bilinear form $b : \hcirc \times \hcirc \rightarrow \mathbb{R}$ where
% \begin{linenomath}
% $$ b(\hat{u}_0, \tilde{u}_0) = \duality{D_u^2 \mathcal{F}(u)\hat{u}_0}{\tilde{u}_0} \; \forall \hat{u}_0, \tilde{u}_0 \in \hcirc$$
% \end{linenomath}

\end{document}